\documentclass[12pt,a4paper,twoside,openright]{report}
\usepackage{amsmath,amsthm,amsfonts,amscd,amssymb,eucal,bbm}

\def\RR{{\mathbb R}}
\def\CC{{\mathbb C}}
\def\NN{{\mathbb N}}
\def\ZZ{{\mathbb Z}}
\def\diff{{{\rm Diff}^+(S^1)}}
\def\mob{{\rm M\ddot{o}b}}
\def\mo{{M\"obius }}
\def\pcw{{\rm Pcw}^1(\mob)}
\def\pcwm{{\rm Pcw}^1({\mathfrak m})}
\def\A{{\mathcal A}}
\def\B{{\mathcal B}}
\def\C{{\mathcal C}}
\def\D{{\mathcal D}}

\def\H{{\mathcal H}}
\def\I{{\mathcal I}}
\def\K{{\mathcal K}}
\def\M{{\mathcal M}}
\def\N{{\mathcal N}}

\def\P{{\mathcal P}}
\def\R{{\mathcal R}}

\def\T{{\mathcal T}}
\def\U{{\mathcal U}}
\def\V{{\mathcal V}}

\def\Z{{\mathcal Z}}
\newcommand{\fin}{{\mathcal D}\!_{f\!in}}
\newcommand{\resume}[1]{
\begin{flushleft}\bf \Large \vspace{-2mm} Chapter Summary
\end{flushleft} \noindent {\small #1}\newpage}
\newtheorem{theorem}{Theorem}[section]
\newtheorem{definition}[theorem]{Definition}
\newtheorem{corollary}[theorem]{Corollary}
\newtheorem{proposition}[theorem]{Proposition}
\newtheorem{lemma}[theorem]{Lemma}

\newcounter{paragrafus}[section]
\newcommand{\paragrf}{\addtocounter{paragrafus}{1}
\paragraph*{\S \arabic{paragrafus}.}}

\begin{document}
\begin{titlepage}
\hphantom{a} \vspace{0.4cm}
\begin{flushleft}
{\Huge  CONFORMAL COVARIANCE}
\end{flushleft}
\begin{flushleft}
{\Huge  AND RELATED PROPERTIES}
\end{flushleft}
\begin{flushleft}
{\Huge OF CHIRAL QFT}
\end{flushleft}
\begin{flushleft}
\vspace{3cm} {\Large Mih\'aly Weiner} \\
{\large \vspace{10cm} PhD thesis in Mathematics \\
supervisor: Roberto Longo \\
\vspace{-4mm}\hrulefill \smallskip \\
Universit\`a di Roma ``Tor Vergata'', Dip.\! di Matematica, 2005}.
\end{flushleft}
\end{titlepage}
\thispagestyle{empty} \cleardoublepage 
\thispagestyle{empty}
\hphantom{a} \vspace{5cm}
\begin{flushright}
{\large
{\it Though this be madness, yet there is method in 't.}\\
\noindent (William Shakespeare: {\it Hamlet}, Act 2., Scene I\!I.)
}\end{flushright}
\newpage
\thispagestyle{empty} \cleardoublepage \begin{flushleft} 
{\LARGE
PREFACE}
\end{flushleft} \smallskip
\setcounter{page}{1}

\noindent This PhD thesis in Mathematical Physics contains both a
general overview and some original research results about low
dimensional conformal Quantum Field Theory (QFT). The closer area
of the work is chiral conformal QFT in the so-called algebraic
setting; i.e.\!\! the study of \mo covariant and diffeomorphism
(or, as it is sometimes called: \!{\it conformal}) covariant local
nets of von Neumann algebras on the circle. From the mathematical 
point of view, it involves mainly functional analyses, operator algebras, 
and representation theory (and in particular, the representation theory 
of infinite dimensional Lie groups). 

The preliminary chapters give an overview, summarize some of the
most important known facts and give some examples for such local
nets. This is then followed by the presentation of the new
results.

The questions considered are all about conformal covariance.
First, what algebraic property of the net could ensure the
existence of an extension of the \mo symmetry to full
diffeomorphism symmetry? Second, if exists, is this extension of
the symmetry unique? Third, what (new, so far unknown) properties
of the net are automatically guaranteed by the mere existence of
conformal covariance?

In relation with the existence, the author pinpoints a certain
(but purely algebraic) condition on the ``relative position'' of
three local algebras. (This condition is implied by diffeomorphism
covariance and even without assuming diffeomorphism covariance it
can be proved to hold, assuming complete rationality, for
example.) Then it is shown that the \mo symmetry of a regular net
(i.e.\! a net which is $n$-regular for every $n\in\NN^+$; a
property which is weaker than strong additivity and expected to
hold always in case of conformal covariance) can be extended to
certain nonsmooth transformations if and only if this condition is
satisfied, and that in this case the extension is unique. Although
the existence of an even further extension to the full
diffeomorphism group is not proved, for the author, as it will be
explained, it indicates that the introduced condition together
with regularity indeed ensure the existence of conformal
covariance.

Two different proofs will be given to show that there is at most
one way the \mo symmetry of a net can be extended to $\diff$. One
of them is close to the setting of Quantum Fields (and in particular, 
to Vertex Operator Algebras, as it is formulated at the level of Fourier
modes) in its spirit: it is based on the observation that there are only 
few possibilities for commutation relations between fields of dimension 
$2$, whose integral is the conformal energy. (So it is in some sense a
version of the L\"uscher-Mack theorem.)

This idea will be also used to show that every (\mo covariant) subnet of 
a conformal net is conformal covariant. In fact it will be proved that 
given a subnet, the stress-energy tensor can be decomposed into two in
which one part belongs to the subnet while the other part to the
coset of the subnet.

The other argument needs $4$-regularity but it is closer to the
idea of local algebras and it is a simple application of those
nonsmooth symmetries discussed before. This argument, too, implies
the automatic conformal covariance of a subnet, given that it is
strongly additive, for example.

Further applications of these results, such as the construction of
new examples of nets admitting no $\diff$ symmetry, will be also
discussed. Also, a model independent proof will be given for the
commutation of vacuum preserving gauge symmetries with
diffeomorphism symmetries.

In connection with nonsmooth symmetries a new technique is
developed: the smearing of the stress-energy tensor with nonsmooth
functions. This technique will be also used to prove that every
locally normal representation of a conformal net is of positive
energy.

Some of these results form the base of two articles. The first one
is a joint work and has been already published in Commun.\!
Math.\! Phys., while the second one is a single-author work
accepted by the same journal.

\tableofcontents

\pagestyle{myheadings}

\chapter{Introduction}

\markboth{CHAPTER \arabic{chapter}. INTRODUCTION} {CHAPTER
\arabic{chapter}. INTRODUCTION}

\section{Why doing QFT in low dimensions}

Quantum Physics, despite the mathematical, philosophical and
interpretational problems, is truly amazing and proved to be
extremely fruitful at the practical level (e.g.\! applications to
semiconductors) as well as at the theoretical level. Of course the
expression ``Quantum Physics'' is wide, and covers a variety of
theories and models. Usually it means only that the {\it
uncertainty principle} is taken into account and so a Quantum
Probabilistic description is used. More specifically we talk about
Quantum Mechanics (QM), Quantum Thermodynamics, the Quantum Theory
of Fields (QFT), etc.

Quantum Mechanics at its birth --- just like almost all new
theories --- contained a large number of heuristic and (therefore)
not too precise, sometimes awkward looking calculations
(involving Dirac's delta ``function'', ``noncommutative
variables'', decomposition into the ``sum of generalized
eigenvectors'', and many other things). Finally these
constructions were given a clear meaning and a completely new
mathematics was born: namely, the theory of Hilbert spaces and
unbounded operators.

One may compare this with the enormous inspiration that Newton's
mechanics gave to differential and integral calculus. The physical
motivation had probably a key role in the development of the new
mathematics. Of course in the beginning the concept of the
differential was somewhat heuristic and involved considerations
with ``infinitesimal'' quantities. It was only much later that with
the invention of limits it became something mathematically
rigorous.

In case of QM a good part of the new mathematics became clear
almost instantaneously. Indeed, it was in $1932$ that von Neumann
published his work \cite{neumann} about the mathematical
foundation of QM; only some $7$ years later that Heisenberg and
Schr\"odinger independently from each other proposed what at the
time was called ``matrix-'' and ``wavemechanics''.

Similarly, at its birth QFT suffered from all sorts of problems.
An important step towards the solution was done with the invention
of {\it renormalization}, first introduced by Schwinger, Feynman,
Tomonaga and Dyson in the late $40'$s for the special case of
Quantum Electrodynamics (QED). It lead to a perturbation expansion
whose lowest order terms are in excellent agreement with
experimental physics \cite{weinberg}, e.g.\! the remarkable
prediction of the electron's anomalous magnetic momentum.
Nevertheless, from the theoretical point of view it gave rise to
new complications that till now, after more than fifty years have
no satisfactory solution.

So the current state is the following. Using renormalization and a
perturbative approximation (containing many mathematical
obscurities) these theories can be used to give results that are
comparable with experimental values. Yet at the mathematical level
we do not know whether these theories {\it exist} at all; so for
examples that the formulas involved really converge (even if their
value with some heuristic calculations or computer simulations was
approximated).

Of course, in order to prove mathematically the existence of a
certain theory, first one needs to {\it define} what, in general,
a QFT should be. There are essentially two (accepted) axiomatic
frameworks. The first of the two (in the historical order) is
essentially about fixing the notion of a quantum field and
postulating a set of required properties which are usually called
the {\it Wightmann-axioms}, see e.g.\! the book \cite{wightmann}.
The second one is based on the concept of the {\it local structure
of observables} which is described by a net of operator algebras.
It is usually called the framework of Algebraic Quantum Field
Theory, for a good introductory book see \cite{Haag}.

Works done in these frameworks have greatly improved our
understanding of some general structural properties, such as for
example the existence of PCT symmetry, the correspondence between
spin and statistics and the existence of superselection sectors.
However, while a solid base, an axiomatic framework in which
questions can be precisely formulated is important to have, in
itself it is not enough for solving all problems. The frameworks
are of course free from those ``mathematical uncertainties'': the
problem appears in a different way. Namely, now the difficulty
lies in reconstructing within these frameworks such (only
heuristically described) theories as QED or QCD (\,=\, ``Quantum
Chromodynamics''; the quantum field theory of strong
nuclear interactions).

In fact, the only type of models (so far) constructed in these
rigorous manners are the so-called {\it free} ones. From the
physical point of view, they have the most simple structure possible, 
as they describe non-interacting fields and particles.  Although
these models deserve both some physical and mathematical interests,
they provide a very small variety of examples.

Let us immediately make a correction: this is the current
situation concerning models ``living'' in the  the usual
$4$-dimensional spacetime. However, it turned out that in $2$
dimensions a great variety of interesting models can be
constructed! A key difference --- at least for massless 
theories\footnote{In the massive case the conformal structure
does not have a direct influence. Nevertheless, the simple geometry of 
the $2$-dimensional spacetime still enables us to exhibit interesting
models. In fact, in $2$ dimensions there examples of massive
models with nonzero interactions (!), 
see e.g.\! \cite{glimm-jaffe} and also \cite{Bu+PhD} for recent 
developments. This is in some sense is even ``better'' 
than what we have in the conformal case.} --- lies in the {\it conformal} 
structure of the $2$-dimensional Minkowski-spacetime (of which more will 
be said in the next section) and the geometry of lightcones. While in 
higher dimensions the future-like lightcone of a spacetime point is
connected, in $2$ dimensions it consists of two disjoint
halflines. To illustrate the importance of this with a classical
(i.e.\! not quantum physical) example, consider in coordinates the wave
equation of the free massless Klein-Gordon field in the $2$-dimensional
Minkowski-spacetime:
\begin{equation}
\square \Phi (t,x)\equiv (\frac{d^2}{dt^2}-\frac{d^2}{dx^2})
\Phi(t,x)=0.
\end{equation}
As it is well known, the general solution splits into two parts
``living'' on the {\it left} and {\it right} lightcones:
\begin{equation}
\Phi (t,x)=\Phi_L(t-x) + \Phi_R(t+x).
\end{equation}
These one-dimensional parts are called the {\it chiral}
components. In a similar way, a massless QFT in $2$-dimensions 
gives rise to two chiral QFTs. This thesis is concerned 
with such chiral components in the setting of Algebraic Quantum Field 
Theory.

This factorization in higher dimensions does not happen: there are
no nontrivial examples of fields in $4$-dimensions that could be
factorized in any meaningful sense into lower dimensional ones.
What is the importance then of considering QFTs in low spacetime
dimensions, when certain things change so radically with the
dimension? The direct physical motivation is connected to $4$
spacetime dimensions.

Of course, the question whether it is, or is not worth to make
research in a certain direction is also up to one's belief.
Nevertheless, there are several ``justifications''; i.e.\!
arguments indicating the importance of the work done in this area.

$1.$ These models provide a ``playground'' on which one can
collect experience --- coming from calculations carried out in
concrete models --- about such relevant mathematical objects as
nets of operator algebras. Later this could prove decisive in
constructing $4$-dimensional models.

$2.$ These models turned out to give examples for various
interesting mathematical phenomena, e.g.\! the so-called
half-sided modular inclusion of von Neumann algebras. Through
these models one can directly use his or her physical inspiration
to do great mathematics. In some sense, to many mathematical
problems (especially in operator algebras) they are like vision to
geometry. While it is certainly possible to do geometry completely
formally (without actually having any kind of vision of what we
are doing), it is clear how much {\it seeing} the problem can help
in resolving it.

$3.$ These models may turn out to have a more direct application
to ``real physics''. As it was discovered \cite{BPZ}, some of these models 
correspond to the most important $2$-dimensional statistical mechanical 
models. Also, there are arguments showing that 
starting from an (assumed) model in $4$-dimensions, by some kind of
projection or restriction we can extract one-dimensional parts. By
considering such a one-dimensional part {\it only}, of course we
probably loose information. Still, we may be able to predict {\it
some} properties of the original $4$-dimensional theory just by
using our knowledge about chiral components.

\section{About conformal covariance}

Consider a physical system. If it is {\it autonomous} then, for
example, a translation or a certain rotation, should take a
possible process of this system into (another) possible process.
In general, in the mathematical description it means that a model
of an autonomous physical system should carry a representation (an
action) of the symmetry group of spacetime.

This idea proved to be fundamental in physics. Already in
classical models it can be used, for example, to find
restrictions on possible equations (that could govern the
behaviour of some physical quantities) written in turns of
coordinates. Its importance became even more significant when
through Noether's theorem the relation was discovered between (the
generators of) the symmetries (of flat spacetime) and preserved
quantities like total energy-impulse of the system.

In quantum physics by Wigner's theorem \cite{wigner} continuous
symmetries appear as projective unitary representations acting on
the Hilbert space of the physical model. By the mentioned
connection between energy-impulse and symmetries and by the
fundamental concept of the {\it positivity of energy} these
representations are required to satisfy a certain positivity
condition. The classification of such positive energy
representations of the symmetry group of (the flat) spacetime gave
similar tools to find restrictions on possible quantum physical
models as for example the mentioned concept of ``coordinate
independence of equations'' in classical physics. For example, an
immediate application lead to the understanding of the possible
values of the spin of particles.

So far, everything said was about autonomous systems, in general,
and the word ``conformal'' appeared only in the title of the
section. Let us see now the particular case of {\it massless}
systems.

In general a {\it dilation} is not a symmetry of a physical
system. Consider for example a classical system of masspoints with
gravitation being the only interaction between them. It is not
hard to see that {\it rescaling} a solution of the equations of
motion will not give again a solution of the same equations, but
rather the ones obtained by a similar rescaling of the masses of
the masspoints. A dilation carries a possible process of the
system into a possible process of another system; so in this case
it is not a symmetry.

However, there are physical systems in which no specific masses
appear. The most common example in classical physics is the free
electromagnetic field (``free''\,=\, no charges or currents). An
electromagnetic wave is a solution of the equations of the free
electromagnetic system and indeed, rescaling an electromagnetic
wave gives another electromagnetic wave; i.e.\! again, a solution
of the same equations. One may say that the Maxwell equations with
zero charges and currents are {\it scale invariant}.

For such systems the symmetry group is larger and so it gives
stronger restrictions; in the same time it also reveals more. 
To put it in another way: in general, using more symmetries, more 
can be established.

Actually, even more symmetries should be taken account. In fact,
it is expected that for a massless system all spacetime
transformations preserving the {\it causal} structure (the
structure of past, future and spacelike separation) should be a
symmetry transformation, at least in a certain slightly weaker
sense\footnote{In Algebraic Quantum Field Theory this ``slightly
weaker sense'' means that such transformations are indeed
represented by unitary (or antiunitary) operators whose adjoint
actions indeed preserve the local structure of the theory (more
precisely: act in a covariant manner on the local algebras), but
in general do not preserve the vacuum state.}. Such symmetries are
called {\it conformal}.

An initial enthusiasm about conformal symmetries could be taken away by
the observation that in three or higher dimensional spacetime there are no
conformal transformations apart from the ones generated by the usual
symmetry transformations and the dilations.  Allowing local
transformations (that may not extend to the full 
spacetime), too, in addition one also finds the so-called {\it special 
conformal} transformations (and of course everything that they, together 
with the already listed transformations, generate). They indeed have some
interesting applications, but still, they only ``enlarge'' the symmetry
group by a finite many additional parameters.

The situation in dimension $2$ however, is completely different. Consider
the $2$-dimensional Minkowski-spacetime with coordinates $(t,x)\in \RR^2$
(where the speed of light is $1$) and take two orientation preserving
diffeomorphisms of the real line $\varphi_L,\varphi_R\in {\rm
Diff}^+(\RR)$. It is almost evident that the transformation defined by the
map $(t,x)\mapsto (t',x')$ where \begin{equation} t'-x' =
\varphi_L(t-x),\;\;\;t'+x' = \varphi_R(t+x) \end{equation} preserves the
causal structure. (See what happens to the worldline of a left and what
happens to a wordline of a right light signal under this transformation.)
So now, apart from the usual symmetry transformations we do not only have
dilations but an infinite variety of diffeomorphism as conformal
transformations! To put it in another way: in a $2$-dimensional conformal
theory we should have an action\footnote{It turns out that not all
diffeomorphisms give symmetries; certain groth conditions at infinity must
be satisfied.} of the group ${\rm Diff}^+(\RR)\times{\rm 
Diff}^+(\RR)$.

In the ``normal'' (i.e.\! $4$-dimensional) Minkowski-spacetime, the
number of parameters in the geometrical symmetry group of a system
is always finite. Assuming no characteristic masses, the symmetry
group ``grows'' as it becomes enriched for example by the dilations, 
but ``not much'': it will still have just a finite number of parameters.
Contrary to this situation, as we have seen, in $2$-dimensions we
can have infinitely many geometrical symmetry parameters; in each
light-direction we can have a full diffeomorphism symmetry.

As it was already explained, a chiral theory is a one-dimensional
part of a two dimensional field theory (the part ``living'' on a
worldline of a light signal). Under some conditions such a theory
can be ``carried over'' to the one-point compactification of the
line; that is, to the circle $S^1$. So by what was explained for
such chiral QFT we should expect a $\diff$ symmetry to appear.

Indeed, this is the case in many (although not all) examples. In
fact, one of the questions considered in this thesis: can we find a
condition that assures the existence of the diffeomorphism
symmetry? Even more, having found the right condition, can we
somehow reconstruct the action of the full conformal group on a
given chiral model in case it was originally ``hidden for us''?

To the opinion of the author this question is not only of
mathematical interest. To give an example, consider a 
massless field theory in $4$ dimensions. Compared with 
the massive case we have in ``more'' the dilational covariance but 
certainly not an additional infinite number of geometrical symmetry 
parameters. On the other hand, by taking an appropriate restriction to a 
certain lower dimensional hyperplane we may get a lower dimensional field
theory possessing an infinite number of conformal symmetry
parameters. However, {\it a priori} we do not know whether this is
indeed the case and even if it is so there seems to be no direct
way of describing the corresponding action. The point is that these
symmetries are not coming from the original QFT (as it did not
have them), they are not inherited; they are not obtained simply
by taking restrictions.

Of course, the really interesting question would be the following.
Suppose we manage to construct such geometrical symmetries for a
chiral theory obtained by restriction. They are not present in the
original $4$-dimensional field theory, they are not among the
``visible'' geometrical symmetries. There are many interesting
things which we can calculate for a conformal chiral QFT
relatively easily, exactly by making use of the conformal symmetry.
However, can we use the conformal symmetry which only appears in
the restriction to get some information about the original field
theory?

This problem is out of the scope of this thesis, though. This work
is concerned ``only'' with chiral QFT in itself and not in
relation to QFT in higher spacetime dimensions. It concentrates on
the following questions. $1.$ Can we find conditions that could
guarantee the existence of full diffeomorphism symmetry? $2.$ Is
the action of the diffeomorphism group unique? In fact, can we
find a construction resulting this action (if it exists at all)?
$3.$ What properties of a conformal chiral QFT can we derive just
by using diffeomorphism symmetry? (More precisely: what {\it new}
properties can we derive --- as there are already many such
results.)

\section{About this thesis}

The aim of this section is twofold. First, to explain the
structure of this thesis and to give a short summary of its main
results. Second, to clarify what is the author's own result, what
is result of a joint work with someone else, what is just a
recollection of known results, and finally, to acknowledge 
contributions from others.

After this introductory chapter where right now the reader is, two
preliminary chapters follow. They serve to give an overview of
some well known facts concerning the material. The first concerns
general structural properties and some examples of \mo covariant
local nets. Apart from the known results, given with the
appropriate references, with Theorem \ref{localcore} it will be
shown that for any \mo covariant local net $(\A,U)$ the set
$\A(I)\Omega$ contains a core for every positive power of $L_0$
(where $I,\Omega$ and $L_0$ stand for what they usually stand
for). This affirmation is important for many (technical) reasons,
and although the proof relies on standard arguments similar to the
ones used in the demonstration of the {\it Reeh-Schlieder}
property, to the knowledge of the author, it has been not stated
before.

The second chapter is a lengthy introduction to diffeomorphism
covariance, representation theory of $\diff$ and to the Virasoro
algebra. In many cases, statements are not only given with citing
the original articles as references, but with sketched or
indicated proofs. Moreover, it will be demonstrated that any
positive energy representation of $\diff$ is (equivalent to) a
direct sum of irreducible representations (Theorem
\ref{pos.diffrep=+irr.}). Again, the proof is almost just a copy
of other standard arguments (for example of those similar
affirmations about positive energy representations of ${\rm
SL}(2,\RR)$) but nevertheless, to the knowledge of the author it
has not appeared in previously published works.

The remaining three chapters contain the research results achieved
during the PhD studies. All of these three chapters begin with a
``chapter summary'' which is a similar abstract to the one given
here, except that it is specific to that chapter. It is only for
helping readability.

The first of these chapters is about the (possible) reconstruction
of diffeomorphism symmetry, or to put it in another way, about the
(algebraic) characterization of the existence of diffeomorphism
symmetry. It consists of results so far not yet published.

First, the group of (once) differentiable piecewise \mo
transformations is introduced. Then with Theorem
\ref{existence.pcw.symm.} it is shown that the \mo symmetry of a
regular net $\A$ (i.e.\! a net which is $n$-regular for every
$n\in\NN^+$) has an extension to this group of piecewise \mo
transformations if and only if a certain (but purely algebraic)
condition is satisfied regarding the ``relative position'' of
$\A(I_1),\A(I_2)$ and $\A(I)$, where $I_1,I_2 \subset I$ are the
intervals obtained by removing a point from the nonempty, nondense
open interval $I \subset S^1$.

The group of piecewise \mo transformations is indeed a large
symmetry group. For example, it is shown that any diffeomorphism
can be uniformly approximated by elements of this group; even in a
stronger, local sense (Corollary \ref{pcwdensity}). Thus, although
the symmetry is extended only to these piecewise \mo
transformations and not to the group $\diff$, assuming regularity
this algebraic condition is conjectured to imply diffeomorphism
covariance.

The relation of this condition with others is investigated. First,
it is confirmed that conformal covariance implies this condition.
However, the point is that it can be checked in many cases without
assuming conformal covariance. In particular it is shown to be
satisfied for a large class of regular nets including the
completely rational ones (Theorem \ref{3cond.->pcw}); indicating
that for example in case of complete rationality diffeomorphism
symmetry should be automatic.

The second of the last three chapters contains the essential
estimates of the joint work \cite{CaWe} of the author of this
thesis and S.\! Carpi. With Theorem \ref{theo.main.nonsmooth} and
Prop.\! \ref{f->T(f)continuity} it is shown that the stress-energy
tensor of a positive energy representation of $\diff$ can be
continuously (in the sense of strong resolvents) extended to
functions of finite $\|\cdot\|_{\frac{3}{2}}$ norm. The fruits of
the extension theorem are harvested only in the last chapter, here
the focus is on the extension itself. The crucial point of the
argument is a certain estimate involving the contraction semigroup
associated to the conformal Hamiltonian.

The last chapter contains several applications of the results of
the previous chapters as well as some new arguments. By the
extension theorem the stress-energy tensor can be directly
evaluated on piecewise \mo vector fields. On the other hand, the
piecewise \mo symmetry of a net is completely determined by the
\mo symmetry and the local algebras. This can be used to give an
elegant and short demonstration for the uniqueness of the
diffeomorphism symmetry, given that the net is $4$-regular
(Theorem \ref{4reg->unique}).

This was already stated in the mentioned joint work, although with
a slightly different proof. There the problem regarding the
density of piecewise \mo transformations was cleverly avoided as
by the minimality of Virasoro nets a single non-trivial piecewise
\mo transformation turned out to suffice for the proof. However,
in order to make that argument work, a separate proof was needed
for the case of Virasoro nets.

In this thesis minimality is not used, and nor there is a need for
a separate argument treating the case of Virasoro nets. Rather,
the question of density is directly approached. For example, it is
proved that for every $\gamma\in\diff$ a sequence of piecewise \mo
transformations $\gamma_n\; (n\in\NN)$ converging uniformly to
$\gamma$ can be constructed having the property that
$U(\gamma_n)\to U(\gamma)$ strongly for any regular conformal
local net $(\A,U)$ (Theorem \ref{density.in.right.sense}). With
the same technique it is also shown that a \mo covariant subnet of
a conformal net is automatically a conformal subnet, too, assuming
again some regularity condition (Theorem \ref{4reg->covsubnet}).

In this thesis both the uniqueness and the conformality of \mo
covariant subnets are demonstrated by a completely different
argument, too: in fact even without assuming any kind of
regularity condition and thus in full generality (Theorem
\ref{theo:uniq2} and Theorem \ref{auto.conf.subnet.2a}). This
argument has not been published elsewhere. It is less related to
the idea of local algebras: in spirit it is closer Quantum Fields
(and in particular to Vertex Operator Algebras, as it is formulated
at the level of Fourier modes). It is based on the 
observation that there are only few possibilities for commutation 
relations between fields of dimension $2$ whose integral is the conformal 
energy. It is, in some sense, a version\footnote{Although found 
independently, the argument employed --- as it was pointed out to the 
author by Karl-Henning Rehren --- is essentially a Fourier-transformed 
version of an old idea of L\"uscher and Mack (see \cite[Sect.\! 
I\!I\!I.]{lmack} for their argument). What seems to be really new 
however, is the application in the algebraic setting, showing the 
mentioned uniqueness property and the conformal covariance
of \mo covariant subnets. Note that for such application one
needs a careful analyses involving the use of certain energy bounds
and in particular the relation of these bounds with the local structure.} 
of the L\"uscher-Mack theorem.

About the conformality of \mo covariant subnets an even stronger statement
is proved: namely, that the stress-energy tensor can be decomposed
into two in which one part belongs to the subnet while the other
part to the coset of the subnet (Theorem \ref{theo:T=TB+TB'}).

After solving these problems some applications are discussed. They
are mainly the same ones that appear in the mentioned joint work,
too. Namely, with Corollary \ref{inf.tens.product->nodiff} it is
shown how new examples (even strongly additive ones!) of nets
admitting no $\diff$ symmetry can be constructed. Also, a model
independent proof is given for the commutation of vacuum
preserving gauge symmetries with diffeomorphism symmetries
(Corollary \ref{com_int}).

Finally, as a further application of the developed technique of
evaluating the stress-energy tensor on nonsmooth functions, it is
proved that every locally normal representation of a conformal net
is of positive energy (Theorem \ref{mainresult}). This result is
part of the single-author work \cite{weiner1} accepted by
Commun.\! Math.\! Phys.

The thesis is also ``equipped'' with two appendices. The first one
is about the \mo group and its representation theory (so in fact
about the representation theory of $\widetilde{{\rm SL}(2,\RR)}$).
It is ``all-inclusive'': everything is demonstrated rather then
just referred to some articles. The statements presented can be all found 
in basic textbooks. Nevertheless, the author found it convenient to 
collect together all these facts and reprove them directly in the 
formalism of local nets in which they are then applied.

The second appendix is about the ``geometry'' of piecewise \mo
transformations. It contains only elementary (but not necessarily
short and trivial) geometrical calculations, with no von Neumann
algebras or even functional analysis involved. So, in some sense a
good high-school student should be able to reproduce all that can
be found there. Nevertheless, as these facts are used when the \mo
symmetry of local net is extended to the piecewise \mo group, it
is included for ``correctness''.

The author would like to thank Roberto Longo for introducing him to the
area of chiral conformal nets, for proposing problems (in particular,
about the existence and uniqueness of diffeomorphism symmetry) and for the
numerous inspirating discussions made during the four years of his PhD 
studies. He would also like to thank Sebastiano Carpi for the fruitful
research made together. He is moreover grateful to all members of 
the PhD committee for the careful reading of (an earlier version of) 
his thesis, for finding and correcting several unprecise arguments,  
and for the useful comments given. In particular, he would like to
thank Karl-Henning Rehren for pointing out the similarity between one 
of his arguments and that of the L\"uscher-Mack theorem.

\chapter{Preliminaries I: local nets}

\markboth{CHAPTER \arabic{chapter}. LOCAL NETS}{CHAPTER
\arabic{chapter}. LOCAL NETS}

\section{Axioms for \mo covariant local nets}

\paragrf
We shall not dwell much more on the physical background. An
excellent introduction and general arguments about why and how QFT
is to be described as a net of operator algebras can be found in
\cite{Haag}, although the book does not discuss the special case
of chiral theories.

A chiral (massless) QFT, as it was shortly explained in the
introduction, is a one-dimensional part of a two-dimensional
theory; it ``lives'' on the worldline of a light signal. Under
some conditions it can be carried over to the one-point
compactification of the line which is modelled by the unit circle
$S^1 =\{z\in \CC : \,|z|=1 \}$. These conditions, together with
the connection between the theory on the line and on the circle,
will be shortly discussed in Sect.\! \ref{sec:furtherprop}.
Although the setting on the line is closer to the physical
motivation, from the mathematical point of view it is more
convenient to consider and axiomatize chiral (massless) QFT on the
circle.

\paragrf
Throughout the thesis $\I$ will stand for the set of open,
nonempty and nondense intervals of the unit circle $S^1$. (The
word ``arc'' would be more appropriate but by now ``interval'' is a
convention. It also unites the terminology of field theory on the
circle and that one on the line.)

A {\bf \mo covariant local net of von Neumann algebras on $\mathbf
{S^1}$} is a map $\A$ which assigns to every $I \in \I$ a von
Neumann algebra $\A(I)$ acting on a fixed complex, Hilbert space
$\H_\A$ (``the vacuum Hilbert space of the theory''), together
with a given strongly continuous representation $U$ of $\mob
\simeq {\rm PSL}(2,\RR)$, the group of \mo
transformations\footnote{diffeomorphisms of $S^{1}$ of the form $z
\mapsto \frac{az+b}{\overline{b}z+\overline{a}}$ with $a,b\in
\CC$, $|a|^2-|b|^2=1$. See more about the group $\mob$ and its
representation theory in appendix \ref{chap:A}.} of the unit
circle $S^1$ satisfying for all $I_1,I_2,I \in \I$ and $\varphi
\in \mob$ the following properties.

\begin{itemize}
\item[(i)] {\bf Isotony.} $I_1 \subset I_2 \,\Rightarrow\, \A(I_1)
\subset \A (I_2).$

\item[(ii)] {\bf Locality.} $I_1 \cap I_2 = \emptyset
\,\Rightarrow\, [\A(I_1),\A(I_2)]=0.$

\item[(iii)] {\bf Covariance.}
$U(\varphi){\A}(I)U(\varphi)^{-1}={\A}(\varphi(I)).$

\item[(iv)] {\bf Positivity of the energy.} The representation $U$
is of positive energy type: the conformal Hamiltonian $L_0$,
defined by $U(\rho_\alpha)=e^{i\alpha L_0}$ where $\rho_\alpha \in
\mob$ is the anticlockwise rotation by an angle of $\alpha$, is
positive. Thus in particular ${\rm Sp}(L_0) \subset \NN$.

\item[(v)] {\bf Existence and uniqueness of the vacuum.} Up to
phase there exists a unique $\Omega \in \H_\A$ called the ``vacuum
vector'' which is invariant for $U$. Equivalently: the eigenspace
of $L_0$ associated to $0$ is one-dimensional.

\item[(vi)] {\bf Cyclicity of the vacuum.} ${\Omega}$ is cyclic
for the von Neumann algebra $\A(S^1)\equiv (\cup_{I\in\I} \A(I))''$.
\end{itemize}

\paragrf
Suppose we drop the last axiom. Then by restriction to the Hilbert
space $\overline{\A(S^1)\Omega}$ we get a \mo covariant net
satisfying all of the above axioms. Thus point (vi) says that the
Hilbert space is not ``unnecessarily'' too big. We shall now see
that it follows that the vacuum vector is cyclic even for a single
local algebra $\A(I),\; I\in\I$. This is a well-known general
phenomenon not specific to chiral QFT, although in order to prove
it, usually one needs to require some kind of {\it additivity} for
the net. In our case, due to the \mo symmetry, we need no
additional assumptions.

All through the thesis $\rho,\tau,\delta$ stand for the
one-parameter group of rotations, translations and dilations in
$\mob$, respectively. Their precise definition can be found in
appendix \ref{app:mobgroup} but the reader may also look at the short
explanation in the beginning of the next section. The only thing
to mention here is that the so-called infinite point in this
thesis is chosen to be the point $1\in S^1$, which is probably not
the most conventional choice, but has certain advantages when it
comes to coverings of $\mob$. Regarding the scalar product
$\langle\cdot,\cdot\rangle$ on the Hilbert space $\H_\A$ the
physicist convention is used according to which it is linear in
the second and antilinear in its first variable.

\begin{theorem}\!{\rm \cite{FrG, FJ} {\bf ``Reeh-Schlieder property''.}}
$\Omega$ is cyclic and separating for each local algebra $\A(I)$,
$I\in\I$.
\end{theorem}
\begin{proof}
By {\it locality} it is enough to prove cyclicity. It is also
clear that by {\it covariance} if cyclicity holds for a certain
interval then it holds for all as $\mob$ acts transitively on the
set of intervals $\I$. So we may assume that $I$ contains the
point $1 \in S^1$.

Suppose $\Psi \in \H_\A$ is an orthogonal vector to the set
$\A(I)\Omega$. Let us choose a smaller open nonempty interval
$K\subset I$ which has $1$ as one of its endpoint and whose
closure is still contained in $I$. For an operator $B \in \A(K)$
consider the function $f_B$ defined on $\RR \times \RR$ by
\begin{equation}
f_B(\alpha,a)\equiv \langle\Psi,\, U(\rho_\alpha) U(\tau_a)
A\Omega\rangle = \langle\Psi,\, U(\rho_\alpha\circ\tau_a) A U
(\rho_\alpha\circ\tau_a)^* \Omega\rangle.
\end{equation}
It is well known about positive energy representations of $\mob$
that not only the self-adjoint generator of the rotations, but
also that of the translations is positive (see Theorem
\ref{pos.mob^n} in the appendix). Therefore $f_B$ has a unique
continuous extension $F_B$ to the domain $(\RR + i \RR^+_0) \times
(\RR + i \RR^+_0) $ which is holomorphic on $(\RR + i \RR^+)
\times (\RR + i \RR^+)$. On the other hand, by the geometric
position of $K$, there exists an $\epsilon >0$ such that if
$|\alpha|,|a|<\epsilon$ then $\rho_\alpha (\tau_a(K)) \subset I$
and thus $f_B(\alpha,a)=0$ as by {\it covariance}
\begin{equation}
U(\rho_\alpha \tau_a) B U (\rho_\alpha \tau_a)^* \in
\A(\rho_\alpha \tau_a(K)) \subset \A(I).
\end{equation}
Hence the function $F_B$ is everywhere zero.

Suppose $L\in\I$ is an arbitrary interval. Depending on the sign
of its parameter a translation either increases or decreases the
interval $K$. By choosing $a\in\RR$ right, the length of
$\tau_a(K)$ can be made exactly equal to the length of $K$. Thus
there exist two parameters $\alpha_0,a_0 \in \RR$ such that
$\rho_{\alpha_0} (\tau_{a_0}(K))=L$. By this, and the previous
considerations about the function $F_B$ one can easily show
that $\left(\A(I)\Omega\right)^\perp =
\left(\A(L)\Omega\right)^\perp$ for every interval $L\in\I$. This
means that the subspace $\overline{\A(L)\Omega}$ is independent of
$L\in\I$. Hence the orthogonal projection $P$ to this subspace
belongs to $\A(S^1)'$ as $P \in \A(L)'$ for each $L\in\I$.
Therefore the range of $P$ is invariant for $\A(S^1)$ and as it
contains the vector $\Omega$, by the {\it cyclicity} axiom it must
be the full Hilbert space $\H_\A$ and so $\overline{\A(I)\Omega} =
{\rm Ran}(P)=\H_\A$.
\end{proof}

\paragrf
We shall now see that the set $\A(I)\Omega$ is not only dense for
each interval $I\in\I$ but it even contains a core for $L_0^k$
where $L_0$ is the conformal Hamiltonian and $k$ is any positive
number. For $k=1$ this was shown in the appendix of
\cite{Carpi99}. But to the knowledge of the author, for $k>1$ this has 
not been explicitly stated in the literature, although it can be proved
by standard arguments that were already used in the demonstration
of the Reeh-Schlieder theorem.
\begin{definition}
For an interval $I\in\I$ let $\A^\infty(I)$ be the star-algebra of
smooth elements (with respect to the rotations) in $\A(I)$; i.e.
$$
\A^\infty(I)\equiv \left\{A\in\A(I): \alpha\mapsto
U(\rho_\alpha)AU(\rho_\alpha)^*\;\;{\rm
is\;op.\!-\!norm\;smooth}\right\}.
$$
\end{definition}
It is almost evident that $\overline{\A^\infty(I)\Omega}=\H_\A$
for every interval $I\in\I$. In fact, we may choose a sequence of
nonnegative smooth functions on the real line
$\{\varphi_n:n\in\NN\}$ such that $\forall n\in\NN:\int\varphi_n
=1$ and their support contracts to the point $0$. Then for
$K\in\I,\,\overline{K}\subset I$ and $B\in\A(K)$ we have that the
norm-convergent integral
\begin{equation}
B_{\varphi_n} \equiv \int U(\rho_\alpha)BU(\rho_\alpha)^*
{\varphi_n}(\alpha) d\alpha
\end{equation}
for $n$ large enough belongs to $\A^\infty(I)$. Then by the
density of $\A(K)$ and the fact that $B_{\varphi_n}\Omega \to
B\Omega$ as $n\to\infty$ we can conclude that $\A^\infty(I)$ is a
dense subspace of $\H_\A$.
\begin{theorem}
\label{localcore} Let $k$ be a positive number. Then for every
interval $I\in\I$ the subspace $\A^\infty(I)\Omega$ is a core for
$L_0^k$.
\end{theorem}
\begin{proof}
We may choose a $K\in\I$ so that $\overline{K}\subset I$. We know
that both $\A^\infty(K)\Omega$ and $\A^\infty(K)\Omega$ is a dense
subset of $\D(L_0^k)$. Suppose $\Psi\in\H_\A$ is orthogonal to the
set $L_0^k \A^\infty(I)\Omega$ and $B\in\A(K)$. Then just like in
the proof of the Reeh-Schlieder theorem, there exists an
$\epsilon>0$ such that if $|\alpha|<\epsilon$ then
\begin{equation}
f_B(\alpha)\equiv \langle \Psi,\, e^{i\alpha L_0} L_0^k B
\Omega\rangle
\end{equation}
is equal to zero. Then, similarly to how we proceeded in the proof
of the mentioned theorem, we conclude that $f_B$ is everywhere
zero. Using this one can easily show that for $n\neq0$ the vector
$\Psi$ must be orthogonal to $E_n B \Omega$, where $E_n$ is the
orthogonal projection to the eigenspace of $L_0$ associated to the
eigenvalue $n\in {\rm Sp}(L_0)\subset \NN$. Thus by the
arbitrariness of $B\in \A^\infty(K)$, the uniqueness of the vacuum
and the density of $\A^\infty(K)\Omega$ in $\H_\A$ we find that
$\Psi$ must be a multiple of the vacuum vector.

So let $\Phi$ be a vector in the domain of $L_0^k$. Then
$L_0^k\Phi$ is orthogonal to the vacuum vector and thus by our
previous arguments it belongs to
$\overline{L_0^k\A^\infty(I)\Omega}=\{\Omega\}^\perp$\!. Hence
there exists a sequence $\{\Phi_n : n\in\NN\} \subset
\A^\infty(I)\Omega$ such that $L_0^k\Phi_n \to L_0^k\Phi$ as
$n\to\infty$ and evidently the sequence can be chosen so that
$\forall n\in\NN: \,\langle\Omega,\,\Phi_n\rangle=0$, which,
considering the spectral properties of $L_0$ implies that as $n\to
\infty$ we have $\Phi_n \to \Phi$. Then by the arbitrariness of
$\Phi\in \D(L_0^k)$ we get exactly what we wanted.
\end{proof}

\section{Connection with modular theory of von Neumann algebras}
\label{sec:modular}

\paragrf
We have seen that given a \mo covariant local net, the vacuum
vector $\Omega$ is cyclic and separating and for the local algebra
$\A(I)$ for every $I\in\I$. Thus one can consider the modular
one-parameter group $t\mapsto \Delta^{it}_{\A(I),\Omega}$ and the
modular conjugation $J_{\A(I),\Omega}$ associated to the von
Neumann algebra $\A(I)$ and the cyclic and separating vector
$\Omega$. For a reference on the modular theory of von Neumann 
algebras see e.g.\!  the book of Takesaki \cite{takesaki}.

One of the most important recognitions in the study of local
covariant nets was that these modular objects have geometrical
meaning. In order to describe their geometrical meaning, we need to
fix some one-parameter subgroups of $\mob$. The one-parameter subgroup of 
rotations $\alpha \mapsto \rho_\alpha \in \mob$ was already introduced. 
The one-parameter groups of $\mob$ generated by the vector fields
\begin{equation}
t(z) \equiv \frac{1}{2}
-\frac{1}{4}(z+z^{-1}),\;\;d(z)\equiv\frac{i}{2}(z-z^{-1})
\end{equation}
will be denoted by $a\mapsto \tau_a$ and $s\mapsto \delta_s$,
respectively. The first one is called the one-parameter group of
{\bf translations} while the second one is the one-parameter group
of {\bf dilations}. Their name becomes clear if one identifies the
circle with the one-point compactification of the real line in a
certain way with the point $1\in S^1$ being the infinite point,
see more on these in Sect.\! \ref{app:mobgroup} of the appendix.
The dilations {\bf scale} the translations: $\delta_s \tau_a
\delta_{-s} = \tau_{e^s a}$.

The introduced dilations will be also called the ``dilations
associated to the upper circle $S^1_+\equiv \{z\in S^1: {\rm
Im}(z)>0\}$'' as up to parametrization this is the only
one-parameter subgroup of $\mob$ preserving $S^1_+$. If $I\in\I$
then there exists a $\varphi\in\mob$ such that $\varphi(S^1_+)=I$.
The one parameter group $s\mapsto \delta^I_s\equiv
\varphi\circ\delta_s\circ\varphi^{-1}$ does not depend on the
choice of $\varphi$ and up to parametrization it is the only
one-parameter subgroup preserving the interval $I$. It is
therefore called the dilations associated to the interval $I$. For
a more detailed introduction on these notions see again the
mentioned section of the appendix.

\paragrf
The following result is of key importance in finding out the
``geometrical'' meaning of these modular objects.
\begin{theorem}\!{\rm \cite{borchers}}
Let $\H$ be a Hilbert space, $\M \subset B(\H)$ a von Neumann
algebra, $\Omega\in\H$ a cyclic and separating vector for $\M$ and
$a\mapsto U(a)\in B(\H)$ a strongly continuous one parameter group
such that
\begin{itemize}
\item ${\rm Ad}(U(a))(\M) \subset \M$ for every $a\leq 0$ \item
the self-adjoint generator of $a\mapsto U(a)$ is nonnegative \item
$U(a)\Omega = \Omega$ for every $a\in\RR$.
\end{itemize}
Then  $\Delta^{it}U(a) \Delta^{-it} = U(e^{2\pi
t}a)$ and $J U(a)J = U(-a)$ for every $(t,a)\in\RR^2$, where
$\Delta\equiv \Delta_{\M,\Omega}$ is the modular operator and
$J\equiv J_{\M,\Omega}$ is the modular conjugation associated to
$(\M,\Omega)$.
\end{theorem}
Let us try to apply the above theorem for local nets. A
translation with negative parameter shrinks the upper half circle.
So given a \mo covariant local net $(\A,U)$ with vacuum vector
$\Omega$ we have
\begin{equation}
{\rm Ad}(U(\tau_a))(\A(S^1_+)) \subset \A(S^1_+)\;\;{\rm for\,
every}\;a\leq 0,
\end{equation}
and of course $U(\tau_a)\Omega=\Omega$ for every $a\in\RR$.
Moreover, since $U$ is a positive energy representation of $\mob$,
it follows that the self-adjoint generator of the translations
$a\mapsto U(a)$ is positive (see Prop.\! \ref{pos.mob^n}). So by
the above theorem the adjoint action of
$\Delta^{it}_{\A(S^1_+),\Omega}$ on $U(\tau_a)$ is exactly the
same as that of $U(\delta_s)$ with the parameter of dilation
$s=2\pi t$.

That the modular one-parameter group ``scales'' a certain
one-parameter subgroup of the geometrical symmetry group, is a
general feature, not specific to low dimensional QFT. What is nice
about \mo covariant local nets, that in their case without any
further assumption\footnote{On the $4$-dimensional
Minkowski-spacetime, assuming that the local nets are generated by
Wightmann-fields, Bisognano and Wichmann \cite{biwi} proved the
geometrical nature of the modular objects associated to wedges.}
one can prove that $\Delta^{it}_{\A(S^1_+),\Omega}$ does not only
scale the translations {\it like} $U(\delta_s)$: the two operators
actually coincide. Moreover, the modular conjugation has a
geometrical meaning, too. To describe this latter one, for an
$I\in\I$ we shall introduce the (geometrical) {\bf conjugation
associated to ${\mathbf I}$}. First we define
$\iota\equiv\iota_{S^1_+}$ to be the conjugation $z\mapsto
\overline{z}$, and then set
\begin{equation}
\iota_I \equiv \varphi \circ \iota \circ \varphi^{-1}
\end{equation}
where $\varphi(S^1_+)=I$. It is easy to check that the definition
is ``good'', see the similar argument at Def.\! \ref{def:dil.I} in
the appendix. With these conventions the following result can be
established.
\begin{theorem}\!{\rm \cite{BGL, FrG} {\bf ``Bisognano-Wichmann 
property''}}
Let $(\A,U)$ be a \mo covariant local net with vacuum vector
$\Omega$. Then
$$\Delta^{it}_{\A(I),\Omega} = U(\delta^I_{2 \pi t})$$
for every $I\in\I$ and $t\in\RR$. Moreover, for every $K\in\I$
$$
J_{\A(I),\Omega} \,\A(K)\, J_{\A(I),\Omega} = \A(\iota_I(K)).
$$
\end{theorem}
There are several consequences; some are immediate and some are
less immediate. Most of them were already stated in the mentioned
aricles (so where there is no reference given, one should look for
\cite{BGL, FrG}). To begin with: as $\mob$ is generated by the set of
dilations (associated to different intervals), it follows that the
representation $U$ is completely determined by the local algebras
and the vacuum vector (via modular structure).

By modular theory the modular conjugation of an algebra takes the
algebra into its commutant. So $J_{\A(I),\Omega} \,\A(I)\,
J_{\A(I),\Omega}= \A(I)'$ and hence by the above theorem we
find\footnote{Actually the real order of demonstration is rather
vice versa: first one proves Haag-duality by the theorem of
Takesaki that is later mentioned on this page, and then deduces
the geometrical action of the modular conjugation.} that a \mo
covariant local net always satisfies the property called {\bf
Haag-duality}:
\begin{equation}
\A(I)'=\A(I^c)
\end{equation}
for every $I\in \I$ where the following notation is used (and will
be used all through this thesis):
\begin{equation}
I^c \equiv {\rm \; the\;interior\;of\;} (S^1\setminus I).
\end{equation}
(Note that while the complement of $I$ is not, $I^c$ is again an
element of $\I$).

If a vector is globally invariant under the one-parameter group
$s\mapsto U(\delta^I_s)$ then it is actually invariant under
$U(\varphi)$ for every $\varphi\in\mob$ (see \cite{GuLo96}) and
hence by the uniqueness of the vacuum it is a multiple of the
vacuum vector. So by the {\it Bisognano-Wichmann property}, the
modular group of $\A(I)$ with respect to $\Omega$, is ergodic,
showing that $\A(I)$ is a type ${\rm I\!I\!I}_1$ factor. This
property is called the {\bf factoriality} of the local algebras.

Actually, in the definition used earlier in the literature, a \mo
covariant local net was not required to satisfy all axioms listed
in this thesis. The property ``missing'' was the {\it uniqueness}
of the vacuum.

In \cite{GuLo96} it is proved that the vacuum is unique if and
only if the local algebras of the net are factors and that it is
further equivalent with the {\bf irreducibility} of the net, that
is, with the property that
\begin{equation}
(\cup_{I\in\I}\A(I))' = \CC \mathbbm 1.
\end{equation}
To distinguish between the two different set of axioms, when
uniqueness is not required, we shall talk about a ``non necessarily
irreducible'' \mo local net.

Let us see more implications of the {\it Bisognano-Wichmann
property}. By the theorem of Takesaki concerning modular invariant 
subalgebras and expectations (see e.g.\!
\cite[Vol.\! I\!I, Theorem 4.2]{takesaki}), if $\N\subset \M$ is an
inclusion of von Neumann algebras, $\Omega$ is a cyclic and
separating vector for $\M$, then
\begin{equation}
\left.\begin{matrix}
\overline{\N\Omega}=\overline{\M\Omega}& \vspace{2mm}\\
\forall t\in\RR:\;\Delta_{\M,\Omega}^{it}\,\N\,\Delta_{\M,\Omega}^{it}
&=& \N &\end{matrix}\right\}\;\;\Rightarrow\;\; \N=\M.
\end{equation}
Taking this into account, the {\it Bisognano-Wichmann property}
turns out \cite{jors} to also imply the {\bf continuity} of the
local algebras of a \mo covariant local net. Indeed, suppose $I_n
\;\;(n\in\NN)$ is a growing  sequence of intervals with
$\cup_{n\in\NN} I_n =I\in \I$. Then by isotony
\begin{equation}
\N\equiv (\cup_{n\in\NN} \A(I_n))'' \subset \A(I) \equiv \M
\end{equation}
and by the {\it Reeh-Schlieder property} the vacuum vector
$\Omega$ is cyclic and separating for both $\M$ and $\N$. But it
is also clear that by the {\it Bisognano-Wichmann property} $\N$
is globally invariant under the adjoint action of
$\Delta^{it}_{\M,\Omega} \;\;(t\in\RR)$ and so actually
$(\cup_{n\in\NN} \A(I_n))'' = \A(I)$. Likewise, if $I_n
\;\;(n\in\NN)$ is a decreasing sequence of intervals with $I\in
\I$ being the interior of $\cap_{n\in\NN} I_n $ then it follows
that $\cap_{n\in\NN} \A(I_n)=\A(I)$. (This second form of the
continuity we can get easily from the first form by
considering {\it Haag-duality}, for example.)

Another important structural property which can be shown by
similar arguments is the {\bf additivity} of the net. This means
that if $\{I_n \in\I: n\in\NN\}$ is a collection of intervals with
$\cup_{n\in\NN} I_n$ completely covering an interval $I\in\I$ then
$\A(I)\subset (\cup_{n\in\NN} \A(I_n))''$. We shall not overview
the exact proof, but it is based on similar ideas that were used
above for showing the continuity of the net.

\paragrf
So far we have discussed properties that automatically follow. In
relation to the last of them (additivity), we shall now introduce
some {\it additional} properties that are known to hold in many,
but not all models.

By {\it additivity} (and isotony), if $I_1,I_2,I \in \I$ are such
that $I_1 \cup I_2 = I$ then $\A(I_1)\vee\A(I_2)=\A(I)$. As it was
mentioned, the local algebras are continuous, one may wonder
whether $\A(I_1)\vee\A(I_2)=\A(I)$ holds even if $I_1\cup I_2$
does not cover $I$ but the closer $\overline{I_1\cup I_2}$ does.
(So $I_1$ and $I_2$ are the connected components obtained by
removing a certain point from $I$.) In many (but not all)
physically interesting models this stronger version of the
additivity property, called {\bf strong additivity}, indeed holds.
However, there are counterexamples.

For an $n=2,3,..$ the net $\A$ is said to be {\bf n-regular}, if
whenever we remove $n$ points from the circle the algebras
associated to the remaining intervals generate the whole of
$\A(S^1)={\rm B}(\H_\A)$. By isotony $n$-regularity is a stronger
property then $m$-regularity if $n>m$, and by {\it Haag-duality}
(and {\it factoriality}) every \mo covariant local net is at least
$2$-regular. Strong additivity is of course stronger than
$n$-regularity for any $n$. Finally, if a net is $n$-regular for
every $n\in\NN^+$, then we shall just simply call it {\bf
regular}.

Later on we shall discuss examples of \mo covariant local nets
showing the existence of strongly additive nets, as well as nets
that are {\it not} strongly additive, $3$-regular but not
$4$-regular nets, and finally, a net which is not even
$3$-regular. So all these properties are indeed ``additional''.

\paragrf
Let $(\A,U)$ be a \mo covariant local net and suppose $I,K\in \I$
are two intervals with a common endpoint and with $K\subset I$.
Depending on which ``side'' they have a common endpoint, the
dilations associated to $I$ for positive (or negative) parameters
shrink the interval $K$ and so by the {\it Bisognano-Wichmann
property} the modular group associated to $\A(I)$ and $\Omega$,
for positive (negative) parameters maps $\A(K)$ into a smaller
algebra. So for example by setting $\M\equiv \A(S^1_+)$ and
$\N\equiv \A(S^1_{+,R})$ where the right upper half-circle
$S^1_{+,R} \equiv \{z\in S^1_+: {\rm Re}(z)>0\}$, we have that
$\Delta^{it}_{\M,\Omega}=U(\delta_{2\pi t})$ and so
\begin{equation}
\Delta^{it}_{\M,\Omega} \,\N\, \Delta^{-it}_{\M,\Omega} \subset
\N\;\;{\rm for\;all}\; t \geq 0.
\end{equation}
This kind of ``relative position'' of the algebras $\M$, $\N$ and
the vector $\Omega$, was given a special name.
\begin{definition}
Let $\N\subset \M$ an inclusion of von Neumann algebras on a
Hilbert space $\H$ with a common cyclic and separating vector
$\Omega\in\H$. Then $(\Omega,\N\subset \M)$ is called a {\bf
``$\mathbf +$'' (``$\mathbf -$'') half-sided modular inclusion}
(in abbreviation: a $\pm$hsm inclusion) if
$$
\Delta^{it}_{\M,\Omega} \,\N \,\Delta^{-it}_{\M,\Omega} \subset
\N\;\;{\rm for\; all}\; t \geq 0\;\; ({\rm for\; all}\; t \leq 0).
$$
If further $\Omega$ is cyclic even for $\N'\cap\M$ then
$(\Omega,\N\subset \M)$ is called a {\bf standard} $\pm$hsm
inclusion.
\end{definition}
Note that by the mentioned theorem of Takesaki if
$\Delta^{it}_{\M,\Omega} \N \Delta^{-it}_{\M,\Omega} \subset \N$
for all $t \in\RR$ then actually $\N=\M$. So in some sense this is
the most natural, nontrivial situation the inclusion 
$\N\subset \M$ can have regarding the modular group associated to $\M$ and 
the common and cyclic and separating vector $\Omega$. It turns out that in
this case the modular group of the smaller algebra has a certain
commutation relation with the modular group of the bigger one.
\begin{theorem}\!{\rm \cite{wiesbrock1,arazsido}}
Let $(\Omega,\N\subset \M)$ be a $\pm$hsm. Then, there is a unique
$(s,a) \mapsto V(s,a)$ a strongly continuous unitary
representation of the $2$-dimensional Lie group of transformations
of $\RR$ of the form
$$
\RR^2 \ni (s,a): x\mapsto e^s x + a
$$
such that omitting the vector $\Omega$ in the indices, for all
$t\in \RR$
$$
V(2\pi t,0)\equiv \Delta^{it}_\M,\;\;\;V(2\pi t, e^{2\pi
t}-1)\equiv \Delta^{i t}_\N.
$$
Moreover, for this representation the self-adjoint generator $T$
of $a\mapsto V(0,a)$ is positive: $T\geq 0$, and $V(0,1)\,\M\,
V(0,1)^* = \N$.
\end{theorem}
The reason for citing two articles is the following. The original
argument of Wiesbrock contained an error, which was then noticed
and corrected by Araki and Zsid\'o.

An inclusion of von Neumann algebras $\N\subset \M$ is said to be
{\bf normal} if $\N^{cc}=\N$ where $\R^c \equiv \R' \cap \M$ is
the relative commutant. The inclusion is said to be {\bf conormal}
if $\M$ is generated by $\N$ and $\N^c$ (which is equivalent with
the normality of $\M'\subset \N'$).

Using the above theorem, Wiesbrock proved \cite{wiesbrock2} that
given a $+$hsm conormal inclusion $(\Omega,\N\subset \M)$ one can
construct a unique non-necessarily irreducible, strongly additive
\mo covariant local net $(\A,U)$ with $\A(S^1_+)=\M$,
$\A(S^1_{+,R})=\N$ and $\Omega$ being the vacuum vector.

In \cite{GLW} it is proved that a standard hsm inclusion $(\Omega,
\N\subset \M)$ is always normal and conormal and that the center
of $\N$ and $\M$ coincides. So putting everything together, we
have the following scenario.
\begin{itemize}
\item Given a \mo covariant local net $(\A,U)$ with vacuum vector
$\Omega$, the triple $(\Omega,\A(S^1_{+,R})\subset \A(S^1_+))$ is
a $+$hsm inclusion of factors.

\item Given a $+$hsm inclusion of von Neumann factors
$(\Omega,\N\subset\M)$ one can construct a unique strongly
additive \mo covariant local net with $\A(S^1_+)=\M$,
$\A(S^1_{+,R})=\N$ and $\Omega$ being the vacuum vector.
\end{itemize}
Moreover, it is clear, that with the isomorphism defined in the
natural way, the above two constructions give a one-to-one
correspondence between isomorphism classes of strongly additive
\mo covariant local nets and isomorphism classes of $+$hsm
inclusions of von Neumann factors.

What about the \mo covariant local nets that are not necessarily
strongly additive: how can they be described by modular
inclusions? In \cite{GLW} the notion of a {\bf $\mathbf \pm$hsm
factorization} is introduced. It is a quadruple $(\Omega,
\N_0,\N_1,\N_2)$ where $\{\N_j:j\in\ZZ_3\}$ are pairwise commuting
von Neumann algebras on a common Hilbert space and $\Omega$ is a
common cyclic and separating vector such that $(\Omega,
\N_j\subset \N_{j+1})$ is a $\pm$hsm inclusion for each
$j\in\ZZ_3$.

Suppose we are given a \mo covariant local net $(\A,U)$ with
vacuum vector $\Omega$. Fix three (different) points on the circle
and let $I_0,I_1,I_2\in \I$ be the intervals obtained by the
removal of these three points, with the numbering corresponding to
the clockwise order. Then by the {\it Bisognano-Wichmann
property}, $(\Omega, \A(I_0),\A(I_1),\A(I_2))$ is a $+$hsm
factorization. Further, by the {\it factoriality}, all three von
Neumann algebras involved are factors.

So given a \mo covariant local net we can construct such
structures. However, the real point is that by \cite[Theorem
1.2]{GLW}, given a $+$hsm factorization $(\Omega, \N_0,\N_1,\N_2)$
with $\N_j \;(j\in\ZZ_3)$ being factors, and $I_0,I_1,I_2\in\I$
three intervals as before, there exists a unique \mo covariant
local net $(\A,U)$ with $\A(I_j)=\N_j\;(j\in\ZZ_3)$, and $\Omega$
being the vacuum vector. So this is how a general, not necessarily
strongly additive, \mo covariant local net can be described by
modular data.

\section{Further structural properties and representation theory}
\label{sec:furtherprop}

\paragrf
So far everything was considered on the circle, although by the
physical motivation (``the restriction of a theory to the
worldline of a light signal'') one should work with nets on the
real line. Also, by restriction, it is only translation-dilation
covariance, and not \mo covariance that we directly receive. (The
translation-dilation group is the group of transformations of
$\RR$ of the form $x\mapsto e^s x + a$ with $(s,a)\in \RR^2$.)

Denote by  $\I_\RR$ the set of nonempty, nondense open intervals
of $\RR$. A map $I\mapsto \A(I)$ assigning to each $I\in\I_\RR$ a
von Neumann algebra $\A(I)$ acting on a fixed, complex Hilbert
space, is called a local net of von Neumann algebras on the real
line, if it satisfies the usual isotony and locality requirement.
We shall further say that given a local net $\A$ on $\RR$ and a
strongly continuous unitary representation $U$ of the
translation-dilation group (acting on the Hilbert space of $\A$),
the pair $(\A,U)$ is a translation-dilation covariant net on $\RR$
if the following properties hold:
\begin{itemize}
\item ({\it covariance}) $U(\gamma)\A(I)U(\gamma)^* =
\A(\gamma(I))$ for every element $\gamma$ of the
translation-dilation group and $I\in\I_\RR$,

\item ({\it existence and uniqueness of vacuum}) there exists an
up to phase unique unit vector (called: the vacuum vector), which
is invariant for $U$,

\item ({\it cyclicity of the vacuum}) the vacuum vector is cyclic
for the von Neumann algebra $\A(\RR)\equiv
(\cup_{I\in\I_\RR} \A(I))''$.
\end{itemize}
Identifying $\RR$ with $S^1\setminus \{1\}$ through the
Cayley-transformation of Eq.\! (\ref{Cayley}), the translations of
$\RR$ become the one-parameter group $a\mapsto \tau_a$ and the
dilations the one-parameter group $s\mapsto \delta_s$ (which is
actually the explanation of the name of these one-parameter
groups). Moreover, as it was already mentioned (see
\cite{GuLo96}), if a vector is globally invariant under the
one-parameter group $s\mapsto U(\delta^I_s)$ then it is actually
invariant under $U(\varphi)$ for every $\varphi\in\mob$.

From what was said it is clear that given a \mo covariant local
net $(\A,U)$, the map $\I_\RR \ni I \mapsto \A_\RR(I)\equiv \A(I)$
where the identification is defined by this Cayley transformation,
is a local net on $\RR$, and that after defining $U_\RR$ by a
similar restriction, we have that $(\A_\RR,U_\RR)$ (called the
``restriction onto $\RR$'') is a translation-dilation covariant
local net on $\RR$. Of course $(U,\A)$ can be completely restored
from its restriction to $\RR$, since $U$ is completely determined
by the local algebras of $\A_\RR$ and the vacuum vector via
modular theory (for the geometrical part of this statement consider for
example Lemma \ref{IRgenerators}.)

\begin{theorem}\!{\rm \cite[Theorem 1.4]{GLW}}
A translation dilation covariant local net $(\A_\RR,U_\RR)$ is the
restriction of a \mo covariant local net (in which case it is the
restriction of a unique \mo covariant local net) if and only if it
satisfies
\begin{itemize}
\item the Reeh-Schlieder property: the vacuum vector is cyclic and
separating for each local algebra $\A_\RR(I),\;(I\in\I_\RR)$,
\item and the Bisognano-Wichmann property for $\RR^+$: omitting
the vacuum vector in the index,
$\Delta_{\A_\RR(\RR^+)}^{it}=U(\{x\mapsto e^{2\pi t}x\})$ for
every $t\in\RR$.
\end{itemize}
\end{theorem}
A translation-dilation covariant local net obtained by some kind
of restriction of a QFT generated by  Wightmann-fields should of
course satisfy {\it additivity} (this is clear at the level of
test functions on which the fields are evaluated) and thus also
the {\it Reeh-Schlieder property}, but by \cite{biwi} it should
also satisfy the Bisognano-Wichmann property for half-lines. Thus
the mentioned theorem essentially confirms that it is indeed
physically meaningful to consider chiral nets as \mo covariant
local nets on the {\it circle} instead of translation-dilation
covariant local nets on the {\it real line}.

\paragrf
A {\bf (\mo covariant) subnet} of the \mo covariant net $(\A,U)$
is a map $I\mapsto \B(I)$ assigning a von Neumann algebras acting
on the Hilbert space of $\A$ to each $I\in\I$ such that for all
$I_1,I_2,I \in \I$ and $\varphi \in \mob$
\begin{itemize}
\item[(i)] $\B(I) \subset \A(I)$ \item[(ii)] $I_1 \subset I_2
\Rightarrow \B(I_1) \subset \B(I_2)$ \item[(iii)] $\varphi \in
\mob \Rightarrow U(\varphi)\B(I)U(\varphi)^* = \B(\varphi(I))$.
\end{itemize}
We shall use the simple notation $\B \subset \A$ for subnets.

A subnet $\B \subset \A$ which is really smaller then $\A$ is not
a \mo covariant net in the precise sense of the definition because
the vacuum will not be cyclic for $\B$. However, this
inconvenience can be overcome by restriction to the so-called
``vacuum Hilbert space of $\B$''; that is to
\begin{equation}
\H_\B\equiv\overline{\B(S^1)\Omega},
\end{equation}
where $\B(S^1)\equiv (\cup_{I \in \I}\B(I))''$ and $\Omega$ is the
vacuum vector.

It is evident that $\H_\B$ is invariant for $U$. The map $I
\mapsto \B(I)|_{\H_\B}$ together with the restriction of $U$ onto
$\H_\B$ is a \mo covariant net. Rather direct consequences of the
definition and of the properties of \mo covariant nets (such as
for example the {\it Reeh-Schlieder} and {\it Haag property}) are:
\begin{itemize}
\item[(i)] for any $I \in \I$ the restriction map from $\B(I)$ to
$\B(I)|_{\H_\B}$ is an isomorphism between von Neumann algebras,
\item[(ii)] the map $A \mapsto [\H_\B] A |_{\H_\B}$ where
$[\H_\B]$ is the orthogonal projection onto $\H_\B$ and $A \in
\A(I)$ for a fixed $I \in \I$ defines a faithful normal
conditional expectation from $\A(I)$ to $\B(I)$, after identifying
$\B(I)$ with $\B(I)|_{\H_\B}$ using point (i), \item[(iii)]
$\B(S^1) \cap \A(I) = \B(I)$ for all $I \in \I$.
\end{itemize}

\paragrf
A {\bf locally normal representation} $\pi$ (or for shortness,
just simply representation) of a \mo covariant local net $(\A,U)$
consists of a Hilbert space $\H_\pi$ and a normal representation
$\pi_I$ of the von Neumann algebra $\A(I)$ on $\H_\pi$ for each $I
\subset \I$ such that the collection of representations $\{\pi_I :
I \in \I \}$ is consistent with {\it isotony}:
\begin{equation}
I \subset K \Rightarrow \pi_K |_{\A(I)} = \pi_I.
\end{equation}
It follows easily from the axioms and the known properties of
local nets listed in the last section that if $I \cap K =
\emptyset$ then $[\pi_I(\A(I)),\pi_K(\A(K))]=0$, if ${\mathcal S}
\subset \I$ is a covering of $K \in \I$ then $(\cup_{I\in
{\mathcal S}} \pi_I(\A(I)))''\supset \pi_K(\A(K))$, if 
${\mathcal S} \subset\I$ is a covering of $S^1$ then 
$(\cup_{I\in {\mathcal S}}\pi_I(\A(I)))''= (\cup_{I\in\I}\pi_I(\A(I)))'' 
\equiv\pi(\A)$ and finally, that for each $I \in \I$ the representation 
$\pi_I$ is faithful. The representation $\pi$ is called {\it 
irreducible}, if $\pi(\A)' = \CC \mathbbm 1$. Such an irreducible 
representation is sometimes also called a {\bf charged sector} of $\A$.

In this section, although it was not pointed out, we have already
meat with representations of local nets. Suppose $(\A,U)$ is \mo
covariant local net and $\B \subset \A$ is a subnet. As it was
explained, the restriction of $\B$ onto its vacuum Hilbert space
$\H_\B$, is a \mo covariant local net. The restriction map from
$\B(I)$ to $\B(I)|_{\H_\B}$ is an isomorphism, so in particular
also its inverse, which we shall here denote by $\pi^\B_I$, is an
isomorphism. It is easy to see then that the collection
$\pi^\B\equiv\{\pi^\B_I:I\in\I\}$ is a locally normal
representation of the \mo covariant local net obtained by the
restriction of the subnet $\B$ onto its vacuum Hilbert space
$\H_\B$.

A locally normal representation $\pi$ of $(\A,U)$ on the Hilbert
space $\H_\pi$ is said to be {\bf covariant}, if there exists a
strongly continuous unitary representation $\tilde{U}_\pi$ of
$\widetilde{\mob}$ (i.e.\! the universal covering of $\mob$)
acting on $\H_\pi$ such that for every $I\in\I$ and
$\tilde{\gamma} \in \widetilde{\mob}$
\begin{equation}
{\rm Ad}(\tilde{U}_\pi(\tilde{\gamma})) \circ \pi_I =
\pi_{\gamma(I)}\circ {\rm Ad}(U(\gamma))|_{\A(I)},
\end{equation}
where $\gamma = p(\tilde{\gamma})$ with
$p:\widetilde{\mob}\rightarrow \mob$ being the canonical covering
map. In other words, $\pi$ is covariant if the \mo symmetry of the
net is implementable in the representation $\pi$. A covariant
representation is said to be of {\bf positive energy}, if the
implementation of the \mo symmetry in the representation
$\tilde{U}_\pi$ can be chosen to be a positive energy
representation of $\widetilde{\mob}$.

In the usual way, one can define for locally normal
representations, too, what is meant by {\it equivalence,
invariant subspace} and {\it direct sum}. It was a fundamental
recognition of Doplicher, Haag and Roberts \cite{DHR1,DHR2} (made
in the $4$-dimensional setting), that furthermore, a certain {\it
product} operation can be introduced, too.

We shall not enter into details as it will be not needed for the
results of this thesis. We shall only overview some of the most
relevant points in the setting of \mo covariant nets.

Suppose $\pi$ is a locally normal representation of $(\A,U)$ and
$I$ is an element of $\I$. Then, since $\A(I^c)$ is a type ${\rm
I\!I\!I}_1$ factor, there must exists a unitary $W^\pi_{I^c}$
implementing $\pi_{I^c}$.

By {\it Haag-duality} it is easy to see that
\begin{equation}
\varrho^\pi_I \equiv {\rm Ad}(W^\pi_{I^c})^{-1}\circ \pi_I
\end{equation}
is an {\it endomorphism} of $\A(I)$. Of course it depends on the
choice of $W^\pi_{I^c}$, but a change in the implementing operator
``modifies'' this endomorphism by an inner automorphism only.

For a local net $\A$ one can introduce in a canonical way its {\it
universal algebra} $C^*(\A)$ (see e.g.\! \cite{fredenhagen90} for
precise definitions); it is a $C^*$ algebra with given injective
embeddings $\iota_I : \A(I) \mapsto C^*(\A)$ such that every
representation of $\A$ factors through a representation of
$C^*(\A)$ in a unique way. (It is similar to the regular
representation of groups.)

The endomorphism defined previously can be interpreted as a {\it
localized endomorphism} of $C^*(\A)$ (with $I$ being the interval
of localization); that is, an endomorphism such that if $K\in\I,
K\supset I$ then $\varrho^\pi_I(\iota_K(\A(K)))\subset
\iota_K(\A(K))$ while if $K \cap I =\emptyset$ then
$\varrho^\pi_I|_{\iota(\A(K))}={\rm id}_{\iota(\A(K))}$.

Then one can identify an equivalence class of representations with
a certain equivalence class of localized endomorphisms (with the
equivalence defined by the inner automorphisms). In this way,
since endomorphisms can be composed, we get a new operation on
(the equivalence classes) of representations. Moreover, using
these localized endomorphisms, for each representation something
which is called ``a statistical operator'' can be defined. In
turn, this statistical operator allows us to define the so-called
{\it statistical phase} and the {\it statistical dimension} of a
representation.

As it was already said, these concepts (although they are
fundamental) will be not explained in this overview. Apart from
the mentioned original research articles, an excellent
introduction to the theory of charged sectors (in the
$n$-dimensional Minkowski spacetime) can be found in \cite{Haag}.
There it is also explained, in what (physical) sense such
representations describe charges, and why the statistical phase is
indeed related to the ``statistics'' of particles.

The only point that is to be stressed here, is the following. If
the dimension of the spacetime is bigger or equal to $3$, then ---
through the statistical operator --- a representation of a net always
gives rise to a representation of the permutation group $S_n$ for
every $n\in\NN^+$. This leads to the conclusion that (if the
statistics is finite, then) the statistical phase can be either
$1$ or $-1$ (corresponding to particles following {\it Bose-}, and
particles following {\it Fermi-statistics}), while the statistical
dimension must be a positive integer.

In $2$-dimensions and in the case of chiral nets, the situation is
completely different. Due to simple geometrical reasons, instead
of the permutation group, one has representations of the {\it braid
group} $\mathbb B_n$ for every $n\in \NN^+$. These representations
do not necessarily factor through a representation of the
permutation group. As a result, the statistical phase is not
necessarily $\pm 1$ but in general a unit complex number, and
neither the statistical dimension has to be an integer. So, in low
dimensions, apart from the usual bosons and fermions, we have more
exotical particles, too.

\paragrf
If $\pi$ is a covariant representation of the \mo covariant local
net $(\A,U)$, then the Jones (minimal) index\footnote{For a 
general reference on index theory see e.g.\! the book \cite{takesaki}. 
For the special case of type I\!I\!I algebras, see also the lectures notes 
\cite{kosaki} of Kosaki.} of the inclusion
\begin{equation}
\pi_I(\A(I)) \subset \pi_{I^c}(\A(I^c))
\end{equation}
is independent of the interval $I\in\I$. (In the vacuum sector,
that is, with $\pi$ being the identical representation, by {\it
Haag-duality} this inclusion is actually an equality. So this
index measures ``how much Haag-duality is broken''.) It was a
fundamental recognition that this index is related to the
statistical dimension $d_\pi$ of the representation $\pi$. By the
works \cite{longo1,GuLo96}, if $\pi$ is of positive energy, then
$d_\pi$ is finite if and only if the mentioned Jones index is
finite and in this case
\begin{equation}
d_\pi  = [\pi_{I^c}(\A(I^c)):\pi_I(\A(I))]^{\frac{1}{2}}.
\end{equation}
Actually, as a consequence of \cite[Prop.\! 2.14]{GuLo96} and
\cite[Corollary 4.4]{BCL}, the condition on the positivity of
energy in the representation $\pi$ can be dropped. More precisely,
if $\pi$ is a covariant representation with either finite index or
finite statistics, then both are finite, $\pi$ is automatically of
positive energy and the above relation between index and
statistics holds. Moreover, in this case, by \cite[Prop.\!
2.2]{GuLo96} the representation $\tilde{U}_\pi$, implementing the
\mo symmetry in the representation $\pi$, is unique. Hence if
$\pi$ is further irreducible, then
$\tilde{U}_\pi(\tilde{\rho}_{2\pi})$ is a multiple of the
identity. By the {\it Conformal Spin and Statistics Theorem}
\cite[Theorem 3.13]{GuLo96}, this multiple is exactly the
statistical phase of the sector.

\paragrf
As we shall later see, there are local nets possessing positive
energy sectors with infinite index/statistics, as well as nets
with extremely well-behaved representation theory: having only a
finite number of inequivalent sectors, each with finite index. As
it was rather indicated then explained, sectors can be composed
resulting new representations. It was not mentioned, but for a
covariant sector $\pi$ of finite statistical dimension one can also
introduce its {\it conjugate sector} $\overline{\pi}$, which is a 
covariant sector defined up to equivalence. 
The conjugate of the conjugate is equivalent with
the original sector and $d_{\overline{\pi}}=d_\pi$. Moreover, 
the product of two finite dimensional covariant sectors is finite 
dimensional, and can be decomposed into a finite sum of finite dimensional 
covariant sectors, see \cite{GuLo96}. 

So suppose our net is ``well-enough'' behaved: it has a finite 
number of (equivalence classes of) covariant sectors, each 
of which is of finite statistical dimension.
Let $[\pi_j] \;\;(j\in \{0,..,m\}\equiv L)$ be the equivalence
classes of the covariant sectors of such a net with
$\pi_0={\rm id}$ being the vacuum representation. For a $j\in L$ define
$\overline{j}$ by $[\pi_{\overline{j}}]\equiv \overline{[\pi_j]}$.
Then the hyper-matrix 
\begin{equation}
(N^a_{b,c})_{(a,b,c)\in L^3}
\end{equation}
defined by the multiplicities appearing in
the decomposition of the products, that is
\begin{equation}
[\pi_a][\pi_b] \simeq \mathop{\oplus}_{c\in L}N^c_{a,b}[\pi_c],
\end{equation}
has only nonnegative integers in its entries, and
satisfies certain identities. Some of these are just direct
consequences of the commutativity and associativity of the product
and some are consequences of the definition of statistical
dimension and conjugate sector. More precisely, with $d_k$
standing for the statistical dimension of $\pi_k$, we have that
(see e.g.\! \cite[Lemma 4.46]{FrG}) for all $a,b,c \in L$
\begin{eqnarray}\nonumber
&1.& N^a_{b,c}= N^a_{c,b} = N^{\overline{c}}_{\overline{a},b} =
N^{\overline{c}}_{a,\overline{b}}\\
\nonumber &2.& \sum_{k,l\in L} N^k_{a,b} N^l_{k,c} = \sum_{k,l\in
L} N^k_{a,l} N^l_{b,c} \\
\nonumber &3.& \sum_{k\in L} N^k_{a,b} d_k = d_a d_b \\
&4.& N^0_{a,b} = \delta_{a,\overline{b}}.
\end{eqnarray}
The hyper-matrix $N$ contains all informations about the
category of finite dimensional covariant representations 
of the net. It completely determines the rules governing the algebraic 
structure of such representations. These rules are also called the 
``fusion rules''. Note that the third listed property is due to the fact 
that the statistical dimension is indeed a kind of dimension: for direct 
sums it gets summed up, for products it is to be multiplied together. This
property completely determines the dimensions of the sectors through the 
fusion rules.

How do we know that a net is ``well-enough'' behaved? In the work
\cite{KLM}, three properties were pinpointed. One is strong
additivity, which was already introduced in the last section. The
second one is the so-called {\it split} property which we shall
introduce now.

Two disjoint intervals $I_1,I_2\in \I$ are said to be {\it
distant}, if their distance is positive; i.e.\! if
$\overline{I_1}\cap \overline{I_2} = \emptyset$. The net $\A$ is
said to be {\bf split} if the map
\begin{equation}
\A(I_1)\A(I_2) \ni AB \mapsto A \otimes B \in \A(I)\otimes \A(K)
\end{equation}
extends to an isomorphism between the von Neumann algebra $\A(I_1)
\vee \A(I_2)$ and the von Neumann algebra $\A(I)\otimes \A(K)$,
for every pair of distant intervals $I_1,I_2\in\I$. This property
is known to hold if the ${\rm Tr}(e^{-\beta}L_0)<\infty$ for all
$\beta>0$, where $L_0$ is the conformal Hamiltonian, cf.\!
\cite{BuD'ALongo,D'ALoRa}.

Finally, consider the so-called $2$-interval inclusion
\begin{equation}
\A(I_1)\vee \A(I_2) \subset (\A(I_3)\vee \A(I_4))'
\end{equation}
where $I_3,I_4 \in \I$ are the two distant intervals obtained by
the removal of the closure of the union of the two distant
intervals $I_1,I_2\in\I$. If the net is strongly additive and
split, and moreover the Jones index of this inclusion is finite
for a certain pair of distant intervals $I_1,I_2\in\I$, then by
\cite[Prop.\! 5]{KLM} it is finite for all distant pairs and its
value is independent from the choice of the distant intervals.
Thus it is an invariant of the net $\A$; it is usually called the
``mu index of the net'' and denoted by $\mu_\A$. (Note that the
independence from the choice of the intervals does not follow from
\mo covariance; although the \mo group acts transitively on the
set of intervals, it does not act transitively on the set of {\it
pairs} of distant intervals.)

A net is said to be {\bf completely rational}, if it is strongly
additive, split, and has a finite $\mu$-index. The main point of
the mentioned work \cite{KLM} where complete rationality is
introduced, is that under this condition one can prove that the
$2$-interval inclusion is isomorphic to the so-called {\it
Longo-Rehren inclusion}, cf.\! \cite{longorehren}. Through this
identification then it is shown that the net has finitely many
(equivalence classes) of sectors, each with finite index and that
in fact, using the previous notations
\begin{equation}
\mu_\A = \sum_{k\in L} d^2_k.
\end{equation}
(So actually there is even a bound on the number of inequivalent
sectors and their dimensions.) Moreover, in this case every
representation of the net is a direct sum of irreducible
representations, the conjugate sector exists for every covariant 
sector as it is automatically of finite dimension,
and so we have all previously described structures regarding the 
fusion rules.

\section{A first example: the $U(1)$ current model and its derivatives}
\label{sec:U(1)}

\paragrf
Given a Hilbert space $\K$ the {\bf bosonic Fock-space} is defined as
\begin{equation}
\Gamma_+(\K) \equiv \CC \oplus \overline{\mathop{\oplus}_{n\in
\NN^+} (\mathop{\otimes}_{\rm sym.}^n \K)}
\end{equation}
where $\Omega_\K \equiv 1 \in \CC \hookrightarrow
\Gamma_+(\K)$ (with the embedding being the natural one coming from
the above direct sum formula) is called the vacuum vector. The dense 
set of vectors given by the above direct sum without taking closure
is called the {\it finite particle space}.

The point is of course not the construction of a new Hilbert space, but 
that of a certain system of operators. Namely, in this Hilbert space 
one can canonically construct a set of self-adjoint operators
$\{A(u) : u\in \K\}$ acting on $\Gamma_+(\K)$ (see e.g.\!
\cite[Sect.\! 5.2]{fockspace}) satisfying the following relations:
\begin{itemize}
\item the finite particle space is a common invariant core
for $A(u) \; (u\in \K)$, \item on this core $u\mapsto A(u)$ is
real-linear and $[A(u),A(v)] = - i\, {\rm Im}((u,v))\mathbbm 1$ \item
the Weyl-operators $W(u)\equiv e^{iA(u)}$ $(u\in\K)$ satisfy the
relation
\begin{equation}\label{Weyl.rel}
W(u)W(v)= e^{\frac{i\,{\rm Im}((u,v))}{2}}W(u+v)
\end{equation}
\item
if $\D\subset \K$ is dense then the set $\{W(u) : u\in \D\}$ is
irreducible in the Fock-space; moreover, polynomials of
$A(u)\;(u\in\D)$ applied to $\Omega_\K$ give a dense set in
$\Gamma_+(\K)$
\item if $V$ is a unitary operator acting on $\K$ then by setting
\begin{equation}
V_{\Gamma_+}\equiv \mathbbm 1_\CC \oplus (\mathop{\oplus}_{n\in
\NN^+} (\mathop{\otimes}_{\rm sym.}^n \V))
\end{equation}
for the unitary acting on $\Gamma_+(\K)$ we have
\begin{equation}\label{VGamma.on.Weyl}
V_{\Gamma_+} W(u) V_{\Gamma_+}^* = W((Vu)).
\end{equation}
\end{itemize}
Let us consider now the positive energy representation $V_1$
of $\mob$ with lowest weight equal to $1$ and its Hilbert space
$\H_1$ with scalar product $(\cdot,\cdot)_1$ (see Sect.\!
\ref{app:mobrep} of the appendix), and let us apply 
the above construction to get a set of self-adjoint operators
$\{A(\Phi) : \Phi\in \H_1\}$ and a vacuum vector which we shall
just simple denote by $\Omega$. 

For an $f\in C^\infty(S^1,\RR)$ with 
Fourier coefficients
$\hat{f}_n \equiv 
\frac{1}{2\pi}\int_0^{2\pi}f(e^{i\theta})e^{-in\theta}d\theta$
$(n\in\ZZ)$, set
\begin{equation}
[f] \equiv \{n\mapsto \sqrt{2(n+1)} \hat{f}_{n+1}\} \in l^2(\NN).
\end{equation}
In this thesis the Hilbert space $\H_1$ is identified with 
$l^2(\NN)$ (see again the mentioned section of the appendix for details), 
so we shall think of $[f]$ as an element of $\H_1$. Through this 
identification we find the following relation of the above introduced
element and the representation $U_1$: for every
$\varphi \in \mob$ we have 
\begin{equation}
U_1(\varphi)[f] =
[f\circ\varphi^{-1}].
\end{equation}
Moreover, by an easy calculation
\begin{equation}\label{Im[f][g]}
{\rm Im}( ([f],[g])_1) =
\frac{1}{2\pi}\int_0^{2\pi}f'(e^{i\theta})g(e^{i\theta})d\theta.
\end{equation}
So if we set $U\equiv (U_1)_{\Gamma_+}$ and
\begin{equation}
\A_{U(1)}(I)\equiv\{e^{iJ(f)}:f\in C^\infty(S^1,\RR),
f|_{I^c}=0\}\;\;\;(I\in\I)
\end{equation}
where the {\bf current} $J(f)\equiv A([f])$, then by Eq.\!
(\ref{Im[f][g]}) and by the already mentioned Weyl-relations we
see that $\A_{U(1)}$ satisfies locality (and obviously also
isotony). Moreover, by Eq.\! (\ref{VGamma.on.Weyl}), for all $f\in
C^\infty(S^1,\RR)$ and $\varphi\in \mob$
\begin{equation}\label{UJU^*}
U(\varphi) J(f) U(\varphi)^* = J(f\circ\varphi^{-1})
\end{equation}
and so the action of $U$ on $\A_{U(1)}$ is covariant. It
is easy to see that $\Omega$ is indeed a vacuum vector for
$(\A_{U(1)},U)$ and that the uniqueness of the vacuum is also
satisfied. (The tensor products do not contain the trivial
representation since the lowest eigenvalue of the self-adjoint
generator of the rotations is bigger then $1$ in them.) Thus we
have a \mo covariant local net: the so-called $U(1)$ current model.

For a complex function $f\in C^\infty(S^1,\CC)$ we define $J(f)$
to be the closure of $J({\rm Re}(f))+ i J({\rm Im}(f))$. Then the
Fourier modes of the current can be introduced as
\begin{equation}
J_n \equiv J(\{z\mapsto z^n\})|_{\fin}
\end{equation}
where $\fin$ stands for the dense subspace of {\it finite energy
vectors} of $U$; i.e.\! the algebraic span of the eigenvectors of
$L_0$, where $L_0$ is the self-adjoint generator of rotations in
$U$. (The finite energy space $\fin$ is contained in the finite
particle space.) By Eq.\! (\ref{UJU^*}) one finds that $J_n\fin
\subset \fin$ as in fact $[L_0,J_n]=-n J_n$. Thus $J_n$
``decreases the energy by $n$''; in particular, by the spectral
properties of $L_0$ we have that
\begin{equation}
J_n \Omega = 0 \;\;{\rm for\;all}\; n > 0 \;{\rm integer.}
\end{equation}
In fact, as the constant $1$ function $[1]=0 \in \H_1$, we have
that $J_0 = 0$ so also $J_0\Omega = 0$. Moreover, by the
previously listed properties, on $\fin$ 
\begin{equation}\label{[Jn,Jm]}
[J_n,J_m]=n \delta_{-n,m} \mathbbm 1
\end{equation}
for all $n,m\in\ZZ$, and by the self-adjointness of
$J(f)$ for a real function $f\in C^\infty(S^1,\RR)$,
\begin{equation}\label{J^*=J}
J_n^*\supset J_{-n}.
\end{equation}

\paragrf
Let us discuss now some properties of the $U(1)$ current model. 
First of all, let us mention that by a simple argument 
\cite{BS-M}, using essentially just the Weyl-relations, shows that 
this net is strongly additive. 

Let us see now the spectral properties of $L_0$. From what was so far 
said it is more or less evident that ${\rm Sp}(L_0)=\NN$ and the 
eigenspace associated to $n\in \NN$ is spanned by the vectors
\begin{equation}
\{J_{-m_1}..J_{-m_k}\Omega: m_1+..+m_k = n\}.
\end{equation}
Actually, by the commutation relation of the Fourier modes of the
current and the fact that for positive $n$ the operator $J_n$
annihilates $\Omega$, one can even see that the span of
\begin{equation}
\{J_{-m_1}..J_{-m_k}\Omega: m_1+..+m_k = n; \; m_1\geq ..\geq m_k
> 0\}
\end{equation}
still gives the whole eigenspae. 
So the dimension of ${\rm Ker}(L_0-n\mathbbm
1)$ is less than equal --- actually it turns out that {\it exactly} equal
--- to  the number of partitions $p(n)$ of $n$. As it is well known, if 
$0<q<1$ then $\sum_{n\in\NN} q^n p(n)$ is convergent and in fact 
equal to $\mathop{\Pi}_{n\in\NN}(1/(1-q^n))$. Hence
$e^{-\beta L_0}$ has a finite trace for every $\beta>0$. 
This shows that the $U(1)$ current model is also split.

\paragrf
For a fixed
$\gamma:S^1 \rightarrow S^1$, the map
\begin{equation}
e^{iJ(f)} \mapsto e^{iJ(f\circ\gamma^{-1})}
\end{equation}
respects the Weyl-relations. We know that if $\gamma$ is in
particular an element of $\mob$, then this map is implementable.
Actually, using the theory of representations of the Weyl-algebra,
it is possible to prove that this map is implementable {\it for
every orientation preserving diffeomorphism}, see e.g.\!
\cite[Theorem 9.3.1]{loop}. The implementing operator $U(\gamma)$, up to
phase, is unique, since the currents form an irreducible set. So
we find that in the projective sense
$U(\gamma_1)U(\gamma_2)=U(\gamma_1\circ \gamma_2)$. To put it in
another way, the representation of $\mob$, given with the $U(1)$
model, can be extended to be a projective representation of
$\diff$ acting on the net in a covariant manner. Actually, one can
prove that the extension is strongly continuous with respect to the
natural $C^\infty$ topology of the set of diffeomorphisms. Moreover, if
$\gamma|_{I}={\rm id}_I$ for a certain interval $I\in\I$, and
$f|_I = 0$, then by the above equation $U(\gamma)$ commutes with
$J(f)$ and hence by {\it Haag-duality} $U(\gamma)\in
\A_{U(1)}(I^c)$. We say that the representation of $\diff$ is {\it
compatible} with the local structure of the net.

What we have thus seen is that the $U(1)$ current model is {\bf
diffeomorphism} covariant. In the next chapter precise definitions
will be given and diffeomorphism covariance will be discussed at
length.

Let us move on to representation theory. For a fixed $q\in\RR$ the
map
\begin{equation}
e^{iJ(f)}\mapsto
e^{q\hat{f}_0} e^{iJ(f)}
\end{equation}
respects the Weyl-relations. One may wonder whether it is
implementable. If $q\neq 0$ the answer is no; instead  the map
gives rise to a locally normal, irreducible representation of the
$U(1)$ current model. Actually it turns out that if $q_1\neq q_2$
then these representations are inequivalent and that all
irreducible representations are of this type, see \cite{BMT1,BMT2}.
So in particular the $\mu$-index of the $U(1)$ current model is
infinite.

\paragrf
With the Cayley transformation of Eq.\! (\ref{Cayley}) one may 
identify $S^1\setminus\{1\}$ with the real line $\RR$. Through this 
identification we may define the current on the real line; in particular 
for an $
f\in C^\infty(\RR,\RR)
$ the self-adjoint operator
$
J_\RR(f)
$. The net defined by the formula (in which ``$f^{(n)}$'' stands for the 
${\rm n}^{\rm th}$ derivative)
\begin{equation}
\A^{(n)}_\RR(I) \equiv \{ e^{i J_\RR(f^{(n)})} :  f\in C^\infty(\RR,\RR), 
f|_{\RR \setminus I} = 0
\}\;\;\;(I\in\I_\RR)
\end{equation}
together with a corresponding restriction of $U$, 
is easily shown to be a translation-dilation covariant net
on the real line for which {\it additivity} holds. Moreover, 
simple argument shows  that for every $x\in \RR$ 
\begin{equation}
\A^{(n)}_\RR((x,\infty))=
\A^{(0)}_\RR((x,\infty))\equiv
\A_\RR((x,\infty)).
\end{equation}
(For example, consider that all these algebras are invariant under 
dilations centered at $x$, and that such a dilation has a modular origin
with respect to $\A_\RR$.)
This, by what whas said before (and in particular by what 
was explained in Sect.\! \ref{sec:furtherprop} about the relation between
\mo covariant nets and nets on the real line), shows that 
$\A^{(n)}_\RR$, together with a corresponding restriction of $U$, extends 
to a unique \mo covariant local net $(\A^{(n)}_{U(1)}, U^{(n)})$, which is 
called the {\bf ${\rm \bf n}^{\rm \bf th}$ derivative of the current}.

These models are in some sense ``pathological'', their real 
significance is exactly that they provide counter-examples to many 
properties that otherwise one could (wrongly) believe to be general 
features of \mo covariant local nets. In \cite{GLW} it was proved, that
$(\A^{(1)}_{U(1)}, U^{(1)})$ is $3$- but {\it not} 
$4$-regular, and for $n>1$, $(\A^{(n)}_{U(1)}, U^{(n)})$ is not 
even $3$-regular. In the mentioned article it was also proved, that
these models admit non-covariant sectors. Another interesting 
thing: it turned out \cite{koester03a, GLW} that these models have no 
diffeomorphism symmetry (this will be explained 
with somewhat more details after diffeomorphism covariance will be 
introduced in the next chapter.)

\section{More on examples}
\label{sec:more.ex}

\paragrf
Let $G$ be a finite dimensional connected Lie group with Lie algebra
$\mathfrak g$. The {\bf loop group $\mathbf{LG}$} is defined to be 
$C^\infty(S^1,G)$, with the group operation being the 
pointwise multiplication. It is an infinite dimensional Lie group with Lie 
algebra identified with $C^\infty(S^1,\mathfrak g)$, where the commutator 
is the pointwise commutator. A reference book on loop groups and on their
representation theory is \cite{loop}.

The Lie algebra of $U(1)$ can be naturally identified with $\RR$ and so  
we may think of a real-valued smooth function $f$ on the circle,
as an element of the Lie algebra of the loop group $LU(1)$. Thus,
in case of the $U(1)$ current model, by the Weyl-relations, the map
\begin{equation}
LU(1) \ni {\rm Exp}(f) \mapsto e^{iJ(f)}
\end{equation}
defines a projective unitary representation of the connected part of 
the loop group $LU(1)$.
This representation is strongly continuous, and in this representation, 
the action of $\widetilde{\mob}$ on the loop group, 
\begin{equation}
\mob \ni \tilde{\varphi} 
\mapsto 
\{\sigma \mapsto \sigma\circ p(\tilde{\varphi}^{-1})\},
\end{equation}
where $p:\widetilde{\mob} \rightarrow \mob$ is the universal covering,
is implementable by a unique, inner, positive energy representation. 
(``Inner'' = the representing operators belong to the von Neumann algebra 
generated by the representation of the loop group.) Thus we say that this 
representation is a {\it positive energy representation} of $LU(1)$.

One way of generalizing the $U(1)$ current model is to replace $U(1)$ by
other, non-necessarily commutative Lie groups. Such generalization is
usually called a {\bf loop group model}. The most important ingredients of
the construction are the positive energy representations of the given loop
group. The previous example served exactly to
give a sense of what is such a representation; in this overview
there will be no formal definitions regarding positive energy 
representations. (For the detailed treatment one should look for the
mentioned book.)

The most studied case is that of $LSU(N)$ with $N$ being an integer bigger 
or equal to $2$. Let us see a little bit more in detail, how does the 
construction go in this case.

There are two important real invariants, that, given a positive energy
irreducible representation of $LSU(N)$ one can introduce. One of them is 
the lowest value of the energy, which, by definition, is a nonnegative 
number. The other one is the so-called {\it level}. It is an invariant of 
the representation, which, in some sense, measures, ``how much'' the
representation is not true (but projective). To explain how it appears,
consider the scalar product $(\cdot,\cdot)$ on $\mathfrak{su}(N)$ defined 
by the formula
\begin{equation}
(a,b)\equiv \frac{-1}{2N}{\rm Tr}({\rm ad}(a){\rm ad}(b))
\end{equation}
and choose a base $(a^k \in \mathfrak{su}(N))_{(k=1,..,N)}$ orthonormed
with respect to the above scalar product.
 
The functions $z\mapsto z^n a^k$ $(n\in\ZZ, k\in\{1,..,N\})$ form a
(non-algebraic) base in the complexification of the Lie algebra of 
$LSU(N)$. Via infinitesimal generators, a projective representation 
of $LSU(N)$ gives rise to a representation of a central extension of
$\mathfrak{lsu}(N)$. Thus, given a projective representation, via 
infinitesimal generators, one can introduce the {\it Fourier modes of the 
currents}; the operators $J^k_n$ $(n\in\ZZ,k\in\{1,..,N\})$ corresponding 
to the above introduced base.
(Although we deal with a projective representation, and so the 
generators are only defined up to additive constants, in this case there 
is a certain canonical way to fix these additive constants. We shall 
not discuss this more in detail. Also, one has to address the question of 
domains, etc.\! which again,  we shall not do in this short overview.)
Then by the unitariness of the representation, the Fourier modes of the
(self-adjoint) generators satisfy the {\bf hermicity condition}
\begin{equation} 
(J^k_n)^* \supset J^k_{-n}
\end{equation}
and the commutator $[J^j_n,J^k_m]$, up to an additive constant, must be 
equal to $\sum_{l=1}^N -i f^{j,k}_l J^l_{n+m}$, where the {\it structure 
constants} $f^{j,k}_l$ of $\mathfrak{su}(N)$ corresponding to our base,
are defined by the equation
\begin{equation}
[a^j,a^k] = \sum_{l=1}^N f^{j,k}_l a^l.
\end{equation}
Actually, it turns out that there must exists a positive integer 
$K\in\NN^+$ such that
\begin{equation}
[J^j_n,J^k_m] = \sum_{l=1}^N -i f^{j,k}_l J^l_{n+m} + K n 
\delta_{-n,m}\delta_{j,k}.
\end{equation}
It is not hard to see that this integer is independent from the chosen 
orthonormed base $(a^k \in \mathfrak{su}(N))_{(k=1,..,N)}$, and thus it is 
an invariant of the representation. This is the so-called level of the 
representation.

For a function $f\in C^\infty(S^1,\RR)$ 
with Fourier coefficients $\hat{f}_n\;(n\in\ZZ)$ and $l\in \{1,..,N\}$,
the formula 
\begin{equation}
J^l(f) \equiv \overline{\sum_{n\in\ZZ} \hat{f}_n J^l_n}
\end{equation} 
is well defined (i.e.\! the sum converges and gives a closable operator)
and gives a self-adjoint operator. We may think of $J^l$ as an
operator valued functionals. It is usually called a {\bf current}.

It is known, that for each $K\in \NN^+$ there exists a unique positive
energy irreducible representation of $LSU(N)$ with the level equal to $K$
and lowest energy equal to zero. In this representation the eigenspace
corresponding to the zero energy is one-dimensional, and so up to phase it
contains a unique unit vector $\Omega$. Moreover, the net 
of algebras $I\mapsto \A_{SU(N)_K}(I)$ $(I\in\I)$, given by 
\begin{equation}
\A_{SU(N)_K}(I)\equiv \{e^{iJ^l(f)}: l\in\{1,..,N\}, f\in
C^\infty(S^1,\RR), f|_{I^c}=0\}'', 
\end{equation} 
together with the zero energy vector $\Omega$, define a \mo covariant
local net, which is actually also diffeomorphism covariant. It is
called the $\mathbf{SU(N)}$ {\bf model at level} $\mathbf K$.  This net
is strongly additive and split; essentially for the same arguments, that
were outlined for the case of the $U(1)$ current model.  There is however,
a significant difference between these models and the $U(1)$ current
model. These nets are {\bf completely rational}. In fact, their $\mu$
index and representation theory, so in particular, the fusion algebra, is
known for all values of $N=2,3,..$ and $K=1,2,..$, see  
\cite{wasserman,Xu}.

\paragrf
There are some ``abstract'' methods, with which, starting from some 
example, we can generate more examples. First of all, given two \mo 
covariant local nets, one can define, in an obvious way, their tensor 
product, which is again a \mo covariant local net. Actually, this
construction can be even generalized to infinite sequences of \mo 
covariant local nets.
 
Let $(\A_n, U_n)$, $(n\in\NN)$ be a sequence of 
\mo covariant local nets with $\Omega_n$, $(n\in\NN)$ being the 
corresponding sequence of vacuum vectors. Consider the infinite tensor 
product Hilbert space defined with respect to the infinite sequence of
vacuum vectors,
\begin{equation}
\H_\A\equiv \bigotimes_{n\in\NN}^{(\Omega_n)}\H_{\A_n}.
\end{equation}
It is fairly easy to show, that on this Hilbert space, the net defined 
by the map
\begin{equation}
\A(I) \equiv \otimes_{n\in\NN} \A_n(I) \;\;\; (I\in\I)
\end{equation}
together with the representation $U\equiv \otimes_{n\in\NN}U_n$ is a
\mo covariant local net. (Note that the above infinite tensor product 
representation indeed exists as $U_n$ leaves invariant $\Omega_n$ for 
every $n\in\NN$.)

It is not too hard to see that the tensor producty of a {\it finite} 
number of diffeomorphism covariant and/or completely rational nets is 
again a diffeomorphism covariant and/or completely rational net. For 
example, if $\mu_j <\infty$ is the $\mu$-index of $\A_j$ 
$(j\in\{1,..,N\})$, then one can show, that the $\mu$-index of the 
correponding tensor product net is simply $\mu = \Pi_{j=1}^N 
\mu_j$. This formula already indicates, that some properties may not 
``survive'' the infinite tensor product construction; in fact in this 
thesis it will be shown that the infinite tensor product of a sequence of 
nontrivial nets is never diffeomorphism covariant. Some properties of 
course survive; for example the tensor product of a sequence of strongly 
additive (resp.\! $n$-regular) nets is easily shown to be strongly 
additive (resp.\! $n$-regular).

Another abstract construction was already described previously: given a 
\mo covariant subnet $\B\subset \A$, by restriction to the vacuum Hilbert 
space of $\B$, we get a \mo covariant local net. However, there is also 
a way, with which, starting from a subnet, we can get a new subnet 
(and thus a new example). Consider the net of algebras defined by
\begin{equation}
C(I)\equiv \A(I) \cap \B' \;\;\;(I\in\I),
\end{equation}
where, as usual, $\B\equiv \B(S^1)\equiv (\cup_{I\in\I}\B(I))''$. It is 
not hard to see, that it is again a \mo covariant subnet. It is called
the {\bf coset} of the subnet $\B$.

In the previous point the $SU(N)$ models were introduced. If $N_1 < N_2$, 
then we have a natural embedding $SU(N_1)\subset SU(N_2)$, which, in turn,
gives the natural embedding 
\begin{equation}
\A_{SU(N_1)_K} \subset \A_{SU(N_2)_K}
\end{equation}
in the sense of \mo covariant subnets. The coset of this, and other 
similar, naturally arising embeddings in loop group models, are 
known to provide interesting new examples for \mo covariant local nets.

Finally, we shall discuss one more abstract construction. Suppose $(\A,U)$
is a \mo covariant local net which is not strongly additive. Consider its 
real restriction, $(\A_\RR,U_\RR)$, and set
\begin{equation}
\A^d_\RR(I)\equiv \left(\mathop{\cup}_{K\in\I_\RR, K\cap I = \emptyset} 
\A_\RR(K)\right)' \;\;\; (I\in\I_\RR).
\end{equation}
It is not hard to see, that $(\A^d_\RR,U_\RR)$ is still a 
translation-dilation covariant local net on the line, which fulfills
the conditions explained in Sect.\! \ref{sec:furtherprop}, ensuring
the existence of a unique extension to a \mo covariant local net. (In 
fact, note that $\A^d_\RR(I)=\A_\RR(I)$ for any half-interval 
$I=(x,\infty)$, while if $I$ is a bounded interval, then the dual contains 
the original algebra.) The resulting \mo covariant local net is called the 
{\bf dual} of $(\A,U)$. It is not isomorphic to the original net, as it 
is {\bf strongly additive}, see \cite{GLW}. 
(On the other hand, it is evident, that the dual of a strongly additive 
net is itself.)

In the last section, the derivatives of the $U(1)$ current models were 
introduced. It was mentioned, that they not strongly additive, as in fact, 
they are not even $4$-regular. Thus by applying the dual construction, we 
get a new net. However --- as it was noticed there --- by the way these 
models were constructed, it is clear, that they can be brought to the 
Hilbert space of the $U(1)$ current model in such a way, that they will 
have the same vacuum vector as the $U(1)$ current model, and on the real 
line, they will even have the same local algebra associated to any 
half-interval $(x,\infty)$. Hence it is clear, that the dual of any of the 
derivatives is isomorphic to the $U(1)$ current model. Thus in this case 
we did not manage to generate interesting new examples.

In what follows, some knowledge about diffeomorphism covariant local nets 
(e.g.\! central charge, Virasoro nets) will be assumed, although they will 
be only introduced in the next chapter. (This overview, in some sense, is 
not causal.)

The Virasoro net with central charge $c>1$ is not strongly additive, thus 
its dual is not just itself. Unlike in the previous example, in this case
(to the knowledge of the author) there is little known about the dual. It 
is not known whether it is isomorphic with one of the 
already known examples. In fact, it is not even clear, whether it will be 
diffeomorphism covariant, and/or split. 

\paragrf
A PhD thesis on the topic of local nets on the circle would be incomplete 
without mentioning some of the classification results in this area. We 
shall not go into details, as such results are not used in this thesis. 

Diffeomorphism covariant local nets on the circle 
with central charge $c<1$ are completely classified, see \cite{KL}. It 
turns out, that they are all completely rational. The classification 
consists of two ``infinite sequences'' (of which one of them is the 
sequence of Virasoro nets), and four ``exceptional models''. 

Diffeomorphism covariant local nets on the circle with central charge 
$c=1$ are partially classified, see \cite{Carpi2,Xu.strong}. This means, 
that the classification is made assuming a certain property (so it could 
happen, that there actually more models). An example for the $c=1$ case 
is the $U(1)$ current model. All models with $c=1$ are strongly additive, 
but as it can be seen by the mentioned example, they are not all 
completely rational. 

In the region $c>1$, little is known. There many known examples, for 
example all $SU(N)$ model, and of course every Virasoro net with central 
charge $c>1$. From these latter examples we see, that such models are not 
even necessarily strongly additive.

\chapter{Preliminaries I\!I: diffeomorphism covariance}
\label{chap:diffcov}

\markboth{CHAPTER \arabic{chapter}. DIFFEOMORPHISM
COVARIANCE}{CHAPTER \arabic{chapter}. DIFFEOMORPHISM COVARIANCE}

\section[The Lie group $\diff$ and the Virasoro algebra]
{The Lie group $\mathbf {{\rm \bf Diff}^+(S^1)}$ and the Virasoro algebra}
\label{sec:diffgroup}

\paragrf
We denote by ${\rm Vect}(S^1)$ the set of (smooth) vector fields
on $S^1$. We shall think of a vector field symbolically written as
$f(e^{i\vartheta})\frac{d}{d\vartheta}\in  {\rm Vect}(S^1)$ as the
corresponding real function $f$. Thus by this identification ${\rm
Vect}(S^1) \equiv C^\infty(S^1,\RR)$. Moreover, we use the notation $f'$
(calling it simply the derivative) for the function on the circle
obtained by differentiating with respect to the angle:
\begin{equation}
f'(z)\equiv\frac{d}{d\theta}f(e^{i\theta})|_{e^{i\alpha}=z}.
\end{equation}
Since $S^1$ is compact, a smooth function $f:S^1 \rightarrow \RR$
always gives rise to a flow on the circle. We denote by $a\mapsto
{\rm Exp}(af)$ the corresponding one-parameter group of
(orientation preserving) diffeomorphisms.

A $\gamma:S^1\rightarrow S^1$ diffeomorphism has an action on
vector fields; this will be denoted by $\gamma_*$. With the above
identification we have that for a vector field (function) $f$
\begin{equation}
(\gamma_*f)(z)= -i \overline{\gamma(z)}
\frac{d}{d\theta}\gamma(e^{i\theta})|_{e^{i\theta}=z}
f(\gamma^{-1}(z)).
\end{equation}
We denote by $\diff$ the group of orientation preserving (smooth)
diffeomorphisms of the circle. It is an infinite dimensional Lie
group whose Lie algebra is identified with the real topological
vector space ${\rm Vect}(S^1)$ of smooth real vector fields on
$S^1$ with the usual $C^\infty$ topology \cite[Sect.\! 6]{Milnor}
with the negative of the usual bracket of vector fields. This
choice of the sign is ``compulsory'' if we want the ``abstract''
exponential --- defined for Lie algebras of Lie groups
--- to coincide with the exponential ${\rm Exp}$ defined for
vector fields through generated flows. Thus in this thesis 
for two vector fields (functions) on the circle  
$f$ and $g$, we set
\begin{equation}
\label{-bracket} [f,g]\equiv f'g-fg'.
\end{equation}
An important property of the group $\diff$ is the following. For
references see the already cited survey \cite{Milnor}.
\begin{theorem}
The group $\diff$ is simple (in the algebraic sense).
\end{theorem}
This is a truly remarkable result. It does not only assure that
$\diff$ has no nontrivial {\it closed} normal subgroups (which
would be easier to prove). It states that there are no nontrivial
normal subgroups regardless to their topological nature. We shall
frequently use the following simple consequence. In what follows
we say that a vector field is {\bf localized} (somewhere), if its
support is contained in an open nondense interval of the circle.
\begin{corollary}\label{generators.of.diff}
$\diff$ is generated by exponentials: every element can be written
as a finite product of exponentials. In fact, every element can be
written as a finite product of exponentials of {\rm localized}
vector fields (in general, however, with vector fields having no common 
interval of localization).
\end{corollary}
\begin{proof}
The subgroup generated by the local exponentials is a normal
subgroup which also contains elements different from the identity.
Hence by simplicity it is the whole group.
\end{proof}

\paragrf
We shall use some finite dimensional Lie subgroups of $\diff$. The
case of the subgroup $\mob$ was already discussed. We shall now
introduce the group $\mob^{(n)}$ for each positive integer $n$ as
the subgroup of $\diff$ containing all elements $\gamma \in \diff$
for which there exists a \mo transformation $\phi \in \mob$
satisfying
\begin{equation}
\gamma(z)^n = \phi(z^n) \;\; (\forall z \in S^1).
\end{equation}
The group $\mob^{(n)}$ gives a natural $n$-covering of $\mob$.
This group has been already considered and successfully used for
the analysis of conformal nets, see e.g.\! \cite{loke,SchWb} and
\cite{LoXu}.

In $\mob^{(n)}$, just like in $\mob$, one introduces the
one-parameter subgroup of translations $a \mapsto \tau^{(n)}_a$,
which is defined by the usual procedure of lifting: it is the
unique continuous one-parameter subgroup satisfying
$\tau^{(n)}_a(z)^n = \tau_a(z^n)$. Alternatively, one may define
it directly with its generating vector field
$t^{(n)}(z)=\frac{1}{2n} - \frac{1}{4n}(z^n+z^{-n})$. Similarly
one introduces the notion of rotations $\alpha \mapsto
\rho^{(n)}_\alpha$ and of dilations $s \mapsto \delta^{(n)}_s$. Of
course the ``$n$-rotations'', apart from a rescaling of the
parameter, will simply coincide with the ``true'' rotations:
\begin{equation}
\label{n-rotation} \rho^{(n)}_\alpha=\rho_{\alpha/n}.
\end{equation}

As $\mob^{(n)}$ covers $\mob$ in a natural way, its universal
cover is canonically identified with $\widetilde{\mob} \simeq
\widetilde{{\rm SL}(2,\RR)}$. Note that if $p^{(n)} :
\widetilde{\mob} \rightarrow \mob^{(n)}$ is the natural covering
map then $p^{(n)}(\tilde{\rho}^{(m)}_\alpha)=\rho^{(n)}_\alpha$
for all $n,m\in\NN^+$, where the sign ``$\tilde{\hphantom{a}}$'',
as always throughout the rest of this article, stands for the
appropriate lifts of one-parameter groups to $\widetilde{\mob}$.

\paragrf
Let us have now a closer look at the Lie algebra of $\diff$.
First, we choose a suitable base in the complexification of this
algebra, of which we can think of as complex valued functions on
the circle. We shall set $l_n$ $(n\in\ZZ)$ for the function
defined by the formula
\begin{equation}
l_n(z)= -i z^n.
\end{equation}
These functions form a (non algebraic) base. By Eq.\!
(\ref{-bracket}) we find that
\begin{equation}\label{deWitt}
[l_n,l_m] = (n-m) l_{n+m}.
\end{equation}
The Lie algebra formed by the (algebraic) span of
$\{l_n:n\in\ZZ\}$ is called the {\it De Witt algebra}. Note that
$l_n=-\overline{l_{-n}}$ and that for any $n\in\NN^+$ the span of
$\frac{1}{n}l_0,\frac{1}{n}l_n,\frac{1}{n}l_{-n}$ form a Lie
subalgebra which is isomorphic to $\mathfrak{sl}(2,\CC)$:
\begin{equation}
[\frac{1}{n}l_{\pm n},\frac{1}{n}l_0]=\pm \frac{1}{n}l_{\pm
n},\;[\frac{1}{n}l_n,\frac{1}{n}l_{-n}] = 2\frac{1}{n}l_0.
\end{equation}
The Lie subgroup of $\diff$, generated by the real parts of the
vector fields $\frac{1}{n}l_0,\frac{1}{n}l_n,\frac{1}{n}l_{-n}$ is
exactly $\mob^{(n)}$.

\paragrf
As we are interested by {\it projective} (and not necessarily
true) representations of $\diff$, what is of real interest for us
is not the De Witt algebra. For a projective representation $U$, it
is possible, that the phases cannot be fixed in a way that
would make it a true representation; that is, with any choice of
the phases in general $U(g)U(h)\neq U(gh)$ but $U(g)U(h)=
\omega(g,h) U(gh)$ where $\omega(g,h)$ is a unit complex number.
At the level of Lie algebra and generators of representation, this
means that the generators can be introduced only up to additive
constants and after fixing these constants, 
in general $[\eta(a),\eta(b)] \neq
\eta([a,b])$ but $[\eta(a),\eta(b)] = \eta([a,b]) + r(a,b)
\mathbbm 1$ for some $r(a,b)\in \CC$ (where $\eta(x)$ is the
generator for $U$ corresponding to the Lie algebra element $x$).

For the explained reason we shall investigate what we can say at
the purely algebraic level if, in stead of Eq.\! (\ref{deWitt}) we
only require the following relations:
\begin{equation}
[\tilde{l}_n,\tilde{l}_m] = (n-m) \tilde{l}_{n+m} +
r_{n,m}\tilde{z}
\end{equation}
where $r_{n,m}$ is a number and $\tilde{z}$ is a {\it central}
element (so $r_{n,m}\tilde{z}$ behaves like the additive constant
$r(a,b) \mathbbm 1$ in our previous consideration):
\begin{equation}
[\tilde{z},\tilde{l}_n]=0.
\end{equation}
Of course the possible dependence of $r_{n,m}$ on $n$ and $m$ is
strongly restricted by the requirement that the elements
$\{\tilde{l}_n:n\in\ZZ\}\cup \{\tilde{z}\}$ should span a Lie
algebra. For example, the antisymmetricalness of the commutator
implies that $r_{m,n} = -r_{n,m}$. The Jacobi identity gives
further equations.

As it was explained, for us the introduced elements $\tilde{l}_n$
have a meaning only up to an additive central element. So let us
choose the following fixing of this freedom:
\begin{eqnarray} \nonumber
L_n &\equiv& \frac{1}{n}[\tilde{l}_n,\tilde{l}_0] = \tilde{l}_n +
\frac{1}{n}r_{n,0}\,\tilde{z} \;\;\;(n\neq 0), \\
\nonumber L_0 &\equiv& \frac{1}{2}[\tilde{l}_1,\tilde{l}_{-1}] =
\tilde{l}_0 + \frac{1}{2} r_{1,-1}\,\tilde{z}, \\
Z &\equiv& 2 [L_2,L_{-2}] - 4 [L_1,L_{-1}] = (2 r_{2,-2} - 4
r_{1,-1})\,\tilde{z}.
\end{eqnarray}
\begin{proposition}\label{deWitt->Vir}
With the above definition we have the following relations:
\begin{eqnarray*}
[Z,L_n] &=& 0, \\
\hphantom{a}[L_n,L_m] &=& (n-m)L_{n+m} + \frac{1}{12}\,Z\,
(n^3-n)\,\delta_{-n,m}.
\end{eqnarray*}
\end{proposition}
\begin{proof}
The first equation follows trivially. Also, the second equation
clearly holds up to an additive element from the center. Hence the
commutator, when it appears inside another commutator, can be
indeed replaced by the expression appearing in the statement. Thus
if $(n+m) \neq 0$ then, using the Jacobi identity, the above
remark and the fact that essentially by definition $[L_n,L_0]=n
L_n$, we find that
\begin{eqnarray}\nonumber
(n-m)L_{n+m} &=&
\frac{n-m}{n+m}[L_{n+m},L_0]=\frac{1}{n+m}[[L_n,L_m],L_0]\\
&=& \frac{1}{n+m}([[L_n,L_0],L_m] + [L_n,[L_m,L_0]]) \\
\nonumber &=& [L_n,L_m].
\end{eqnarray}
We are left with the case when $n=-m$. It is clear that it is
enough to consider the case when $n$ is positive. For $n=1$ we
have that $(n^3-n)=0$ and essentially by definition $[L_1,L_{-1}]
= 2L_0$, thus the relation holds. For $n=2$, again by definition,
we have that
\begin{eqnarray}\nonumber
[L_2,L_{-2}] &=& [L_2,L_{-2}] - 2[L_1,L_{-1}]+4L_0 \\
&=& 4L_0 + \frac{1}{2}Z = 4L_0 + \frac{1}{12}(2^3-2)Z.
\end{eqnarray}
Finally, for $n\geq 2$ one can proceed by induction using the
Jacobi identity and the (already justified) fact that
$L_{n+1}=\frac{1}{n-1}[L_n,L_1]$ and $L_{-n-1}=
\frac{1}{n-1}[L_{-1},L_{-n}]$.
\end{proof}
The Lie algebra spanned by the operators $\{L_n:n\in\ZZ\}$ and $Z$
is called the {\bf Virasoro algebra}. What the above proposition
essentially says is that up to isomorphism there is only one
non-trivial central extension of the De Witt algebra, and it is
the Virasoro algebra. By redefining the central element, the
normalizing factor $\frac{1}{12}$ could be changed (without
changing the isomorphism class of the algebra). However, from the
point of view of positive energy representations of $\diff$, as we
shall see in the next section, the suitable choice is exactly the one 
made here.

\section{Positive energy representations}
\label{sec:diffposenergy} 

\paragrf
In what follows $\mathbb T$ stands for the set of complex numbers
of unit length. We consider the {\bf projective unitary group}
$\U(\H)/\mathbb T$ of a Hilbert space $\H$ as a topological group
with the topology being induced by the strong operator topology on
$\U(\H)$. A {\bf strongly continuous projective representation}
$U$ of a topological group $G$ on a Hilbert space $\H$ is a
continuous homomorphism from $G$ to $\U(\H)/\mathbb T$. Sometimes
we think of $U(g)$ $(g \in G)$ as a unitary operator. Although
there are more than one ways to fix the phases, note that
expressions like Ad$(U(g))$ or $U(g) \in \M$ for a von Neumann
algebra $\M \subset {\rm B}(\H)$ are unambiguous.

Unlike in case of true representations, in general we cannot consider
the direct sum of projective representations. On the other hand,
we can still introduce the concept of invariant subspaces (it is
independent from the choice of phases) and thus, that of
irreducible representations.

Note also that the self-adjoint generator of a one-parameter group
of strongly continuous projective unitaries $t \mapsto Z(t)$ is
well defined up to a real additive constant. This means that there
exists a self-adjoint operator $A$ such that Ad$(Z(t))=$
Ad$(e^{iAt})$ for all $t \in \RR$, and if $A'$ is another
self-adjoint with the same property then $A'=A+r\mathbbm 1$ for
some $r \in \RR$.

A strongly continuous projective unitary representation $U$ of
$\diff$ is called a {\bf positive energy representation} of
$\diff$ if the self-adjoint generator of the one-parameter group
associated to the rotations $\alpha \mapsto U(\rho_\alpha)$ is
bounded from below. Note that this fact does not depend on the
fixing of the additive constant in the definition of the
self-adjoint generator.

The restriction of $U$ to $\mob \subset \diff$ always lifts to a
unique strongly continuous unitary representation of the universal
covering group $\widetilde{\mob}$ of $\mob$ (see Sect.\! \ref{app:mobrep} 
in the appendix). By Prop.\! \ref{pos.mob^n} the self-adjoint generator of 
rotations in this lift is not just bounded below but actually {\it 
positive}. This explains the name ``positive energy representations''. 
From now on we always think of the self-adjoint generator of rotations 
with the additive constant fixed in this way, and we shall denote it
by $L_0$. Note that as $U(\rho_{2\pi})$ is a multiple of the identity, the 
spectrum of $L_0$ is contained in $h+\NN$ where $h\geq 
0$ is the lowest point of ${\rm Sp}(L_0)$.

\paragrf
By the previous observation the (algebraic) span of $L_0$ is a dense
subspace in the Hilbert space of the representation. The vectors
in this set are called the {\bf finite energy vectors} of the
representation.

If $U$ is a positive energy representation of $\diff$, then its
restriction $U^{(n)}$ to $\mob^{(n)}$ lifts to a unique positive
energy representation $\tilde{U}^{(n)}$ of $\widetilde{\mob}$.
Note that by Eq.\! \ref{n-rotation}, the self-adjoint generator of $\alpha 
\mapsto \tilde{U}^{(n)}(\tilde{\rho}^{(n)}_\alpha)$, up to an additive
constant must be equal to $\frac{1}{n}L_0$. Thus by Corollary
\ref{dom.of.Epm} the domain of $L_0$ is contained in every
self-adjoint generator of $\tilde{U}^{(n)}$. Hence the definition
\begin{eqnarray}\label{def.Ln}\nonumber
L_n &\equiv& n\left(H^{(n)} - 2 T^{(n)} - i D^{(n)}\right)|_\fin\\
L_{-n} &\equiv& n\left(H^{(n)} - 2 T^{(n)} + i
D^{(n)}\right)|_\fin
\end{eqnarray}
is meaningful for every $n\in\NN^+$, where $\fin$ is the set of
finite energy vectors for $U$ and
\begin{eqnarray}\nonumber
H^{(n)} &=& -i \frac{d}{d\alpha}
\tilde{U}^{(n)}(\tilde{\rho}^{(n)}_\alpha)|_{\alpha = 0}, \\
\nonumber T^{(n)} &=& -i \frac{d}{da}
\tilde{U}^{(n)}(\tilde{\tau}^{(n)}_a)|_{a = 0}, \\
D^{(n)} &=& -i \frac{d}{ds}
\tilde{U}^{(n)}(\tilde{\delta}^{(n)}_s)|_{s = 0}.
\end{eqnarray}
Of course it is clear that the operator $L_n$, in some sense,
is the generator corresponding to the complex vector field 
$l_n$ introduced in the last section. Indeed, as
$t^{(n)}(z)= \frac{1}{2n}-\frac{1}{4n}(z^n+z^{-n})$ and 
$d^{(n)}(z)=\frac{i}{2n}(z^n-z^{-n})$, where the vector 
fields $t^{(n)}$ and $d^{(n)}$ generate the $n$-rotations 
and the $n$-dilations, respectively, and by considering, that
the generator of the $n$-rotations is the constant vector field 
$\frac{1}{n}$, we find that $L_{\pm n}$ should correspond to
\begin{equation}
z\mapsto -in\left(\frac{1}{n}- 2t^{(n)} \mp i\,  d^{(n)}(z)\right) = 
-iz^{\pm n} =l_{\pm n}(z).
\end{equation}
(We use the physicist convention, according to which the
generator of a one-parameter group of linear operators $t\mapsto 
X(t)$ is minus $i$ times its derivative at time equal zero.)
Rather than introducing $L_n$ as ``the generator corresponding to $l_n$'',
we made this complicated-looking construction because in this way we fixed 
the additive constants and ensured the existence of a common domain.

By the last section, we expect that the introduced operators
satisfy the Virasoro commutation relations with central
element being just a real number. Of course, by their definition
and by Corollary \ref{comm.on.Dinfty}, for example we have already 
that $[L_n,L_0]= n L_n$ and so in particular $L_n\fin \subset \fin$.
However, using only relations coming from the group $\mob^{(n)}$ will
not give information on commutators like $[L_1,L_2]$.

Using the full group $\diff$ and not just the considered subgroups, it 
can be proved that there exists   
a unique real number $c\in\RR$, called the {\bf central charge} of   
the representation $U$ such that
\begin{equation}\label{Vir.rel}
[L_n,L_m] = (n-m)L_{n+m} + \frac{c}{12}\, (n^3-n)\,\delta_{-n,m}.
\end{equation}
The proof can be found in \cite{loke} (although there an additional
condition is assumed, but that can be dropped; see the
remarks made before \cite[Prop.\! A.1]{Carpi2}).

By their commutation relation, the operators $\{L_n:n\in\ZZ\}$ 
are called the {\bf representation of the Virasoro algebra associated to} 
$\mathbf U$. Note that we usually think of $L_n \;(n\in\ZZ)$ as operators 
defined on $\fin$, with the exception of $L_0$ of which sometimes we 
think as a positive self-adjoint operator and thus write expressions like
$\D(L_0)$ for its domain.

$L_0$ is self-adjoint and by Eq.\! (\ref{def.Ln}) in general
\begin{equation}\label{Ln^*=L-n}
L_n^* \supset L_{-n}.
\end{equation}
Moreover, we have
\begin{equation}
L_0 \geq 0.
\end{equation}
Thus we say that the operators $\{L_n:n\in\ZZ\}$ form a {\bf
positive energy unitary representation} of the Virasoro algebra.
The ``unitariness'' is the condition that $L_n^* \supset L_{-n}$;
in some sense it is due to the fact that $\diff$ is represented by
(projective) unitary operators.

The self-adjoint generator of the ``$n$-rotations'' is
\begin{equation}
H^{(n)}= \frac{1}{2}[\frac{1}{n}L_n,\frac{1}{n}L_{-n}] =
\frac{1}{n}L_0 + \frac{c}{24}(n-\frac{1}{n}).
\end{equation}
Hence by the positivity of $H^{(n)}$ (recall that since $H^{(n)}$ is 
bounded from below it is actually positive; see Prop.\! \ref{pos.mob^n})
and the arbitrariness of $n\in\NN^+$ we find that the {\it central
charge must be a positive number}. Also, using that
$\tilde{U}^{(n)}$ is a positive energy representation of
$\widetilde{\mob}$ we find that its self-adjoint generators
satisfy certain bounds with respect to $H^{(n)}$. In fact, taking
account that $H^{(n)}\geq \frac{c}{24}(n-\frac{1}{n})$, by
Corollary \ref{dom.of.Epm} we find (note that it is the operator
$\frac{1}{n}L_n$ that corresponds to a restriction of the one
there denoted by $E_-$) the following {\it energy bound} for
nonnegative $m$:
\begin{equation}
\| L_m \Psi \| \leq \|\sqrt{(L^2_0 + \frac{c}{12}(m^2)L_0}\,\Psi\|
\end{equation}
for every $\Psi \in \fin$. Moreover, as
\begin{eqnarray}\label{|L-n|^2}\nonumber
\|L_{-m}\Psi\|^2 &=& \langle L_{-m}\Psi,\,L_{-m}\Phi\rangle =
\langle \Psi,\,L_m L_{-m}\Psi\rangle\\ \nonumber &=& \langle
\Psi,\,[L_m ,L_{-m}] \Psi\rangle +
\langle \Psi,\,L_{-m} L_m \Psi\rangle \\
\nonumber &=& \langle \Psi,\,(2mL_0 + \frac{c}{12}(m^3-m)\mathbbm
1)\Psi\rangle + \langle L_m\Psi,\, L_m \Psi\rangle \\
&=& 2m \| L^{\frac{1}{2}}_0 \Psi\| +
\frac{c}{12}(m^3-m)\|\Psi\|^2+ \|L_m \Psi\|^2,
\end{eqnarray}
for a general (not necessarily positive) $n\in\ZZ$ we can derive
for example the following bound:
\begin{equation}\label{e.bound.1b}
\| L_n \Psi \| \leq \sqrt{1+\frac{c}{12}}\, (1+|n|^{\frac{3}{2}})
\,\|(\mathbbm 1 + L_0)\Psi\|
\end{equation}
for every $\Psi\in \fin$. (See more on energy bounds in
\cite{GoWa,BS-M}).

As a consequence of the energy bounds, if $f\in C^\infty(S^1,\CC)$
with Fourier coefficients $\hat{f}_n \equiv
\frac{1}{2\pi}\int_0^{2\pi} f(e^{i\theta})e^{-in\theta} d\theta
\;\;\;(n\in\ZZ)$ then the operator $\sum_{n\in\ZZ}\hat{f}_n L_n$
is well defined (i.e.\! the sum converges on $\fin$) and by Eq.\!
(\ref{Ln^*=L-n})
\begin{equation}\label{fL*.subset.overline.f.L-}
(\sum_{n\in\ZZ}\hat{f}_n L_n)^* \supset
\sum_{n\in\ZZ}\hat{\overline{f}}_n L_n,
\end{equation}
thus it is actually closable. Its closure is usually denoted by
$T(f)$. The operator valued functional $T$ is called the {\bf
stress-energy tensor} associated to the representation $U$. By the
energy bound of Eq.\! (\ref{e.bound.1b}), $\D(L_0)\subset \D(T)$
and for every $\Psi\in\D(L_0)$
\begin{equation}\label{e.bound.2}
\| T(f)\Psi \| \leq \sqrt{1+\frac{c}{12}}\, \|f\|_{\frac{3}{2}}
\|(\mathbbm 1 + L_0)\Psi\|
\end{equation}
where the $\|\cdot\|_{\frac{3}{2}}$ norm is defined by the formula
\begin{equation}
\|f\|_{\frac{3}{2}}\equiv \sum_{n\in\ZZ}\hat{f}_n (1 +
|n|^{\frac{3}{2}}).
\end{equation}
By its definition $T(1) = L_0$ and in general if $f,g$ are
polynomials of the variables $z,z^{-1}\in S^1$, then $T(f),T(g)$
are closures of finite sums of Virasoro operators. Thus in this
case $T(l)\fin \subset \fin\;(l=f,g)$ and moreover by Eq.\!
(\ref{Vir.rel}) we have that on $\fin$
\begin{equation}\label{[T(f),T(g)]}
[T(g),T(f)] = i T(gf'-g'f) + i \frac{c}{12}(g,f) \mathbbm 1
\end{equation}
where 
\begin{equation}
(g,f)\equiv -\frac{1}{2\pi}\int_0^{2\pi} 
g(e^{i\theta})(f'(e^{i\theta})+f'''(e^{i\theta}))
d\theta.
\end{equation}
If we set $g$ to be the constant $1$ function, then
$(g,f)= -(f,g) = 0$. Let now $f$ be any smooth function (not
necessarily a polynomial). Then we can find a sequence of
polynomials $f_n\; (n\in\NN)$ such that $f_n\to f$ and $f'_n\to
f'$ in the $\|\cdot\|_\frac{3}{2}$ norm. By the energy bound of
Eq.\! (\ref{e.bound.2}) this implies that on $\fin$ both
$T(f_n)\to T(f)\Psi$ and $T(f'_n) \to T(f')\Psi$ and so by the
previous equation if $\Psi\in\fin$ then
\begin{eqnarray}\nonumber
L_0 T(f_n)\Psi &=& [T(1),T(f_n)]\Psi + T(f_n)L_0\Psi \\
&=& i T(f'_n)\Psi + T(f_n)L_0 \Psi \to (i T(f') + T (f)L_0)\Psi.
\end{eqnarray}
Since also $T(f_n)\Psi \to T(f)\Psi$ the above equation implies
that $T(f)\Psi\in\D(L_0)$ and that
\begin{equation}
[L_0,T(f)]\Psi = iT(f')\Psi.
\end{equation}
To summarize, we have the following situation. By definition
$\fin$ is a core for the closed operators $T(f),T(\overline{f})$
and by Eq.\! (\ref{fL*.subset.overline.f.L-}), $T(f)^*\supset
T(\overline{f})$. The dense subspace $\fin$ is an invariant core
for the nonnegative self-adjoint operator $L_0$ and $T(f)\fin
\subset \D(L_0)$ so the commutator $[L_0,T(f)]$ is well defined.
Moreover, by the energy bound Eq.\! (\ref{e.bound.2}) and the
previous equation, setting $r\equiv \sqrt{1+\frac{c}{12}}\,
(\|f\|_{\frac{3}{2}}+ \|f'\|_{\frac{3}{2}})$ we have that both
$\|T(f)\Psi\|\leq r \|(\mathbbm 1+L_0)\Psi\|$ and
$\|[L_0,T(f)]\Psi\| \leq r \|(\mathbbm 1+L_0)\Psi\|$ for every
$\Psi\in\fin$. Thus by an application of Nelson's commutator theorem 
\cite[Prop.\! 2]{Ne}
\begin{equation}
T(f)^* =\overline{T(f)^*|_\fin} =  T(\overline{f}).
\end{equation}
In particular, if $f$ is real then $T(f)$ is self-adjoint and
moreover by the energy bounds it is essentially self-adjoint on
any core of $L_0$.

If $f,f_n$ $(n\in\NN)$ are real smooth functions with
$f_n\;(n\in\NN)$ converging to $f$ in the
$\|\cdot\|_{\frac{3}{2}}$ sense, then by the energy bound of Eq.\!
(\ref{e.bound.2}) and the previous remark, we have that the sequence
$T(f_n)\;(n\in\NN)$ converges to $T(f)$ on a common core and thus
also in the strong resolvent sense. In particular, we have that
$e^{iT(f_n)}\to e^{iT(f)}$ strongly as $n\to \infty$ (see more on
convergence of self-adjoint operators in \cite[Sect.
VIII.7]{RSI}).

Essentially by definition, if $f$ is such that $\forall
t\in\RR:\;{\rm Exp}(ft)\in\mob^{(n)}$ for a certain positive
integer $n$, then in the projective sense $U({\rm
Exp}(f))=e^{iT(f)}$. If both $f,g$ are such functions (with
possibly different values of $n$), then as $k\to \infty$
\begin{equation}
(e^{i T(f/k)}e^{iT(g/k)})^k \to e^{i T(f+g)}
\end{equation}
strongly, as the Trotter product formula can be indeed employed by
the existence of a common core (see e.g.\! \cite[Theorem
VIII.31]{RSI}). On the other hand, if $f$ and $g$ are sufficiently
small, a similar product formula holds at the level of
diffeomorphisms (see \cite[Theorem 4.1]{flows} and the comment
in \cite{loke} about the $C^\infty$ convergence in the 
$C^\infty$ case). In this way one can
prove that if $e^{iT(h)}=U({\rm Exp}(h))$ for $h=f,g$, then the
same holds for $h=f+g$. By taking sums we get all polynomials of
$z,z^{-1}$. Finally, if $f_n \to f$ in the $C^\infty$ sense, then
of course $\|f_n-f\|_{\frac{3}{2}}\to 0$ which shows that
$e^{iT(f_n)}\to e^{iT(f)}$. Thus one gets that for every $f\in 
C^\infty(S^1,\RR)$, in the projective sense
\begin{equation}
e^{iT(f)}=U({\rm Exp}(f)).
\end{equation}
The above equation can be considered as the {\it
meaning} of the stress-energy tensor. Note that by Corollary
\ref{generators.of.diff}, this also shows that the following three
objects: the representation $U$, the representation of the
Virasoro algebra $L$ and the stress-energy tensor $T$ completely
determine each other.

\paragrf
Let us discuss some further important consequences of the
fundamental property of the stress-energy tensor mentioned above.
First of all, if $f$ is a smooth function and $\gamma\in \diff$
then
\begin{equation}
U(\gamma)T(f)U(\gamma)^* = T(\gamma_*f)+r_{\gamma,f}\mathbbm 1
\end{equation}
for a certain $r_{\gamma,f}\in \CC$. Indeed, if $f$ is real then
$\gamma\circ{\rm Exp}(ft)\circ\gamma^{-1}= {\rm Exp}(t\gamma_* f)$
for all $t\in\RR$ and hence the generators of the one-parameter
groups of unitaries $t\mapsto U(\gamma)e^{itT(f)}U(\gamma)^*$ and
$t\mapsto e^{itT(\gamma_* f)}$ can only differ in an additive
constant. (If $f$ is not real, then first we can write it as a sum
of its real part plus its imaginary part, etc.) This in turn
implies another important property which is general to any
positive energy representation of $\diff$. Although we shall give
a self-contained argument, the reader should note that this was
proved in \cite{BS-M}. We begin first with an observation which is
useful in its own.
\begin{lemma}\label{poly.e.bound}
Let $T$ be a stress-energy tensor with $L_0\equiv T(1)\geq 0$ and
$\D^\infty\equiv \cap_{n\in\NN}\D(L_0^n)$. Suppose $f_1,..f_n$ are
smooth functions. Then there exists an $r>0$ such that
$$
\|T(f_1)..T(f_n)\Psi\| \leq r \|(\mathbbm 1 + L^n_0)\Psi\|
$$
for all $\Psi\in\D^\infty$.
\end{lemma}
\begin{proof}
For $n=1$ we have already stated this property. One can proceed by
induction, the argument is straightforward with no difficulties,
but for clarity here only the case $n=2$ will be concerned. By the
energy bounds and by the commutation relations with $L_0$ we have
that
\begin{eqnarray}\nonumber
\|T(f_1)T(f_2)\Psi\| &\leq& r_1 \|(\mathbbm 1 + L_0) T(f_2)\Psi\|
\\
\nonumber &\le& r_1 \| T(f_2)(\mathbbm 1 + L_0)\Psi\| +r_1\|[L_0,T(f_2)]\Psi\|\\
\nonumber &\leq & r_1 r_2 \|(\mathbbm 1 + L_0)^2\Psi\| +
r_1\|T(f'_2)\Psi\| \\
\nonumber &\leq& 2r_1 r_2\|(\mathbbm 1 + L^2_0)\Psi\| + r_1r_3
\|(\mathbbm 1 + L_0)\Psi\| \\
&\leq& 2(r_1 r_2+r_3)\|(\mathbbm 1 + L^2_0)\Psi\|
\end{eqnarray}
where $\Psi\in\D^\infty$, and the constants $r_1,r_2$ and $r_3$ are
coming from the energy bounds on $T(f_1),T(f_2)$ and $T(f'_2)$
(and are independent from $\Psi\in\D^\infty$).
\end{proof}
\begin{corollary}\label{UD^infty=D^infty}
Let $U$ be a positive energy representation of $\diff$ with $L_0$
and $\D^\infty$ as before. Then $U(\gamma)\D(L_0^n) =
\D(L_0^n)\;(n\in\NN)$ and hence $U(\gamma)\D^\infty = \D^\infty$
for every $\gamma\in\diff$.
\end{corollary}
\begin{proof}
By the transformation property of the stress-energy tensor we have
that 
\begin{equation}
U(\gamma)^* \,L^n_0 \, U(\gamma) = (T(\gamma^{-1}_* 1)+r\mathbbm
1)^n
\end{equation}
for a certain $r\in\RR$. By the previous polynomial energy
bound if $\Psi\in\D(L^n_0)$ then it is also in the domain of
$T(\gamma^{-1}_* 1)^k$ for $k=1,..,n$ and hence in the domain of
$(T(\gamma^{-1}_* 1)+r\mathbbm 1)^n$, too. Therefore
$U(\gamma)\Psi\in \D(L_0^n)$ and thus $U(\gamma)\D(L_0^n)\subset
\D(L_0^n)$. Exchanging $\gamma$ with $\gamma^{-1}$ the same
argument shows that actually $U(\gamma)\D(L_0^n)= \D(L_0^n)$
\end{proof}

\paragrf
We have seen that once we have a positive energy representation of
$\diff$, we also have a positive energy unitary representation of
the Virasoro algebra. So now in turn we shall investigate such
representations of the Virasoro algebra with a certain central
charge $c\in\RR$. Our aim is to classify such representations and 
then to ``integrate them'', in order to get positive energy 
representations of $\diff$. 

By a positive energy unitary representation of the Virasoro algebra we 
shall mean a set of operators $\{L_n: n\in\ZZ\}$ defined on a common 
invariant dense domain $\D$ in a Hilbert space $\H$, satisfying the
following properties:
\begin{itemize}
\item $L_0$ is a positive diagonalizable operator on $\D$ \item
$L_n^* \supset L_{-n}$ \item $[L_n,L_m] = (n-m)L_{n+m} +
\frac{c}{12}\, (n^3-n)\,\delta_{-n,m}$.
\end{itemize}
The eigenspaces of $L_0$ are sometimes called the {\bf energy
levels}. As $L_0 \geq 0$ it has a lowest eigenvalue; the
corresponding eigenspace is the ``$0^{\rm th}$ energy level''.

By the commutation relations if $\Phi\in\D$ is from the eigenspace
of $L_0$ associated to the eigenvalue $q$, then $L_0 (L_n \Phi) =
(q-n)(L_n\Phi)$ and so $L_n \Phi$ belongs to the eigenspace of
$L_0$ associated to the eigenvalue $q-n$. Hence we say that $L_n$
``decreases the energy by $n$''. It follows that if $h \geq 0$ is
the lowest eigenvalue of $L_0$ and $\Phi$ is an eigenvector of
$L_0$ with eigenvalue $h$ then $\Phi$ is annihilated by all $L_n$
(that is, $L_n\Phi = 0$) for $n>0$. Hence using the commutation
relations and the unitariness condition --- exactly like in Eq.\!
(\ref{|L-n|^2}) --- we find that $\|L_{-n}\Phi\|^2=(2nh+
\frac{c}{12}(n^3-n))\|\Phi\|^2$ for every $n\in\NN^+$. This, by
the positivity of the norm and the arbitrariness of $n\in\NN^+$,
implies that $c$ cannot be negative (which we have already noted
in case the representation is coming from a positive energy
representation of $\diff$).

Since $L\geq 0$ and diagonalizable on $\D$, it must have a lowest
eigenvalue $h\geq 0$ and there must be a least one-dimensional
eigenspace associated to $h$. We shall fix a unit vector $\Psi_h$
in this space. We can always consider the minimal invariant
subspace for the representation containing this vector; that is,
the subspace
\begin{equation}
M_{\Phi_h}\equiv {\rm Span}\{L_{n_1}L_{n_2}..L_{n_j}\Phi_h:
j\in\NN,\, n_1,n_2,..,n_j\in\ZZ \}
\end{equation}
where also $j=0$ is allowed (so the vector $\Phi_h$,
with no operator acting on it, is also included). It is an easy
observation that due to the commutation relations and the fact
that $L_n\Phi_h=0$, in reality, to generate $M_{\Phi_h}$ it
suffices to use $L$ indexed with negative values only, and that in
fact one can even require an ordering. Thus
\begin{equation} \label{strong.generates}
M_{\Phi_h}={\rm Span}\{L_{n_1}L_{n_2}..L_{n_j}\Phi_h: j\in\NN,\,
n_1\leq n_2\leq .. \leq n_j < 0 \}.
\end{equation}
The simplest possible situation is when actually $M_{\Phi_h}$ is
the full domain of the representation; $\D=M_{\Phi_h}$. In this case
the representation is irreducible and it is called an {\bf
irreducible lowest weight representation} (with the ``lowest
weight'' being the value of $h$). 

Indeed, suppose $N\subset M_{\Phi_h}$ is an invariant subspace. 
Any vector of $M_{\Phi_h}$ is a finite linear combination of eigenvectors 
of $L_0$. Thus if $\Psi \in N$ then by the invariance under $L_0$ it 
follows that all eigenvectors appearing as ``components'' of $\Psi$ belong 
to $N$. If $\Phi_h \in N$ then of course $N=M_{\Phi_h}$. On the other
hand, by the previous argument if $\Phi_h \notin N$ then it
follows that $N$ is orthogonal to the vector $\Phi_h$. Actually it
is not hard to see, that by the unitariness condition and the
invariance of $N$ all vectors of the type
$L_{n_1}L_{n_2}..L_{n_j}\Phi_h$ are orthogonal to $N$ and hence
$N=\{0\}$.

Let us see what we can say about the energy levels. At the $0^{\rm
th}$ level, up to multiples there is a unique vector: it is
$\Phi_h$. By Eq.\! (\ref{strong.generates}), the $1^{\rm th}$ energy 
level, which is orthogonal to the $0^{\rm th}$ level, is
spanned by $L_{-1}\Phi_h$. The $2^{\rm nd}$ energy level, which is
orthogonal to both the $0^{\rm th}$ and the $1^{\rm th}$ level, is
spanned by $L_{-1}L_{-1}\Phi_h$ and $L_{-2}\Phi_h$. In general,
one can see that to span the $n^{\rm th}$ energy level (at most)
we need as many vectors as many partitions of $n$ we can find.
Thus, similarly to the case of the $U(1)$ current model, we find that for 
such a representation $e^{-\beta L_0}$ is {\bf trace class} for every
$\beta>0$.

We have seen that the length of the vector $L_n \Phi_h$, that is
the scalar product of this vector with itself, is completely
determined by the value of $n, h$ and $c$. Actually one can easily
convince him- or herself, that in general, by the commutation
relations, unitariness, and the facts that $L_0\Phi_h= h\Phi_h$ and
$L_n\Phi_h=0$ for all $n>0$, the scalar product of
$L_{n_1}L_{n_2}..L_{n_j}\Phi_h$ with
$L_{m_1}L_{m_2}..L_{m_k}\Phi_h$ is completely determined by the
sequences $(n_1,..,n_j),\,(m_1,..,m_k)$ and the value of $h$ and
$c$. Thus one can see that the representation, up to unitary
equivalence is completely determined by the value of $c$ and $h$.

The question which is to be answered: for what values of $(c,h)$
there exists such an irreducible lowest weight representation? (We
have seen, that if it exists, then up to unitary equivalence, it
exists uniquely.) That not all values are possible, is sure. First
of all, $h$ was required to be positive. We have also seen, that
$c$ must be nonnegative. But are there more limitations?

As we have already noticed, the second energy level is spanned by
the vectors $\Psi_{1,1} = L_{-1}L_{-1}\Phi_h$ and $\Psi_2 =
L_{-2}\Phi_h$. So if we assume that we indeed have such a unitary
representation, then by the positivity of the scalar product,
using the Cauchy-Schwartz inequality, we find that
\begin{equation}
\left(\begin{matrix} \langle \Psi_{1,1},\,\Psi_{1,1}\rangle && \langle 
\Psi_{1,1},\,\Psi_2\rangle\\
\langle \Psi_{2},\,\Psi_{1,1}\rangle && \langle
\Psi_{2},\,\Psi_{2}\rangle
\end{matrix}\right)
\end{equation}
must be a positive semidefinite matrix. As it was already
explained, we can calculate each matrix element in the above
expression, in fact we have already calculated that
$\|L_{-2}\Phi_h\|^2 = (4h+ \frac{c}{2})$. By calculating the
remaining elements we find that the determinant
\begin{equation}
\left|\begin{matrix} 4h(2h+1) && 6h\\
6h && (4h+ \frac{c}{2})
\end{matrix}\right| \geq 0,
\end{equation}
which immediately gives an inequality between $h$ and $c$. Similar
inequalities using (higher order) determinants can be written using
generating vectors of the $n^{\rm th}$ energy level. These determinants 
give further restriction on the possible values of $c$ and $h$.

We have seen that we can find couples of $(c,h)$ for which surely
{\it does not} exist an irreducible unitary lowest weight
representation. Can we find a couple for which it does? Of course
for example for $c=h=0$ it does, in fact it is not to hard to see
that the corresponding representation is trivial: the
representation space is one-dimensional, and all Virasoro
operators are represented as multiplication by zero. However, can
we find a nontrivial example?

One of the last observations of the previous chapter was that the
$U(1)$ current model is diffeomorphism covariant; in fact in this
overview this was the motivating reason to investigate positive
energy representations of $\diff$. So studying that model surely
gives a unitary lowest weight representation of the Virasoro
algebra. Indeed, on the finite energy vectors of the $U(1)$ current
model the {\it normal product}
\begin{equation}
:\!J^2\!:_n \equiv \sum_{k+n \geq -k} J_{-k}J_{k+n}
+ \sum_{k+n < -k} J_{k+n}J_{-k}
\end{equation}
is meaningful: on any given finite energy vector the expression
results a {\it finite} sum (any finite energy vector is
annihilated by all terms appearing in the expression of $:\!J^2\!:,$
apart from a finite number of them). Given the commutation
relation Eq.\! (\ref{[Jn,Jm]}) of the currents, it is an exercise
to calculate that
\begin{equation}
[\frac{1}{2}:\!J^2\!\!:_n,\frac{1}{2}:\!J^2\!\!:_m] =
(n-m)\frac{1}{2}:\!J^2\!\!:_n + \frac{1}{12}(n^3-n)\delta_{-n,m},
\end{equation}
that is $\frac{1}{2}:\!J^2\!\!:$ satisfies the Virasoro relations
with central charge equal to $1$. Moreover, since the current by
Eq.\! (\ref{J^*=J}) is hermitian, we easily find that
$\frac{1}{2}:\!J^2\!\!:_n^* \supset \frac{1}{2}:\!J^2\!\!:_{-n}$.
Finally, by its definition $\frac{1}{2}:\!J^2\!\!:_0$ is a sum of
positive expressions, and evidently on the vacuum vector $\Omega$
of the $U(1)$ current model it acts as multiplication by zero. So
let $M_{\Omega}$ be the minimal invariant subspace for
$\frac{1}{2}:\!J^2\!\!:$ containing $\Omega$. Then by what was
said, $\frac{1}{2}:\!J^2\!\!:$ on $M_{\Omega}$ gives an
irreducible unitary representation of the Virasoro algebra with
lowest weight equal to $0$ and central charge equal to $1$.

To summarize: on one hand, as it was explained, we can find some
limitations on the possible values of $c$ and $h$, on the other
hand, we can find some models providing examples for admissible
pairs of $(c,h)$. It is a truly remarkable achievement though that
the set of admissible pairs was {\it completely characterized}.
Namely, apart from the case $c=h=0$ (in which case the
representation is trivial) there exists an irreducible
unitary lowest weight representation of the Virasoro algebra if
and only if $(c,h)$ belongs to the {\it continuous} part
$(1+\RR^+_0) \times \RR^+_0$ or to the {\it discrete} sequence
$\{(c(m),h_{p,q}(m)): m\in\NN,\, p=1,..,m+1;\, q=1,..,p\}$ where
\begin{equation}
c(m)=1-\frac{6}{(m+2)(m+3)}
\end{equation}
and
\begin{equation}
h_{p,q}(m)=\frac{((m+3)p-(m+2)q)^2-1}{4(m+2)(m+3)}.
\end{equation}
For references, and an explanation how these result was achieved
(and by who), see the book \cite{Kac}. Note that the continuous part gives 
only values $c\geq 1$ for the central charge, while the discrete part 
gives a discrete set for $c$ in the interval $[\frac{1}{2},1)$ 
accumulating at the point $1$. Note also that the lowest possible 
nontrivial value for $c$ is $\frac{1}{2}$.

\paragrf
After finding all irreducible unitary lowest weight representation
of the Virasoro algebra, the next question we should ask: can we
integrate them to obtain a positive energy representation of
$\diff$? Let us first note, that given such a representation on a
dense subset of a Hilbert space $\fin\subset \H$, setting $h\equiv
\frac{1}{n}L_0 + \frac{c}{24}(n-\frac{1}{n})$ (where $c$ is the
central charge of the representation) and $e_\pm\equiv
\frac{1}{n}L_{\pm n}$ for a positive integer $n$, the operators
$h,e_-$ and $e_+$ form a positive energy representation of the Lie
algebra of $\widetilde{\mob}$ and thus (more precisely, by Prop.\!
\ref{energybound:mob}) we have that
\begin{equation}
\| e_\pm \Psi\|\leq \| (\mathbbm 1 + h)\Psi\|
\end{equation}
for all $\Psi\in \fin$. So if $f$ is a smooth function on the
circle with Fourier coefficients $\{\hat{f}_n:\in\ZZ\}$, then the
sum $\sum_{n\in\ZZ}\hat{f}_n L_n$ is well defined (converges) and
gives a closable operator. Setting in general $T(f)$ for its
closure, we have that $T(f)^*\supset
T(\overline{f})$, $\D(L_0)\subset \D(T(f))$ and $T(f)\fin\subset
\D(L_0)$ for every $f\in C^\infty(S^1,\CC)$. Moreover on $\fin$
we have that $[L_0,T(f)]=i T(f')$, and there exists a constant $r_f
> 0$ such that $\|T(f)\Psi\|\leq r_f \| (\mathbbm 1 + L_0)\Psi\|$
for all $\Psi \in \fin$. Thus by Nelson's commutator theorem $T(f)^* =
T(\overline{f})$.

All this is verified exactly in the way as it was shown in the
first part of this section when the stress-energy tensor was
introduced; although there $T$ was defined by a positive energy
unitary representation of the Virasoro algebra {\it associated} to
a positive energy representation of $\diff$. However, note that
the fact that the positive energy unitary representation of the
Virasoro algebra is coming from a representation of $\diff$ was
never used directly: only through the energy bound which we have
in this case, too. So, also Eq.\! (\ref{[T(f),T(g)]}) holds,
showing that the functional $f \mapsto iT(f)$ can be viewed as a
representation of the Lie algebra of $\diff$ (with a cocycle).
Hence we may expect the formula
\begin{equation}\label{U(exp)U(exp)=exp(T)exp(T)}
U({\rm Exp}(f_1))U({\rm Exp}(f_2))\ldots U({\rm Exp}(f_n))\equiv
e^{iT(f_1)}e^{iT(f_2)}\ldots e^{iT(f_n)}
\end{equation}
to define a positive energy representation of $\diff$. (Recall
that the group $\diff$ is generated by the exponentials.) However,
there are some difficulties. The most serious is the following:
unlike in case of finite dimensional Lie groups, the exponential
map between ${\rm Vect(S^1)}$ and $\diff$ is not a local
diffeomorphism. In fact, there is not even a neighbourhood of the
identity of $\diff$ in which the exponentials would form a dense
set: every element can be written as a {\it product} of
exponentials but there are only ``few'' elements that {\it are}
exponentials.

After so much tension built up, it is time to reveal to the reader
that yet this problem was solved; the proof for the existence of the
corresponding positive energy representation of $\diff$ can be 
found in \cite{GoWa}. In this short overview however the argument of
Goodman and Wallach will be not discussed.

\paragrf
It is easy to prove that the positive energy representation
$U_{c,h}$ obtained by integrating the unique irreducible unitary
lowest weight representation of the Virasoro algebra with lowest
weight $h$ and central charge $c$, is irreducible. There are two
questions remaining. First, is every irreducible positive energy
representation of $\diff$ is of this type? Second, what can we say
about the representations that are not irreducible? To answer, the
following lemma will be useful.
\begin{lemma}
The group generated by $\cup_{n\in\NN^+} \mob^{(n)}$ is dense in
$\diff$.
\end{lemma}
\begin{proof}
Taking account of the mentioned Trotter product-formula for the
exponential of a sum of vector fields, the statement is rather 
clear at the Lie algebra level: the span of the functions 
$l_n\;(n\in\ZZ)$ is dense in $C^\infty(S^1,\CC)$. All we need is an
``exponential'' version of this density, since as it was mentioned in 
the previous section, every element in $\diff$ is a finite product of
exponentials. Thus, we want to know if ${\rm Exp}$ is continuous:
so that for a sequence of real functions $f_n\; (n\in\ZZ)$ converging
to $f$ in the $C^\infty$ sense, wether it follows that ${\rm
Exp}(f_n)\to {\rm Exp}(f)$ in the topology of $\diff$. The answer
is yes, for references see for example \cite{Milnor}.
\end{proof}
\begin{corollary}\label{closure.of.M.is.inv}
Let $U$ be a positive energy representation of $\diff$ with
associated representation of the Virasoro algebra
$\{L_n:n\in\ZZ\}$ and finite energy vectors $\fin$. Then if
$M\subset \fin$ is invariant for $\{L_n:n\in\ZZ\}$ then its
closure is invariant for $U$.
\end{corollary}
\begin{proof}
By Prop.\! \ref{energybound:mob} the finite energy vectors are
analytic for the self-adjoint generators of $\tilde{U}^{(n)}$. It
follows that if $\gamma\in\mob^{(n)}$ then
$U(\gamma)\overline{M}\subset \overline{M}$. Hence the result
follows by the previous density-lemma.
\end{proof}
We can now answer those two questions. We have seen that given an
irreducible positive energy representation $U$ of $\diff$ we can
construct a positive energy unitary representation of the Virasoro
algebra $\{L_n:n\in\ZZ\}$ with a certain central charge
$c\in\RR^+_0$. As $L_0$ is nonnegative and $e^{i2\pi L_0}\in
\CC\mathbbm 1$, its spectrum has a minimum point $h\geq 0$ and
there exists a unit vector $\Phi_h$ such that it is an eigenvector
of $L_0$ with eigenvalue $h$. Let $M_{\Phi_h}$ be the minimal
invariant subspace for the Virasoro operators containing the
vector $\Phi_h$. Then, by Corollary \ref{closure.of.M.is.inv} the
non zero-dimensional closed subspace $\overline{M_{\Phi_h}}$ is
invariant for $U$ and hence it is the full Hilbert space.
Moreover, we have also seen that the operators $\{L_n:n\in\ZZ\}$
on $M_{\Phi_h}$ form a unitary lowest weight representation of the
Virasoro algebra which is irreducible. Thus by the result of
Goodman and Wallach it integrates to an irreducible positive
energy representation of $\diff$, which must be the original
representation $U$ as it has the same associated representation of
the Virasoro algebra. Let us also note, that as the values of $c$
and $h$ are constructed in an unambiguous way from the
representation, two positive energy irreducible representations
are equivalent if and only if they have the same central charge
and the same lowest weight.

Thus we have proved that a positive energy irreducible
representation $U$ of $\diff$ is equivalent to $U_{c,h}$ for a
certain (unique) admissible pair $(c,h)$ of the central charge and
the lowest weight. (The reader should note that this was first
stated in \cite[Theorem A.2]{Carpi2}.)

Let us move onto the second question. Can we establish that a
reducible representation is necessarily a direct sum (or integral)
of irreducible ones? What does it mean, at all: direct sum of {\it
projective} representations? What is rather clear is the {\it
tensor product} of projective representations, but not their
direct sum.

Suppose $U_1,U_2$ are projective unitary representations of a
group $G$ on the Hilbert spaces $\H_1$ and $\H_2$, respectively.
We can try to fix the phases of $U_1$ and $U_2$ in such a way,
that, after this arrangement $g\mapsto U_1(g)\oplus U_2(g)$ will
give a projective unitary representation of $G$ on $\H_1\oplus
\H_2$. It is not clear that we can make such fixing of the phases,
but that is not a problem: in that case we shall simply say that
the direct sum of $U_1$ and $U_2$ does not exist. The real
problem is that different fixing of the phases of $U_1$ and $U_2$
can lead into different (even inequivalent!) projective
representations. If it is so, the direct sum of $U_1$ and $U_2$ is
not well-defined and thus in general the decomposition of a
reducible representation into irreducible ones is also not well
defined.
\begin{lemma}
Suppose $G$ is a simple noncommutative group. Then the direct sum
of some of its projective representations, if exists, then exists
uniquely.
\end{lemma}
\begin{proof}
Suppose we have two representations $U_1$ and $U_2$ with some
concrete fixing of the phases with which the direct sum can be
formed and suppose $g\mapsto \lambda_j(g)U_j(g)\;(j=1,2)$ is
another choice of the phases with which, again, the direct sum can
be formed. Then we can easily show that there must exists an
$\omega : G\times G \rightarrow \mathbb T$ such that
$\lambda_j(g)\lambda_j(h)= \omega(g,h)\lambda_j(gh)$ for every
$g,h\in G$ and $j=1,2$. Hence $\frac{\lambda_1}{\lambda_2}$ is a
character on $G$ and thus by the conditions of the lemma
$\lambda_1 = \lambda_2$ which shows that the two direct sums
result the same (not only equivalent!) projective representation.
\end{proof}
As it was mentioned, $\diff$ is a simple noncommutative group. So
by the previous lemma, if a certain closed subspace $\K\subset \H$
is invariant for a projective unitary representation $U$ of
$\diff$ given on the Hilbert space $\H$, then we can indeed write
that $U\equiv U|_{\K}\oplus U|_{\K^\perp}$. So given a positive
energy representation of $\diff$, we can try to decompose it. The
next theorem is probably not an ``original'' one although the
author has not yet seen it stated in the literature.
\begin{theorem}\label{pos.diffrep=+irr.}
A positive energy representation $U$ of\, $\diff$ is (equivalent)
to a direct sum of irreducible ones. Moreover, every
representation appearing in the direct sum decomposition has the
same central charge as $U$, and the difference of the lowest
weights of two irreducible components is always an integer.
\end{theorem}
\begin{proof}
The affirmations about the values of the lowest weight and the
central charge are trivial. All we have to do is to prove the
existence of a decomposition into irreducible representations.

Let $\{L_n: n\in\ZZ\}$ be the associated representation of the
Virasoro algebra and $h$ is the minimum point of the spectrum of
$L_0$. We can choose a complete orthonormed system $\Phi_h^\alpha
\;(\alpha \in I_h)$ in ${\rm Ker}(L_0 -h\mathbbm 1)$, where $I_h$
is some index set. Let $M_{\Phi_h^\alpha}$ be the minimal
invariant subspace for the Virasoro operators containing the
vector $\Phi_h^\alpha$. Then, by Corollary
\ref{closure.of.M.is.inv} for every $\alpha \in I_h$ the closed
subspace $\overline{M_{\Phi_h^\alpha}}$ is invariant for $U$.
Moreover, just like as we argued before, we can show that the
restriction of $U$ onto $\overline{M_{\Phi_h^\alpha}}$ is
irreducible and that the subspaces $M_{\Phi_h^\alpha}\;(\alpha\in
I_h)$ are pairwise orthogonal. Let $\H_1\subset \H$ be the
orthogonal in the Hilbert space $\H$ of the representation to
$\oplus_{\alpha\in I_h} \overline{M_{\Phi_h^\alpha}}$. It is again
a closed invariant subspace for $U$, but here the lowest point of
the spectrum is at least $h+1$. Thus by choosing a complete
orthonormed system $\Phi_{h+1}^\alpha \;(\alpha \in I_{h+1})$ in
${\rm Ker}((L_0 -h\mathbbm 1)|_{\H_1})$ we can repeat all what was
said. Note that
\begin{equation}
{\rm Ker}((L_0 -h\mathbbm 1)|_{\H_1}) = {\rm Ker}(L_0 -h\mathbbm
1) \cap (\oplus_{\alpha\in I_h} M_{\Phi_h^\alpha})^\perp
\end{equation}
as $L_0$ is diagonal on $M_{\Phi_h^\alpha}$ for every $\alpha\in
I_h$. Then we take the orthogonal of these subspaces in $\H_1$ and
we define $\H_2$, etc. Continuing like this (it is clear how to
make this ``continuing like this'' more precise by Zorn's lemma)
we get for every $n\in\NN$ an index set $I_{h+n}$ and a collection
of closed subspaces $\overline{M_{\Phi_{h_n}^\alpha}}\;(\alpha\in
I_{h+n})$. From the construction it is clear that the subspaces
$\overline{M_{\Phi_{h_n}^\alpha}}\;(n\in \NN,\,\alpha\in I_{h+n})$
are pairwise orthogonal minimal invariant closed subspaces and
that their sum (even their algebraic sum) contains every
eigenvector of $L_0$ and hence
\begin{equation}
\overline{\mathop{\oplus}_{n\in\NN} \mathop{\oplus}_{\alpha \in
I_{h+n}} \overline{M_{\Phi_{h_n}^\alpha}}} =\H.
\end{equation}
\end{proof}
{\it Remark.} We have not answered the question: what are the
irreducible positive energy representations of $\diff$ that can be
added? We shall not show in this overview, but actually the
converse of the above theorem is also true. Namely, if a
collection of irreducible positive energy representations of
$\diff$ have all the same central charge and the pairwise
differences in their lowest weight are integers, then we can form
their direct sum.

Actually, once showing this, we can prove the above theorem {\it
without} using Corollary \ref{closure.of.M.is.inv} about the
invariance. In fact, we did not use that lemma for the
decomposition of the associated representation of the Virasoro
algebra into irreducible ones. After that we could have said, that
in (the closure of) each subspace we can integrate the
representation of the Virasoro algebra and we get a representation
of $\diff$. If we know that these can be added, then by adding
them we get a positive energy representation of $\diff$ on the
full Hilbert space possessing the same associated representation
of the Virasoro algebra as the original representation of $\diff$
--- so it {\it is} the original one. In this way, the decomposition at
the level of the Virasoro algebra can be passed onto the level of the 
representation of the group, without the need of showing that if a 
subspace is invariant for the Virasoro operators then its closure is 
invariant for the representation.

\paragrf
Let us make a final remark. In general, if $U$ is a projective
representation of a group $G$ and $g,h\in G$, it is meaningful to
ask whether $U(g)$ and $U(h)$ commute or not. Indeed, although
$U(g)$ and $U(h)$, as operators, depend on the choice of phases,
whether they commute or not is unambiguous. (One can also think
about the von Neumann algebras $\{U(g), U(g)^*\}''$ and $\{U(h),
U(h)^*\}''$, which are obviously well defined, and ask if one is
in the commutant of the other or not.) Note that the commutation
of $U(g)$ and $U(h)$, unlike in case of a true representation,
does not necessarily follow from the commutation of $g$ and $h$.

Unless its central charge is zero (and thus, as it was already
noticed, the whole representation is trivial) a positive energy
representation $U$ of $\diff$ is {\it genuinely} projective: from
the commutation relation of the associated Virasoro algebra, there
is no arrangement of the phase factors that would make it a true
representation. Thus, even if $\gamma,\varphi\in \diff$ are
commuting elements. it does not necessarily follow that
$U(\gamma),U(\varphi)$ are commuting operators. In order to
generate local nets though, we shall need to know, that
nevertheless, in certain cases it does follow. The following easy
observation, not particular to $\diff$ but rather a general fact
about projective representations, will suffice.

{\it Observation.} Suppose $U$ is a projective projective
representation of a group $G$ and $g \in G$. Set $C(g)$ for the
{\it centralizer} of $g$ in $G$; i.e.\! the subgroup of $G$ formed
by the elements commuting with $g$. Then by setting
\begin{equation}
U(g)U(h)U(g)^*U(h)^* = \lambda(h)\mathbbm 1 \;\;\;(h\in C(g))
\end{equation}
one can see that $\lambda$ is a well-defined (i.e.\! independent
from the choice of phases) character on $C(g)$. Hence if
$h=h_1h_2h_1^{-1}h_2^{-1}$ for some $h_1,h_2\in C(g)$ then $U(g)$
commutes with $U(h)$.

\section[Virasoro nets and $\diff$ covariance]
{Virasoro nets and $\mathbf {{\rm \bf Diff}^+(S^1)}$ covariance}
\label{sec:diffcov}

\paragrf
We say that a strongly continuous projective unitary
representation $U$ of $\diff$ {\bf extends} the unitary (true)
representation $V$ of $\mob$ if in the projective sense
$U(\varphi)=V(\varphi)$ for all $\varphi\in\mob$. Note that by
considering $U$ only on $\mob$ one can reconstruct the original
{\it true} representation $V$ (due to the fact that there are no
continuous nontrivial characters on $\mob$). Thus we do not loose
any information if we ``throw out'' the true representation $V$
and keep only the projective representation $U$. Note also that if
$U$ extends a positive energy representation then it also must be
of positive energy. Let us now proceed to the notion of
diffeomorphism covariance.

\begin{definition}
\label{diffcov:def} A \mo covariant local net $(\A,U)$ is {\bf
conformal (or diffeomorphism)  covariant} if there is a strongly
continuous projective unitary representation of $\diff$ on $\H_\A$
extending $U$ (and by an abuse of notation we shall denote this
extension, too, by $U$), such that for all $\gamma \in \diff$ and
$I \in \I$
\begin{itemize}
\item $U(\gamma)\A(I)U(\gamma)^* = \A(\gamma(J)),$ \item
$\gamma|_I={\rm{id}}_I \Rightarrow
\rm{Ad}(U(\gamma))|_{\A(I)}=\rm{id}_{\A(I)}$.
\end{itemize}
\end{definition}
{\it Remarks.} There are several important things to note here.
First, it seems that there is a slight ambiguity in the
definition. When next time we speak about a conformal net do we
mean: with a {\it given} representation of $\diff$, or with the
representation of the \mo group only (with diffeomorphism
covariance meaning only the existence of such an extension of the
symmetry)? However, it will be proved that {\bf there is at most
one way} the representation of the \mo group can be extended to
$\diff$ with the given conditions. (In fact, this will be one of
the ``highlights'' of this thesis.) So there is no ambiguity. This
also explains the mentioned ``abuse'' in the definition.

Another important observation is that although we have required
two properties from the extension, it turns out that the second
one implies the first, as it was already noted in \cite[Prop.\!
3.4]{Carpi2}. Let us see why. As a consequence of {\it Haag
duality}, if a diffeomorphism is localized in the interval $I$ ---
i.e. it acts trivially (identically) elsewhere --- then, by the
second listed property, the corresponding unitary is also
localized in $I$ in the sense that it belongs to $\A(I)$.
(Actually, by the continuity of the net it is enough to assume
localization in $\overline{I}$).

Assuming the second property, we want to show the first, namely,
the covariance of the net under the action of the diffeomorphism
group. So suppose $\gamma\in\diff$ and $I\in\I$. Of course we know
that in particular covariance holds for \mo transformations. Hence
we can safely assume that $\gamma(I)= I$. With some elementary
(but not trivial) geometrical argument one can show that for every
$\epsilon >0$ there exists a diffeomorphism $\gamma_{I_\epsilon}$
which is localized in the interval $I_\epsilon= I$ ``plus-minus
$\epsilon$'', such that its action on $I$ coincides with that of
$\gamma$. Then $\gamma_{I_\epsilon}^{-1}\circ\gamma$ is localized
in $\overline{I^c}$ and hence $U(\gamma_{I_\epsilon})^*U(\gamma)$
commutes with $\A(I)$. It follows that
\begin{equation}
U(\gamma)\A(I)U(\gamma)^*=U(\gamma_{I_\epsilon})\A(I)U(\gamma_{I_\epsilon})^*\subset
\A(I_\epsilon),
\end{equation}
since by our assumption $U(\gamma_{I_\epsilon})\in\A(I_\epsilon)$.
As the above equation is true for all $\epsilon>0$, by the
continuity of the net we conclude that
$U(\gamma)\A(I)U(\gamma)^*=\A(I)$.

This second property is usually called {\bf compatibility with the
local structure} of the net $\A$. By what was explained,
diffeomorphism covariance means that the \mo symmetry of the net
can be extended to $\diff$ in a way which is compatible with the
local structure of the net.

\paragrf
So far we have seen the following examples of \mo covariant local
nets: the $U(1)$ current model and its derivatives. As it was
noted, the first one is diffeomorphism covariant. What about the
derivatives? It turns out that they do not admit diffeomorphism
symmetry \cite{GLW,koester03a}. One way to see this is the
following. Using diffeomorphism covariance it is possible to show
that \mo symmetry of a conformal net $(\A,U)$ is implementable in
every locally normal representation of $\A$. As it was already
mentioned though, the derivatives of the $U(1)$ current model
admit locally normal representations in which the \mo symmetry is
not implementable.

One can of course say that the derivatives are not ``nice''
models; they are not even $4$-regular. In fact, it is expected
that under some conditions diffeomorphism covariance is automatic.

In the last section all positive energy irreducible representation
of the group $\diff$ were described. We shall now see that some of
them can be used to give new examples of local nets. Let
$U_c\equiv U_{c,0}$ be one of those irreducible positive energy
representations of $\diff$ with lowest weight $h=0$. (See the last
section for the possible values of $c$). Denote by $\Omega$ the
(up to phase) unique unit eigenvector of $L_0$ with zero
eigenvalue. As it was explained, the distance between two
eigenvalues of $L_0$ is always an integer and hence ${\rm
Sp}(L_0)\subset \NN$. This implies that the phase factors of $U_c$
on the elements of $\mob$ can be arranged in a unique way so that
it will become a true representation for $\mob$.

Let $\{L_n:n\in\ZZ\}$ be the associated representation of the
Virasoro algebra. By Eq.\! (\ref{|L-n|^2}) $\|L_{-1}\Omega\|^2 =
0$ and so $\Omega$ is annihilated by both $L_1,L_0$ and $L_{-1}$.
Hence if $\varphi \in \mob$ then $U_c(\varphi)\Omega = \Omega$ (if
the phases are fixed in the above explained way). We set
\begin{equation}
\A_{{\rm Vir},c}(I)\equiv
\{U(\gamma_1\gamma_2\gamma_1^{-1}\gamma_2^{-1}):
\gamma_1,\gamma_2\in {\rm Diff}^+_{I^c}(S^1)\}''\;\;\;(I\in\I)
\end{equation}
where for a $K\in\I$
\begin{equation}
{\rm Diff}^+_{K}(S^1)\equiv \{\gamma\in\diff: \gamma|_K = {\rm
id}_K\}
\end{equation}
that is, the {\it stabilizer} of the interval $K$, or to put it in
another way, the diffeomorphisms localized in $K^c$ (\,=\,acting
nontrivially only in the interval $K^c$).

It is clear that the net $\A_{{\rm Vir},c}(I)$ satisfies isotony
and by the last observation of the previous section also locality.
Moreover, as it was mentioned, the group $\diff$ is
(algebraically) simple. As the subgroup generated by the set of
elements
\begin{equation}
\{\gamma_1\gamma_2\gamma_1^{-1}\gamma_2^{-1}:I\in\I,\,
\gamma_1,\gamma_2\in {\rm Diff}^+_{I^c}(S^1)\}
\end{equation}
is normal and evidently contains elements different from the
identity, by the simplicity of $\diff$ it must be the whole group.
Hence
\begin{equation}
\{\A_{{\rm Vir},c}(I):I\in\I\}'\subset U_c(\diff)'=\CC \mathbbm 1
\end{equation}
and so $\Omega$ is cyclic for $\A_{{\rm Vir},c}(S^1)$. Hence $(\A_{{\rm 
Vir},c},U_c)$ is a diffeomorphism local covariant net; it is called the 
{\bf Virasoro net at central charge $\mathbf c$}. (Actually to check the 
compatibility condition one again needs to use the last observation of the 
previous section.) By Haag-duality if $\gamma\in {\rm Diff}^+_{I^c}(S^1)$ 
for a certain $I\in\I$ then $U(\gamma)\in\A(I)$, so it turns out that we
could have defined the net by
\begin{equation}
\A_{{\rm Vir},c}(I)= \{U(\gamma): \gamma\in {\rm
Diff}^+_{I^c}(S^1)\}''\;\;\;(I\in\I).
\end{equation}
(We have not done so because in this case to check locality would
have been not so straightforward.) In fact, by Corollary
\ref{generators.of.diff} we could have defined the same net also
by
\begin{equation}
\A_{{\rm Vir},c}(I)= \{e^{iT_c(f)}: f\in C(S^1,\RR),
f|_{I^c}=0\}''\;\;\;(I\in\I)
\end{equation}
where $T_c$ is the stress-energy tensor of $U_c$.

The Hilbert space of a Virasoro net, by construction (see again
the last section), is separable. Moreover, as it was remarked,
$e^{-\beta L_0}$ is trace class for every $\beta
>0$ so a Virasoro net is actually {\bf split}.

\paragrf
The Virasoro nets are not just {\it some} new models: their real
importance lies in the fact that they are essential building
blocks for any conformal net. Indeed, suppose $(\A,U)$ is a
conformal net. Then by the stress-energy tensor $T$ of $U$ we can
define
\begin{equation}
\T_\A(I)\equiv \{e^{iT(f)}: f\in C(S^1,\RR),
f|_{I^c}=0\}''\;\;\;(I\in\I),
\end{equation}
which is obviously a \mo covariant subnet of $\A$. Actually,
$\T_\A$ is evidently covariant under the action of the full
diffeomorphism group, too. Moreover, its vacuum Hilbert space is
invariant for $U$ and the restriction of $U$ to this space is
irreducible. So by the classification of the positive energy
irreducible representations of $\diff$ (see the last section, or
directly \cite[Theorem A.1]{Carpi2}), the restriction of this
subnet to its vacuum Hilbert space is isomorphic to the Virasoro
net with central charge $c$, where $c$ is the central charge of
$U$. Thus any conformal local net can be viewed as an extension of
a Virasoro net.

It is interesting to note that while every diffeomorphism
covariant net carries a natural embedding of a Virasoro net, there
is no nontrivial net that could be embedded in a Virasoro net
(apart from itself). More precisely, a Virasoro net has no
nontrivial \mo covariant subnets, cf.\! \cite{Carpi98}. So in
some sense the Virasoro nets are {\it minimal}.

\paragrf
What other structural properties are known about the Virasoro
nets? As for their representations: it seems natural that the map
for a fixed $I\in\I$
\begin{equation}
\pi^c_{h,I}:\,e^{iT_c(f)} \mapsto e^{iT_{c,h}(f)} \;\;\;(f\in
C(S^1,\RR), f|_{I^c}=0),
\end{equation}
--- where $T_{c,h}$ is the stress-energy tensor of a positive energy
irreducible representation of $\diff$ with the same central
charge, but with lowest weight $h$ possible different from zero
--- should extend to a normal representation $\pi^c_{h,I}$ of
$\A_{{\rm Vir},c}$. (The central charge, by Eq.\!
(\ref{[T(f),T(g)]}) is locally detectable, but the lowest weight
should not be.) If it is so, then the collection $\pi^c_h\equiv
\{\pi^c_{h,I}:I\in\I\}$ is evidently a locally normal
representation of our Virasoro net.

What are the values of $c$ and $h$ (among those admissible ones,
for which there exists a positive energy representation), for
which $\pi^c_h$ is indeed known to be a locally normal
representation? The answer is that for almost all. More precisely,
for all with the exception of $1<c<2$, $h<(c-1)/24$.

One way local normality can be proved was already indicated in the
previous chapter. Namely, whenever we have a subnet of a net
$\B\subset \A$, it can be considered as a representation of the
restriction of $\B$ to its vacuum Hilbert space. Then such a
representation, and all subrepresentations of such a
representation are automatically locally normal. Moreover, if
$\pi$ is a representation of $\A$, then it can be considered as
representation of $\B$, too. Of course, if it was locally normal
as a representation of $\A$, it will be locally normal for $\B$,
too. Again, all subrepresentations of such a representation are
automatically locally normal.

Without aiming completeness, let us see at least one example; 
let us show the local normality of $\pi^c_h$ for the special 
case of $c=h=1$. For other values see the explanation in \cite[Sect.\! 
2]{Carpi2} with the references given there.

We have seen that the $U(1)$ current model is diffeomorphism
covariant and in fact we have even seen that the normal product of
the current with itself gives a positive energy representation of
the Virasoro algebra with central charge $c=1$. We have not proved
that it indeed corresponds to a positive energy representation $U$
of $\diff$ which makes the $U(1)$ current model diffeomorphism
covariant, but in fact it is easy to show this; so just suppose
that this is the case. Then, by what was explained we have the
natural embedding
\begin{equation}
\A_{{\rm Vir},1} \subset \A_{U(1)}.
\end{equation}
Since $L_{-1}\Omega = 0$, by the what was explained in the last
section, in a Virasoro net, the spectrum of $L_0$ does not contain
the value $1$. On the other hand, as $J_{-1}\Omega \neq 0$ (in fact,
by the commutation relation of the current, the hermicity
condition and the fact that $J_1\Omega = 0$ one can easily compute
that $\|J_{-1}\Omega\|=1$), it follows that that $\pi^1_1$ appears 
as a subrepresentation of the representation given by the above
embedding, showing its local normality.

So, as it was explained, for many values of $(c,h)$ we know the
existence of a locally normal representation $\pi^c_h$. As for the
converse direction, by \cite[Prop.\! 2.1]{Carpi2} we know that an
irreducible representation of a Virasoro net $\A_{{\rm Vir},c}$
must be (equivalent with) one of those $\pi^c_h$ for a certain
value of $h$ --- given that the representation is of positive
energy. Actually in this thesis the condition about positivity
will be eliminated. So in fact all irreducible representations are
of the discussed type.

By looking at the admissible values of $c$ and $h$ and by what was
said, it is clear that $\A_{{\rm Vir},c}$, in case $c<1$, has
finitely many, and in case $c\geq 1$ has infinitely many
inequivalent irreducible representations. Thus the $\mu$-index of
$\A_{{\rm Vir},c}$ is infinite if $c\geq 1$ and we have a strong
indication that it should be finite if $c<1$; in fact this is exactly
the case, cf.\! \cite{KL}.

Finally, we have not discussed strong additivity. Without further
explanation on how the result was achieved, let it be mentioned
that $\A_{{\rm Vir},c}$ is known to be strongly additive for
$c\leq 1$ (see \cite{KL} for the case $c<1$ and \cite{Xu.strong} for
the case $c=1$) and not strongly additive for 
$c>1$, cf.\! \cite{BS-M}. Thus $c=1$ is an
interesting ``border'' point: if $c<1$ then the Virasoro net is
completely rational. If $c=1$ then its $\mu$-index is infinite but
at least it is still strongly additive. Finally, if $c>1$ then it
is not even strongly additive.

\chapter[Towards the construction of diffeomorphism symmetry]
{Towards the construction of diffeomorphism symmetry}

\resume{
The aim of this chapter is to pinpoint some (algebraic) conditions
on a \mo covariant local net that could ensure the existence of
$\diff$ symmetry. After some motivating arguments, the group of
(once) differentiable piecewise \mo transformations will be
introduced. It will be then shown that the \mo symmetry of a
regular net $\A$ (i.e.\! a net which is $n$-regular for every
$n\in\NN^+$) has an extension to this group of piecewise \mo
transformations if and only if a certain (but purely algebraic)
condition is satisfied regarding the ``relative position'' of
$\A(I_1),\A(I_2)$ and $\A(I)$, where $I_1,I_2 \subset I$ are the
intervals obtained by removing a point from $I\in \I$.

It will be shown that any diffeomorphism can be (for example
uniformly) approximated by them; even in a stronger, local sense.
(Uniform topology serves simply to illustrate that there are
``many'' such transformations
--- even among the local ones. In a later chapter a more relevant
statement will be proved about their density). Thus, although the
symmetry is extended only to these piecewise \mo transformations
and not to the group $\diff$, assuming regularity, this algebraic
condition is conjectured to imply diffeomorphism covariance.

The relation of this condition with others will be investigated.
First, it will be confirmed that conformal covariance implies this
condition. However, the point is that it can be checked in many
cases without assuming conformal covariance. In particular it will
be shown to be satisfied for every completely rational net (which
are in particular also regular); indicating that in case of
complete rationality diffeomorphism symmetry should be automatic.
}

\markboth{CHAPTER \arabic{chapter}. CONSTRUCTING NEW
SYMMETRIES}{CHAPTER \arabic{chapter}. CONSTRUCTING NEW SYMMETRIES}

\section{Motivations for using piecewise \mo transformations}

\paragrf
Suppose we have a \mo covariant local net of von Neumann algebras
on the circle $(\A,U)$. We do not know whether our net is
conformal. We would like to find a way to construct
--- if it is possible --- an extension of the symmetry to a
larger geometrical group; possibly to $\diff$.

In order to go beyond the \mo symmetry, we consider a (not
necessarily \mo\!\!) transformation $\kappa:S^1 \rightarrow S^1$
and we try to implement it in our model. In particular we try to
find a unitary $V(\kappa)$ satisfying
\begin{equation}
V(\kappa) \, \A(I) \, V(\kappa)^* = \A(\kappa(I))
\end{equation}
for all $I\in \I$. Of course for a $\kappa \in \mob$ we have such
a unitary. The problem is that in general we do not even seem to
have a ``candidate'' for the implementing operator.

Let us suppose for a moment that our net is actually a conformal
net with $U$ being the representation of the full diffeomorphism
group. If the diffeomorphism $\kappa \in \diff$ acts like a \mo
transformation on a certain interval $I\in \I$, that is $\exists
\varphi \in \mob: \, \kappa|_I=\varphi|_I$, then
$(\varphi^{-1}\circ\kappa)|_I={\rm id}|_I$ and so
$U(\varphi^{-1})U(\kappa)\in \A(I^c)$. Putting this into the form
\begin{equation}
\label{pcwidea1}
{\rm Ad}\left(U(\kappa)\right)(X)=
{\rm Ad}\left(U(\varphi)\right)(X)\;\;\;(\forall X\in \A(I))
\end{equation}
we see that the adjoint action of $U(\kappa)$ on elements of
$\A(I)$ is completely determined by the representation of the \mo
group.

Let $\{I\}_{n=1}^N$ be the set of disjoint open intervals that we
get by removing $N>1$ points of the circle; so $I_j \cap I_k
=\emptyset$ when $j\neq k$ and $\overline{\cup_{n=1}^N I_n} =
S^1$. Suppose $\kappa:S^1 \rightarrow S^1$ is a transformation of
the circle that acts on each interval $I_n$ like a certain \mo
transformation $\varphi_n\in \mob$ $(n=1,..,N)$. In other words
$\kappa$ is a transformation which is obtained by putting together
the action of $N$ (possibly different) \mo transformations on $N$
different pieces.

Of course if do not want $\kappa$ to be simply a \mo
transformation, we cannot require smoothness. But for this time
let us not care about such questions. What is important for us is
that by equation (\ref{pcwidea1}) when trying to implement
$\kappa$ we should look for a unitary $V(\kappa)$ whose adjoint
action $\alpha_\kappa \equiv {\rm Ad}(V(\kappa))$ satisfies
\begin{equation}
\label{pcwidea2}
\alpha_\kappa (X) = {\rm Ad}\left(U(\varphi_n)\right)(X)
\;\;\; (\forall X\in \A(I_n))
\end{equation}
for $n=1,..,N$. But if our net is at least $N$-regular, then such
an $\alpha_\kappa$ --- if it exists at all --- is completely
determined by the above equation! So if exists, then up to phase
there exists a unique unitary $V(\kappa)$ satisfying our
expectations. As it is uniquely determined (in the projective
sense) by the \mo transformations, by an abuse of notations we
shall use the symbol $U(\kappa)$ and we shall forget the letter
``$V$''.

\paragrf
It is clear that $\kappa$ must satisfy some conditions
more then just --- let us say --- continuity. Indeed,
take $\kappa$ to be the continuous transformation which
acts on the upper half circle as a dilation with a certain
parameter $s\neq 0$, and acts identically on the lower half
circle. Consider the (projective) unitary operator associated
to $\kappa$ in the sense of equation (\ref{pcwidea2}). Its
adjoint action on $\A(S^1_-)$ should be trivial. So by Haag
duality this unitary should be localized in $\A(S^1_+)$. On
the other hand, by the Bisognano-Wichmann property its adjoint
action on $\A(S^1_+)$ should be a modular automorphism with
parameter $t=-\frac{s}{2\pi}$. But unless the net is trivial,
this is impossible: we know that $\A(S^1_+)$, if not trivial,
then it is a type ${\rm I\!I\!I}_1$ factor whose modular
automorphism with nonzero parameter is never inner.

It is easy to check that in this little (counter)example although
$\kappa$ was continuous, it was not differentiable. We shall now
see that starting with a conformal net $(\A,U)$, nonsmooth
symmetries naturally appear, but the ones we shall construct will
all turn out to be (once) differentiable. So suppose that $(\A,U)$
is conformal, and take a $\gamma\in\diff$ such that
$\gamma(S^1_+)=S^1_+$. Thus ${\rm Ad}\left(U(\gamma)\right)$ is an
automorphism of the von Neumann algebra $\A(S^1_+)$ and so we can
consider its standard implementing unitary\footnote{Recall that if
$\M$ is a von Neumann algebra with a cyclic and separating vector
$\Omega$ and $\alpha \in {\rm Aut}(\M)$, then there exists a
unique unitary $V$ whose adjoint action on $\M$ is $\alpha$ and
leaves $\Omega$ in the natural cone $\P^\natural_{\M,\Omega}$; see
e.g.\! \cite[Vol.\! I\!I, Theorem 1.15]{takesaki}. $W$ automatically 
commutes with the modular conjugation $J_{\M,\Omega}$.} $W$ with respect 
to the mentioned algebra and the vacuum vector. Omitting the vacuum
vector $\Omega$ from the subscript of the modular conjugation, we
have that
\begin{equation}
\label{pcwidea3}
J_{\A(S^1_+)} W = W J_{\A(S^1_+)}.
\end{equation}
But we know that the modular conjugation $J_{\A(S^1_+)}$ acts
geometrically, more concretely like the reflection $\iota: z\mapsto
\overline{z}$. So by equation (\ref{pcwidea3}) if $I$ is a
subinterval of $\A(S^1_+)$ or $\A(S^1_-)$, then
\begin{equation}
W \A(I) W^* = \left\{
\begin{matrix} \A(\gamma(I)) &{\rm if}& I \subset S^1_+
\\ \A((\iota\circ\gamma\circ\iota)(I)) &{\rm if}& I \subset S^1_-
\end{matrix} \right.
\end{equation}
which by assuming strong additivity, implies that
the adjoint action of $W$ is geometrical and corresponds
to the transformation
\begin{equation}
z\mapsto \left\{\begin{matrix} \gamma(z) &{\rm if}& z \in S^1_+
\\ (\iota\circ\gamma\circ\iota)(z) &{\rm otherwise.}&
\end{matrix}\right.
\end{equation}
This transformation is not smooth; it is easy to think of a
$\gamma$ for which it will have no second derivative at the points
$\pm 1$. On the other hand, evidently it will be always (once)
differentiable. (If $f$ is the restriction of a smooth function
onto $\RR^+_0$ with $f(0)=0$, then by extending it to the
negatives with the formula $f(-x)= -f(x)$ we get an odd continuous
function. Thus for the $n^{\rm th}$ derivative at $x\neq 0$ we
have the formula $f^{(n)}(-x) = (-1)^{(n+1)}f^{(n)}(x)$ which
shows the existence of the first derivative even at $x=0$ and
explains why in general we do not get a twice differentiable
function.)

\paragrf
Motivated by what was explained so far, we shall introduce the
notion of piecewise \mo transformations as follows.
\begin{definition}
A $\gamma: S^1 \rightarrow S^1$ bijection is said to be a {\bf
piecewise \mo transformation} if it is (once) differentiable and
there exist a finite number of intervals $I_k\in \I$ $(k=1,..,N)$
and \mo transformations $\phi_k \in \mob$ $(k=1,..,N)$ such that
\begin{itemize}
\item
$\gamma|_{I_k} = \phi_k|_{I_k}$ for $k=1,..,N$,
\item
$\cup_{k=1}^N \overline{I_k} = S^1$.
\end{itemize}
We shall denote by {\bf Pcw$^{\mathbf 1}$(M\"ob)} the set of all
piecewise \mo transformations. As a \mo transformation is in
particular differentiable, by definition $\mob\subset\pcw$.
\end{definition}

One can prove various geometrical properties of $\pcw$. First of all,
with respect to the composition it is a group and every element of it
is an orientation preserving transformation (Prop.\! \ref{pcwgroup}).
Second, if $\{z_n\}_{n=1}^N$ and $\{z'_n\}_{n=1}^N$ are two sets of
$N$ (different) points in $S^1$ where $N$ is any positive integer and the
numbering of the points corresponds to the anticlockwise order, then there
exists a $\gamma \in \pcw$
such that
\begin{equation}
\label{ntransitivity} \gamma(z_n)=z'_n \;\;{\rm for}\; n=1,..,N.
\end{equation}
In fact, an even stronger version of this statement is true. For
two points of the circle  $p,q\in S^1$ denote by $[p,q]$ the
closed interval whose first point (in the anticlockwise direction)
is $p$ and last point is $q$. (If $p=q$ then set $[p,q]$ to be the
set containing the single point $p$.) Then if both $[z_1,z_N]$ and
$[z'_1,z'_N]$ are contained in a certain interval $I\in \I$, then
the transformation $\gamma$ in equation (\ref{ntransitivity}) can
be chosen to be localized in $I$; i.e.\! can be chosen so that
$\gamma|_{I^c}={\rm id}_{I^c}$ (Theorem \ref{localntransitivity}).

The proofs of these statements are not completely trivial.
Nevertheless in some sense they are ``elementary'' and have
nothing to do either with the theory of von Neumann algebras or
with Quantum Field Theory. Thus the author of this thesis
postpones the proofs and refers to the appendix where these
questions are treated in detail.

There is a rather trivial but nevertheless
suggestive (in relation with the construction of
diffeomorphism symmetry for a local net) consequence
of the above mentioned Theorem \ref{localntransitivity},
which is probably worth to state on its own.

\begin{corollary}
\label{pcwdensity} Let $\gamma \in \diff$. Then there exists a
sequence of piecewise \mo transformations
$\{\gamma_n\}_{n=1}^\infty \in \pcw$ converging uniformly to
$\gamma$. Moreover, if $\gamma$ is localized in the interval $I\in
\I$ then the sequence can be chosen so that every element of it is
localized in $I$.
\end{corollary}

\section{Existence of piecewise \mo
symmetry}

\paragrf
Recall that two intervals $I_1,I_2 \in \I$ are said to be {\it
distant}, if $\overline{I_1} \cap \overline{I_2} = \emptyset$.
Recall also that for each interval $I \in \I$ there is an
associated one-parameter group $s\mapsto \delta^I_s \in \mob$
which is called the one-parameter group of dilations associated to
$I$.
\begin{definition}
\label{def:kappa}
Let $I_1,I_2 \in \I$ be two distant intervals. The interior of the
complement of $I_1 \cup I_2$ is again the union of two distant
elements of $\I$ which we shall call $K_1$ and $K_2$. Then for
each $s\in \RR$ we define $\mathbf{\kappa^{I_1,I_2}_s}$ to be the
transformation given by the formula
$$
z \mapsto \left\{\begin{matrix} \delta^{I_j}_{+s}(z) &{\rm if}&
z\in I_j \;(j=1,2) \vspace{2mm} \\ \delta^{K_j}_{-s}(z) &{\rm if}&
z\in \overline{K_j} \;(j=1,2).
\end{matrix}\right.
$$
\end{definition}
In the appendix it is shown that $\kappa^{I_1,I_2}_s$ is once
differentiable (but if $s\neq 0$ it is not twice differentiable:
$(\kappa^{I_1,I_2}_s)''$ does not exists exactly at the four
points that are the endpoints of the intervals $I_1$ and $I_2$)
and so $s\mapsto \kappa^{I_1,I_2}_s$ is a one-parameter group of
piecewise \mo transformations (Lemma \ref{kappainpcw}). The
importance of such transformations is that $\pcw$ is actually
generated by the elements $\kappa^{I_1,I_2}_s$ (with $s$ ranging
over all real numbers $I_1,I_2$ all distant intervals of $\I$) and
the \mo group (Theorem \ref{pcwgenerators}).

\paragrf {\it Observation.}
Let $(\A,U)$ be a \mo covariant local net which is at least
$4$-regular. Suppose $s\in \RR$ and $I_1,I_2\in \I$ are two
distant intervals with $S^1\setminus\overline{(I_1\cup I_2)}=K_1
\cup K_2$ where $K_1,K_2\in I$. Consider an automorphism $\alpha$
of the algebra of all bounded operators of $\H_A$, whose action on
the four algebras $\{\A(X): X=I_1,I_2,K_1,K_2\}$ satisfies
\begin{equation}
\label{k-kov.action} \alpha|_{\A(X)}= \left\{\begin{matrix} {\rm
Ad}\!\left(U(\delta^{X}_{+s})\right)\!|_{\A(X)} &\vspace{2mm}{\rm when}& X=I_1,I_2 \\
{\rm Ad}\!\left(U(\delta^{X}_{-s})\right)\!|_{\A(X)} & {\rm when}&
\;X=K_1,K_2. \end{matrix}\right.
\end{equation}
Then if there exists such an automorphism, it must be unique.
Indeed, by $4$-regularity the algebras associated to the intervals
$I_1,I_2,K_1$ and $K_2$ generate the algebra of all bounded
operators of $\H_\A$ and thus the action of the automorphism by
the above equation is fixed on a generating set.
\begin{corollary} \label{uniq.pcw.symm.}
Let $(\A,U)$ be a \mo covariant local net which is at least
$4$-regular. Suppose $U$ extends to a projective unitary
representation $\tilde{U}$ of the group $\pcw$ which is compatible
with the local structure of the net $\A$; i.e.\! for all
$\gamma\in\pcw$ and $I\in\I$
\begin{equation}
\gamma|_{I^c}={\rm id}_{I^c} \;\Rightarrow\;
\tilde{U}(\gamma)\in\A(I).
\end{equation}
Then such an extension is unique (and so we may omit the symbol
$\tilde{\hphantom{U}}$ and by an abuse of notation denote the
extension, too, by $U$).
\end{corollary}
\begin{proof}
It is a consequence of the observation made above and the fact
that the elements in question, together with the \mo
transformations generate the whole group of piecewise \mo
transformations (Theorem \ref{pcwgenerators}).
\end{proof}

\paragrf
We shall now see that if the piecewise \mo transformations
$\{\kappa^{I_1,I_2}_s\}$ are implementable, then --- given that
the net is {\it regular} (i.e.\! $n$-regular for any positive
integer $n$)--- such an extension of $U$ always exists: i.e.\! the
group relations are automatically satisfied.\footnote{By using
only a set of defining group relations rather than checking the
``full'' multiplication table of $\pcw$ it is possible that one
can weaken the regularity condition. In fact the author spent some
time trying to give a proof using only $6$-regularity rather than
complete regularity, but then decided to not do so as it would
make everything much more complicated without resulting a
statement of really more use. In any case, regularity is expected
to hold in every diffeomorphism covariant net and to the knowledge
of the author there are no known examples or general situations
where the existence of diffeomorphism covariance is unknown and
regularity cannot be easily checked while $6$-regularity can be.}
Moreover, it turns out that such an extension acts always in a
covariant way on the net, which is not at all a trivial statement.
Take for example, the adjoint action of $U(\kappa^{I_1,I_2}_s)$.
If $K\in I$ is a subinterval of $I_1$, then by definition
\begin{equation}
\label{k-geo.action}
U(\kappa^{I_1,I_2}_s)\A(K)U(\kappa^{I_1,I_2}_s)^* =
\A(\kappa^{I_1,I_2}_s(K))
\end{equation}
as $\kappa^{I_1,I_2}_s(K) = \delta^{I_1}_s(K)$ and the restriction
of the adjoint action of the unitary in question to
$\A(K)\subset\A(I)$ by Eq.\! (\ref{k-kov.action}) coincides with
that of $U(\delta^{I_1}_s)$. But if $K$ does not lie in any of the
``pieces'' because for example it contains one of the endpoints of
$I_1$, then Eq.\! (\ref{k-geo.action}) --- as we did not assume
strong additivity --- does not immediately follow from the
definitions.
\begin{proposition} \label{prop.pcw}
Let $(\A,U)$ be a \mo covariant local net which is regular and
admits an implementation of the transformation
$\delta^{I_1,I_2}_s$ in the sense of Eq. (\ref{k-kov.action}) for
every $s\in\RR$ and $I_1,I_2$ distant elements of $\I$. Then $U$,
as a projective unitary representation, extends to $\pcw$ in a
unique manner which is compatible with the local structure of
$\A$; i.e.\!
$$
\gamma|_{I^c}={\rm id}_{I^c} \;\;\Rightarrow\;\;
U(\gamma)\in\A(I)\;\;\;\; (\gamma\in\pcw, I\in\I).
$$
Moreover, this extension acts in a covariant way on $\A$:
$$
U(\gamma) \A(I) U(\gamma)^* = \A(\gamma(I)) \;\;\;\;
(\gamma\in\pcw, I\in\I).
$$
\end{proposition}
\begin{proof}
By Corollary \ref{uniq.pcw.symm.}, the extension, if exists, is
unique. By the condition of the proposition we are given with the
unitaries $\{U(\delta^{I_1,I_2}_s)\}$ and for each such unitary we
know its action on certain local algebras.

Since by Theorem \ref{pcwgenerators} the transformations
$\{\delta^{I_1,I_2}_s\}$ generate the group of piecewise \mo
transformations, if $\gamma \in \diff$ then there exists a
decomposition $\gamma=\gamma_1\circ\gamma_2\circ .. \gamma_n$
where each transformation $\gamma_j$ $(j=1,2..,n)$ is either a \mo
transformation or of the mentioned kind. Thus we can try to define
the unitary associated to $\gamma$ by the formula
$U(\gamma)=U(\gamma_1)U(\gamma_2)..U(\gamma_n)$. However, we need
to assure that $U(\gamma)$ is well-defined (in the projective
sense) by this procedure and that the thus extended $U$ is indeed
a (projective) representation.

For all this what we need to check is that
$U(\gamma_1)U(\gamma_2)..U(\gamma_n)$ is a multiple of the
identity whenever $\gamma_1\circ\gamma_2\circ..\gamma_n=id_{S^1}$.
But it is rather clear, that we can find a set of disjoint
intervals $I_r\in\I$ $(r=1,2..,k)$ whose union equals to $S^1$
minus a finite set of points and satisfy for every
$r\in\{1,2,..,k\}$ the following requirement: $\gamma_n$ is smooth
on $I_r$, $\gamma_{n-1}$ is smooth on $\gamma_n(I_r)$,
$\gamma_{n-2}$ is smooth on $\gamma_{n-1}(\gamma_n(I_r))$, etc.
Then it follows that the restriction of the adjoint action of
$U(\gamma_1)U(\gamma_2)..U(\gamma_n)$ to $\A(I_r)$ is trivial for
every $r\in\{1,2,..,k\}$. Thus by the regularity of $\A$ the
unitary $U(\gamma_1)U(\gamma_2)..U(\gamma_n)$ must be a multiple
of the identity.

We shall now show the compatibility with the net. That is, the
property that if $\gamma\in\pcw$ is such that it acts identically
on the complement of a certain interval $I\in\I$ then the
restriction of the adjoint action of $U(\gamma)$ is trivial on
$\A(I^c)$.

By the consideration made before it is clear that we can find a
finite set of points such that ${\rm Ad}(U(\gamma))$ is trivial on
the algebras associated to the subintervals of $I^c$ obtained by
the removal of these points. This in itself though is not enough
as we did not assume strong additivity. Suppose $K\subset I^c$ is
one of those nonempty open subintervals of $I^c$ for which ${\rm
Ad}(U(\gamma))$ acts trivially on $\A(K)$. We may choose $K$ so
that it has a common endpoint with $I$. Then we can find an
interval $F\in\I, F\subset I^c$ which is distant from $K$ such
that also $F$ has a common endpoint with $I$ (of course on the
``other side'' with respect to where $K$ is). Then for every
$s\in\RR$ the transformation $\delta^I_{s}\kappa^{K,F}_s$ is
identical on $I$ and thus commutes with $\gamma$, showing that in
the projective sense
\begin{equation} \label{eq.pcwimpl.1}
U(\delta^I_{s}\kappa^{K,F}_s) U(\gamma)
U(\delta^I_{s}\kappa^{K,F}_s)^* = U(\gamma).
\end{equation}
On the other hand $U(\delta^I_{s}\kappa^{K,F}_s) =
U(\delta^I_{s})U(\kappa^{K,F}_s)$ and the adjoint action of
$U(\kappa^{K,F}_s)$ on $\A(K)$ coincides with that of
$U(\delta^K_s)$. Then, since by \mo covariance we have
$U(\delta^K_s)\A(K)U(\delta^K_s)^*=\A(K)$, it follows that
\begin{equation}
{\rm Ad}\left(U(\delta^I_{s}\kappa^{K,F}_s)\right)(\A(K))= {\rm
Ad}\left(U(\delta^I_{s})\right)(\A(K)) = \A(\delta^I_{s}(K))
\end{equation}
which, together with Eq.\! (\ref{eq.pcwimpl.1}), shows that the
action of ${\rm Ad}(U(\gamma))$ is trivial on
$\A(\delta^I_{s}(K))$. But the parameter $s\in\RR$ was arbitrary
and by the continuity of the net
$\cup_{s\in\RR}\A(\delta^I_{s}(K))$ generates $\A(I^c)$. Thus in
fact we have verified that ${\rm Ad}(U(\gamma))$ acts trivially on
$\A(I^c)$ showing that the extension of $U$ to the piecewise \mo
group is compatible with the local structure of $\A$.

We are left only with the part regarding the covariance of the
action of the extension. Let $\gamma\in\pcw$ and $I\in\I$. We can
safely assume that $\gamma(I)=I$; even if it was originally not so
we can modify $\gamma$ by a \mo transformation so that it holds
(and of course a \mo transformation acts in a covariant way). For
the moment we shall further assume that $\gamma$ acts identically on a
certain $K\in\I$ where $K\subset I^c$. We would like to show that
$U(\gamma)\A(I)U(\gamma)^*=\A(I)$. By Haag-duality and the
compatibility which we have just proved it is clear that
$U(\gamma)\A(I)U(\gamma)^*\subset\A(K^c)$. Then by the same trick
which we have just employed for the demonstration of the
compatibility, for any $s\in\RR$ we can find another piecewise \mo
transformation $\tilde{\gamma}$ such that $\tilde{\gamma}|_I =
\gamma|_I$ and $\tilde{\gamma}$ acts identically on
$\A(\delta^{I}_s(K))$. Then it follows that
\begin{equation}
U(\gamma)\A(I) U(\gamma)^* =
U(\tilde{\gamma})\A(I)\U(\tilde{\gamma})^*\subset
\A(\delta^I_s(K^c))
\end{equation}
for every parameter $s\in\RR$. Thus by the additivity of the net
$\A$ we can conclude that $U(\gamma)\A(I)\U(\gamma)^*=\A(I)$.

Finally, we need to remove the condition on the existence of a
nonempty open interval $K\subset I^c$ on which $\gamma$ is
identical. Of course we can choose an interval $K\in\I, K\subset
I^c$ such that $\gamma|_K$ is smooth and $K$ has a common endpoint
with $I$. Then there exists an $s_1\in\RR$ such that
$\gamma(K)=\delta^I_{s_1}(K)$ and therefore the restriction of
$\delta^I_{-s_1}\circ \gamma$ to $K$ is a dilation associated to
$K$ with a certain parameter $s_2\in\RR$. We may choose another
open nonempty interval $F\subset I^c$ which is distant from $K$
and which also has a common endpoint with $I$. Then by what was
said
\begin{eqnarray} \nonumber
\kappa^{K,F}_{-s_2}\circ\delta^I_{-s_1}\circ\gamma|_K&=&{\rm id}_K
\;\;\;{\rm and} \vspace{2mm} \\
\kappa^{K,F}_{-s_2}(\delta^I_{-s_1}(\gamma(I)))&=&I.
\end{eqnarray}
Thus the adjoint action of $U(\kappa^{K,F}_{-s_2})
U(\delta^I_{-s_1}) U(\gamma)$ preserves the algebra $\A(I)$. But
as $\kappa^{K,F}_{-s_2}|_I =\delta^I_{s_2}|_I$ by the geometrical
position of $K,F\subset I^c$ and by the definition of
$\kappa^{K,F}_{-s_2}$, we have that also the adjoint action of
$U(\kappa^{K,F}_{-s_2})$ (and of course that of
$U(\delta^I_{s_2})$, too) preserve the algebra in question showing
that $U(\gamma)\A(I)U(\gamma)=\A(I)$. This concludes the proof of
the covariance of the action.
\end{proof}
\paragrf
When considering diffeomorphism symmetry, apart from compatibility
with local structure we also require continuity of the
representation. Having seen the that the one-parameter groups of
the type $s\mapsto \kappa^{I_1,I_2}_s$ generate the extension to
the piecewise \mo group, we shall fix the notion of continuity
exactly by these. Afterwards we shall also fix what we mean by
piecewise \mo symmetry.
\begin{definition}
A projective unitary representation $U$ of $\pcw$ is called
continuous if its restriction to \mo group as well as its
restriction to any of the one-parameter groups $s\mapsto
\kappa^{I_1,I_2}_s$, where $s\in\RR$ and $I_1,I_2$ are distant
elements of $\I$, is strongly continuous.
\end{definition}
\begin{definition}
A \mo covariant local net on the circle $(\A,U)$ is called {\bf
piecewise \mo covariant} if $U$ extends to a continuous projective
unitary representation of the group $\pcw$ which is compatible
with the local structure of the net and acts in a covariant manner
on it; i.e.\! for all $\gamma\in\pcw$ and $I\in\I$
\begin{itemize}
\item[$1$.]\, $\gamma|_{I^c}={\rm id}_{I^c}
\;\;\Rightarrow\;\;U(\gamma)\in\A(I)$ \item[$2$.]\,
$U(\gamma)\A(I) U(\gamma)^* = \A(\gamma(I))$.
\end{itemize}
\end{definition}

\paragrf
Assuming regularity, we shall now establish the sufficient and
necessary algebraic condition (on the local algebras) of the
existence of piecewise \mo symmetry. As in general a \mo covariant
local net can be given by a +hsm factorization involving three von
Neumann factors (see in the preliminary Sect.\!
\ref{sec:modular}), one wishes to find a condition that concerns
the ``relative position'' of three local algebras. In the
particular case of strong additivity, the condition should be
simplified and expressed in a way that involves only two local
algebras, in agreement with the fact that a strongly additive \mo
covariant local net can be characterized as a +hsm inclusion of
two von Neumann factors (see again in the preliminary Sect.\!
\ref{sec:modular}).

\begin{theorem}
\label{existence.pcw.symm.} Let $(\A,U)$ be a \mo covariant local
net and $I_1,I_2 \subset I$ the two connected subintervals
obtained by removing a point from the interval $I\in \I$. Then
\begin{itemize}
\item if $\A$ is regular then it is piecewise \mo covariant if and
only if there exists a continuous one-parameter group $t\mapsto
\alpha_t \in {\rm Aut}(\A(I))$ preserving both $\A(I_1)$ and
$\A(I_2)$ such that the restriction of $\alpha_t$ to $\A(I_k)$
$(k=1,2)$ is a modular automorphism with parameter $(-1)^k t$;

\item if $\A$ is strongly additive then it is piecewise \mo
covariant if and only if there exists an unbounded 
expectation\footnote{Note that in the mathematical literature
it is usually called a ``normal, semifinte, operator valued weight''.
In this thesis the name ``unbounded expectation'' is used, 
to stress that a (usual) expectation is a special case of such
an unbounded one.} from $\A(I)$ to $\A(I_1)$.
\end{itemize}
In any of these cases with the given algebraic condition
satisfied, the extension of the symmetry from the \mo group to
$\pcw$ is unique.
\end{theorem}
\begin{proof}
Uniqueness is ensured by Corollary \ref{uniq.pcw.symm.}. First we
shall treat the more general case of a regular net. Let $J$ be the
modular conjugation associated to the algebra $\A(I)$ and the
vacuum vector and set $I^*_1$ and $I^*_2$ for the intervals
obtained by applying the geometrical part of the conjugation to
$I_1$ and $I_2$; i.e.\! $J\A(I_k)J=\A(I^*_k)$ for $k=1,2$.

Suppose the net is piecewise \mo covariant. Then the one-parameter
group $t\mapsto {\rm Ad}(U(\kappa^{I^*_1,I_2}_{2\pi t}))|_{\A(I)}$
trivially satisfies the required properties. In the converse
direction, suppose $t\mapsto \alpha_t \in {\rm Aut}(\A(I))$ is the
continuous one-parameter group appearing in the condition of the
theorem. By the Bisognano-Wichmann property and the cocycle
theorem (see e.g.\! \cite[Vol.\!  I\!I, Theorem 3.19]{takesaki}) 
of Connes it is clear that we may assume 
$t\mapsto\alpha_{(-1)^k t}|_{\A(I_k)}$ to be not just any
modular automorphism group of $\A(I_k)$ $(k=1,2)$ but in fact the
one coming from the vacuum-state. Then setting
$U(\kappa^{I^*_1,I_2}_s)$ $(s\in\RR)$ for the standard
implementing unitary (with respect to the vacuum) of the
automorphism $\alpha_{\frac{s}{2\pi}}$, we have that $s\mapsto
U(\kappa^{I^*_1,I_2}_s)$ is a strongly continuous one-parameter
group of unitaries whose adjoint action on $\A(I_k)$, by
definition, coincides with that of $U(\delta^{I_k}_{(-1)^k s})$.
Moreover, as $JU(\kappa^{I^*_1,I_2}_s)J =U(\kappa^{I^*_1,I_2}_s)$
we also have the ``right'' action on the algebras $\A(I^*_1)$ and
$\A(I^*_2)$. Nevertheless we still cannot apply Prop.\!
\ref{prop.pcw}. If $K\in\I$ and $K_1,K_2\subset K$ are the two
intervals obtained by removing a point from $K$, then --- if the
numbering of the intervals is made in the right way --- there
exists a \mo transformation sending $(I_1,I_2\subset I)$ into
$(K_1,K_2\subset K)$ and then of course we can repeat everything
with $I$ replaced by $K$. The problem is that we cannot get any
pair of distant intervals as $(K^*_1,K_2)$ with some choice of $K$
and the removed inner point. The \mo group does not act
transitively on the set of pairs of distant intervals.

However, we have already managed to extend the symmetry to {\it
some} piecewise \mo transformations, as we already have the
operator $U(\kappa^{I^*_1,I_2}_s)$. It is not claimed that these
transformations generate the group $\pcw$. Nevertheless, they are
already many enough to move any pair of distant intervals into any
other pair.

Suppose $F_1,F_2,F_3,F_4\in\I$ are four disjoint intervals
numbered in an anticlockwise manner whose union equals to $S^1$
minus four points. Again, by the Bisognano Wichmann property and
Connes' cocycle theorem it is clear that whether the one-parameter
group of transformations $s\mapsto\kappa^{F_1,F_3}_s$ is
implementable in the sense of Eq.\! (\ref{k-kov.action}),  in a
strongly continuous way, is an algebraic property; i.e. a property
which depends only on the relative position of the four algebras
$(\A(F_1),\A(F_2),\A(F_3)$ and $\A(F_4))$ but does not depend on
the vacuum state. We know that this property holds for
$(\A(I_1),\A(I_2),\A(I^*_1),\A(I^*_2) )$. All we need to show is
that by ``moving'' this quartet of algebras repeatedly with the
adjoint action of $U(\gamma)$ with $\gamma$ being either a \mo
transformation or $\kappa^{I^*_1,I_2}_s$ for some $s\in\RR$, we
can obtain the quartet associated to $(F_1,F_2,F_3,F_4)$.

By \mo covariance we can safely assume that $F_1=I$ and
$F_2=I^*_2$. Then $F_3,F_4\subset I^*_1$ and in fact their union
is $I^*_1$ minus a point. Moreover, there exists a \mo
transformation $\phi\in\mob$ such that $\phi(I_1)=I=F_1$ and
$\phi(I_2)=I^*_2=F_2$. Then it follows that
$\phi(I^*_1),\phi(I^*_2\subset I^*_1$ and that also in this case
their union is $I^*_1$ minus a point. So it is easy to see that
there must exist a certain parameter $s\in \RR$ such that
$\delta^{I^*_1}_s(\phi(I^*_1))=F_3$ and
$\delta^{I^*_1}_s(\phi(I^*_1))=F_3$. Then the action of ${\rm Ad}
\left(U(\kappa^{I^*_1,I_2}_s)U(\phi)\right)$ sends
\begin{eqnarray} \nonumber
\A(I_1)&\mapsto&
\alpha_{\frac{s}{2\pi}}(\A(I))=\A(I)=\A(F_1), \\
\nonumber \A(I_1)&\mapsto& U(\kappa^{I^*_1,I_2}_s) \A(I^*_2)
U(\kappa^{I^*_1,I_2}_s)^*
=\A(I^*_2)=\A(F_2)  \\
\nonumber \A(I^*_1)&\mapsto& U(\delta^{I^*_1}_s) \A(\phi(I^*_1))
U(\delta^{I^*_1}_s)^*
=\A(\delta^{I^*_1}_s(\phi(I^*_1))) = \A(F_3)  \\
\A(I^*_2)&\mapsto& U(\delta^{I^*_1}_s) \A(\phi(I^*_2))
U(\delta^{I^*_1}_s)^* = \A(\delta^{I^*_1}_s(\phi(I^*_2)))= \A(F_4)
\end{eqnarray}
where we have used that $U(\kappa^{I^*_1,I_2}_s)$, by its
construction, preserves both $\A(I)$ and $\A(I^*_2)$, and its adjoint
action on $\A(I^*_1)$ coincides with that of
$U(\delta^{I^*_1}_s)$. This shows that $s\mapsto\kappa^{F_1,F_3}$
is implementable in the sense of Eq.\! (\ref{k-kov.action}). Thus
by Prop.\! \ref{prop.pcw}, $U$ extends to a projective
representation of $\pcw$ which is compatible with the net and
whose action on the net is covariant. Thus, since the
one-parameter groups in question are continuous, the net is
piecewise \mo covariant.

We are left with the part of the theorem regarding the strongly
additive case. The existence of an
unbounded expectation is equivalent with the fact that there is a
modular automorphism group of the bigger algebra (in general
associated to a weight and not necessarily to a state) that leaves
the smaller algebra globally invariant and its restriction to the
smaller algebra is again modular, see e.g.\! \cite[Vol.\! I\!I, Theorem 
4.12]{takesaki}. 
Again, by Connes' cocycle theorem, it then follows that any modular 
automorphism group of the smaller algebra extends to a continuous 
one-parameter automorphism group of the bigger.

If $\A$ is strongly additive and there exists an unbounded
expectation from $\A(I)$ to $\A(I_1)$ then the modular group of
$\A(I_1)$ associated to the vacuum state extends to a continuous
one-parameter automorphism group of $\A(I)$. Using the standard
the standard implementing unitary with respect to $\A(I)$ and the
vacuum vector, it is easy to see that what we get is a strongly
continuous one-parameter group $t\mapsto W_t$ whose adjoint action
on $\A(I)$ coincides with that of
$U(\delta^{I_1}_{\frac{t}{2\pi}})$, while on $\A(I^*_1)=J\A(I_1)J$
with that of $U(\delta^{I^*_1}_{\frac{-t}{2\pi}})$. Set 
\begin{equation}
K\equiv \{{\rm the\;interior\;of\;} \overline{I_1 \cup I^*_1}\} \in \I. 
\end{equation}
By strong additivity $\A(I_1)$ and $\A(I^*_1)$ generates $\A(K)$ and
hence the adjoint action of $W_t$ preserves $\A(K)$. Then by what
was said we can apply the first part of the statement (which we
have already proved) to $(I_1,I^*_1\subset K)$ and thus we can
conclude that our net is piecewise \mo covariant.

As for the converse direction, suppose our net is strongly
additive and piecewise \mo covariant. Then the restriction of
\begin{equation}
t\mapsto {\rm Ad}(U(\kappa^{I_1,\tilde{I_2}}_{\frac{t}{2\pi}}))
\end{equation}
to $\A(I_2)$ is an inverse-parameterized modular automorphism
group. As $\A(I_2)$ is a factor in its standard form it follows
that the restriction to $\A(I^c_2)=\A(I_2)'$ is a modular
automorphism group. But of course, by its very definition, also
the restriction to $\A(I_1)$ is a modular automorphism group.
Therefore there exists an unbounded expectation from $\A(I^c_2)$
to $\A(I_1)$ and then of course, by \mo symmetry also from $\A(I)$
to $\A(I_1)$.
\end{proof}

\section{Piecewise \mo symmetry: conclusions}

\paragrf
Using piecewise \mo transformations one can move any $n$
consecutive points into $n$ other consecutive points. So, for
example one can prove that the $\mu$-index of a piecewise \mo
covariant net is well-defined. Using only \mo transformations we
cannot move any pair of distant intervals into any other pair of
distant intervals. (Even without diffeomorphism covariance, by
\cite[Prop.\! 5]{KLM} the $\mu$-index is independent of the choice
of intervals but only with the assumption of split property and
strong additivity.) However, the author does not think that such a
use of piecewise \mo transformations (i.e.\! to show something
which follows obviously from diffeomorphism covariance) is of real
importance. Piecewise \mo symmetry is introduced to derive
diffeomorphism symmetry. Assuming that sooner or later this will
be indeed carried out, there is of course no point in using
piecewise \mo symmetry for something that diffeomorphism symmetry
can be easily used for.

The piecewise \mo symmetry was easily shown to be unique. By
establishing a connection between piecewise \mo symmetry and
diffeomorphism symmetry (at least in the case when both are
assumed to exist), this can serve to show the uniqueness of the
latter one. Also, it is almost obvious that a strongly additive
\mo covariant subnet of a piecewise \mo covariant net is
automatically piecewise \mo covariant. Indeed, a piecewise \mo
transformation is ``put together'' by \mo transformations, and the
action of a \mo transformation --- by definition --- preserves the
subnet. So again, in case both piecewise \mo symmetry and
diffeomorphism symmetry exist, this observation can be used to
show the conformal covariance of a subnet.

However, for the moment we shall suspend discussing how these new,
nonsmooth symmetries can be exploited. The mentioned ideas will be
developed in a later chapter. In this section the only question we
shall investigate is: {\it when} can we establish that a \mo
covariant local net admits piecewise \mo symmetry. Of course in
some sense we have a complete characterization (Theorem
\ref{existence.pcw.symm.}) of the existence of such symmetry. But
it is worth to investigate the relation of the condition there
given with such ``more traditional'' assumptions as for example
complete rationality. Does complete rationality imply piecewise
\mo covariance? Another interesting question is the following.
Consider a Virasoro net with central charge $c>1$. As it was
mentioned, it is not strongly additive. Thus by taking its dual
(see Sect.\! \ref{sec:more.ex}), we obtain a new (\mo
covariant) net. Unfortunately, at the moment there is no direct
description of this new net and there seems to be no easy way of
establishing the existence of diffeomorphism symmetry. So can we
prove at least that it admits piecewise \mo symmetry? We shall
answer both questions positively.

\paragrf
Let us first see what can we say if we assume conformal
covariance.
\begin{lemma} \label{diffcov->pcwcov}
Let $(\A,U)$ be a conformal local net and $I_1,I_2$ two intervals
obtained by removing a point from $I\in\I$. Then there exists a
continuous one-parameter group $t\mapsto \alpha_t \in {\rm
Aut}(\A(I))$ preserving both $\A(I_1)$ and $\A(I_2)$ such that the
restriction of $\alpha_t$ to $\A(I_k)$ $(k=1,2)$ is a modular
automorphism with parameter $(-1)^k t$.
\end{lemma}
\begin{proof}
We can assume that $I=S^1_+$ and that the point removed is its
``middle'' point. Then the one-parameter group of $2$-dilations
preserves both $I_1$, $I_2$ and $I$. Thus the corresponding
one-parameter unitary group preserves the three local algebras
$\A(I_1), \A(I_2)$ and $\A(I)$. Moreover, the actions on $\A(I_1)$
and $\A(I_2)$, with right parametrizations, are modular, since the
geometrical part of the action on each piece is conjugate to a
dilation; see also \cite[Prop.\! A.1]{LoXu}. Thus the automorphism
group implemented by this one-parameter group satisfies all that
was required.
\end{proof}
\begin{corollary}
Let $(\A,U)$ be a conformal local net. Then its dual is piecewise
\mo covariant.
\end{corollary}
\begin{proof}
Let $I_1,I_2$ be two intervals obtained by removing a point
$z_\infty\in \I$ from the interval $I\in\I$. By choosing the
infinite point to be $z_\infty$, we may assume that $\A^d(I_1) =
\A(I_1)$ and $\A^d(I_1) = \A(I_1)$, where $\A^d$ is dual net. By
the previous lemma there exists a continuous one-parameter group
$t\mapsto \alpha_t \in {\rm Aut}(\A(I))$ preserving both $\A(I_1)$
and $\A(I_2)$ such that the restriction of $\alpha_t$ to $\A(I_k)
= \A^d(I_k)$ $(k=1,2)$ is a modular automorphism with parameter
$(-1)^k t$. But by the strong additivity of the dual net $\A^d(I)
= \A^d(I_1)\vee \A^d(I_2)$, hence $\alpha_t$ preserves $\A^d(I)$,
too. Thus by Theorem \ref{existence.pcw.symm.} the statement is
proven.
\end{proof}
Consider the three requirements of complete rationality. There are
no known examples of diffeomorphism covariant nets not satisfying
split property. Strong additivity does not always hold, but the
weaker assumption of regularity again seems to be a general
feature of diffeomorphism covariant nets. Finally, once assuming
regularity, the two-interval inclusion is irreducible. If this
inclusion is expected (i.e.\! there exists an expectation from the
bigger algebra to the smaller) then by covariance so is the
inclusion obtained by going over the commutants. Thus by a theorem
of Longo (see \cite[I\!I, Theorem 4.4]{longo1}) it follows that the
two-interval inclusion is of finite index. The contrary, of course,
is true essentially by definition: if the two-interval inclusion is of
finite index then it is expected. However, as we know not all
diffeomorphism covariant nets have a finite $\mu$-index. So
instead let us consider the weaker assumption that the
two-interval inclusion is expected {\it in the unbounded sense}.
By a similar argument that was employed in the first lemma of this
section, using $2$-dilations it is easy to show that this
condition is always satisfied in case of diffeomorphism
covariance. (One has to use what was already mentioned in the last
section, namely, that the existence of an unbounded expectation is 
equivalent with the fact that there is a modular automorphism group of the 
bigger algebra --- in general associated to a weight and not necessarily 
to a state --- that leaves the smaller algebra globally invariant and its 
restriction to the smaller algebra is again modular.)
\begin{theorem}\label{3cond.->pcw}
Let $(\A,V)$ be a \mo covariant local net. Suppose
\begin{itemize}
\item $\A$ is split, \item $\A$ is regular, \item the two-interval
inclusion in $\A$ is expected at least in the unbounded sense.
\end{itemize}
Then it follows that $(\A,V)$ is piecewise \mo covariant.
\end{theorem}
\begin{proof}
Let $I_1,I_2,I_3,I_4\in\I$ be four disjoint intervals, numbered in
the anticlockwise direction, whose union is $S^1$ minus the four
endpoints. By the isomorphism $\A(I_1)\vee \A(I_3) \simeq
\A(I_1)\otimes\A(I_3)$ we can consider the tensor product
$\sigma^1 \otimes \sigma^3$ of the modular automorphism groups of
$\A(I_1)$ and $\A(I_3)$ with the state being the vacuum state.
Similarly we have $\sigma^2 \otimes \sigma^4$. Since $\sigma^1
\otimes \sigma^3$ is a modular automorphism group of $\A(I_1)\vee
\A(I_3)$, by the third listed property it extends to a modular
automorphism group $\tilde{\sigma}^{1,3}$ of $(\A(I_2)\vee
\A(I_4))'$. This algebra is in its standard form, and therefore
$\tilde{\sigma^{13}}$ can be implemented by a strongly continuous
one-parameter unitary group $\tilde{U}^{1,3}$. Then the adjoint
action of $\tilde{U}^{1,3}$ preserves the algebra $\A(I_2)\vee
\A(I_4)$ and in fact by the inverse of the cocycle theorem (see
e.g.\! \cite[Vol.\! I\!I, Theorem 3.21]{takesaki}), it is easy to 
see that its action on
$\A(I_2)\vee \A(I_4)$ is modular with inverted parametrization.
Thus, using the Connes' cocycle theorem there exists a strongly
continuous unitary cocycle $W \in \A(I_2)\vee \A(I_4)$ such that
the adjoint action of $\tilde{U}^{1,3}(t)W(t)$ on $\A(I_2)\vee
\A(I_4)$ is exactly $\sigma^2_{-t} \otimes \sigma^4_{-t}$. So the
following is true about the adjoint action $\alpha_t \equiv {\rm
Ad}\tilde{U}^{1,3}(t)W(t)$:
\begin{equation}
\alpha_t|_{\A(I_j)} = \left\{\begin{matrix} \sigma^j_t &{\rm if}&
j=1,3 \\ \sigma^j_{-t} &{\rm if}& j=2,4. \end{matrix}\right.
\end{equation}
Then by Prop.\! \ref{prop.pcw} $\alpha_t$ preserves the algebra
$\A(I)$ where $I$ is the interior of $\overline{I_1\cup I_2}$,
which in turn, by Theorem \ref{existence.pcw.symm.} implies that
the net is piecewise \mo covariant.
\end{proof}

\chapter{Stress-energy tensor on nonsmooth functions}

\resume{
In this chapter it will be shown that the stress-energy tensor,
associated to a positive energy representation of $\diff$, can be
naturally extended from the smooth functions to functions of
finite $\|\cdot\|_{\frac{3}{2}}$ norm. To this extension Nelson's
commutator theorem cannot be directly used to show that the
operators involved are self-adjoint. However, an estimate will be
found involving the contraction semigroup associated to the
conformal Hamiltonian, which allows to demonstrate
self-adjointness.

The possibility of evaluating the stress-energy tensor on
nonsmooth functions has various use. In particular, one can
directly construct the nonsmooth symmetries discussed in the
previous chapter by evaluating the stress-energy tensor on
nonsmooth functions. In the next chapter this will serve both for
showing the uniqueness of the diffeomorphism symmetry and for the
automatic conformal covariance of subnets.}

\markboth{CHAPTER \arabic{chapter}. TREATING NON-SMOOTH FUNCTIONS}
{CHAPTER \arabic{chapter}. TREATING NON-SMOOTH FUNCTIONS}

\section{Why dealing with the nonsmooth case and how}
\label{sec:nonsmooth1}

\paragrf
In the previous chapter it was shown by abstract arguments that in
a large class of models, the \mo symmetry of the net can be extended to 
the piecewise \mo group. One may wonder whether in case of conformal
covariance the construction of nonsmooth symmetries could be
directly done by evaluating the stress-energy tensor on some
nonsmooth functions. Since the estimates in the paper of Goodmann
and Wallach \cite{GoWa} are not sufficient for our purpose, we
need a more detailed analysis. We shall use ideas coming from the
paper \cite{Ne} of E. Nelson. Essentially all presented results of 
this chapter are extracted from the joint work \cite{CaWe} of the
author of this thesis with S.\! Carpi.

Evaluating the stress-energy tensor on nonsmooth functions can
give advantages in various calculations and provides new tools for
proofs. In particular, since the piecewise \mo symmetry of the net
is unique, the direct construction of this symmetry from the
stress-energy tensor gives a way to demonstrate the uniqueness of
the diffeomorphism symmetry. It also turns out to be the key for
showing that a locally normal representation of a conformal net is
always of positive energy. All these results will be derived in
the next chapter, in this chapter we shall consider only the
problem of evaluating the stress-energy tensor on nonsmooth
functions.

\paragrf
Recall that for a real continuous function $f$ on the circle with Fourier
coefficients
\begin{equation}
\hat{f}_n \equiv \frac{1}{2\pi}\int_0^{2\pi} f(e^{i\theta})
e^{-i\theta}d\theta \;\; (n\in\ZZ)
\end{equation}
we have set
\begin{equation}
\|f\|_{\frac{3}{2}} \equiv \sum_{n\in\ZZ} |\hat{f}_n|
(1+|n|^{\frac{3}{2}}) \, \in \RR^+_0 \cup \{\infty\}.
\end{equation}
Suppose that $\{L_n: n\in\ZZ\}$ is the representation of the
Virasoro operators associated to the positive energy
representation $U$ of $\diff$. If $f$ is a smooth function, then
--- as we have seen in the preliminaries (Sect.\!
\ref{sec:diffposenergy}) --- the stress-energy tensor $T$ of $U$
evaluated on $f$ can be defined as the closure of the operator
$\sum_{n\in\ZZ}\hat{f}_n L_n$. In the next section we shall show
that if $\|f\|_{\frac{3}{2}}<\infty$ then even if $f$ is
nonsmooth, the above sum is meaningful (i.e.\! convergent) and
results an essentially self-adjoint operator. To do so we shall
find (by a rather straightforward calculation) an
$\epsilon$-independent bound on the norm of the commutator
$[\mathop{\sum}_{n \in {\mathbb Z}}\hat{f}_n L_n, e^{-\epsilon
L_0}]$ where $\epsilon
>0$.

\paragrf
Let us return for a second to piecewise \mo transformations and
introduce the notion of piecewise \mo vector fields.
\begin{definition}\label{def:pcwm}
A function $f:S^1\rightarrow \RR$ is said to be a {\bf piecewise
\mo vector field} if it is (once) differentiable and there exist a
finite number of intervals $I_k\in \I$ $(k=1,..,N)$ and \mo vector
fields $g_k \in {\mathfrak m}$ $(k=1,..,N)$ such that
\begin{itemize}
\item $f|_{I_k} = g_k|_{I_k}$ for $k=1,..,N$, \item $\cup_{k=1}^N
\overline{I_k} = S^1$.
\end{itemize}
We shall denote by $\mathbf{{\rm \bf Pcw}^1(\mathfrak m)}$ the set
of all piecewise \mo vector fields. As a \mo vector field is in
particular differentiable, by definition ${\mathfrak m}
\subset\pcwm$.
\end{definition}
The following lemma is useful for checking the finiteness of the
relevant norm.
\begin{lemma}
\label{finite1.5norm} Let $f$ be a (once) differentiable function
on the circle. Suppose that there exists a finite set of intervals
$I_k\in\I$ and smooth functions $g_k$ on the circle $(k=1,..,N)$
such that $\overline{\cup_{k=1}^N I_k}=S^1$ and
$f|_{I_k}=g_k|_{I_k}$. Then $\|f\|_{\frac{3}{2}} < \infty$.
\end{lemma}
\begin{proof}
The conditions imply that $f''$, which is everywhere defined apart
from a finite set of points, has Fourier coefficients
$\hat{(f'')}_n=-n^2 \hat{f}_n$ and is of bounded variation.
Therefore $|n^2\hat{f}_n| \leq |\frac{{\rm Var}(f'')}{n}|$ (see
\cite[Sect.\! I.4]{katznelson}), from which the claim follows
easily.
\end{proof}
\begin{corollary}
\label{pcwmfin1.5norm} Every piecewise \mo vector field has finite
$\|\cdot\|_{\frac{3}{2}}$ norm.
\end{corollary}
It is also important to investigate the density of the smooth
functions among the functions of finite $\|\cdot\|_{\frac{3}{2}}$
norm.
\begin{lemma}
\label{smoothdensity} Let $f\in C(S^1,\RR)$ with
$\|f\|_{\frac{3}{2}}<\infty$. Then there exists a sequence of
smooth functions $f_n \in C^\infty(S^1,\RR)\; (n\in \NN)$
converging to $f$ in the $\|\cdot\|_{\frac{3}{2}}$ norm. Moreover,
if ${\rm Supp}(f) \subset I$ for a certain interval $I\in\I$
and/or $f\geq 0$ then the sequence can be chosen so that
${\rm Supp}(f_n)\subset I$ and/or $f_n\geq 0$ for every $n\in\NN$.
\end{lemma}
\begin{proof}
The proof follows standard arguments relying on convolution with
smooth functions. Clearly, it is enough to prove the part of the
statement concerning the case when $f$ is localized and
nonnegative. Let $\varphi_k \,\, (k\in\NN)$ be a sequence of
positive smooth functions on $S^1$ with support shrinking to the
point $1 \in S^1$ such that for all $k \in \NN$
$\frac{1}{2\pi}\int_0^{2\pi} \varphi_k(e^{i\alpha}) \,d\alpha=1$.
Then the convolution $f*\varphi_k$ defined by the formula
\begin{equation}
(f*\varphi_k)(z) \equiv \frac{1}{2\pi}\int_0^{2\pi}
f(e^{i\theta})\varphi_k(e^{-i\theta}z) d\theta
\end{equation}
is a nonnegative real smooth function whose support, for $k$ large
enough, is contained in $I$. Moreover we have
\begin{equation}
\|f*\varphi_k-f\|_{\frac{3}{2}} = \sum_{n\in \ZZ}
|(\hat{(\varphi_k)}_n-1) \hat{f}_n|\,(1+|n|^{\frac{3}{2}}) \to 0
\end{equation}
as $k \to \infty$ since $|(\hat{\varphi_k})_n| \leq 1$,
$\|f\|_{\frac{3}{2}}<\infty$ and
$\lim_{k\to\infty}(\hat{\varphi_k})_n = 1$.
\end{proof}
Let $C^1(S^1,\RR)$ be the space of real valued continuously 
differentiable functions on the circle. We shall think of this 
space as a normed space, equipped with the norm defined by the formula 
\begin{equation}
\|f\|_1 \equiv {\rm max}\{|f|\}+{\rm max}\{|f'|\}
\end{equation}
for every $f\in C^1(S^1,\RR)$. It is evident, that every continuous real 
functions with finite $\|\cdot\|_{\frac{3}{2}}$ norm belongs to this 
space, as in fact
\begin{equation}
\|f\|_1 \leq \|f\|_{\frac{3}{2}}.
\end{equation}
Moreover, it is clear that every element of $C^1(S^1,\RR)$ is 
actually globally Lipschitz, and thus, by the compactness of $S^1$, 
one can consider its exponential, defined by the generated flow.
The following is an important observation
about the continuity of the exponential.

\begin{lemma} 
\label{f.3/2.exp.uniform}
Let $f_n\in C^1(S^1,\RR)$ be a sequence converging 
to the function $f\in C^1(S^1,\RR)$. 
Then the sequence ${\rm Exp}(f_n)\; (n\in \NN)$ 
converges uniformly to ${\rm Exp}(f)$. 
\end{lemma}
\begin{proof}\hspace{-2mm}\footnote{The argument here given
was indicated to the author by S.\! Carpi.}
In stead of working on the circle, we can obviously transfer the
problem to the real line, with $f,f_n (n\in \NN)$ being periodic 
functions. Then for any fixed point $x\in\RR$, the real function 
defined by the map $t\mapsto h_n(t) \equiv {\rm Exp}(tf)(x)-{\rm 
Exp}(tf_n)(x)$, is differentiable. In fact, 
\begin{eqnarray} \nonumber 
|h'_n(t)| &=& |f({\rm Exp}(tf)(x))-f_n({\rm Exp}(tf_n)(x))| \\
\nonumber &\leq& |f({\rm Exp}(tf)(x)) - f({\rm Exp}(tf_n)(x))| + 
|(f-f_n)({\rm Exp}(tf_n)(x))| \\
&\leq& |h_1(t)|\, \|f\|_1  + \|f-f_n\|_1. 
\end{eqnarray}
Since $h_n(0)=0$, we have that
$|h_n(t)| = |\int_0^t h'_n(s)ds| \leq \int_0^t |h'_n(s)|ds$, which by
applying the above estimate, results, that if $0\leq t \leq 1$ then
\begin{equation}
|h_n(t)| \leq \|f-f_n\|_1 + \int_0^t |h_n(s)| \|f\|_1 \,ds.
\end{equation}
Thus by an application of the Gronwall-lemma (see e.g.\! 
\cite[Sect.\! 8.4]{diff})
we find that $|h_n(1)| \leq \|f-f_n\|_1 \, e^{\|f\|_1}$. The statement 
then follows trivially.
\end{proof}
\paragrf
In the next section, as it was mentioned, we shall prove that
$T(f)$ can be defined as a self-adjoint operator (which is
essentially self-adjoint on the finite energy vectors) even when
$f$ is not smooth but at least $\|f\|_{\frac{3}{2}}<\infty$
(Theorem \ref{theo.main.nonsmooth}). We shall also prove that the
thus extended $T$ has the following continuity property: if $f_n
\to f$ in the $\|\cdot\|_{\frac{3}{2}}$ norm as $n\to \infty$,
then $T(f_n) \to T(f)$ in the strong resolvent sense (Prop.\!
\ref{f->T(f)continuity}). Taking these two statements granted for
the rest of the section (we shall not use in the next section
anything of what will follow now), we have the following picture:
\begin{itemize}
\item if $f,f_n$ $(n \in \NN)$ are real smooth functions on the
circle and $f_n\rightarrow f$ in the $\|\cdot\|_{\frac{3}{2}}$
sense then $T(f_n)$ converges to $T(f)$ in the strong resolvent
sense,
\item if $f_n$ $(n \in \NN)$ is a Cauchy sequence of real
smooth functions with respect to the $\|\cdot\|_{\frac{3}{2}}$
norm then $T(f_n)$ converges to a self-adjoint operator in the
strong resolvent sense, which is essentially self-adjoint on the
finite energy vectors,
\item the real smooth functions form a
dense set among the real continuous functions with finite
$\|\cdot\|_{\frac{3}{2}}$ norm.
\end{itemize}
Thus one can really think of this extension of $T$ as an extension
by continuity. As a last thing in this section, we shall state some 
important properties of this extension of the stress-energy tensor 
in relation with conformal nets.
\begin{proposition}
\label{aff} Let $(\A,U)$ be a conformal local net with
stress-energy tensor $T$ and further let $f\in C(S^1,\RR)$ of
finite $\|\cdot\|_{\frac{3}{2}}$ norm. Then for every interval
$I\in\I$
\begin{eqnarray*}
&(1).& e^{iT(f)}\A(I)e^{-iT(f)} = \A({\rm Exp}(f)(I)),\\
&(2).& f_{I^c} = 0 \,\;\Rightarrow \; T(f)\; {\rm
is\;affiliated\;to}\;\A(I).
\end{eqnarray*}
Moreover, if $\gamma\in\diff$
and $\|\gamma_* f\|_{\frac{3}{2}} < \infty$ where $\gamma_*$
stands for the action of $\gamma$ on vector fields, then up to
phase factors $U(\gamma)\, e^{iT(f)}\, U(\gamma)^*= e^{iT(\gamma_*
f)}$.
\end{proposition}
\begin{proof}
Of course for a smooth function by definition we have all the
required properties. As for the first assertion, by Lemma
\ref{smoothdensity} we can take a sequence $\{f_n:n\in\NN\}$ of
smooth functions converging to $f$ in the
$\|\cdot\|_{\frac{3}{2}}$ norm. Then, by Lemma \ref{f.3/2.exp.uniform},
${\rm Exp}(f_n) \to {\rm Exp(f)}$ uniformly. 
Thus, since by Prop.\! \ref{f->T(f)continuity}
$T(f_n)\to T(f)$, the result is obtained by the strong convergence
of $e^{iT(f_n)}$ to $e^{iT(f)}$ (which follows from the strong
resolvent convergence, see e.g.\! \cite[Sect. VIII.7]{RSI}) and the
continuity of the net $\A$.

The proof of the rest of the statement is similar. For the second
assertion, by the continuity of the net it is enough to prove that
if ${\rm Supp}(f) \subset I$ (and not only in its closure) then
$T(f)$ is affiliated to $\A(I)$. Then by Lemma
\ref{smoothdensity}, there exists a sequence of smooth functions
$f_n$ $(n \in \NN)$ converging to $f$ in the
$\|\cdot\|_{\frac{3}{2}}$ norm whose support is contained in $I$.
Thus again, the result follows from the strong resolvent
convergence of $T(f_n)$ to $T(f)$. Similarly, the proof of the
last part of the statement is obtained by taking appropriate
limits.
\end{proof}
Let $I_1,I_2 \in \I$ be two distant intervals. In the appendix
(Lemma \ref{aI1I2z}) it is proven that the map $s\mapsto
\kappa^{I_1,I_2}_s(z)$ is differentiable at $s=0$ for every
$z\in S^1$ and that denoting by $a^{I_1,I_2}$ its derivative at
$s=0$ we have that $a^{I_1,I_2}\in\pcwm$. As an implication of the
previous proposition we have the following.
\begin{corollary}
\label{U(kappa)=exp(iT(a))} Let $(\A,U)$ be a conformal local net
and $I_1,I_2\in\I$ two distant intervals. Then the unitary
$e^{isT(a^{I_1,I_2})}$ is an implementation of the transformation
$\kappa^{I_1,I_2}_s$ in the sense of Eq.\! (\ref{k-kov.action}).
\end{corollary}
\begin{proof}
We shall show that the restriction of the adjoint action of
$e^{isT(a^{I_1,I_2})}$ to $\A(I_1)$ coincides with that of
$U(\delta^{I_1}_s)$. The case of the other three remaining
intervals can be dealt with in a similar manner.

By the previous proposition the adjoint action of both unitaries
preserves $\A(I_1)$. Set $d^{I_1}$ for the generating \mo vector
field of $s\mapsto\delta^{I_1}_s$. Then $a^{I_1,I_2}-d^{I_1}$ is
zero on $I_1$, thus by the previous proposition (the closure of) 
$T(a^{I_1,I_2})-T(d^{I_1})$ is affiliated to $\A(I_1^c)$.
Then by an application of the Trotter product-formula (since both
the sum and the operators involved have a common core: the set of
finite energy vectors, the strong convergence of the formula
indeed follows, see e.g.\! \cite[Theorem VIII.31]{RSI}) one can
easily finish the proof.
\end{proof}

\section{Proof of essentially self-adjointness}
\label{sec:nonsmoothestimate}

Throughout this section $U$ will be a (fixed) positive energy
representation of $\diff$ with associated representation of the
Virasoro algebra $\{L_n: n\in \ZZ\}$ defined on the dense subspace
$\fin$ of its finite energy vectors. We shall denote by $T$ the
stress-energy tensor corresponding to $U$.

Recall the energy bound given by equation \ref{e.bound.1b},
according to which by setting $r=\sqrt{1+\frac{c}{12}}$ we have
\begin{equation}
\label{energyestimate} \nonumber \|L_n v\| \leq r
(1+|n|^{\frac{3}{2}})\,\,\|(\mathbbm 1 +L_0)v\|
\end{equation}
for every $v \in \fin$.

\begin{proposition}
\label{AonD} If $a_n \in {\mathbb C} \,\, (n \in {\mathbb Z})$ is
such that $\mathop{\sum}_{n \in {\mathbb Z}}|a_n|
(1+|n|^{\frac{3}{2}}) < \infty$ then
\begin{itemize}
\item[(i)] the operator $A\equiv\mathop{\sum}_{n \in {\mathbb
Z}}a_n L_n$ with domain $\fin$ is well defined, (i.e. the sum
strongly converges on the domain);

\item[(ii)] $\|Av\| \leq r \left(\mathop{\sum}_{n \in {\mathbb
Z}}|a_n| (1+|n|^{\frac{3}{2}})\right) \|(\mathbbm 1 + L_0) v\|$
where $r$ is the constant appearing in the energy bound given by
equation (\ref{energyestimate}) and $v\in \fin$;

\item[(iii)] $A^*$ is an extension of the operator
$A^+\equiv\mathop{\sum}_{n \in {\mathbb Z}}\overline{a}_{-n} L_n$.
(This again is understood as an operator with domain $\fin$.)
\end{itemize}
\end{proposition}
\begin{proof}
Since
\begin{equation}
\label{domain} \mathop{\sum}_{n \in {\mathbb Z}} \|a_n L_n v\|
\leq r \left( \mathop{\sum}_{n \in {\mathbb Z}} |a_n|
(1+|n|^{\frac{3}{2}}) \right) \|(\mathbbm 1 +L_0)v\| < \infty
\end{equation}
claim (i) and (ii) holds. Finally, the last claim follows, since
$L_n^* \supset L_{-n}$ for every integer $n$.
\end{proof}

We now consider for every $\epsilon >0$ the operator
$R_{n,\epsilon}=[L_n,e^{-\epsilon L_0}]$ with domain $\fin$. (Note
that $e^{-\epsilon L_0}\fin \subset \fin$ so it is indeed well
defined on the finite energy vectors.) The following proposition
gives an estimate on the norm of this commutator which is
independent of $\epsilon$.
\begin{proposition}
\label{commutatornorm}
There exists a constant $q > 0$ independent of $\epsilon$ and $n$ such
that $\|R_{n,\epsilon}\|^2=\|[L_n,e^{-\epsilon L_0}]\|^2 \leq q |n|^3$.
\end{proposition}
\begin{proof}
For $n=0$ the statement is trivially true as $L_0$ commutes with
any bounded function of itself. As the operators in question all
preserve the finite energy vectors and moreover $L_n = L_{-n}^*$
and $e^{-\epsilon L_0}$ is self-adjoint it follows that
$R_{n,\epsilon} \subset -R_{-n,\epsilon}^*$. Thus it suffices to
demonstrate the statement for negative values of $n$.

Let therefore be $n<0$, $v \in \fin$ and for every $k \in {\rm
Sp}(L_0)$ let $v_k$ be the component of the vector $v$ in the
eigenspace of $L_0$ associated to the eigenvalue $k$. To not to
get confused about positive and negative constants, in the
calculations we shall use the positive $m:=-n$ rather than the
negative $n$. As $L_n$ raises the eigenvalue of $L_0$ by $m$, we
have that
\begin{equation}
\label{Rv} R_{n,\epsilon}v_k=[L_n,e^{-\epsilon L_0}]v_k =
(e^{-\epsilon k}-e^{-\epsilon (k+m)})\,L_n v_k.
\end{equation}
The mapping
\begin{equation}
f_m: \epsilon \mapsto e^{-\epsilon k}-e^{-\epsilon (k+m)}
\end{equation}
is a positive smooth function on $\RR^+$ which goes to zero both
when $\epsilon \rightarrow 0$ and when $\epsilon \rightarrow
\infty$. Therefore $f_m$ has a maximum on $\RR^+$. But the only
solution of the equation $f'_m(\epsilon)=0$ is
\begin{equation}
\epsilon_m=-\frac{1}{m} \ln \left(\frac{k}{k+m}\right)
\end{equation}
which then must correspond to the point of maximum. Then by direct
substitution we get
\begin{equation}
\label{f_m}
\mathop{\rm sup}_{\epsilon \in \RR^+}{|f_m(\epsilon)|^2}=
f_m(\epsilon_m)^2= \left(\frac{k}{k+m}\right)^{\frac{2k}{m}}
\left(\frac{m}{k+m}\right)^2 \leq \left( \frac{m}{k+m} \right)^2.
\end{equation}
We can now return to the question of the norm of the commutator.
Equation (\ref{Rv}) shows that the vectors $R_{n,\epsilon}v_k\,\,
(k \in {\rm Sp}(L_0))$ are in particular pairwise orthogonal.
Using this and the fact that only for finitely many values of $k$
the vector $v_k \neq 0$ we find that
\begin{eqnarray}
\nonumber \|R_{n,\epsilon}v\|^2 &=& \|R_{n,\epsilon}
\!\!\mathop{\sum}_{k \in {\rm Sp}(L_0)} \!\! v_k\|^2=
\!\!\mathop{\sum}_{k\in {\rm Sp}(L_0)} \|R_{n,\epsilon}v_k\|^2
=\!\!\mathop{\sum}_{k \in {\rm Sp}(L_0)}
\!|f_m(\epsilon)|^2\|L_nv_k\|^2
\\ \nonumber
&\leq& \!\mathop{\sum}_{k \in {\rm Sp}(L_0)} \mathop{\rm
sup}_{\epsilon \in \RR^+}\{|f_m(\epsilon)|^2\} \, \|L_n v_k\|^2
\\ \nonumber
&\leq&
\!\mathop{\sum}_{k \in {\rm Sp}(L_0)} \left(\frac{m}{k+m}
\right)^2 r^2 (k^2+km^2+m^3) \,\|v_k\|^2
\\ &\leq&
\!\mathop{\sum}_{k \in {\rm Sp}(L_0)} r^2 (m^2+m^3+m^3)
\,\|v_k\|^2 \leq  3 r^2 |n|^3 \,\|v\|^2,
\end{eqnarray}
where we have used the inequality (\ref{f_m}) and the constant $r$
is the one appearing in equation (\ref{energyestimate}).
\end{proof}

\begin{theorem}
\label{theo.main.nonsmooth} If $a_n \in {\mathbb C} \,\, (n \in
{\mathbb Z})$ is such that $\mathop{\sum}_{n \in {\mathbb Z}}|a_n|
|n|^{\frac{3}{2}} < \infty$ then $A$ is closable and
$\overline{A}=(A^+)^*$, where $A\equiv\mathop{\sum}_{n \in
{\mathbb Z}}a_n L_n$ and $A^+\equiv\mathop{\sum}_{n \in {\mathbb
Z}}\overline{a}_{-n} L_n$ considered as operators with domain
$\fin$. In particular, if $a_n=\overline{a_{-n}}$ for all $n \in
{\mathbb Z}$, then $A$ is essentially self-adjoint on $\fin$.
\end{theorem}
\begin{proof}
By Prop.\! \ref{AonD} the operators $A$, $A^+$,
$R_{A,\epsilon}\equiv[A,e^{-\epsilon L_0}]$ and
$R_{A^+,\epsilon}\equiv[A^+,e^{-\epsilon L_0}]$ are indeed well
defined on the finite energy vectors. By using Prop.\!
\ref{commutatornorm} with its constant $q$ appearing in it and the
condition on the sequence $a_n \,\, (n \in {\mathbb Z})$,
\begin{equation}
\mathop{\sum}_{n \in {\mathbb Z}}\|a_n R_{n,\epsilon}\| \leq
\mathop{\sum}_{n \in {\mathbb Z}}|a_n| q^{\frac{1}{2}} |n|^{\frac{3}{2}}
< \infty.
\end{equation}
This means that $\|R_{A,\epsilon}\|$ is bounded by a constant
independent of $\epsilon$. Obviously, the same is true for
$\|R_{A^+,\epsilon}\|$. Thus for example, taking account of the
fact that $A^+\subset A^*$ (see Prop.\! \ref{AonD}) we have that
$R_{A,\epsilon}^*$ is the closure of both $-R_{A^+,\epsilon}$ and
of $-[A^*,e^{-\epsilon L_0}]$.

If $v_k$ is an eigenvector of $L_0$ with eigenvalue $k \in {\rm
Sp}(L_0)$ then as $\epsilon$ tends to zero,
\begin{equation}
R_{{A^+},\epsilon}v_k=(e^{-\epsilon k}{\mathbbm 1}-e^{-\epsilon
L_0})A^+v_k
\rightarrow 0.
\end{equation}
Hence the operators $R_{A^+,\epsilon}$ on $\fin$ strongly converge
to zero as $\epsilon \to 0$. Then since their norm is bounded by a
constant independent of $\epsilon$,  also the everywhere defined
bounded operators $R_{A,\epsilon}^*$ converge strongly to zero as
$\epsilon \to 0$.

Finally, remarking on Prop.\! \ref{AonD} let us observe two more
things. First, that the operator $A^+$ is closable (as its adjoint
is an extension of the densely defined operator $A$). Second, the
smooth energy vectors $\D^\infty=\cap_{n-1}^\infty \D(L_0^n)$ are
in the domain of its closure (as $A^+$ is energy bounded).

From here the proof of the theorem continues exactly as in
\cite{Ne}, but for self-containment let us revise the concluding
argument. Suppose $x$ is a vector in the domain of $A^*$. Then,
since $e^{-\epsilon L_0}x \in \D^\infty \subset
\D(\overline{A^+})$, we have that
\begin{equation}
\label{eq:Nelson} \overline{A^+} e^{-\epsilon L_0} x = A^*
e^{-\epsilon L_0}x = e^{-\epsilon L_0} A^* x  - R^*_{A,\epsilon}x.
\end{equation}
As $\epsilon$ tends to zero of course $e^{-\epsilon L_0}x
\rightarrow x$, but now Eq.\! (\ref{eq:Nelson}) shows that also
$A^+ e^{-\epsilon L_0} x \rightarrow A^*x$ strongly. Therefore
$A^*=\overline{A^+}$.
\end{proof}

With this we have proved the main theorem of this section on which
the previous section was essentially based. The result ensures,
that if the continuous function $f: S^1 \rightarrow \RR$ with
Fourier coefficients $\hat{f}_n$ $(n \in \ZZ)$ has a finite
$\|\cdot\|_{\frac{3}{2}}$ norm then $\mathop{\sum}_{n \in {\mathbb
Z}} \hat{f}_n L_n$ is an essentially self-adjoint operator on
$\fin$. As in the case of smooth functions, we denote by $T(f)$
the corresponding self-adjoint operator obtained by taking
closure. We finish by showing the continuity property of the
stress-energy tensor $T$ which too, was mentioned (and used) in
the previous section without proof.
\begin{proposition}
\label{f->T(f)continuity} For every $f \in C(S^1,\RR)$ with
$\|f\|_{\frac{3}{2}}<\infty$ and $v \in \D(L_0)$
\begin{equation}
\|T(f)v\| \leq r \|f\|_{\frac{3}{2}}\, \|(\mathbbm 1+L_0)v\|
\end{equation}
where $r$ is the positive constant appearing in the energy bound
Eq.\! (\ref{energyestimate}). Moreover, if $f$ and $f_n \,\, (n
\in \NN)$ are continuous real functions on $S^1$ of finite
$\|\cdot\|_{\frac{3}{2}}$ norm, and in this norm $f_n \to f$ as
$n\to \infty$, then $T(f_n) \to T(f)$ in the strong resolvent
sense. In particular, $e^{iT(f_n)} \to e^{iT(f)}$ strongly.
\end{proposition}
\begin{proof}
The claimed inequality is an immediate consequence of the
inequality in Eq.\! (\ref{energyestimate}), the definition of
$T(\cdot)$ and that of the $\|\cdot\|_{\frac{3}{2}}$ norm. It
follows then that $T(f_n)v$ converges to $T(f)v$ for every $v \in
\fin$ --- in fact, for every $v\in\D(L_0)$. Since $\fin$ is a
common core for these self-adjoint operators, the conclusion
follows (see e.g. \cite[Sect.\! VIII.7]{RSI}).
\end{proof}

\chapter[Applications and further results]
{\vspace{1cm}Applications and further results}

\resume{
In this chapter, by a direct application of the piecewise \mo
transformations discussed before, a short proof will be given to
demonstrate the uniqueness of the diffeomorphism symmetry assuming
$4$-regularity. The argument relies on the density of piecewise
\mo transformations, so it will be also shown that for every
$\gamma\in\diff$ a sequence of piecewise \mo transformations
$\gamma_n\; (n\in\NN)$ converging uniformly to $\gamma$ can be
constructed with the property that $U(\gamma_n)\to U(\gamma)$
strongly for any regular conformal local net $(\A,U)$ (so they are
dense in the ``right'' sense). Still with the same method it will
be shown that a \mo covariant subnet of a conformal net is always
conformal, assuming again some regularity condition.

Then a completely different argument will be presented to prove
again these two results. This second argument will be somewhat
longer but on the other hand it will work without any regularity
assumption. It will also show that given a subnet, the
stress-energy tensor can be split into a part belonging to the
subnet and another part belonging to the coset of the subnet.

Afterwards some applications, such as the construction of new
examples of \mo covariant nets admitting no $\diff$ symmetry, will
be discussed. Finally with the technique of evaluating the
stress-energy tensor on nonsmooth function it will be shown that
every locally normal representation of a diffeomorphism covariant
net is of positive energy.
}

\markboth{CHAPTER \arabic{chapter}. APPLICATIONS AND FURTHER
RESULTS}{CHAPTER \arabic{chapter}. APPLICATIONS AND FURTHER
RESULTS}

\section[Uniqueness of the $\diff$ symmetry]
{Uniqueness of the $\mathbf {{\rm \bf Diff}^+(S^1)}$ symmetry}
\label{sec:uniq2}

\paragrf
Let $(\A,V)$ be a \mo covariant local net and suppose $V$ can be
extended to be a projective representation of $\diff$ in a
compatible way with the local structure of $\A$. Such an extension
makes the net diffeomorphism covariant (see the remark after
Def.\! \ref{diffcov:def} on the automatic covariance of the
action). The problem we shall consider is whether such an
extension of the \mo symmetry is unique.

Two different proofs will be presented for the uniqueness of the
diffeomorphism symmetry. The first one will be a direct
application of the previous considerations regarding local
algebras and nonsmooth symmetries, but requires $4$-regularity.
The second one will be an argument that, in its spirit, is closer
to Quantum Fields (and in particular, to Vertex Operator Algebras as it is 
carried out at the level of Fourier modes)
than to local algebras. It relies on
exploiting restrictions on possible commutation relations coming
from the fact that the stress-energy tensor is a so-called
dimension $2$ field whose integral is the total conformal energy. 
Although found independently, it is  --- as it was pointed out to the
author by Karl-Henning Rehren ---  essentially a Fourier-transformed
version of an old idea of L\"uscher and Mack (see \cite[Sect.\!
I\!I\!I.]{lmack} for their argument). However, as it was mentioned in
the introduction, what seems to be really new is the application
in the algebraic setting discussed in this and the next section.
Note also --- as it will be clear later on --- that for such application
a detailed analyses is needed involving the use of energy bounds and in 
particular the relation between local structure and  energy bounds.

The first proof is contained in the joint work \cite{CaWe} of the
author with S.\! Carpi, although with some differences which will
be later explained. The second one is the author's own but also
takes ideas from the mentioned work; in particular from the way
the uniqueness was there proved for the special case of Virasoro
nets.

\paragrf
Throughout this section $(\A,V)$ will be a \mo covariant local net
and $U$ and $\tilde{U}$ two strongly continuous projective unitary
representations of $\diff$ both extending $V$ and making $\A$
conformal covariant. Of course by definition the restrictions of
these two representations to $\mob$, coincide (as both
restrictions result $V$). Our aim is to prove that $\tilde{U}$
must coincide with $U$ on the whole group of $\diff$.

The two representations $U$ and $\tilde{U}$ give rise to two
stress-energy tensors $T$ and $\tilde{T}$ and also two associated
representations of the Virasoro algebra $\{L_n: n\in \ZZ\}$ and
$\{\tilde{L}_n: n\in \ZZ\}$. By our assumptions, $L_j=\tilde{L}_j$
for $j=-1,0,1$.

\paragrf
Let $I_1,I_2\in \I$ be two distant intervals. Recall the
the one-parameter group of piecewise \mo
transformations $s\mapsto\delta^{I_1,I_2}_s$ 
introduced by Def.\! \ref{def:kappa},
and its generating
(piecewise \mo) vector field $a^{I_1,I_2}$ (see Lemma \ref{aII}).
\begin{lemma}
If the net $\A$ is at least $4$-regular then
$T(a^{I_1,I_2})=\tilde{T}(a^{I_1,I_2})$.
\end{lemma}
\begin{proof}
By Corollary \ref{U(kappa)=exp(iT(a))}, in the sense of Eq.\!
(\ref{k-kov.action}) both $s\mapsto e^{isT(a^{I_1,I_2})}$ and
$s\mapsto e^{is\tilde{T}(a^{I_1,I_2})}$ are implementations of the
one-parameter group of transformations
$s\mapsto\delta^{I_1,I_2}_s$ in $(\A,V)$, and by the observation
made before Corollary \ref{uniq.pcw.symm.}, if the net is
$4$-regular, then in the projective sense, such an implementation is 
unique. Thus up to an additive constant $T(a^{I_1,I_2})$
and $\tilde{T}(a^{I_1,I_2})$ must coincide. But since the vacuum
expectation of both of them is zero, this additive constant is actually
zero.
\end{proof}
{\it Remark.} The operator
$V(\delta^{I_1}_1)^*e^{iT(a^{I_1,I_2})}$ is localized in $I^c$ and
it equals to $V(\delta^{I_1}_1)^*e^{i\tilde{T}(a^{I_1,I_2})}$.
This shows that the local intersection of the two Virasoro subnets
(defined by $U$ and $\tilde{U}$) is not trivial. Then by the
minimality of the Virasoro net (cf.\! \cite{Carpi98}) it
follows that the two subnets coincide. Together with a separate
argument treating the case of the Virasoro net, this gives a proof
of the uniqueness for the general case of a $4$-regular conformal
net. In the mentioned joint work \cite{CaWe} that is how the
result is achieved. Thus the problem regarding the density of
piecewise \mo transformations was cleverly avoided as, by
minimality, a single non-trivial piecewise \mo transformation was
sufficient for the proof. Here the argument given will be
different. Minimality will be not used, and nor will be need for a
separate argument treating the case of the Virasoro net. Rather,
we shall directly address the question of density.
\begin{theorem}
\label{4reg->unique} Let $(\A,V)$ be a \mo covariant local net
which is at least $4$-regular and let $U$ and $\tilde{U}$ be two
strongly continuous projective unitary representations of $\diff$
both extending $V$ and making $\A$ a conformal covariant local
net. Then $U=\tilde{U}$.
\end{theorem}
\begin{proof}
If $I_1,I_2\in\I$ are distant intervals then by the previous lemma
we have that $T(a^{I_1,I_2})=\tilde{T}(a^{I_1,I_2})$. The two
stress-energy tensors also coincide on any \mo vector field.
Moreover, it is clear that if the stress-energy tensors coincide
on some functions then in fact they coincide on any linear
combination of those functions (the set of finite energy vectors
is a common core so there is no problem with taking sums). Thus
$T(f)=\tilde{T}(f)$ for every piecewise \mo vector field, since,
as it is proved in the appendix (Theorem \ref{pcwmgenerators}),
the functions $\{a^{I_1,I_2}\}$ (where $I_1,I_2\in\I$ run through
all pairs of distant intervals) together with the \mo vector
fields span the space of piecewise \mo vector fields. But it is also
proved in the appendix (Theorem \ref{pcwmdensity}), that the
piecewise \mo vector fields are dense among the functions of
finite $\|\cdot\|_{\frac{3}{2}}$ norm. Thus by Prop.\!
\ref{f->T(f)continuity} we can conclude the equality of the two
stress-energy tensors. As it was mentioned, $\diff$ is generated
by exponentials; hence by the equality of the stress-energy
tensors $U=\tilde{U}$.
\end{proof}
{\it Remark.} We have used the density of the piecewise \mo vector
fields, but it is worth discussing the density of piecewise \mo
transformations. By Lemma \ref{diffcov->pcwcov} and Theorem
\ref{existence.pcw.symm.} if a conformal net $(\A,U)$ is regular
it is also piecewise \mo covariant and its piecewise \mo symmetry
is unique. Thus we can safely write expressions like $U(\gamma)$
where $\gamma\in\pcw$.
\begin{theorem}\label{density.in.right.sense}
Let $\gamma\in\diff$. Then there is a sequence of piecewise \mo 
transformations $\gamma_n\; (n\in\NN)$ converging
uniformly to $\gamma$ such that $U(\gamma_n)\to U(\gamma)$
strongly for any conformal, regular, local net $(\A,U)$.
\end{theorem}
\begin{proof}
We can assume that $\gamma={\rm Exp}(f)$ for a certain smooth
function $f$, since $\diff$ is generated by the exponentials. Consider 
this statement for all such ${\rm Exp}(f)$ with $f$ being of 
finite $\|\cdot\|_{\frac{3}{2}}$ norm. If it is true
for $f,g$ and any multiples of them, then it is also true for
$f+g$. Indeed, since $f$ and $g$ are globally Lipschitz, and since
the circle is compact, by an application of \cite[Theorem 4.1]{flows} 
we find that $({\rm Exp}(f/n)\circ {\rm
Exp}(g/n))^n$ converges uniformly to ${\rm Exp}(f+g)$, while by
the Trotter product-formula (as the operators involved have a
common core) we have that $(e^{iT(f/n)}e^{iT(g/n)})^n$ converges
strongly to $e^{iT(f+g)}$. Of course the statement is clearly true
for $\lambda a^{I_1,I_2}$ where $I_1,I_2\in\I$ are two distant
intervals and $\lambda \in \RR$. Also, it is by definition true
for any \mo vector field. Thus by the above argument, taking into
account the already mentioned generating property of the listed
functions, we have justified the statement for any piecewise \mo
vector field.

If $f$ is a smooth function, then by the density of piecewise \mo
vector fields we can take a sequence $f_n$ of piecewise \mo vector
fields converging to $f$ in the $\|\cdot\|_{\frac{3}{2}}$ norm.
Then by Prop.\! \ref{f->T(f)continuity} we have that $e^{iT(f_n)}$
converges strongly to $e^{iT(f_n)}$ while by Lemma 
\ref{f.3/2.exp.uniform}, ${\rm Exp}(f_n)\to {\rm Exp}(f)$ uniformly. 
This, by what was explained, proves the theorem.
\end{proof}
{\it Remark.} Note that this theorem really gives a way to
``create'' diffeomorphism symmetry. The piecewise \mo
transformations are constructed by the local algebras and the
representation of the \mo group. On the other hand, a general
diffeomorphism (and the unitary associated to it) we can obtain by
taking limits. So if we have a net for which it is not known
whether it does or does not have diffeomorphism symmetry, but we
do know that it is piecewise \mo covariant
--- and as we have seen, this can be established in many cases
and concrete examples even when the existence of conformal covariance
is not known --- then we have a ``formula'' for $U(\gamma)$: it
should be the strong limit of $U(\gamma_n)$. For this reason the
author conjectures the piecewise \mo symmetry implies
diffeomorphism symmetry.

\paragrf
We shall now make a completely different proof, which is slightly
less related to the idea of local algebras but has the advantage
that it does not require $4$-regularity. So we still assume that
we have two representations of $\diff$. The notations will be
those that were fixed in the beginning of this section. By our
assumptions $L_j=\tilde{L}_j$ for $j=-1,0,1$. In particular the
dense subset of the finite energy vectors $\fin$ is unambiguous
(i.e.\! it is the same for both representations) and any
polynomial expression of the operators $\{L_n: n\in \ZZ\}$ and
$\{\tilde{L}_n: n\in \ZZ\}$ are well-defined on it. As usual
$\Omega$ will stand for the vacuum vector.
\begin{proposition}
\label{L2Ln} $[L_2, \tilde{L}_{-n}]\Omega =
(2+n)\tilde{L}_{-n+2}\Omega$ for every integer $n\geq 3$.
\end{proposition}
\begin{proof}
The operator $L_{-1}=\tilde{L}_{-1}$ annihilates $\Omega$. As the
vector $L_2 \tilde{L}_{-2}\Omega$ is a vector of zero energy
(i.e.\! an eigenvector of $L_0$ with eigenvalue $0$), we have that
$L_{-1} L_2 \tilde{L}_{-2}\Omega = 0$. Moreover, using the
commutation relations we find that
\begin{eqnarray}
\nonumber [L_2,L_{-1}] \tilde{L}_{-2}\Omega &=& 3 L_1
\tilde{L}_{-2}\Omega \\ \nonumber &=& 3 [L_1,\tilde{L}_{-2}]\Omega
\\ &=& 9 L_{-1}\Omega = 0.
\end{eqnarray}
Therefore, as $L_0 \Omega = 0$,
\begin{eqnarray}
\nonumber [L_2 ,\tilde{L}_{-3}]\Omega &=& L_2\tilde{L}_{-3}\Omega
\\ \nonumber
&=& L_2 L_{-1}\tilde{L}_{-2}\Omega \\
&=& [L_2,L_{-1}] \tilde{L}_{-2}\Omega + L_{-1} L_2
\tilde{L}_{-2}\Omega = 0
\end{eqnarray}
and of course also $(2+3)\tilde{L}_{-3+2}\Omega = 5 L_{-1}\Omega =
0$. Thus the statement for $n=3$ is true. We shall continue by
induction. Assume that the statement holds for a certain integer
$n=k\geq 3$. Then for $n=k+1$ we have
\begin{equation}
\tilde{L}_{-n}\Omega = \frac{1}{k-1} [L_{-1},L_{-k}]\Omega =
\frac{1}{k-1} L_{-1} L_{-k}\Omega
\end{equation}
and thus
\begin{eqnarray}
\nonumber [L_2, \tilde{L}_{-n}] \Omega &=&
\frac{1}{k-1} L_2 L_{-1} \tilde{L}_{-k} \Omega \\
\nonumber &=& \frac{1}{k-1}\left([L_2,L_{-1}] L_{-k} +
L_{-1} L_2 \tilde{L}_{-k}\right)\Omega \\
\nonumber &=& \frac{1}{k-1} \left( 3 L_1 \tilde{L}_{-k} +
(2+k)L_{-1}\tilde{L}_{-k+2} \right)\Omega \\ \nonumber &=&
\frac{1}{k-1} \left( 3 [L_1,\tilde{L}_{-k}] +
(2+k)[L_{-1},\tilde{L}_{-k+2}] \right)\Omega \\ \nonumber &=&
\frac{3 (1+k) + (2+k)(k-3)}{k-1} \tilde{L}_{-k+1} \Omega \\
&=& (3+k)\tilde{L}_{-k+1} \Omega =(2+n)\tilde{L}_{-n+2} \Omega
\end{eqnarray}
where in the third line we have used the inductive assumption.
\end{proof}
Recall that the function $l_n$ for each integer $n\in\ZZ$ was
defined by the formula $l_n(z)=-iz^n$ and that on the finite
energy vectors $T(il_n)=L_n$ (and also
$\tilde{T}(il_n)=\tilde{L}_n$). For a smooth function $f$ the set
of smooth vectors $\D^\infty$ is an invariant core for both $T(f)$
and $\tilde{T}(f)$. Thus any polynomial expression containing the
two stress-energy tensors evaluated on some smooth functions is
well-defined on $\D^\infty$ (and hence also on the smaller
subspace $\fin$).
\begin{lemma}
Let $f$ be a smooth function on the circle with $\hat{f}_{-2}=0$
where $\{\hat{f}_n: n\in\ZZ\}$ are the Fourier coefficients of the
function $f$. Then $$
[T(l_2),\tilde{T}(f)]\Omega=[\tilde{T}(l_2),\tilde{T}(f)]\Omega.
$$
\end{lemma}
\begin{proof}
It is a simple consequence of our previous proposition. As $f$ is
smooth its Fourier coefficients are fast decreasing and by a
simple use of the energy bounds it is easy to see that the
operators $T(l_2)$ and $\tilde{T}(l_2)$ can be brought behind the
summation when applied on sums like $\sum_{n\in\ZZ}\hat{f}_n L_n
\Omega$. Thus by Prop.\! \ref{L2Ln},
\begin{eqnarray}
\nonumber i[T(l_2),\tilde{T}(f)]\Omega &=& i[T(l_2),
\sum_{n\in\ZZ} \hat{f}_n \tilde{L}_n] \Omega = i[T(l_2),
\sum_{n=3}^{\infty} \hat{f}_{-n} \tilde{L}_{-n}] \Omega \\
\nonumber &=& \sum_{n=3}^{\infty}  \hat{f}_{-n}[ L_2,
\tilde{L}_{-n}] \Omega =
\sum_{n=3}^{\infty}  \hat{f}_{-n} (2+n) \tilde{L}_{-n+2}\Omega \\
\nonumber &=& \sum_{n=3}^{\infty}  \hat{f}_{-n} [\tilde{L}_2,
\tilde{L}_{-n}] \Omega = \sum_{n\in\ZZ} \hat{f}_{-n} [\tilde{L}_2,
\tilde{L}_{-n}] \Omega \\ &=& i[\tilde{T}(l_2),\tilde{T}(f)] \Omega
\end{eqnarray}
as in the summations $\hat{f}_{-n} \tilde{L}_{-n}\Omega = 0$
whenever $n<3$.
\end{proof}
\begin{proposition}
\label{[L2,T(f)]loc} Let $f$ be a smooth function on the circle
with $\hat{f}_{-2}=0$ and suppose that ${\rm Supp}(f)\subset I$
for a certain $I\in \I$. Then on the set of finite energy vectors
$$ [T(l_2),\tilde{T}(f)]=[\tilde{T}(l_2),\tilde{T}(f)].
$$
\end{proposition}
\begin{proof}
We can choose a bigger interval $K\in\I, \overline{I}\subset K$
and a smooth function $\varphi$ which is constant $1$ on the
interval $I$ and zero outside $K$. Then on the smooth energy
vectors $T(l_2) = T(\varphi l_2)+T((1-\varphi)l_2)$, and since
$T((1-\varphi)l_2)$ is affiliated to $\A(I^c)$ while
$\tilde{T}(f)$ is affiliated to $\A(I)$, (still on the smooth
energy vectors) we find that
\begin{equation}
[T(l_2),\tilde{T}(f)]=[T(\varphi l_2),\tilde{T}(f)]
+[T((1-\varphi)l_2), \tilde{T}(f)] =[T(\varphi l_2),\tilde{T}(f)]
\end{equation}
and similarly $[\tilde{T}(l_2),\tilde{T}(f)]= [\tilde{T}(\varphi
l_2),\tilde{T}(f)]$. Thus if $A'\in\A^\infty(K^c)$ then by our
previous lemma
\begin{eqnarray}
\nonumber [T(l_2),\tilde{T}(f)] A'\Omega &=& [T(\varphi
l_2),\tilde{T}(f)] A' \Omega \\ \nonumber &=& A' [T(\varphi
l_2),\tilde{T}(f)] \Omega
\\ \nonumber &=& A' [T(l_2),\tilde{T}(f)] \Omega \\ \nonumber &=&
A'[\tilde{T}(l_2),\tilde{T}(f)]\Omega \\ \nonumber &=&
A'[\tilde{T}(\varphi l_2),\tilde{T}(f)] \Omega \\ &=&
[\tilde{T}(\varphi l_2),\tilde{T}(f)] A'\Omega  =
[\tilde{T}(l_2),\tilde{T}(f)] A' \Omega
\end{eqnarray}
where we have also used that the operators $T(\varphi l_2)$ and
$\tilde{T}(\varphi l_2)$ are affiliated to $\A(K)$ and thus they
commute with $A'$.

We have proved the desired equality on a dense set of vectors but
not exactly on $\fin$. In order to finish the demonstration we
must take limits. To do so we need something better then just the
``simple'' energy bound on the stress-energy tensor, since a
commutator involves dealing with quadratic expressions. However,
by the polynomial energy bound of Lemma \ref{poly.e.bound} there
exists a positive constant $q$ such that for all $\Psi \in
\D^\infty$
\begin{equation}
\label{energy^2bound} {\rm
max}\left\{\|[\tilde{T}(l_2),\tilde{T}(f)]\Psi\|,\,
\|[T(l_2),\tilde{T}(f)]\Psi\|\right\} \leq q \|(\mathbbm 1 +
L_0^2)\Psi\|.
\end{equation}
As it was shown in the preliminaries (Theorem \ref{localcore}),
the set $\A^\infty(K^c)\Omega$ is a core for every positive power
of $L_0$. Thus by using the constant $q$ of the quadratic energy
bound of Eq.\! (\ref{energy^2bound}) on the commutators and the
fact that on the mentioned set we have already verified the
desired equality, one can easily finish the proof.
\end{proof}
\begin{lemma}
Let $f$ be a smooth function on the circle with $\hat{f}_{-2}=0$.
Then there exist two intervals $I,K \in \I$ and two smooth
functions $g$ and $h$ such that \\
\indent $(1)$. $f=g+h$, \\
\indent $(2)$. ${\rm Supp}(g)\subset I$ and ${\rm Supp}(h)\subset
K$.\\ \indent $(3)$. $\hat{g}_{-2}=\hat{h}_{-2}=0$
\end{lemma}
\begin{proof}
If $f=0$ then the statement is clearly true. If $f$ is not the
constant $0$ function then we can choose two intervals $I,K\in\I$
for which $I\cup K = S^1$ and $(I\setminus \overline{K^c})\cap
{\rm Supp}(f) \neq \emptyset$. We can also take a smooth function
$\varphi_0$ such that it is constant $0$ on $I^c$ and constant $1$
on $K^c$.

It is clear that by setting $g_0\equiv \varphi_0 f$ and $h_0
\equiv (1-\varphi_0)f$ property $(1)$ and $(2)$ are satisfied.
However, it can happen that $(\hat{g_0})_{-2}\neq 0$. So consider
the function $f l_2$. Its support coincides with that of $f$ so by
our conditions it cannot be that every real smooth function
supported in $I\setminus \overline{K^c}$ is orthogonal to it in
the ${\mathcal L}^2$ sense. Thus there exists a real smooth
function $\varphi$ such that
\begin{equation}
\frac{1}{2\pi} \int_0^{2\pi}
\varphi(e^{i\theta})l_2(e^{i\theta})f(e^{i\theta})d\theta =
(\hat{g_0})_{-2}.
\end{equation}
But the right hand side is nothing else than $(\hat{\varphi
f})_{-2}$. Thus by setting
\begin{eqnarray}
\nonumber g &\equiv& g_0 -\varphi f =
(\varphi_0-\varphi)f \\
h &\equiv& f-g = (1 - (\varphi_0-\varphi))f
\end{eqnarray}
we have the desired decomposition.
\end{proof}
\begin{corollary}
\label{[L2,T(f)]} Let $f$ be a smooth function on the circle with
$\hat{f}_{-2}=0$. Then on the set of finite energy vectors $$
[T(l_2),\tilde{T}(f)]=[\tilde{T}(l_2),\tilde{T}(f)].
$$
\end{corollary}
\begin{proof}
According to our previous lemma there exists a decomposition
\begin{equation}
f=g+h
\end{equation}
in which both $g$ and $h$ satisfy the conditions of Prop.\!
\ref{[L2,T(f)]loc}. Therefore on the set of finite energy vectors
\begin{eqnarray}
\nonumber [T(l_2),\tilde{T}(f)] &=& [T(l_2),\tilde{T}(g)] +
[T(l_2),\tilde{T}(h)] \\
&=&[T(l_2),T(g)] + [T(l_2),T(h)] = [T(l_2),T(f)].
\end{eqnarray}
\end{proof}
\begin{theorem}\label{theo:uniq2}
Let $(\A,V)$ be a \mo covariant local net and let $U$ and
$\tilde{U}$ be two strongly continuous projective unitary
representations of $\diff$ both extending $V$ and making $\A$ a
conformal covariant local net. Then it follows that $U=\tilde{U}$.
\end{theorem}
\begin{proof}
The two representations $U$ and $\tilde{U}$ give rise to two
stress-energy tensors $T$ and $\tilde{T}$ and also two
corresponding representations of the Virasoro algebra $\{L_n: n\in
\ZZ\}$ and $\{\tilde{L}_n: n\in \ZZ\}$.

As $(\hat{l_1})_{-2}=0$, by Corollary \ref{[L2,T(f)]} on the set
of finite energy vectors
\begin{equation}
[T(l_2),\tilde{T}(l_1)]=[\tilde{T}(l_2),\tilde{T}(l_1)].
\end{equation}
But on the finite energy vectors the left hand side equals to
$-[L_2,L_1]=-L_3$ while the right hand side is
$-[\tilde{L}_2,\tilde{L}_1]=-\tilde{L}_3$. Hence
$L_3=\tilde{L}_3$. Then by taking commutators with $L_1$ and
$L_{-1}$ we find that $L_n=\tilde{L}_n$ for every positive integer
$n$ and thus in fact for all integers since both representations
are unitary.
\end{proof}

\section{Conformal covariance of subnets}

\paragrf
Let $(\A,U)$ be diffeomorphism covariant local net on the circle.
As usual, we shall denote the stress-energy tensor of $U$ by $T$,
the corresponding representation of the Virasoro algebra by $\{L_n
: n\in\ZZ\}$, its central charge by $c$, the vacuum vector by
$\Omega$ and the dense set of finite energy vectors by $\fin$.
Suppose $\B \subset \A$ is a \mo covariant subnet of $\A$. In
particular it means that $U(\phi)\B U(\phi)^*=\B$ for every \mo
transformation $\phi\in\mob$. However, {\it a priori} we do not
know whether this invariance holds also for a diffeomorphism
$\phi\in \diff\setminus \mob$. In what follows we shall prove that
in fact this follows.

Just like in the previous section, two different arguments will be
presented. The first, again, will be a direct application of those
nonsmooth symmetries. The point here is the following. We know
that the subnet is \mo covariant. Suppose the subnet is strongly
additive. (In fact a weaker assumption will be enough.) Then the
unitary associated to a piecewise \mo transformation leaves the
subnet invariant as it is piecewise ``put'' together of \mo
transformations and the subnet is generated by the local algebras
associated to those pieces. This argument is conceptually very
clear and it is expressed in terms of local algebras. However it
needs some regularity assumption of the subnet: we need to know
that the subnet is generated by its local algebras associated to
intervals that cover the circle only after taking closures. The
other argument, just like in the previous section, does not need
any assumption. It is however less related to the idea of local
algebras; it is again about the restrictions on possible commutation
relations coming from the fact that the stress-energy tensor is a
dimension $2$ field.

\paragrf
We shall shortly discuss the first argument. We shall need the
notion of strong regularity.
\begin{definition}
A \mo covariant local net $(\C,V)$ is called {\bf strongly
$\mathbf n$-regular} if it remains regular even in
representations; i.e.\! if $\{I_k:k=1,..,n\}$ are the intervals
obtained by removing $n$ points of the circle, then
$$
\{\pi_{I_k}(\C(I_k)):k=1,..,n\}'' = \pi(\C)\equiv
\{\pi_{I}(\C(I)):I\in\I\}''
$$
for any locally normal representation $\pi$ of $\C$.
\end{definition}
Obviously a strongly additive net is also strongly regular (i.e.\!
strongly $n$-regular for every $n\in\NN^+$. If $\B\subset \A$ is a
subnet of $\A$ with vacuum vector $\Omega$, then its restriction
to $\overline{\B\Omega}$ is a \mo covariant local net. We shall
say that $\B$ is strongly $n$-regular if this \mo covariant local
net is so. Note that it follows that $\{\B(I_k)):k=1,..,n\}''=\B$
for any intervals $\{I_k:k=1,..,n\}$ obtained by removing $n$
points of the circle, since $\B$ can be considered to be a
representation of $\B|_{\overline{\B\Omega}}$.
\begin{lemma}
Let $(\A,U)$ be a conformal local net with stress-energy tensor
$T$ and $\B \subset \A$ a \mo covariant, strongly $4$-regular
subnet. Further, let $I_1,I_2\in\I$ two distant intervals. Then
$e^{isT(a^{I_1,I_2})} \,\B\, e^{-isT(a^{I_1,I_2})} = \B$ for every
$s\in\RR$.
\end{lemma}
\begin{proof}
It is a trivial observation. By Corollary
\ref{U(kappa)=exp(iT(a))} the adjoint action of
$e^{isT(a^{I_1,I_2})}$ is an implementation of
$\kappa^{I_1,I_2}_s$ in the net. So for example, this adjoint
action on $\B(I_1)\subset\A(I_1)$ coincides with that of
$U(\delta^{I_1}_s)$ and hence preserves it. Of course we may
repeat this argument with any of the four intervals obtained by
removing the endpoints of $I_1$ and $I_2$ from the circle. Thus by
the strong $4$-regularity of $\B$, the statement is true.
\end{proof}
\begin{theorem}\label{4reg->covsubnet}
Let $(\A,U)$ be a conformal local net with stress-energy tensor
$T$ and $\B \subset \A$ a \mo covariant, strongly $4$-regular
subnet. Then $U(\gamma)\,\B\, U(\gamma)^*=\B$.
\end{theorem}
\begin{proof}
We shall prove that if $f$ is a real function of finite
$\|\cdot\|_{\frac{3}{2}}$ norm then $e^{isT(f)}\, \B \,e^{-isT(f)}
=\B$ for every $s\in\RR$. This is clearly enough as $\diff$ is
generated by exponentials.

The statement is true for \mo vector fields and by our previous
lemma for the functions of the type $a^{I_1,I_2}$. Just like in
the proof of Theorem \ref{4reg->unique}, by an application of the
Trotter product-formula we get that then it is also true for any
sum of such functions and hence by Theorem \ref{pcwmgenerators}
for all piecewise \mo vector fields. But by Theorem
\ref{pcwmdensity} for any real continuous function $f$ with
$\|f\|_{\frac{3}{2}}<\infty$ there exists a sequence of piecewise
\mo vector fields $\{f_n:n\in\NN\}$ converging to $f$ in the
$\|\cdot\|_{\frac{3}{2}}$ norm. Then taking into account the
strong convergence of the corresponding unitaries (Prop.\!
\ref{f->T(f)continuity}) one can easily finish the proof.
\end{proof}
{\it Remark.} We could now go on and for example deduce that the
stress-energy tensor can be split  into a sum. However we shall
not do so as in any case we will obtain all these results by the
other argument.

\paragrf
We shall now proceed to the second proof. Since for our purposes
it is important to distinguish between the \mo symmetry and the
conformal symmetry, we shall denote the representation of the \mo
group whose extension is $U$ by $V$. As it was explained in the
preliminaries, the subspace $\K \equiv \H_\B =\overline{\B\Omega}$
is invariant for both $\B$ and $V$ and denoting the restrictions
by $\B\!|_\K$ and $V|_\K$ we have that $(\B\!|_\K,V|_\K)$ is a \mo
covariant net on the Hilbert space $\K$. Note that $\K\!_{f\!in}$,
the set of finite energy vectors for this \mo covariant net is
just simply $\K\cap\fin$.
\begin{definition}
For every $n\in\ZZ$ let $L_n^{\K} \equiv [\K] L_n|_{\K\!_{f\!in}}$
where $[\K]$ is the orthogonal projection from $\H_\A$ to $\K$.
\end{definition}
\begin{lemma} \label{prop.of.Ln^B}
$L_n^{\K}(\K\!_{f\!in}) \subset \K\!_{f\!in}$, $L_n^{\K} \subset
(L_{-n}^{\K})^*$, the operators $L^\K_1,L^\K_{-1}$ and $L^\K_0$
are generators of the representation $V|_\K$ (in particular
$L^\K_0$ is the self-adjoint generator of rotations for $V|_\K$)
and
$$
[L^{\K}_j,L^{\K}_n] = (j-n) L^{\K}_{n+j}
$$
for every $n\in\ZZ$ and $j\in\{-1,0,1\}$. Moreover,
$L^{\K}_k\Omega=0$ for every $k\geq-1$, and there exists an $r>0$
such that for every $\Phi\in\K\!_{f\!in}$ and $n\in\ZZ$
$$
\|L_n \Phi \| \leq r (1+|n|^{\frac{3}{2}}) \|(\mathbbm 1 +
L^{\K}_0) \Phi \|.
$$
\end{lemma}
\begin{proof}
It is a rather trivial consequence of the definition of
$L_n^{\K}$, the well-known energy bounds on the Virasoro operators
of $U$ and the fact that $L_0$ preserves $\K\!_{f\!in}$ (and thus
in particular $\|(\mathbbm 1 + L^{\K}_0)\Phi \| = \| (\mathbbm 1 +
L_0)\Phi\|$ for every $\Phi \in \K\!_{f\!in}$).
\end{proof}
\begin{proposition} \label{prop.of.TK(f)}
Let $f$ be a smooth function on the circle with Fourier
coefficients $\{\hat{f}_n: n\in\ZZ\}$. Then the operator
$T_\K(f)\equiv \overline{\sum_{n\in\ZZ}\hat{f}_n L^{\K}_n}$ is
well-defined (i.e.\! the sum converges on the finite energy
vectors and gives a closable operator), $T_\K(f)^* \supset
T_\K(\overline{f})$ and the set of smooth vectors $\K^\infty
\equiv \cap_{k\in\NN} \D((L^\K_0)^k)$ form an invariant core for
$T_\K(f)$. On this core $[L^\K_0,T_\K(f)] = T_\K(if')$ and for
every positive integer $k$ and smooth functions $f_1,..,f_k$ there
exists a $q>0$ such that
$$
\|T_\K(f_1)..T_\K(f_k)\Phi\| \leq q \|(\mathbbm 1 +
(L^\K_0)^k)\Phi\|
$$
for every $\Phi \in \K^\infty$. Thus $T_\K(f)^* =
T_\K(\overline{f})$.
\end{proposition}
\begin{proof}
Again, everything is essentially trivial (in the sense that it is
an exact copy of the standard arguments used for establishing the
properties of the stress-energy tensor, see the preliminary
Sect.\! \ref{sec:diffposenergy}); the proposition serves only to
clearly state what we know. In particular, for the polynomial
energy bound see Lemma \ref{poly.e.bound} while for the last
assertion the idea (explained in the mentioned section) of
Goodmann and Wallach of using Nelson's commutator theorem.
\end{proof}
{\it Remark.} Note that the subindex of $T_\K$ is indeed ``right''
in the sense that $T_\K$ was defined by the original stress-energy
tensor $T$ and by the subspace $\K$ and not (directly) by the
subnet $\B$.
\begin{proposition} \label{T(t)inA(I)}
If\, ${\rm Supp}(f) \subset I$ for a certain interval $I\in\I$
then the operator $T_\K(f)$ is affiliated to $\B\!|_\K(I)$.
\end{proposition}
\begin{proof}
Suppose indeed ${\rm Supp}(f) \subset I$ for a certain interval
$I\in\I$ and take an element $B \in \B\!|_\K(I^c)$ and other two
smooth elements $C,D\in \B\!|^\infty_\K(I^c)$. Clearly, due to its
definition and the existence of energy bounds, $T_\K(f)\Phi =
[\K]T(f)\Phi$ for every smooth vector $\Phi\in\K^\infty$. Then
identifying $\B\!|_\K(I^c)$ with $\B(I^c)$ by the isomorphism
given by the restriction we have that
\begin{eqnarray} \nonumber
\langle BC\Omega,\, T_\K(f)^*D\Omega \rangle &=&\langle
BC\Omega,\, T_\K(\overline{f})D\Omega\rangle
\\ \nonumber &=& \langle BC\Omega,\, T(\overline{f})D\Omega\rangle \\ \nonumber
&=& \langle BT(f)C\Omega,\, D\Omega\rangle \\ &=& \langle BT_\K(f)
C\Omega,\, D\Omega\rangle
\end{eqnarray}
where we have used that $T(f)$ is affiliated to $\A(I)\subset
\B(I^c)'$. By Theorem \ref{localcore} the set
$\B\!|^\infty_\K(I^c)\Omega$ forms a core for $L^\K_0$ and thus it
is also a core for $T_\K(f)$ and $T_\K(f)^*=T_\K(\overline{f})$.
Hence by the arbitrariness of $D\in \B\!|^\infty_\K(I^c)$ the
above equation shows that $BC\Omega$ is in the domain of
$(T_\K(f)^*)^* = T_\K(f)$ and $T_\K(f) BC \Omega = B T_\K(f) C
\Omega$. Then, this time by the arbitrariness of $C\in
\B\!|^\infty_\K(I^c)$, we can conclude that $B T_\K(f) \subset
T_\K(f) B$ and thus --- since $B$ is just any element of
$\B\!|_\K(I^c)$ --- that the operator $T_\K(f)$ is affiliated to
$\B\!|_\K(I^c)'= \B\!|_\K(I)$.
\end{proof}

\paragrf
What we have so far established is that starting from the original
stress-energy tensor by restricting and projecting onto $\K$ we
get a well-behaved local field. That is, a field which satisfies
linear energy bounds, for real test-functions gives an operator
which is essentially self-adjoint on the finite energy vectors and
for a local test function gives an operator which is affiliated to
the corresponding local algebra. Moreover, by the commutation
relation of its Fourier modes $\{L^\K_{n}\}$ with the generators
of the \mo symmetry, given in Lemma \ref{prop.of.Ln^B}, it is a
so-called ``dimension $2$'' field and its $0^{\rm th}$ Fourier
mode is exactly the conformal Hamiltonian. The main observation of
Sect.\! \ref{sec:uniq2}, which enabled the second proof given for
the uniqueness of the diffeomorphism symmetry, was that such a
field can have only a certain type of commutation relation. We can
now repeat the whole argument with almost no change.
\begin{proposition} \label{L2Ln:2}
By setting $c_\B \equiv 2\|L^\K_{-2}\Omega\|^2 \in \RR^+_0$ we
have $$[L^{\K}_2,L^{\K}_{-n}]\Omega = \left((2+n)L^\K_{-n+2} +
\frac{c_\B}{12} (2^3-2)\delta_{2,n}\right)\Omega.$$
\end{proposition}
\begin{proof}
For $n < 2$ the statement is trivial as then all operators
involved in it annihilate the vacuum. For $n=2$ the equation holds
by the definition of $c_\B$ and by the fact that due to the
uniqueness of the vacuum $[L^{\K}_2,L^{\K}_{-2}]\Omega$ must be a
scalar multiple of $\Omega$. Finally, the proof of Prop.\!
\ref{L2Ln} shows that it is also true for $m>2$ as we
carefully avoided the use of the Virasoro algebra relations there, 
too. All we did was the application of the commutation relations with
the generators of the \mo group and the exploitation of the fact
that ``any of the $L_k$ operators'' (with or without tilde or
index ``$\K$'') annihilates the vacuum for $k\geq -1$, and of
course these two properties we have in this case, too.
\end{proof}
\begin{corollary}
Recall that $l_n(z)=-iz^n$ and let $f$ be a smooth function on the
circle with Fourier coefficients $\{\hat{f}_n: n\in\ZZ\}$. Then
$$
[T_\K(il_2),T_\K(f)]\Omega= \left(i T_\K(il_2f'-il_2'f) +
\frac{c_\B}{2}\hat{f}_{-2} \mathbbm 1\right)\Omega.
$$
\end{corollary}
\begin{proposition}
Let $f$ be a smooth function on the circle and suppose that ${\rm
Supp}(f)\subset I$ for a certain $I\in \I$. Then on the set of
finite energy vectors
$$
[T_\K(il_2), T_\K(f)]= i T_\K(il_2f'-il_2'f) +
\frac{c_\B}{2}\hat{f}_{-2} \mathbbm 1.
$$
\end{proposition}
\begin{proof}
The argument is essentially a copy of the proof of Prop.\!
\ref{[L2,T(f)]loc} but for the convenience of the reader let us
briefly recall the main points. We can choose a bigger interval
$K\in\I, \overline{I}\subset K$ and a smooth function $\varphi$
which is constant $1$ on the interval $I$ and constant $0$ outside
$K$. Then $T_\K(l_2) = T_\K(\varphi l_2)+T_\K((1-\varphi)l_2)$,
and since $T_\K((1-\varphi)l_2)$ is affiliated to $\B\!|_\K(I^c)$
while $T_\K(f)$ is affiliated to $\B\!|_\K(I)$, on the smooth
energy vectors
\begin{eqnarray} \nonumber
[T_\K(l_2),T_\K(f)]&=&[T_\K(\varphi l_2),T_\K(f)]
+[T_\K((1-\varphi)l_2), T_\K(f)] \\ &=& [T_\K(\varphi
l_2),T_\K(f)]
\end{eqnarray}
Thus if $B'\in\B\!|^\infty_\K(K^c)$ then by our previous corollary
we find that
\begin{eqnarray}
\nonumber [T_\K(l_2),T_\K(f)] B'\Omega &=& [T_\K(\varphi
l_2),T_\K(f)] B' \Omega \\ \nonumber &=& B' [T_\K(\varphi
l_2),T_\K(f)] \Omega
\\ \nonumber &=& B' [T_\K(l_2),T_\K(f)] \Omega \\ \nonumber &=&
B' \left(i T_\K(l_2f'-l_2'f) -i\frac{c_\B}{2}\hat{f}_{-2} \mathbbm
1\right)\Omega \\ &=& \left(i T_\K( l_2f'- l_2'f) -i
\frac{c_\B}{2}\hat{f}_{-2} \mathbbm 1\right) B'\Omega
\end{eqnarray}
where we have used that since the support of the functions
$\varphi l_2$ and $(l_2f'-l_2'f)$ is contained in $K$,
the operators $T_\K(\varphi l_2)$ and $T_\K(l_2f'-l_2'f)$ are
affiliated to $\B\!|_\K(K)$ and so they commute with $B'$.

We have proved the desired equality on a dense set of vectors but
not exactly on $\fin$. To finish the demonstration we must take
limits, which can be done since the commutators admit a polynomial
energy bound and $\B\!|^\infty_\K(K^c)\Omega$ is a core for any
positive power of $L^\K_0$; see again Prop.\! \ref{[L2,T(f)]loc}
where exactly the same problem was dealt with.
\end{proof}
\begin{corollary}
Let $f$ be a smooth function on the circle (but not necessarily
local). Then on the set of finite energy vectors
$$
[T_\K(il_2),T_\K(f)]= i T_\K(il_2f'-il_2'f) +
\frac{c_\B}{2}\hat{f}_{-2} \mathbbm 1.
$$
\end{corollary}
\begin{proof}
Any smooth function $f$ can be decomposed into a sum $f=g+h$ where
$g$ and $h$ are smooth functions localized in two different
intervals. Thus applying our previous proposition to $g$ and $h$
we have the desired result.
\end{proof}
\begin{theorem}
\label{Automatic.Vir.rel.} $[L^\K_n, L^\K_m] = (n-m)L^\K_{n+m} +
\frac{c_\B}{12}(n^3-n)\delta_{-n,m} \mathbbm 1$ for every integer
$n,m\in\ZZ$.
\end{theorem}
\begin{proof}
For $n=-1,0,1$ the relation reduces to that which has been already
stated in Lemma \ref{prop.of.Ln^B}. For $n = 2$ it is exactly the
previous corollary that by setting $f=il_m$ shows that the claim
holds. For $n\geq 2$ we proceed by induction. Suppose the
commutation relation holds for a certain $n=k\geq 2$. Then for
$n=k+1$ we have $(k-1)[L^\K_n, L^\K_m] = $
\begin{eqnarray} \nonumber
&=& [[L^\K_k,L^\K_1], L^\K_m]
\,=\,[L^\K_k,[L^\K_1,L^\K_m]] - [L^\K_1,[L^\K_k, L^\K_m]] \\
\nonumber &=& (1-m)[L^\K_k, L^\K_{m+1}] - (k-m)[L^\K_1, L^\K_{k+m}] \\
\nonumber &=& (k^2-m k+m-1)L^\K_{k+m+1}
+(1-m)\frac{c_\B}{12}(k^3-k)\delta_{-k,(m+1)} \mathbbm 1 \\ &=&
(k-1)\left((n-m)L^\K_{n+m} + \frac{c_\B}{12}(n^3-n)\delta_{-n,m}
\mathbbm 1\right)
\end{eqnarray}
where in the first line we used the Jacobi-identity, in the second
and the third line we used the inductive assumption, the known
commutation relation with $L^\K_1$ and the fact that the identity
operator falls out from the commutator, and finally, in the last line we
collected the terms taking account of the Kronecker-delta
$\delta_{-k,(m+1)}$ which permits us to exchange $(1-m)$
with $(k+2)$. Dividing by $(k-1)$ (remember that $k\geq 2$ so
indeed we do not divide by zero) we finish the inductive argument
and thus conclude the statement for all $n \geq -1$. Finally, for
negative values, the commutation relation can be easily justified
by taking adjoints and using that $L^\K_{-n} \subset
(L^\K_{n})^*$.
\end{proof}
\begin{corollary} \label{subnet.conf.struct.}
The operators $\{L^\K_n: n\in\ZZ\}$ form a positive energy
representation of the Virasoro algebra. The positive energy
representation $U_\K$ of $\diff$ given by $\{L^\K_n: n\in\ZZ\}$ is
an extension of $V|_\K$ which makes $\B\!|_\K$ diffeomorphism
covariant and so the formula
$$
I\mapsto \T\!_{\B\!|_\K}(I)\equiv \{U_\K({\rm Exp}(f)) :
f_{I^c}=0\}'' \;\;\; (I\in\I)
$$
defines a \mo covariant subnet of $(\B\!|_\K,V|_\K)$ whose
restriction to its vacuum Hilbert space is a Virasoro net.
\end{corollary}
\begin{proof}
The fact that $\{L^\K_n: n\in\ZZ\}$ is a positive energy
representation of the Virasoro algebra and so that it integrates
to a positive energy representation $U_\K$ of $\diff$ has been
already stated by the last theorem. Evidently $U_\K$ will be an
extension of $V|_\K$. Moreover, by Prop.\! \ref{T(t)inA(I)} it is
also clear that the given formula indeed defines a \mo covariant
subnet of $(\B\!|_\K,V|_\K)$ whose restriction to its vacuum
Hilbert space is a Virasoro net.

As it was explained after Def.\! \ref{diffcov:def}, all we have to
show for diffeomorphism covariance is compatibility. Of course by
the mentioned proposition, the compatibility of the action of
exponentials follow. However, although every diffeomorphism is a
finite product of exponentials, to the knowledge of the author the
question whether every local diffeomorphism can be written as a
product of exponentials of vector fields localized in the {\it
same} interval is unsolved. Fortunately, we do not really need to
know this; the proof can be also concluded by \cite[Prop.\!
3.7]{Carpi2}.
\end{proof}
A finite energy vector in $\K$ is in particular also a finite
energy vector in $\H_\A$. Therefore the expression $L_2
L^\K_n\Omega$ is well defined.
\begin{lemma}
$ L_2 L^\K_{-n}\Omega = L^\K_2 L^\K_{-n}\Omega$ for every
$n\in\ZZ$.
\end{lemma}
\begin{proof}
For $n < 2$ the statement is trivially true. Let us consider the
case $n=2$. Then both sides are vectors of zero energy, so by the
uniqueness of the vacuum they are both multiples of $\Omega$. This
multiple for the left hand side is equal to the scalar product
$\langle\Omega,\,L_2L^\K_{-2}\Omega\rangle$ in which $L_2$ can be
safely replaced by $L^\K_2$ as both $\Omega$ and $L^\K_{-2}\Omega$
are vectors of $\K$. This shows the equality of the two sides.

Finally, for $n\geq2$ we can proceed by induction using that on
vectors of $\K$ we have the equality $L_{-1}=L^\K_{-1}$, exactly
like in Prop.\! \ref{L2Ln} (and Prop.\! \ref{L2Ln:2}, so note that
in fact this is third time we employ the same trick).
\end{proof}
\begin{lemma}
$L_m L^\K_{-n}\Omega = L^\K_m L^\K_{-n}\Omega$ for every $n\in\ZZ$
and $m\in\NN$.
\end{lemma}
\begin{proof}\label{LmLn:2}
This time we shall do induction by $m$. For $m=0,1$ the statement
is trivially true, while for $m=2$ it is exactly the previous
lemma. Finally, for $m\geq 2$ exactly the same inductive argument
that was employed in the proof of Theorem \ref{Automatic.Vir.rel.}
applies.
\end{proof}

\paragrf
With the above lemma we have just begun to return to the vacuum
Hilbert space $\H_\A$ of the original net $\A$. For any interval
$I\in\I$ the restriction map from $\B(I)$ to $\B\!|_\K(I)$ is an
isomorphism. Hence locally the restriction map can be inverted and
thus one can lift the subnet $\T\!_{\B\!|_\K} \subset \B\!|_\K$
defined in Corollary \ref{subnet.conf.struct.} to the vacuum
Hilbert space of $\A$. Denoting this lift by $\T\!_\B$ we have the
following chain of \mo covariant subnets: $\T\!_\B\subset
\B\subset \A$. Of course the restriction of the lifted subnet to
its vacuum Hilbert space remains unchanged. Thus we can still look
upon $\T\!_\B$ as a locally normal representation of a Virasoro
net.

As the following regards only normal representations of a Virasoro
net, it will be stated in general. If $\pi$ is a locally normal
representation of the Virasoro net $(\A_{{\rm Vir}\!,c},U_c)$ with
central charge $c\in\RR^+_0$ then by \cite[Theorem 6]{D'AnFreKos}
there exists a unique strongly continuous representation
$\tilde{V}_\pi$ of $\widetilde{\mob}$ on the Hilbert space of the
representation such that $\tilde{V}_\pi(\widetilde{\mob})\subset
\pi(\A_{{\rm  Vir}\!,c})$ and
\begin{equation}
{\rm Ad}(\tilde{V}(\tilde{\gamma}))\circ \pi_I =
\pi_{\gamma(I)}\circ {\rm Ad}(U(\gamma))
\end{equation}
for all $I\in\I$ and $\tilde{\gamma}\in\widetilde{\mob}$ with
$\gamma\equiv p(\tilde{\gamma})$ where
$p:\widetilde{\mob}\rightarrow\mob$ is the covering map. In a later 
section of this thesis it will be actually proved (Theorem
\ref{mainresult}) that $\tilde{V}_\pi$ is always a positive energy
representation.
\begin{proposition}
Let $\pi$ be a locally normal representation of the Virasoro net
$(\A_{{\rm Vir}\!,c},U_c)$ with central charge $c\in\RR^+_0$. Then
if $\pi(\A_{{\rm Vir}\!,c})$ is a factor then in fact $\pi$ is
unitary equivalent to a multiple of one of those irreducible
representations which we get by integrating a positive energy
unitary lowest weight representation of the Virasoro algebra
corresponding to the same central charge.
\end{proposition}
\begin{proof}
By \cite[Lemma 3]{D'AnFreKos}, taking account of the factoriality
of $\pi$ and the fact that $\diff$ is generated by exponentials of
local vector fields (localized in different intervals), it is
clear that there exists a unique strongly continuous projective
representation $U_\pi$ of $\diff$ such that for all $I\in\I$ and
smooth function $f$ localized in $I$
\begin{equation}
\pi_I(U({\rm Exp}(f))) = U_\pi({\rm Exp}(f)).
\end{equation}
By Theorem \ref{mainresult} the representation $U_\pi$ is of
positive energy (although this Theorem is presented later in this
thesis, note that in that section we shall make no use of results
made in this section), and by Theorem \ref{pos.diffrep=+irr.} it
follows that $U_\pi$ is a direct sum of irreducible
representations. As $\pi(\A_{{\rm Vir}\!,c}) = U_\pi(\diff)$ is a
factor, $U_\pi$ can contain only one type of irreducible
representation. However, by \cite[Theorem A.1]{Carpi2}, a positive
energy irreducible representation of $\diff$ is unitary equivalent
to one of those we get by integrating a positive energy unitary
lowest weight representation of the Virasoro algebra. Finally, the
central charge of $U_\pi$ must be equal to $c$ as it can be
calculated by local commutation relations of the stress-energy
tensor, see Eq.\! (\ref{[T(f),T(g)]}).
\end{proof}
For the following proposition note, that if $\pi$ is a locally
normal representation of a \mo covariant local net $(\C,W)$, and
$I\in\I$ then there exists a unitary $X_I$ such that $\pi_I (C)
=X_I CX_I^*$ for all $C\in\C(I)$. The operator $X_I TX_I^*$ for a
closed operator $T$ affiliated to $\C(I)$ is independent from the
choice of $X_I$. Thus the expression $\pi_I(T)\equiv X_I T X_I^*$
(with $T$ not necessarily bounded!) is well defined.
\begin{proposition}
\label{pi(T(f))} Let $\pi$ be a locally normal representation of
the Virasoro net $(\A_{{\rm Vir}\!,c}U_c)$ with central charge
$c\in\RR^+_0$ and stress-energy tensor $T_c$. Further, let
$\tilde{V}_\pi$ be the unique strongly continuous inner
implementation of the \mo symmetry in the representation $\pi$
with $L_0^\pi$ being the self-adjoint generator of rotations. Then
$\D(\pi_I(T_c(f)))\subset \D(L_0^\pi)$ and
$\pi_I(T_c(f))\D^\infty_\pi=\D^\infty_\pi$ for every $I\in\I$ and
smooth function $f$ with ${\rm Supp}(f) \subset I$, where
$\D^\infty_\pi\equiv\cap_{n\in\NN} \D((L_0^\pi)^n)$. Moreover, on
this invariant set $[L^\pi_0,\pi_I(T_c(f))]=\pi_I(T_c(if'))$ and
$$
\|\pi_I(T_c(f))\Phi\| \leq \sqrt{1+\frac{c}{12}}
\|f\|_{\frac{3}{2}}\, \|(\mathbbm 1+ L_0^\pi)\Phi\|
$$
for every $\Phi \in \D(L_0^\pi)$.
\end{proposition}
\begin{proof}
The Hilbert space of $\A_{{\rm Vir},c}$ is separable and so we can
safely consider the direct integral decomposition (see more on at
the comments before Theorem \ref{mainresult} and in \cite[Chapter
14.]{kadison} about generalities regarding direct integral
decompositions) along the center of $\pi(\A_{{\rm Vir},c})$:
\begin{equation}
\pi(\A_{{\rm Vir},c})=\int_{\!\! X}^\oplus \M_x\, d\mu_x
\end{equation}
where $\M_x$ is a von Neumann factor for almost every $x\in X$.
Correspondingly, by \cite[Lemma 8.3.1]{dixmier} (see again the
mentioned comments) the representation $\pi$ can be decomposed
into the direct integral of representations $\pi=\int_X^\oplus
\pi_x d\mu_x$ where $\pi_x$ is a locally normal representation of
our Virasoro net with $\pi_x(\A_{{\rm Vir},c}) = \M_x$ for almost
every $x\in X$. (There is a little abuse of the word
``representation'' as in reality it covers an {\it infinite
collection} of representations. However, evidently one can choose
a countable set of intervals such that a collection of normal
representations associated to these intervals and compatible with
the isotony {\it with respect to a certain, countable set of
inclusions}, uniquely extends to a locally normal representation
of the net. Thus the cited Lemma, which is stated for a single
representation, can be indeed applied.) So the previous
proposition applies, and thus we can write the representation as a
direct integral of irreducible ones. But the irreducible
representations we know well (they are exactly the ones coming
from irreducible positive energy representations of $\diff$, see
the preliminary Sect.\! \ref{sec:diffposenergy}) and so in each
irreducible representation the stated properties are trivially
satisfied. The rest is then an exercise.
\end{proof}
Still in the same setting, with $\pi$ being a locally normal
representation $(\A_{{\rm Vir},c},U_c)$, we shall
introduce the corresponding stress-energy tensor. In what follows
we shall denote by $\D_\pi^\infty$ the set of smooth vectors
defined by the Hamiltonian of the unique inner implementation
$\tilde{V}_\pi$ of the \mo symmetry in $\pi$.
\begin{definition}
For a smooth function $f$, let $T_\pi(f)$ be the closure of
$$
\pi_{I_1}(T(\varphi_1 f))|_{\D_\pi^\infty} +
\pi_{I_1}(T_c(\varphi_1 f))|_{\D_\pi^\infty}
$$
where $I_1,I_2\in\I$ and $1 = \varphi_1 + \varphi_2$ is a
partition of the unity with ${\rm Supp}(\varphi_j)\subset I_j$
$(j=1,2)$.
\end{definition}
\begin{lemma}
The above definition is good, i.e.\! $T_\pi(f)$ is independent of
the choice of the partition of the unity and the operator whose
closure is considered is indeed closable.
$(T_\pi(f))^*=T_\pi(\overline{f})$, the set $\D_\pi^\infty$ is
invariant for $T(f)$ and $T_\pi|_{\D_\pi^\infty}$ is an operator
valued linear functional. If $f$ is localized in $I$ then
$T_\pi(f)=\pi_I(T_c(f))$ and hence $T_\pi(f)$ is affiliated to
$\pi_I(\A_{{\rm Vir},c}(I))$. Moreover,  for every positive
integer $k$ and smooth functions $f_1,..,f_k$ there exists a $q>0$
such that
$$
\|T_\pi(f_1)..T_\pi(f_k)\Psi\| \leq q \|(\mathbbm 1 +
(L^\pi_0))^k\Psi\|
$$
for every $\Psi \in \D_\pi^\infty$.
\end{lemma}
\begin{proof}
Let us consider another partition of the unity $1 =
\tilde{\varphi}_1 + \tilde{\varphi}_2$ with ${\rm
Supp}(\tilde{\varphi}_j)\subset K_j$ $(j=1,2)$ where
$K_1,K_2\in\I$. Then, with everything meant on $\D_\pi^\infty$ and
for shortness setting $T_I(f)\equiv \pi_I(T_c(f))$, we have that
$T_{I_1}(\varphi_1 f) + T_{I_1}(\varphi_1 f) =$
\begin{eqnarray} \nonumber
&=& T_{I_1}(\tilde{\varphi}_1\varphi_1 f) +
T_{I_1}(\tilde{\varphi}_2\varphi_1 f) +
T_{I_2}(\tilde{\varphi}_1\varphi_2 f) +
T_{I_2}(\tilde{\varphi}_2\varphi_2 f) \\ \nonumber &=& T_{I_1\cap
K_1}(\tilde{\varphi}_1\varphi_1 f) + T_{I_1\cap
K_2}(\tilde{\varphi}_2\varphi_1 f) + T_{I_2\cap
K_1}(\tilde{\varphi}_1\varphi_2 f) + T_{I_2\cap
K_2}(\tilde{\varphi}_2\varphi_2 f) \\ \nonumber &=&
T_{K_1}(\tilde{\varphi}_1\varphi_1 f) +
T_{K_2}(\tilde{\varphi}_2\varphi_1 f) +
T_{K_1}(\tilde{\varphi}_1\varphi_2 f) +
T_{K_2}(\tilde{\varphi}_2\varphi_2 f) \\ \nonumber &=&
T_{K_1}(\tilde{\varphi}_1 f) + T_{K_2}(\tilde{\varphi}_2 f).
\end{eqnarray}
The rest of the statement is rather trivial, it is just a copy of
the standard arguments used to derive properties of the
stress-energy tensor, e.g.\! Lemma \ref{poly.e.bound} for the
polynomial energy bound.
\end{proof}
Finally, we shall make some more affirmations. Like the ones just
stated, also these are all essentially trivial: they say that
$T_\pi$ has the usual properties of an energy stress-tensor. (Even
if it does not necessarily integrate to a representation of
$\diff$, but only to its cover.)
\begin{lemma}
For any smooth function $f$ and
$\tilde{\varphi}\in\widetilde{\mob}$,
\begin{equation}
\tilde{V}_\pi(\tilde{\varphi}) T_\pi(f)
\tilde{V}_\pi(\tilde{\varphi})^* = T_\pi(\varphi_*f)
\end{equation}
where the star stands for action on vector fields.
\end{lemma}
\begin{proof} The equality can be easily checked locally (using
that by Corollary \ref{UmobD^infty=D^infty} the smooth vectors for
the rotations form an invariant set for the positive energy
representation of $\widetilde{\mob}$) and then globally using
partitions of unities.
\end{proof}
\begin{lemma} \label{[Tpi,Tpi]}
On the invariant core of $\D_\pi^\infty$
\begin{equation}
[T_\pi(f),T_\pi(g)] = i T_\pi(fg'-f'g) + i \frac{c}{12}(f,g)
\mathbbm 1
\end{equation}
where $(f,g)\equiv
-\frac{1}{2\pi}\int_0^{2\pi}\overline{f(e^{i\theta})}(g'(e^{i\theta})+g'''(e^{i\theta}))
d\theta$.
\end{lemma}
\begin{proof}
Again the equality can be easily checked locally and then globally
using partitions of unities.
\end{proof}
\begin{lemma} With $f$ being a smooth function,
$g$ a real smooth function, setting $\gamma \equiv{\rm Exp}(g)$ we
have
\begin{eqnarray*}
e^{iT(g)} T_\pi(f) e^{-iT(g)} &=& T_\pi(\gamma_*f) + r\mathbbm 1
\\
e^{iT(g)} \D_\pi^\infty &=& \D_\pi^\infty
\end{eqnarray*}
where the star stands for action on vector fields and $r$ is some
complex number. Moreover $\tilde{V}_\pi(\widetilde{{\rm
Exp}}(g))=e^{iT_\pi(g)}$ for every real \mo vector field $g$.
\end{lemma}
\begin{proof}
Again, the first property can be easily checked locally and then
extended to global result using partitions of unities and the
Trotter product-formula (there is no problem with convergence as
we have a common core). The fact that $\tilde{V}_\pi$ can be
obtained by integrating $T_\pi$ follows from the uniqueness of the
inner implementation of the \mo symmetry. Finally the fact that
$e^{iT(g)} \D_\pi^\infty = \D_\pi^\infty$ can be shown exactly as
in the proof of Corollary \ref{UD^infty=D^infty}.
\end{proof}

\paragrf
After this somewhat longer {\it intermezzo} let us return to our
concrete situation. Recall that we had a \mo covariant local net
$(\A,V)$ with conformal Hamiltonian $L_0$ and a \mo covariant
subnet $\B\subset \A$, and that there was a representation $U$ of
$\diff$ given making $(\A,V)$ diffeomorphism covariant. Recall
further how the subnet $\T\!_\B$ was introduced. By our previous
considerations, since this subnet when restricted to its vacuum
Hilbert space gives a Virasoro net, we can introduce the lifted
stress-energy tensor $T^\B$ which is affiliated to $\T\!_\B$ and
has all the properties that were discussed.
\begin{proposition}
\label{prop.of.TB} $\D(T^\B(f)) \subset \D(L_0)$ and
$T^\B(f)(\D^\infty)\subset \D^\infty$ for every smooth function
$f$. Further, there exists a constant $r>0$ independent of $f$
such that
$$
\|T^\B(f))\Phi\| \leq r \|f\|_{\frac{3}{2}}\, \|(\mathbbm 1+
L_0)\Phi\|
$$
and
$$
\|L_0 T^\B(f))\Phi\| \leq r
(\|f\|_{\frac{3}{2}}+\|f'\|_{\frac{3}{2}}) \|(\mathbbm 1+
L^2_0)\Phi\|
$$
for every $\Phi \in \D^\infty$. Moreover $\fin$ is a core for the
closed operator $T^\B(f)$.
\end{proposition}
\begin{proof}
We know that $(\A,V)$ is a \mo covariant local net and
$\T\!_\B\subset \A$ is a \mo covariant subnet. Thus by the result of
\cite{koester02}, there exists a unique pair of commuting positive energy 
representations $(\tilde{V}^{\T\!_\B},\tilde{V}^{\T_{\!\B}'})$ of
$\widetilde{\mob}$ on the vacuum Hilbert space of $\A$ such that
$\tilde{V}^{\T\!_\B}(\widetilde{\mob})\subset \T\!_\B$ and
\begin{itemize}
\item $V(p(\tilde{\gamma}))=
\tilde{V}^{\T\!_\B}(\tilde{\gamma})\tilde{V}^{\T_{\!\B}'}(\tilde{\gamma})$,
\item ${\rm
Ad}\left(\tilde{V}^{\T\!_\B}(\tilde{\gamma})\right)|_{\T\!_\B}
={\rm Ad}\left(V(p(\tilde{\gamma}))\right)|_{\T\!_\B}$
\end{itemize}
for all $\gamma \in \widetilde{\mob}$ with
$p:\widetilde{\mob}\rightarrow \mob$ being the covering map. Note
that $\tilde{V}^{\T\!_\B}$ is the representation that was denoted
by $\tilde{V}_\pi$ in the general setting.

Let us denote by $L_0^{\T_{\!\B}}$ and $L_0^{\T_{\!\B}'}$ the two
commuting conformal Hamiltonians corresponding to
$\tilde{V}^{\T_{\!\B}}$ and $\tilde{V}^{\T_{\!\B}'}$,
respectively. Then $L_0^{\T_{\!\B}}+L_0^{\T_{\!\B}'} = L_0$ and
since we deal with commuting positive operators  $0\leq
(L_0^{\T_{\!\B}})^k \leq (L_0)^k$. Thus, taking account of the
properties of the stress-energy listed before, we have immediately
half of the statement. (Smoothness for $L_0$ imply smoothness for
$L_0^{\T_{\!\B}}$, etc.)

Using that $L_0^{\T_{\!\B}'}$ is positive, commutes (in the sense
of its spectral resolution) with the closed operator $T^\B(f)$,
and that $0\leq (L_0^{\T_{\!\B}})^k \leq (L_0)^k$ for any positive
exponent $k$, it is easy to show that if $\Phi\in\D^\infty$ then
$T^\B_I(f)\Phi \in \D^\infty$ as in fact
\begin{equation}
L_0 T^\B(f)\Phi = T^\B(if')\Phi + T^\B(f)L_0\Phi
\end{equation}
using that $[L_0^{\T_{\!\B}},T^\B(f)] = T^\B(if')$, which was
already shown in the general setting. The second inequality is
obtained by a repeated use of the already proved energy bound,
just like in Lemma \ref{poly.e.bound}. Finally, by Nelson's
commutator theorem (as in the proof of Goodman and Wallach, see
also in the preliminary Sect.\! \ref{sec:diffposenergy}) any core
for $L_0^{\T_{\!\B}}$ is also a core for $T^\B(f)$ and thus the
proof is finished.
\end{proof}
\begin{lemma}
$T^\B$ is a \mo covariant field: $$ U(\varphi) T^\B(f)
U(\varphi)^* = T^\B(\varphi_* f)$$ for all $\varphi \in \mob$ and
$f$ smooth function, where $\varphi_*$ is the action of $\varphi$
on vector fields.
\end{lemma}
\begin{proof}
Follows from the factorization of the \mo symmetry into two parts
in which one commutes with $T^\B$, since for the other, inner part
we have already shown this in the general setting.
\end{proof}
For the next proposition recall that $l_n$ $(n\in\ZZ)$ is the
function defined by the formula $l_n(z)\equiv -iz^n$.
\begin{proposition}
Let $I\in\I$ be an interval and $f$ a smooth function localized in
$I$. Then on $\D^\infty$
$$
[T(l_n),T^\B(f)] = [T^\B(l_n),T^\B(f)]
$$
for every $n\in\NN$.
\end{proposition}
\begin{proof}
On the vacuum vector $T^\B(f)\Omega = \sum_{m\in\ZZ} \hat{f}_m
L^\K_m\Omega$ and so by Lemma \ref{LmLn:2} the claim is justified.
Then by exactly the same argument as in the proof of Prop.\!
\ref{[L2,T(f)]loc}, by appealing to locality, the derived energy
bounds and formulae (with which we can show that expressions like
$T(l_n)T^\B(f)$ admit a quadratic energy bound) and the fact that
$\A^\infty(K)\Omega$ is a core for $L_0^2$ the proof is
straightforward. 
\end{proof}
\begin{corollary}
For every smooth function $f,g$, on $\D^\infty$
$$
[T(g),T^\B(f)] = [T^\B(g),T^\B(f)].
$$
\end{corollary}
\begin{proof}
It follows from the previous proposition by using partitions,
taking adjoints and summing; due to the energy bounds all such
calculations are justified.
\end{proof}
\begin{proposition}
\label{UTBU*} Let $f,g$ be two real smooth functions. Then
$$
U({\rm Exp}(g)) \, T^\B(f)\, U({\rm Exp}(g))^* = e^{iT^\B(g)}\,
T^\B(f)\, e^{-iT^\B(g)}.
$$
\end{proposition}
\begin{proof}
Let $\Phi\in\D^\infty$ and consider the vector valued function on
the real line given by the formula
\begin{equation}
\Psi(t)\equiv U({\rm Exp}(tg)) \, e^{-it T^\B(g)}\, T^\B(f)\,
e^{it T^\B(g)}  \, U({\rm Exp}(tg))^* \Phi
\end{equation}
which is indeed well defined: by Corollary \ref{UD^infty=D^infty}
$U(\gamma)\D^\infty =\D^\infty$ and we know that $e^{-itT^\B(g)}\,
T^\B(f) \, e^{itT^\B(g)} = T^\B(f_t)$ plus a certain number times the
identity where $f_t\equiv{\rm Exp}(-tg)_* f$. In fact, it is not
hard to see then that this function is differentiable and
\begin{equation}
\frac{d}{dt}\Psi(t)= i\, U({\rm Exp}(tg))\,
[T(g)-T^\B(g),T^\B(f_t)]\, U({\rm Exp}(tg))^*\Phi
\end{equation}
which by the previous corollary is zero. Thus $\Psi(t)=\Psi(0)$
and hence the equality is justified as $\D^\infty$ is a core.
\end{proof}
\begin{proposition}
The \mo covariant subnet $\B$ is diffeomorphism covariant:
$U(\gamma) \B U(\gamma)^* = \B$ for all $\gamma\in\diff$. More
specifically, if $f$ is a real smooth function then
$$
{\rm Ad}(U({\rm Exp}(f)))|_\B = {\rm Ad}(e^{iT^\B(f)})|_\B,
$$
the operator $T(f)-T^\B(f)$ is essentially self-adjoint on $\fin$
and its closure is affiliated to $\B'$ and if further $f$ is
supported in the interval $I\in\I$ then it is in fact affiliated
to $\B' \cap \A(I)$.
\end{proposition}
\begin{proof}
We shall use the trick of moving the operators with unknown
commutation relation into disjoint intervals. Let $I,K\in\I$ be
two intervals, $B\in\B(I)$ and $g$ a real smooth function with
support in $K$. It is clear that we can find an interval $F\in\I$
containing $I$ and another smooth function $h$ supported in $F$
such that ${\rm Exp}(h)(I)$ is disjoint from $K$. Further, we can
consider $\B$ as a locally normal representation $\pi$ of
$\B\!|_\K$. Then
\begin{eqnarray}\nonumber
e^{iT^\B(h)}\,B\,e^{-iT^\B(h)} &=& \pi_F\!\left(U_\K({\rm
Exp}(h))\,B|_\K\, U_\K({\rm Exp}(h))^* \right) \\ \nonumber &=&
\pi_{{\rm Exp}(h)(I)}\!\left(U_\K({\rm Exp}(h))\,B|_\K\, U_\K({\rm
Exp}(h))^* \right)\\ &\in& \B({\rm Exp}(h)(I))
\end{eqnarray}
where we have used that by Corollary \ref{subnet.conf.struct.},
$(\B\!|_\K,U_\K)$ is a conformal net. Therefore, setting
$\gamma\equiv {\rm Exp}(g)$
\begin{equation}
e^{iT^\B(h)}\,B\,e^{-iT^\B(h)} = U(\gamma)^* \,
e^{iT^\B(h)}\,B\,e^{-iT^\B(h)}\,  U(\gamma)
\end{equation}
since $U(\gamma)\in \A(K)$ and $K\cap {\rm Exp}(h)(I) =\emptyset$.
So now we can calculate the action of $U(\gamma)$ on the element
$B$. We have that $U(\gamma)\,B\,U(\gamma)^*=$
\begin{eqnarray}\nonumber
 &=& U(\gamma)\,  e^{-iT^\B(h)}\,\left(e^{iT^\B(h)}B
e^{-iT^\B(h)}\right)\,e^{iT^\B(h)}\, U(\gamma)^*  \\ \nonumber &=&
U(\gamma) \, e^{-iT^\B(h)} \, U(\gamma)^*\, \left(e^{iT^\B(h)}\, B\,
e^{-iT^\B(h)}\right)\,U(\gamma)\, e^{iT^\B(h)} \, U(\gamma)^*
\\ &=& {\rm Ad}\left(U(\gamma)\, 
e^{-iT^\B(h)} \, U(\gamma)^*\right) (e^{iT^\B(h)}\, B\,  e^{-iT^\B(h)})
\end{eqnarray}
where we have used our previous equation. We shall now separately
calculate $U(\gamma)\, e^{-iT^\B(h)} \, U(\gamma)^*$. Using Prop.\!
\ref{UTBU*}
\begin{eqnarray}\nonumber
U(\gamma) \, e^{-iT^\B(h)} \, U(\gamma)^* &=& U({\rm Exp}(g))\,
e^{iT^\B(-h)}\, U({\rm Exp}(g)) \\ &=& e^{iT^\B(g)} \, e^{iT^\B(-h)}
\, e^{iT^\B(g)}.
\end{eqnarray}
So finally, putting everything together we get that
\begin{eqnarray}\nonumber
U(\gamma\, B\, U(\gamma)^* &=& {\rm Ad}\left(e^{iT^\B(g)}
\, e^{iT^\B(-h)}\, 
e^{iT^\B(g)} \right) (e^{iT^\B(h)}\, B\,  e^{-iT^\B(h)}) \\
&=& e^{iT^\B(g)} \, B\,  e^{-iT^\B(g)}
\end{eqnarray}
where we have also used that $e^{iT^\B(h)}\,B\,e^{-iT^\B(h)}$ and
$e^{iT^\B(g)}$ are localized in disjoint intervals. Thus we have
proved the desired formula for some $g$ and $B$: namely when both
are localized {\it somewhere} (but not necessarily in the same
interval). Of course, as $\B$ is generated by local elements, the
localization requirement on $B$ can be dropped. Finally, if $g$ is
not local then by using a partition and the Trotter
product-formula (again, there is no problem with convergence: we
have a common core) we can easily justify the equality.

The rest of the statement follows easily. Since $T^\B$ is
affiliated to $\B$, by our formula that we have just derived,
$U(\gamma)\B U(\gamma)^* = \B$ for every exponential and hence for
every $\gamma\in\diff$ as $\diff$ is generated by exponentials.

Finally, the fact that $T(f)-T^\B(f)$ is essentially self-adjoint
on $\fin$ again follows from the energy bounds, the commutator
relations (with $L_0$) and Nelson's commutator theorem. By using
our formula for the action of $U(\gamma)$ on $\B$ which we have
just verified and the Trotter product-formula we get that its
closure is affiliated to $\B'$. If further $f$ is localized in the
interval $I$, then both $T(f)$ and $T^\B(f)$ are affiliated in
$\A(I)$ and its easy to see that it follows that (the closure of)
their difference, too, is affiliated to the same algebra (again,
using that $\fin$ is a common core).
\end{proof}
What we have proved is that given a subnet $\B\subset \A$, the
stress-energy tensor can be ``split'' into a sum of two commuting
parts. One of them is the ``part of the stress-energy tensor in
$\B$''. The remaining part belongs to the so-called {\it coset} of
$\B$, and it is also positive as $T(1)-T^\B(1) = L_0 -
L_0^{\T_{\!\B}} = L_0^{\T_{\!\B}'} \geq 0$.

One can also define the Virasoro operators $L^\B_n\equiv
T^\B(l_n)|_{\fin}$ and reformulate what we have obtained in terms
of ``splitting the representation of the Virasoro algebra''. The
exponential relation of Prop.\! \ref{UTBU*} and the integrated
commutation relation of Lemma \ref{[Tpi,Tpi]} is then turned into
the formula
\begin{equation}
[L_n,L^\B_m] = [L^\B_n,L^\B_m] = (n-m) L^\B_{n+m} +
\frac{c_\B}{12} (n^3-n)\delta_{-n,m} \mathbbm 1
\end{equation}
Although everything has been stated and proved, it is worth to
recollect (and to reformulate) the main results of this section.
We have shown the followings.
\begin{theorem}\label{auto.conf.subnet.2a}
Let $(\A,U)$ be a diffeomorphism covariant local net on $S^1$ and
$\B\subset \A$ a \mo covariant subnet. Then, setting
$\K\equiv{\overline{\B\Omega}}$
\begin{itemize}
\item $U(\gamma)\B\,U\,(\gamma)^*=\B$ for all $\gamma\in\diff$;

\item $(\B\!|_\K, U(\cdot|_\mob)|_\K)$ admits a conformal
structure.
\end{itemize}
\end{theorem}
\begin{definition}
Let $(\A,V)$ be a \mo covariant local net with $\fin$ denoting the
set of finite energy vectors and $L_j\;\;(j=-1,0,1)$ the standard
generators of $V$. The operators $\{\Upsilon_n:n\in\ZZ\}$ with
domain $\fin$ are said to be {\bf Fourier modes of a dimension
$\mathbf 2$ \mo covariant field}, if
$$
[L_j, \Upsilon_n] = (j-n)\Upsilon_{j+n} \;\;\;\; (j\in\{-1,0,1\},
n\in\ZZ).
$$
Moreover $\Upsilon$ is said to be {\bf hermitian}, if
$\Upsilon_n^* \supset \Upsilon_{-n}$ for all $n\in\ZZ$.
\end{definition}
Note that if $\Upsilon$ is hermitian then $\Upsilon_n$ is closable
for every $n\in\ZZ$.
\begin{theorem}
\label{theo:T=TB+TB'} Let $(\A,U)$ be a diffeomorphism covariant
local net on the circle with corresponding representation of the
Virasoro algebra $\{L_n:n\in\ZZ\}$ on the set of finite energy
vectors $\fin$, and central charge $c\in\RR^+_0$. Suppose
$\B\subset \A$ is \mo covariant subnet. Then there exists a unique
decomposition
$$
L_n = L^\B_n + L^{\B'}_n \;\;\;(n\in\ZZ)
$$
such that
\begin{itemize}
\item both collection of operators are Fourier modes of dimension
$2$ \mo covariant hermitian fields,

\item $\overline{L^X_n}$ is affiliated to $X$ where $X=\B,\B'$ and
$n\in\ZZ$.
\end{itemize}
It follows that $\{L^X_n:n\in\ZZ\}$, where $X=\B,\B'$, are
positive energy unitary representations of the Virasoro algebra
with central charges $c_X\in\RR^+_0$, and that the operators
$T^X(f)\equiv\overline{\sum_{n\in\NN}\hat{f}_n L^X_n}$ where
$\hat{f}_n\equiv \frac{1}{2\pi}\int_0^{2\pi}
f(e^{i\theta})e^{-in\theta} d\theta \;\;(n\in\ZZ)$ with $f\in
C^\infty(S^1,\CC)$, are well-defined (i.e.\! the sum exists and
closable). Moreover, we have
\begin{itemize}
\item $[L^\B_n,L^{\B'}_m] = 0$ for all $n,m\in \ZZ$,

\item $0 \leq L^\B_0, L^{\B'}_0 \leq L_0 = L^\B_0 + L^{\B'}_0$ and
similarly, $0\leq c_\B, c_{\B'} \leq c = c_\B + c_{\B'}$,

\item $T^\B(f)$ is affiliated to $\B$ while $T^{\B'}(f)$ to $\B'$,

\item $T^\B(f)$ is affiliated to $\B(I)$ while $T^{\B'}(f)$ to
$\B'\cap \A(I)$ whenever $f$ is supported in the interval $I$,

\item ${\rm Ad}(U({\rm Exp}(f)))|_\B = {\rm Ad}(e^{iT^\B(f)})|_\B$
for all $f\in C^\infty(S^1,\RR)$,

\item $T^\B|_\K$ is the stress-energy tensor that gives conformal
structure to the \mo covariant local net $(\B\!|_\K,
U(\cdot|_\mob)|_\K)$, where $\K=\overline{\B\Omega}$.
\end{itemize}
\end{theorem}
\begin{proof}
We have made an explicit construction resulting a decomposition
$L_n = L^\B_n+L^{\B'}_n$ that satisfies all listed properties and
consequences.

Suppose we have another decomposition $L_n =
\tilde{L}^\B_n+\tilde{L}^{\B'}_n$ satisfying those first listed
two conditions. Then $\Upsilon_n \equiv L^\B_n-L^{\B'}_n$ is still
a Fourier mode of a \mo covariant dimension $2$ hermitian field
whose closure is affiliated to $\B$. On the other hand, as
$L^\B_n+L^{\B'}_n=\tilde{L}^\B_n+\tilde{L}^{\B'}_n$ we have that
$\Upsilon_n = \tilde{L}^{\B'}_n -L^{\B'}_n$ which shows that the
closure of $\Upsilon_n$ is affiliated to $\B'$. It is then easy to
show that on the finite energy vectors $\Upsilon_n$ commutes with
$L_m = L^\B_m + L^{\B'}_m$ for every $n,m\in\ZZ$. Then by \mo
covariance $n \Upsilon_n = [\Upsilon_n,L_0] = 0$ which shows that
$\Upsilon_n = 0$ if $n\neq 0$. Finally, appealing again to \mo
covariance we find that $0 = [L_1,\Upsilon_{-1}] = 2\Upsilon_0$
which concludes the proof.
\end{proof}

\section{Some consequences}

\paragrf
We shall discuss some consequences of the uniqueness of the
diffeomorphism symmetry and the splitting of the stress-energy
tensor given by a subnet. We begin with two rather obvious but
nevertheless important affirmations made in the joint work
\cite{CaWe} of the author with S.\! Carpi.
\begin{corollary}
\label{c.invariant} Let $(\A,V)$ be a \mo covariant net and
suppose the representation $U$ of $\diff$ makes it diffeomorphism
covariant. Then the representation class of $U$, and in particular
its central charge $c\in\RR^+_0$ is an invariant of the M\"obius
covariant net $(\A,V)$. In particular two Virasoro nets as
M\"obius covariant nets are isomorphic if and only if they have
the same central charge.
\end{corollary}
Note that before the uniqueness result, there was no trivial
invariant to distinguish between two Virasoro nets with central
charges bigger than $1$. (Of course when $c<1$, we have several
meaningful invariants, most prominently the $\mu$-index.)

The other interesting thing to note here is the model-independent
proof for the commutation between internal symmetries and
diffeomorphism symmetry.
\begin{definition}
Let $(\A,V)$ be a \mo covariant local net. A unitary $W$ on the
(vacuum) Hilbert space of the net is an {\bf internal symmetry} of
$(\A,V)$ if for every $I \in \I$
\begin{equation}
W \A(I) W^*=\A(I)
\end{equation}
and $W\Omega=\Omega$ where $\Omega$ is the vacuum vector of
$(\A,U)$.
\end{definition}
It is an easy consequence of the Bisognano-Wichmann property that
$W$ commutes with the representation of the \mo group (see
\cite{FrG}). By uniqueness we can now extend this result.
\begin{corollary}
\label{com_int} Let $W$ be an internal symmetry of the
diffeomorphism covariant local net $(\A,U)$. Then $U$ commutes
with $W$.
\end{corollary}
\begin{proof}
Since $W$ commutes with the \mo symmetry, the projective
representation $WVW^*$ of $\diff$ still makes the net $(\A,U)$
diffeomorphism covariant. Hence by uniqueness it must coincide
with $V$. It follows that for every $\gamma \in \diff$ the unitary
$WV(\gamma)W^*V(\gamma)^*$ is a multiple $\lambda(\gamma)$ of the
identity, and in fact it turns out that the complex valued
function $\gamma \mapsto \lambda(\gamma)$ is a character of the
group $\diff$. But as it was already noted, the latter is a simple
noncommutative group, and hence $\lambda$ is trivial.
\end{proof}

\paragrf
We shall now consider subnets and cosets. In \cite{Ko:cosets}
K\"oster proved that a coset can be obtained by local relative
commutants, given that a certain condition is satisfied. According
to the condition, the subnet $\B$, when restricted to its vacuum
Hilbert space, should admit diffeomorphism symmetry and the
corresponding representation of $\diff$ should be possible to
extend to the vacuum Hilbert space of $\A$ in such a way that it
will be compatible with the local structure of $\A$ and
affiliated to $\B$.

In the last section we have seen that given a subnet $\B$ we can
introduce $T^\B$, the part of the stress-energy tensor which
belongs to $\B$. It was not proved though (and maybe it is not
even true) that $T^\B$ integrates to a representation of $\diff$.
(It could integrate to a representation of $\widetilde{\diff}$.)
However, looking at the mentioned article we see that this not
essential for the argument. We shall now make the necessary
rearrangements. Recall that given a locally normal representation
$\pi$ of a Virasoro net we have introduced the stress-energy
tensor $T_\pi$ in the representation.
\begin{lemma}
Let $\pi$ be a locally normal representation of a Virasoro net
$(\A_{\rm Vir}, U)$. Suppose $I_1,..,I_n$ are elements of $\I$ and
$f$ is a smooth function which is zero on $S^1\setminus
\cup_{j=1}^n I_j$. Then $T_\pi(f)$ is affiliated to the von
Neumann algebra $\vee_{j=1}^{n} \pi_{I_j}(\A_{\rm Vir}(I_j))$.
\end{lemma}
\begin{proof}
First of all, it is enough to consider two intervals $I_1,I_2$
with non-dense union in $S^1$. (This case implies the general one
which we can see by considering small enough partitions and/or
appealing to additivity.) Second, the statement is essentially
about Virasoro nets and not about representations, since in this
case the statement involves only operators localized in a certain
(common) interval. So dropping the index ``$\pi$'' from the
stress-energy tensor we shall consider everything in the vacuum
representation.

If $I_1$ and $I_2$ are either distant or overlapping, then the
statement is trivial as then we can find a decomposition $f = f_1
+ f_2$ where $f_1,f_2$ are smooth and ${\rm Supp}(f_j)\subset
\overline{I_j}$ (the fact the we have a closure does not course
any problem, since the local algebras are ``continuous'', see also
the proof of Prop.\! \ref{aff}). So let us consider the case when
$I_1,I_2$ have a common endpoint, which we can safely assume to be
the point $1\in S^1$. Then by our assumption $f(1)=0$. If also
$f'(1) = 0$, then the decomposition 
\begin{equation}
f=f \chi_{I_1} +
f\chi_{I_2}, 
\end{equation}
where $\chi_{I_j}$ is the characteristic function of
$I_j$ $(j=1,2)$, is a decomposition into (once) differentiable
functions. Moreover it is clear that the conditions of Lemma
\ref{finite1.5norm} are satisfied and thus these functions are of
finite $\|\cdot\|_{\frac{3}{2}}$ norm. Then again, by the results
concerning evaluation of the stress-energy tensor on nonsmooth
functions we are finished. Thus it is enough to demonstrate the
statement in a single example when $f'(1)\neq 0$.

We shall use the idea appearing in the proof of \cite[Lemma
11]{Ko:cosets}. For the generating vector field $d$ of the
dilations we have $d(1)=0, d'(1)=-1\neq 0$. Let $f$ be a real
smooth function which coincides with $d$ in some small
neighborhood of $1$ but has its support in $\overline{I_1\cup
I_2}$. Further, with the star standing for actions on vector
field, set $f_s \equiv (\delta_s)_* f$ $(s\in\RR)$.

Let $K\in\I$ distant from both $I_1$ and $I_2$. For $s\to\infty$
the support of $f_s$ is contracting and thus for $s>0$ remains
disjoint from $K$. Hence for any $X\in\A_{\rm Vir}(K)$ and vector
of the Hilbert space $\Psi$ we have
\begin{eqnarray}\nonumber
\lim_{s\to\infty}\langle\Psi,\,T(f_s)X \Omega\rangle &=&
\lim_{s\to\infty}\langle X^*\Psi,\, U(\delta_s)T(f)\Omega\rangle
\\&=& \langle X^*\Psi,\,\Omega\rangle\langle\Omega,\,T(f)\Omega\rangle = 0
\end{eqnarray}
since $T(f)\Omega$ is orthogonal to the vacuum. As a dilation
leaves invariant its generating vector field, in a sufficiently
small neighborhood of $1$ the function $f_s$ will remain equal to
$f$. Hence $T(f-f_s)$ is affiliated to $\A_{\rm Vir}(I_1) \vee
\A_{\rm Vir}(I_1)$. Thus if $Y'\in \{\A_{\rm Vir}(I_j):j=1,2\}'$
and $X_1,X_2\in \A_{\rm Vir}(K)$ then by our previous equation
\begin{eqnarray}\nonumber
\langle X_2\Omega,\,Y'T(f)X_1 \Omega\rangle &=&
\lim_{s\to\infty}\langle X_2\Omega,\,Y' T(f-f_s)X_1\Omega\rangle \\
\nonumber &=&\lim_{s\to\infty}\langle T(f-f_s) X_2\Omega,\,Y'
X_1\Omega\rangle \\ &=& \langle T(f)X_2\Omega,\,Y'X_1
\Omega\rangle.
\end{eqnarray}
The set $\A_{\rm Vir}(K)\Omega$ is a core for $T(f)$ since by
Theorem \ref{localcore} it contains a core for $L_0$. Then by the
above equation it is easy to show that the self-adjoint operator
$T(f)$ commutes with $Y'$ which concludes the proof.
\end{proof}
\begin{lemma}
Let $(\A,U)$ be a conformal local net, $\B\subset \A$ a subnet and
$T^\B$ the part of the stress-energy tensor that belongs to $\B$.
Suppose $f$ is a real smooth function with $f(z_j)=0$ $(j=1,2)$
where $z_1$ and $z_2$ are the endpoints of the interval $I\in\I$.
Then the adjoint action of $e^{iT^\B(f)}$ preserves the algebra
$\A(I)$.
\end{lemma}
\begin{proof}
Suppose we prove this for all functions localized somewhere. Then
by taking partitions and using the Trotter product formula we can
justify the statement for all functions. But once we fix an
interval $K\in \I$, we can use the isomorphism between
$\B(K)|_{\overline{\B\Omega}}$ and $\B(K)$. Therefore for any
diffeomorphism $\gamma$ localized in $K$ we have a unique
decomposition $U(\gamma)= U^\B(\gamma)U^{\B'}(\gamma)$ where
$U^\B(\gamma)\in\B(K)$ and $U^{\B'}(\gamma)\in\B(I)'\cap \A(I)$.
(Note that we did not integrate $T^\B$ into a global
representation of $\diff$: that may be impossible. We have only
integrated it locally.) Then using the (local) representation
$U^\B$ the argument of the proof of \cite[Lemma 2]{Ko:cosets} can
be essentially repeated.
\end{proof}
\begin{corollary}
Let $(\A,U)$ be a conformal local net, $\B\subset \A$ a subnet.
Then $\B'\cap\A(I) = \B(I)'\cap\A(I)$ for every $I\in\I$.
\end{corollary}
\begin{proof}
Our first lemma replaces \cite[Lemma 11]{Ko:cosets} and our second
lemma replaces \cite[Lemma 2]{Ko:cosets}. Thus \cite[Theorem
12]{Ko:cosets} holds, which --- although with a slight difference
--- is exactly what we have claimed.
\end{proof}

\paragrf
Another interesting observation is that the subnet generated by
the stress-energy tensor must be the smallest among the
irreducible subnets.
\begin{corollary}
Let $(\A,U)$ be a conformal local net and $\B\subset \A$ a subnet.
If $\B'\cap \A(I) = \CC \mathbbm 1$ for an interval $I\in\I$ then
it follows that $U(\diff)\subset \B$.
\end{corollary}
\begin{proof}
By \mo covariance it is clear that if the commutant in question is
trivial for a certain interval then it is trivial for elements of
$\I$. So consider the stress-energy tensor $T$ and the part $T^\B$
belonging to $\B$. If $f$ is a function localized in a certain
interval $K\in\I$ then by Theorem \ref{theo:T=TB+TB'} (the closure
of the restriction to the finite energy vectors of) $T(f) -
T^\B(f)$ is affiliated to $\B'\cap \A(K) = \CC \mathbbm 1$, so it
is a complex number. As the expectation of this difference in the
vacuum state is zero, we have that $T(f)= T^\B(f) \in \B$. Thus we
are finished as $\diff$ is generated by exponentials of local
vector fields.
\end{proof}
Suppose we have a chain of subnets. If the local relative
commutants are nontrivial, then at each step we ``loose'' a part
of the stress-energy tensor. As the smallest positive value for
the central charge is $1/2$, this means that at each step the
central charge of the subnet must decrease at least by $1/2$.
Hence we get a restriction on the possible length of this chain.
In turn, if we have an infinite chain of such a subnets then the net 
cannot admit diffeomorphism symmetry. This is exactly the case of a net 
obtained by taking an infinite tensor product (see Sect.\! 
\ref{sec:more.ex} about tensor products). 
Hence in this way we can construct new examples of \mo
covariant local nets that do not admit diffeomorphism symmetry.

The observation that the infinite tensor product net cannot be
diffeomorphism covariant, first appeared in the joint work
\cite{CaWe} of the author with S.\! Carpi. There it is the
uniqueness of the diffeomorphism symmetry, which is exploited to
show that the central charge of the tensor product must be the sum
of the central charges of the nets appearing in the product, and
hence in case of infinite tensor product it should infinite (if
the product was diffeomorphism covariant). Hence the argument
applies in case of taking tensor products of diffeomorphism
covariant nets. However, by the results of the last section we can
consider a more general situation.
\begin{definition}
Let $(\A,U)$ be a conformal local net and $\B\subset \A$ a subnet.
The {\bf central charge of the subnet} $\B$ is the central charge
of the diffeomorphism covariant local net which we get by
restriction to its vacuum Hilbert space.
\end{definition}
\begin{corollary}
Let $(\A,U)$ be a conformal local net and $\B\subset \A$ a subnet.
Then $c=c_\B+c_\C$ where $c$ is the central charge of $\A$ and
$c_\B,c_\C$ are the central charges of $\B$ and the coset
$I\mapsto \C(I)\equiv\B'\cap\A(I)=\B(I)'\cap\A(I)$, respectively.
\end{corollary}
If a conformal net (or subnet) is not trivial, then neither the
corresponding representation of $\diff$ can be trivial. But for a
nontrivial positive energy representation of $\diff$ the lowest
possible value of the central charge is $1/2$. Hence by Theorem
\ref{theo:T=TB+TB'} we make the following observations.
\begin{corollary}
Let $\A=\A_0\supset \A_1 \supset ... \supset \A_n$ be a chain of
subnets in the conformal local net $(\A,U)$ with $\A_{j+1}(I)'\cap
\A_j(I)\neq \CC \mathbbm 1$ for every $j\in\{1,..,n\}$, $I\in\I$. 
Then $c_{\A_j}\geq c_{\A_{j+1}}+\frac{1}{2}$ and hence
$c_A \geq \frac{n}{2}$.
\end{corollary}
\begin{corollary}\label{inf.tens.product->nodiff}
Let $(\A_n,V_n)$, $n\in \NN$ be a sequence of nontrivial \mo
covariant local nets. Then the infinite tensor product net
$\otimes_{n\in\NN}( \A_n,V_n)$ cannot admit diffeomorphism
symmetry.
\end{corollary}
\begin{proof}
$\otimes_{n=1}^{\infty}( \A_n,V_n) \supset \otimes_{n=2}^{\infty}(
\A_n,V_n) \supset ..$ is an infinite chain of subnets with
nontrivial local relative commutants, so by the previous corollary
the full net cannot have a finite central charge.
\end{proof}
Previously the only examples of nets not admitting diffeomorphism
symmetry (at least to the knowledge of the author), were the ones
coming from the derivatives of the $U(1)$ current model (see
Sect.\! \ref{sec:U(1)}). They are not strongly additive, in fact
they are not even $4$-regular. In the presented new way we can
construct examples that are strongly additive. However, on the
other hand --- although it is not proven --- it is likely that
they are not going to be split. (In fact if the infinite sequence
contains an infinite subsequence of isomorphic nets then by
\cite[Theorem 9.2]{WSSI} it is easy to show that they are indeed
not.) So it could be that split property and $4$-regularity are
necessary for a net to admit diffeomorphism symmetry.

\section{Positivity of energy in representations}

\paragrf
The positivity of energy is one of the most important selection
criteria for a model to be ``physical''. In almost all treatments
of Quantum Field Theory it appears as one of the fundamental
axioms. As an axiom, one may say that it is ``automatically''
true, but in a concrete model it is something to be checked. In
particular, to see what are the charged sectors with positive
energy for a model given in its vacuum sector may be difficult (as
calculating the charged sectors can already be a hard problem).

As it was already discussed in the preliminary sections, in the
setting of local nets charged sectors are described as irreducible
representations of the net and their general theory was developed
by Doplicher, Haag and Roberts \cite{DHR1,DHR2}. They proved,
among many other things, that if a covariant sector has a finite
statistical dimension then it is automatically of positive energy.
In \cite{GL} Guido and Longo showed that under some regularity
condition the finiteness of the statistics even implies the
covariance property of the sector.

In case of chiral theories the question of positivity of energy in
a sector with finite statistics is more or less completely
understood \cite{GL,GuLo96,BCL}. If a conformal net is completely
rational, then --- as in this case every sector of it is
automatically covariant and has finite statistics --- we have
positivity of energy in every sector. Yet, although these
conditions cover many of the interesting models, as it was already
discussed there are interesting (not pathological!) models in
which it does not hold and, what is more important in this
context, indeed possessing sectors with infinite statistical
dimension (and yet with positivity of energy). This is clearly in
contrast with the experience coming from massive QFTs (by a
theorem of Buchholz and Fredenhagen \cite{BuFre}, a massive sector
with positive energy is always localizable in a spacelike cone and
has finite statistics).

The first example of a sector with infinite statistics was
constructed by Fredenhagen \cite{Fre}. Rehren gave arguments
\cite{Re} that even the Virasoro model, which is in some sense the
most natural model, should admit sectors with infinite statistical
dimensions when its central charge $c\geq 1$ and that in fact in
this case ``most'' of its sectors should be of infinite
statistics. This was then actually proved \cite{Carpi1} by Carpi
first for the case $c=1$ and then \cite{Carpi2} for many other
values of the central charge, leaving open the question only for
some values of $c$ between $1$ and $2$. Moreover Longo and Xu
proved \cite{LoXu} that if $\A$ is a split conformal net with $\mu
= \infty$ then $(\A \otimes \A)^{\rm flip}$ has at least one
sector with infinite statistical dimension, showing that the case
of infinite statistics is indeed quite general.

Recently D'Antoni, Fredenhagen and K\"oster published a letter
\cite{D'AnFreKos} with a proof that diffeomorphism covariance
itself (in the vacuum sector) is already enough to ensure \mo
covariance in every (not necessary irreducible) representation:
there always exists a unique (projective, strongly continuous)
inner implementation of the \mo symmetry. (In Prop. \ref{n-cover}
we shall generalize this statement to the {\it $n$-\mo} group.)
Thus the concept of the conformal energy, as the self-adjoint
generator of rotations in given representation, is at least
well-defined. (Without the assumption of diffeomorphism covariance
it is in general not true: recall the examples discussed in
Sect.\! \ref{sec:U(1)}.) What remained an open question until now,
whether this energy is automatically positive or not. The author
of this thesis, using the tools previously developed and
explained, gave a proof for the positivity (Theorem
\ref{mainresult}). This will be the main result presented in this
section. Everything that here follows is part of the single-author
work \cite{weiner1} accepted by Commun.\! Math.\! Phys.

\paragrf
The idea behind the proof is simple. The total conformal energy
$L_0$ is the integral of the energy-density; i.e.\! the
stress-energy tensor $T$ evaluated on the constant $1$ function.
So if we take a finite partition of the unity $\{f_n\}_{n=1}^N$ on
the circle, we may write $T(1)$ as the sum $\sum T(f_n)$ where
each element is {\it local}. Thus each term in itself (although
not bounded) can be considered in a given charged sector.
Moreover, it has been recently proved by Fewster and Hollands
\cite{FewHoll} that the stress-energy tensor evaluated on a
nonnegative function is bounded from below. These operators then,
being local elements, remain bounded from below also in the
charged sector. So their sum in the charged sector, which we may
expect to be the generator of rotations in that sector, should
still be bounded from below.

There are several problems with this idea. For example, as the
supports of the functions $\{f_n\}$ must unavoidably ``overlap'',
the corresponding operators will in general not commute. To deal
with sums of non-commuting unbounded operators is not easy. In
particular, while in the vacuum representation --- due to the well
known energy bounds --- we have the natural common domain of the
finite energy vectors, in a charged sector (unless we assume
positivity of energy, which is exactly what we want to prove) we
have no such domain.

To overcome the difficulties we shall modify this idea in two
points. First of all, instead of $L_0=T(1)$, that is, the
generator of rotations, we may work with the generator of
translations --- the positivity of any of them implies the
positivity of the other one. In fact we shall go one step further
by replacing the generator of translations with the generator of
$2${\it-translations}. (This is why, as it has been already
mentioned, we shall consider the $n$-\mo group; particularly in
the case $n=2$.) This has the advantage that the function
representing the corresponding vector field can be written as
$f_1+f_2$, where the two local nonnegative functions $f_1,f_2$
{\it do not} ``overlap''. These functions, at the endpoint of
their support are not smooth (such decomposition is not possible
with smooth functions); they are only once differentiable.
However, as it was discussed all through the previous chapter, the
stress-energy tensor can be evaluated even on nonsmooth functions,
given that their $\|\cdot\|_{\frac{3}{2}}$ norm is finite, which
is exactly the case of $f_1$ and $f_2$. As they are nonnegative
but not smooth, to conclude that $T(f_1)$ and $T(f_2)$ are bounded
from below we cannot use the result stated in \cite{FewHoll}.
However, it turns out to be (Prop.\! \ref{affiliation}) a rather
direct and simple consequence of the construction, thus it will be
deduced independently from the mentioned result, of which we shall
make no explicit use. In fact the author considered this
construction as an argument indicating that if $f \geq 0$ then
$T(f)$ is bounded from below (which by now is of course proven, as
it was already mentioned, in \cite{FewHoll}); see more on this at
the remark after Prop \ref{affiliation} and in the mentioned
article of Fewster and Hollands at the footnote in the proof of
\cite[Theorem 4.1]{FewHoll}.

\paragrf
Recall how the subgroup $\mob^{(n)}\subset \diff$ was defined for
every positive integer $n$ (see the preliminary Sect.\!
\ref{sec:diffgroup}). Recall also that in $\mob^{(n)}$ we defined
the one-parameter subgroup of $n$-translations $a \mapsto
\tau^{(n)}_a$ by the usual procedure of lifting: it was set to be
the unique continuous one-parameter subgroup satisfying
$\tau^{(n)}_a(z)^n = \tau_a(z^n)$. Its generating vector field is
$t^{(n)}(z)=\frac{1}{2n} - \frac{1}{4n}(z^n+z^{-n})$. The notion
of $n$-rotations $\alpha \mapsto \rho^{(n)}_\alpha$ and
$n$-dilations $s \mapsto \delta^{(n)}_s$ were introduced
similarly. The ``$n$-rotations'', however, apart from a re-scaling
of the parameter, simply coincide with the ``true'' rotations, see

As $\mob^{(n)}$ covers $\mob$ in a natural way, its universal
cover is canonically identified with $\widetilde{\mob} \simeq
\widetilde{{\rm SL}(2,\RR)}$. Note that if $p^{(n)} :
\widetilde{\mob} \rightarrow \mob^{(n)}$ is the natural covering
map then $p^{(n)}(\tilde{R}^{(m)}_\alpha)=R^{(n)}_\alpha$ for all
$n,m\in\NN^+$, where the sign ``$\tilde{\hphantom{a}}$'', as
always throughout the rest of this article, stands for the
appropriate lifts of one-parameter groups to $\widetilde{\mob}$. A
strongly continuous projective representation of $\mob^{(n)}$
lifts to a unique strongly continuous representation of
$\widetilde{\mob}$, see more on this in Sect.\!
\ref{app:mobrep} of the appendix.

Regarding the representation theory of the latter group, in
particular the following is known (see Prop.\! \ref{pos.mob^n}).
If $\tilde{V}$ is a strongly continuous unitary representation of
$\widetilde{\mob}$ with $H$ and $P$ being the self-adjoint
generator of rotations and translations in $\tilde{V}$,
respectively, then the following four
conditions are equivalent: \\
\indent $1.\;\,H$ is bounded from below, \\
\indent $2.\;\,P$ is bounded from below, \\
\indent $3.\;\,H \geq 0$, \\
\indent $4.\;\,P \geq 0$. \\
\noindent If any of the above conditions is satisfied, $\tilde{V}$
is called a positive energy representation.

\paragrf
Let us now consider a conformal local net of on the circle
$(\A,U)$ with stress-energy tensor $T$. By equation
(\ref{n-rotation}), $U^{(n)}$, the restriction of the positive
energy representation $U$ of $\diff$ to $\mob^{(n)}$, lifts to a
positive energy representation of $\widetilde{\mob}$. In
particular, the self-adjoint operator $T(t^{(2)})$ must be bounded
from below, since it generates the translations for the
representation $U^{(2)}$. (Note that $T(t^{(2)})$ is bounded from
below but not necessary positive: it is not {\it the} generator
---  it still generates the same projective one-parameter group of
unitaries if you add a real constant to it.) The function
$t^{(2)}(z)= \frac{1}{4}-\frac{1}{8}(z^2+z^{-2})$ is a nonnegative
function with two points of zero: $t^{(2)}(\pm 1)=0$. By direct
calculation of the first derivative: $(t^{(2)})'(\pm 1)=0$, hence
the decomposition
\begin{equation}
\label{decomp} t^{(2)} = t^{(2)}_+ + t^{(2)}_-
\end{equation}
with the functions $t^{(2)}_\pm$ defined by the condition ${\rm
Supp}(t^{(2)}_\pm) = (S^1_\mp)^c$ is a decomposition of $t^{(2)}$
into the sum of two (once) differentiable nonnegative functions
that satisfy the conditions of Lemma \ref{finite1.5norm}.
Therefore, as it was explained, we can consider the self-adjoint
operators $T(t^{(2)}_\pm)$.
\begin{proposition}
\label{affiliation} Let $(\A,U)$ be a conformal net of local
algebras on the circle with stress-energy tensor $T$. Then
$T(t^{(2)}_+)$ is affiliated to $\A(S^1_+)$ and $T(t^{(2)}_-)$ is
affiliated to $\A(S^1_-)$ and so in particular they strongly
commute. Moreover, the operators $T(t^{(2)}_\pm)$ are bounded from
below.
\end{proposition}
\begin{proof}
${\rm Supp}(t^{(2)}_\pm) \subset \overline{S^1_\pm}$ and so by
Prop. \ref{aff} $T(t^{(2)}_\pm)$ is affiliated to $\A(S^1_\pm)$.
So if $P_{[a,b]}$ is a nonzero spectral projection of
$T(t^{(2)}_+)$ and $Q_{[c,d]}$ is a nonzero spectral projection of
$T(t^{(2)}_-)$, then $P_{[a,b]} \in \A(S^1_+),\; Q_{[c,d]} \in
\A(S^1_-)$ and by the algebraic independence of two commuting
factors (see for example \cite[Theorem 5.5.4]{kadison})
$R=P_{[a,b]}Q_{[c,d]} \neq 0$. Of course the range of $R$ is
invariant for (and included in the domain of) $T(t^{(2)}_+) +
T(t^{(2)}_-)$ and the restriction of that operator for this closed
subspace is clearly bigger than $a+c$ and smaller than $b+d$.
Thus, nesting every point of the spectra of $T(t^{(2)}_\pm)$ into
arbitrarily narrow intervals we get that
\begin{equation}
{\rm Sp}\left(T(t^{(2)}_+) + T(t^{(2)}_-)\right) \supset {\rm
Sp}(T(t^{(2)}_+)) + {\rm Sp}(T(t^{(2)}_-)).
\end{equation}
To conclude we only need to observe that by equation
(\ref{decomp}) on the common core of the finite energy vectors
$T(t^{(2)}_+) + T(t^{(2)}_-) = T(t^{(2)})$, and as it was said,
the latter self-adjoint operator is bounded from below.
\end{proof}
\smallskip

{\it Remark.} The author considered this construction to indicate
that if $f \geq 0$ then $T(f)$ is bounded from below, which --- as
it was already mentioned --- by now is a proven fact (cf.\!
\cite{FewHoll}). The point is the following. If $f$ is {\it
strictly} positive then, as a vector field on $S^1$, it is
conjugate to the constant vector field $r$ for some $r>0$. Thus,
using the transformation rule of $T$ under diffeomorphisms, $T(f)$
is conjugate to $T(r)$ plus a constant, and so it is bounded from
below by the positivity of $T(1)=L_0$. The real question is
whether the statement remains true even when $f$ is nonnegative,
but not strictly positive because for example it is {\it local}
(there is an entire interval on which it is zero). One can of
course consider a nonnegative function as a limit of positive
functions, but then one needs to control that the lowest point of
the spectrum does not go to $-\infty$ while taking this limit
(which --- in a slightly different manner --- has been
successfully carried out in the mentioned article). However, even
without considering limits, by the above proposition we find
nontrivial examples of local nonnegative functions $g$ such that
$T(g)$ can easily be checked to be bounded from below. (Take for
example $g=t^{(2)}_\pm$ but of course we may consider conjugates,
sums and multiples by positive constants to generate even more
examples.)
\smallskip

Let us now investigate what we can say about a representation
$\pi$ of the conformal net $(\A,U)$. In \cite{D'AnFreKos} it was
proved that the \mo symmetry is continuously implementable in any
(locally normal) representation $\pi$ by a unique inner way. By
their construction the implementing operators are elements of
$\pi(\T_\A)$ where $\T_\A\subset \A$ is the subnet given by the
stress-energy tensor. Looking at the article, we see that the only
structural properties of the \mo subgroup of $\diff$ that the
proof uses are the following.
\begin{itemize}
\item There exist three continuous one-parameter groups $\Gamma_1,
\Gamma_2$ and $\Gamma_3$ in $\mob$, so that $(s_1,s_2,s_3)\mapsto
\Gamma_1(s_1)\Gamma_2(s_2)\Gamma_3(s_3)$ is a (continuous)
covering of $\mob$ (In the article $\Gamma_1$ is the
translational, $\Gamma_2$ is the dilational and $\Gamma_3$ is the
rotational subgroup; see Prop.\! \ref{iwasawa} for the Iwasawa
decomposition.) \item The Lie algebra of $\mob$ is isomorphic to
${\mathfrak{sl}}(2,\RR)$.
\end{itemize}
These two properties hold not only for the subgroup $\mob$, but
also for $\mob^{(n)}$ where $n$ is any positive integer: the Lie
algebra of $\mob^{(n)}$ is isomorphic to ${\mathfrak{sl}}(2,\RR)$
for all $n\in\NN^+$, and with the rotations, translations and
dilations replaced by $n$-rotations, $n$-translations and
$n$-dilations we still have the required decomposition. Let us
collect into a proposition what we have thus concluded.
\begin{proposition}
\label{n-cover} Let $\pi$ be a locally normal representation of
the conformal local net of von Neumann algebras on the circle
$(\A,U)$. Then for all $n\in\NN^+$ there exists a unique strongly
continuous representation $\tilde{U}^{(n)}_\pi$ of
$\widetilde{\mob}$ such that
$\tilde{U}^{(n)}_\pi(\widetilde{\mob})\subset \pi(\A)$ and for all
$\tilde{\gamma} \in \widetilde{\mob}$ and  $I \in \I$
$$
{\rm Ad}(\tilde{U}^{(n)}(\tilde{\gamma})) \circ \pi_I =
{\pi}_{\gamma^{(n)}(I)} \circ {\rm Ad}(U(\gamma^{(n)})|_{\A(I)}
$$
where $\gamma^{(n)}=p^{(n)}(\tilde{\gamma})$ with $p^{(n)}:
\widetilde{\mob} \rightarrow \mob^{(n)}$ being the natural
covering map. Moreover, this unique representation satisfies
$\widetilde{U}^{(n)}_\pi(\widetilde{\mob})\subset \pi(\T\!_\A)$.
\end{proposition}
We shall now return  to the particular case $n=2$. On one hand,
the action of the $2$-translation $\tau^{(2)}_a$ in the
representation $\pi$ can be implemented by
$\tilde{U}^{(2)}_\pi(\tilde{\tau}^{(2)}_a)$. On the other hand, as
in the projective sense
\begin{equation}
U(\tau^{(2)}_a)= e^{iaT(t^{(2)})}=e^{iaT(t^{(2)}_+)}
e^{iaT(t^{(2)}_-)},
\end{equation}
we may try to implement the same action by $\pi_{S^1_+}(W_+(a))
\pi_{S^1_-}(W_-(a))$, where
\begin{equation}
W_\pm(a)\equiv e^{iaT(t^{(2)}_\pm)} \in \T\!_\A(S^1_\pm).
\end{equation}
\begin{proposition}
\label{localimplementation} The unitary operator in $\pi(\T\!_\A)$
$$W_\pi(a)\equiv \pi_{S^1_+}(W_+(a)) \,\pi_{S^1_-}(W_-(a))
=\pi_{S^1_-}(W_-(a))\, \pi_{S^1_+}(W_+(a))$$ up to phase coincides
with $\tilde{U}^{(2)}_\pi(\tilde{\tau}^{(2)}_a)$.
\end{proposition}
\begin{proof}
It is more or less trivial that the adjoint action of the two
unitaries coincide on both $\pi_{S^1_+}(\A(S^1_+))$ and
$\pi_{S^1_-}(\A(S^1_-))$. There remain two problems to overcome:
\begin{itemize}
\item the algebra generated by these two algebras do not
necessarily contain $\pi(\T\!_\A)$, so it is not clear why the
adjoint action of these two unitaries should coincide on the
mentioned algebra, \item but even if we knew that the actions
coincide, the two unitaries, although both belonging to
$\pi(\T\!_\A)$, for what we know could still ``differ'' in an
inner element.
\end{itemize}
As for the first problem, consider an open interval $I \subset
S^1$ such that it contains the point $-1$ and has $1$ in the
interior of its complement. Note that due to the conditions
imposed on $I$, the sets $K_\pm\equiv I \cup S^1_\pm$ are still
elements of $\I$.
\begin{lemma}
If $a \geq 0$ then $W_+(a) \A(I) W_+(a)^* \subset \A(I)$.
\end{lemma}
\begin{proof}[Proof of the Lemma]
The flow of a vector field given by a nonnegative function on the
circle, moves all points forward (i.e. anticlockwise). Moreover,
the flow cannot move points from the support of the vector field
to outside, and leaves invariant all points outside.

Taking in consideration Prop.\! \ref{aff} and what was said it is
easy to see that ${\rm Exp}(t^{(2)}_+)(I) \subset I$ and
consequently
\begin{equation}
{\rm Ad}\left(e^{i a T(t^{(2)}_+)}\right)(\A(I))\subset \A(I).
\end{equation}
\end{proof}
It follows that if $A \in \A(I)$ and $a \geq 0$ then
\begin{eqnarray} \nonumber
&&\pi_{S^1_+}(W_+(a))\,\pi_I(A)\,\pi_{S^1_+}(W_+(a))^* = \\
&&\pi_{K_+} (W_+(a)\,A\,W_+(a)^*) = \pi_I (W_+(a)\,A\,W_+(a)^*)
\end{eqnarray}
and thus $\rm{Ad}\left(W_\pi(a)\right)(\pi_I(A))=
\rm{Ad}\left(\pi_{S^1_-}(W_-(a))\,\pi_{S^1_+}(W_+(a)) \right)
(\pi_I(A))=$
\begin{eqnarray}
\nonumber
&=&\rm{Ad}\left(\pi_{S^1_-}(W_-(a))\right)(\pi_I(W_+(a)\,A\,W_+(a)^*))
\\ \nonumber
&=&\pi_{K_-}(W_-(a)\,(W_+(a)\,A\,W_+(a)^*)\,W_-(a)^*)
\\
&=&\pi_{K_-}(U(\tau^{(2)}_a)\,A\,U(\tau^{(2)}_a)^*)
=\rm{Ad}\left(\tilde{U}^{(2)}_\pi(\tilde{\tau}^{(2)}_a)\right)(\pi_I(A))
\end{eqnarray}
where in the last equality we have used the fact that for $a \geq
0$ the image of $I$ under the diffeomorphism $\tau^{(2)}_{a}$ is
contained in $K_-$.

We have thus seen that for $a \geq 0$ the adjoint action of
$W_\pi(a)$ and of $\tilde{U}^{(2)}_\pi(\tilde{\tau}^{(2)}_a)$
coincide on $\pi_I(A(I))$. Actually, looking at our argument we
can realize that everything remains true if instead of $I$ we
begin with an open interval $L$ that contains the point $1$ and
has $-1$ in the interior of its complement and we exchange the
``+'' and ``-'' subindices. So in fact we have proved that for $a
\geq 0$ these adjoint actions coincide on both $\pi_I(A(I))$ and
$\pi_{L}(A(L))$ and therefore on the whole algebra $\pi(\A)$,
since we may assume that the union of $I$ and $L$ is the whole
circle. (The choice of the intervals, apart from the conditions
listed, was arbitrary.) Of course the equality of the actions, as
they are obviously one-parameter automorphism groups of $\pi(\A)$,
is true also in case the parameter $a$ is negative. We can now
also confirm that the unitary
\begin{equation}
Z_\pi(a) \equiv W_\pi(a)^* \,
\tilde{U}^{(2)}_\pi(\tilde{\tau}^{(2)}_a)
\end{equation}
lies in $\Z(\pi(\A))\cap \pi(\T\!_\A) \subset \Z(\pi(\T\!_\A))$
where ``$\Z$'' stands for the word ``center''. Thus $a \mapsto
Z_\pi(a)$ is a strongly continuous one-parameter group. (It is
easy to see that as $Z_\pi$ commutes with both $W_\pi$ and
$\tilde{U}^{(2)}_\pi\circ \tilde{\tau}^{(2)}$ it is actually a
one-parameter group.)

We shall now deal with the second mentioned problem. The
$2$-dilations $s \mapsto \delta^{(2)}_s$ scale the
$2$-translations and preserve the intervals $S^1_\pm$. Thus they
also scale the functions $t^{(2)}_\pm$ and so we get some
relations --- both in the vacuum and in the representation $\pi$
--- regarding the unitaries implementing the dilations and
translations and the unitaries that were denoted by $W$ with
different subindices (see Prop.\! \ref{aff}). More concretely,
with everything meant in the projective sense, in the vacuum
Hilbert space we have that the adjoint action of
$U(\delta^{(2)}_s)$ scales the parameter $a$ into $e^sa$ in
$U(\tau^{(2)}_a)$ and in $W_\pm(a)$. In the Hilbert space $\H_\pi$
we have exactly the same scaling of
$\tilde{U}^{(2)}_\pi(\tilde{\tau}^{(2)}_a)$ and of
$\pi_{S^1_\pm}(W_\pm(a))$ by the adjoint action of
$\tilde{U}^{(2)}_\pi(\tilde{\delta}^{(2)}_s)$. Thus (up to phases)
\begin{equation}
Ad\left(\tilde{U}^{(2)}_\pi(\tilde{\delta}^{(2)}_s)\right)(Z_\pi(a))
= Z_\pi(e^s a),
\end{equation}
but on the other hand of course, as $Z_\pi$ is in the center, the
left hand side should be simply equal to $Z_\pi(a)$. So up to
phases $Z_\pi(a)= Z_\pi(e^s a)$  for all values of the parameters
$a$ and $s$ which means that $Z_\pi$ is a multiple of the identity
and hence in the projective sense $W_\pi(a)$ equals to
$\tilde{U}^{(2)}_\pi(\tilde{\tau}^{(2)}_a)$.
\end{proof}
\begin{corollary}
\label{pos:2} The representation $\tilde{U}^{(2)}_\pi$ is of
positive energy.
\end{corollary}
\begin{proof}
The spectrum of the generator of a one-parameter unitary group
remains unchanged in any normal representation. The one-parameter
groups $a\mapsto e^{iaT(t^{(2)}_\pm)}$ are {\it local}. Thus, as
the representation $\pi$ is locally normal, by Prop.
\ref{affiliation} the self-adjoint generator of the one-parameter
group
\begin{equation}
a \mapsto \pi_{S^1_+}\left(e^{iaT(t^{(2)}_+)}\right)
\pi_{S^1_-}\left(e^{iaT(t^{(2)}_-)}\right)
\end{equation}
is bounded from below and by Prop.\! \ref{localimplementation}
this one-parameter group of unitaries, at least in the projective
sense, equals to the one-parameter group $a \mapsto
\tilde{U}^{(2)}_\pi(\tilde{\tau}^{(2)}_a)$. So by Lemma
\ref{pos.mob^n} the representation $\tilde{U}^{(2)}_\pi$ is of
positive energy.
\end{proof}

Let us now take an arbitrary positive integer $n$. By equation
(\ref{n-rotation}) both
$\tilde{U}^{(n)}_\pi(\tilde{\rho}_{n\alpha})$ and
$\tilde{U}^{(2)}_\pi(\tilde{\rho}_{2\alpha})$ implement the same
automorphism of $\pi(\A)$. Since both unitaries are actually
elements of $\pi(\T\!_\A)\subset \pi(\A)$, they must commute and
\begin{equation}
\label{propotion}
C_{\pi}^{(n)}(\alpha)\equiv(\tilde{U}^{(n)}_\pi(\tilde{\rho}_{n\alpha}))^*
\,\, \tilde{U}^{(2)}_\pi(\tilde{\rho}_{2\alpha})
\end{equation}
is a strongly continuous one-parameter group in the von Neumann
algebra $\Z(\pi(\A))\cap \pi(\T\!_\A) \subset \Z(\pi(\T\!_\A))$.

As it was mentioned, (see the preliminary Sect.\!
\ref{sec:diffcov}) the restriction of the subnet $\T\!_\A$ onto
$\H_{\T\!_\A}$
--- unless $\A$ is trivial, in which case
dim$(\H_{\T\!_\A})={\rm dim}(\H_\A)=1$  --- is isomorphic to a
Virasoro net. Thus $\H_{\T\!_\A}$ must be separable (even if the
full Hilbert space $\H_\A$ is not so; recall that we did not
assume separability) as the Hilbert space of a Virasoro net is
separable.

Every von Neumann algebra on a separable Hilbert space has a
strongly dense separable $C^*$ subalgebra. A von Neumann algebra
generated by a finite number of von Neumann algebras with strongly
dense separable $C^*$ subalgebras has a strongly dense $C^*$
subalgebra. Thus considering that for an $I\in\I$ the restriction
map from $\T\!_\A(I)$ to $\T\!_\A(I)|_{\H_{\T\!_\A}}$ is an
isomorphism, one can easily verify that the von Neumann algebra
$\pi(\T\!_\A)$ has a strongly dense $C^*$ subalgebra.

We can thus safely consider the direct integral decomposition of
$\pi(\T\!_\A)$ along its center
\begin{equation}
\label{directint} \pi(\A) = {\int}_{\!\!\! X}^\oplus \pi(\A)(x)
d\mu(x).
\end{equation}
(Even if $\H_\pi$ is not separable, by the mentioned property of
the algebra $\pi(\T\!_\A)$, it can be decomposed into the direct
sum of invariant separable subspaces for $\pi(\T\!_\A)$. Then
writing the direct integral decomposition in each of those
subspaces, the rest of the argument can be carried out without
further changes.) For an introduction on the topic of the direct
integrals see for example \cite[Chapter 14.]{kadison}.

Since the group $\widetilde{\mob}$ is in particular second
countable and locally compact, and the representations
$\{\tilde{U}^{(n)}_\pi\}$ are in $\pi(\T\!_\A)$, the decomposition
(\ref{directint}) also decomposes these representations (cf.\!
\cite[Lemma 8.3.1 and Remark 18.7.6]{dixmier}):
\begin{equation}
\tilde{U}^{(n)}_\pi(\cdot) = {\int}_{\!\!\! X}^\oplus
\tilde{U}^{(n)}_{\pi}(\cdot)(x) d\mu(x)
\end{equation}
where $\tilde{U}^{(n)}_{\pi}(\widetilde{\mob})(x) \subset
\pi(\A)(x)$ and $\tilde{U}^{(n)}_{\pi}(\cdot)(x)$ is a strongly
continuous representation for almost every $x \in X$.

By Lemma \ref{mob.directint} a direct integral of strongly
continuous representations of $\widetilde{\mob}$ is of positive
energy if and only if the the representations appearing in the
integral are so for almost every point of the measure space of the
integration. As $C_{\pi}^{(n)}$ is a strongly continuous
one-parameter group in the center, for almost all $x \in
X:\;\tilde{U}^{(n)}_{\pi} (\tilde{\rho}_{(n\cdot)})(x) =
\tilde{U}^{(2)}_{\pi}(\tilde{R}_{(2\cdot)})(x)$ in the projective
sense. Therefore, since by the mentioned Lemma and Corollary
\ref{pos:2} in $\tilde{U}^{(2)}_{\pi}(\cdot)(x)$ the self-adjoint
generator of rotations is positive, also in
$\tilde{U}^{(n)}_{\pi}(\cdot)(x)$ it must be at least bounded from
below and hence by Lemma \ref{pos.mob^n} it is actually positive.
Thus, by using again Lemma \ref{mob.directint} we arrive to the
following result.
\begin{theorem}
\label{mainresult} Let $\pi$ be a locally normal representation of
the conformal local net of von Neumann algebras on the circle
$(\A,U)$. Then the strongly continuous representation
$\tilde{U}^{(n)}_\pi$ defined by Proposition \ref{n-cover}, is of
positive energy for all $n\in\NN^+$. In particular, the unique
continuous inner implementation of the M\"obius symmetry in the
representation $\pi$ is of positive energy.
\end{theorem}

Carpi proved \cite[Prop.\! 2.1]{Carpi2} that an irreducible
representation of a Virasoro net $\A_{{\rm Vir},c}$ must be one of
those that we get by integrating a positive energy unitary
representation of the Virasoro algebra (corresponding to the same
central charge) under the condition that the representation is of
positive energy. Thus by the above theorem we may draw the
following conclusion.
\begin{corollary}
An irreducible representation of the local net $\A_{{\rm Vir},c}$
must be one of those that we get by integrating a positive energy
unitary representation of the Virasoro algebra corresponding to
the same central charge.
\end{corollary}

\appendix
\chapter{About the \mo group}
\label{chap:A} \markboth{APPENDIX A. ABOUT THE M\"OBIUS GROUP}
{APPENDIX A. ABOUT THE M\"OBIUS GROUP}
\section{Geometrical properties}
\label{app:mobgroup} 

\paragrf
We define $\mob$ to be the set of diffeomorphisms of $S^1$ of the
form $\varphi_{a,b}: z \mapsto
\frac{az+b}{\overline{b}z+\overline{a}}$ with $a,b\in \CC$,
$|a|^2-|b|^2=1$. If $a,b,c,d\in\CC$ such that
$|a|^2-|b|^2=|c|^2-|d|^2=1$, then by setting
\begin{equation} \label{mobgroup.eq.1}
p\equiv (ac+b\overline{d}),\;\; q\equiv (ad + b\overline{c}),
\end{equation}
by direct calculation one can easily find that $|p|^2-|q|^2=1$ and
moreover that $\varphi_{a,b}\circ\varphi_{c,d} = \varphi_{p,q}$.
Thus $\mob$, with the composition, is a group.

For a differentiable transformation $\gamma$ of $S^1$ we 
think of its derivative as a $S^1\rightarrow \RR$ function:
\begin{equation}\label{deriv.of.transf}
\gamma'(z)\equiv -i\overline{\gamma(z)}
\frac{d}{d\theta}\gamma(e^{i\theta})|_{e^{i\theta}=z}.
\end{equation}
Using that $|a|^2-|b|^2=1$ we find that $|\varphi_{a,b}|>0$. Thus
$\mob$ contains only orientation preserving elements and hence
$\mob\subset \diff$.

If $a,b\in \CC$, $|a|^2-|b|^2=1$, then $|- a|^2-|- b|^2=1$ and
$\varphi_{-a,-b} =\varphi_{a,b}$. In the converse direction, if
$\varphi_{a,b} = \varphi_{c,d}$ then by elementary arguments
either $a=c,\,b=d$ or $a=-c,\,b=-d$. Let us denote this
equivalence relation by $\sim$ and set $M \equiv\{(a,b)\in\CC^2:
|a|^2-|b|^2=1\}$. Then the map $(a,b)\mapsto \varphi_{a,b}$ is a
well defined bijection from $M/\!\!\sim$ to $\mob$. Hence it is
easy to see that $\mob$ is actually a three dimensional Lie
subgroup of $\diff$. In fact the $(p,q)$ given by Eq.\!
(\ref{mobgroup.eq.1}) with which $\varphi_{a,b}\circ\varphi_{c,d}
= \varphi_{p,q}$, satisfies
\begin{equation}
\left(\begin{matrix}a & b\\ \overline{b} & \overline{a}
\end{matrix}\right)
\left(\begin{matrix}c & d\\ \overline{d} & \overline{c}
\end{matrix}\right)
=\left(\begin{matrix}p & q\\ \overline{q} & \overline{p}
\end{matrix}\right).
\end{equation}
Thus $\mob\simeq {\rm PSU}(1,1)\simeq {\rm PSL}(2,\RR)$. So one
can think of $\mob$ as a given action of ${\rm PSL}(2,\RR)$ on
$S^1$. We shall now study the geometrical properties of this
action.

\paragrf
We shall think of a vector field symbolically written as
$f(e^{i\vartheta})\frac{d}{d\vartheta}\in  {\rm Vect}(S^1)$ as the
corresponding real function $f$. In case of a vector field we  use
the notation $f'$ (calling it simply the derivative) for the
function on the circle obtained by differentiating with respect to
the angle:
$f'(e^{i\theta})=\frac{d}{d\alpha}f(e^{i\alpha})|_{\alpha=\theta}$.

Since the $S^1$ is compact, an $f:S^1 \rightarrow \RR$ smooth
function always gives rise to a flow on the circle. We denote by
$a\mapsto {\rm Exp}(af)$ this one-parameter group diffeomorphisms.
We say that $f\in {\rm Vect}(S^1)$ is a {\bf M\"obius vector
field}, if $a\mapsto{\rm Exp}(af)\in \mob$. Denoting by $\mathfrak
m$ the set of \mo vector fields, we can naturally identify
$\mathfrak m$ with the Lie algebra of $\mob$ in such a way that
${\rm Exp}$ will indeed correspond to the exponential of Lie
algebra elements. In this case however the Lie bracket will become
the negative of the usual bracket of vector fields. Thus for two
smooth functions $f,g$ we set
\begin{equation}
[f,g]\equiv f'g-fg'.
\end{equation}
In the complexification of ${\mathfrak m}$ the three functions
$l_{-1},l_0$ and $l_1$ form a base, where $l_j(z)\equiv z^j$
$(j=-1,0,1)$. The exponential of the constant $1$ function is the
one-parameter group of (anticlockwise) rotations $\alpha \mapsto
\rho_\alpha$ where $\rho_\alpha(z)\equiv e^{i\alpha}z$.

The \mo vector fields $t,d$ defined by the formulae
\begin{equation}
t(z) = \frac{1}{2} -\frac{1}{4}(z+z^{-1}),\;\;d(z)
=\frac{i}{2}(z-z^{-1})
\end{equation}
are called the translational and dilational vector fields,
respectively. The generated flow $a\mapsto \tau_a\equiv {\rm
Exp}(at)$ is called the {\bf one-parameter group of translations}.
The flow $s\mapsto \delta_s\equiv {\rm Exp}(sd)$ is called the
{\bf one-parameter group of dilations}. It is probably worth to
give them by formulae, too:
\begin{eqnarray}\nonumber \label{formXdilations}
\tau_a(z)&=& \frac{(2+ia)z-ia}{i a z+(2-ia)}, \\
\delta_s(z)&=& \frac{(e^s +1)z+ (e^s -1 )}{(e^s -1 )z+ (e^s +1)}.
\end{eqnarray}
The reason behind the names is the following. Choosing the point
$1\in S^1$ to be the so-called ``infinite point'', we can identify
$S^1\setminus\{1\}$ with the real line through the Cayley 
transformation defined by the map
\begin{equation}\label{Cayley}
x = i \frac{1+z}{1-z}\;\;\Longleftrightarrow z = \frac{x-i}{x+i}.
\end{equation}
The one-parameter groups $a\mapsto \tau_a$ and $s\mapsto \delta_s$
preserve the point $1$ and through the above identification they
go into the transformations $x\mapsto x+a$ and $x\mapsto e^s x$,
respectively.

By what was explained, it is evident that we have the following {\bf
scaling} of the translations by the dilations:
\begin{equation}
\label{scaling} \delta_s \tau_a \delta_{-s} = \tau_{e^s a}.
\end{equation}
This could have been also established by calculating that the
commutator of the generating vector fields resulting $[d,t]=t$.

\paragrf
It is also useful to investigate the transitivity of the action on
the circle. The next few lemmas will be stated without proofs;
they are justified by straightforward calculation.
\begin{lemma} \label{mob.3.trans}
Let $z_1,z_2,z_3\in S^1$ and $z'_1,z'_2,z'_3\in S^1$ two sets of
distinct points with the numbering corresponding to the
anticlockwise direction in both collections. Then there exists a
unique $\varphi \in \mob$ such that $\varphi(z_j)=z'_j$ for
$j=1,2,3$.
\end{lemma}
\begin{corollary}\label{mob.transit.onI}
$\mob$ acts transitively on the set $\I$ of nonempty, nondense
open arcs of $S^1$.
\end{corollary}
The one-parameter group of dilations preserves two points, namely,
the points $\pm 1$. It also follows, that it preserves the
intervals
\begin{equation}
S^1_\pm \equiv \{z\in S^1 : \pm {\rm Im}(z) > 0 \} \in \I.
\end{equation}
\begin{lemma} \label{uniq.dil.of.S+}
Any \mo transformation preserving $S^1_+$ is a dilation $\delta_s$
with a certain parameter $s\in\RR$. Thus up to parametrization
$s\mapsto \delta_s$ is the only one-parameter group in $\mob$
preserving $S^1$.
\end{lemma}
\begin{corollary}\label{uniq.dil.ofI}
If $\varphi_1,\varphi_2 \in\mob$ and
$\varphi_1(S^1_+)=\varphi_2(S^1_+)$ then it follows that $
\varphi_1\circ \delta_s\circ \varphi_1^{-1} = \varphi_2\circ
\delta_s\circ \varphi_2^{-1}$ for all $s\in\RR$.
\end{corollary}
\begin{proof}
By rearranging the equation we see that what we have to prove is
that $\varphi_2\circ\varphi_1$ commutes with $\delta_s$ which
follows since by the previous lemma $\varphi_2\circ\varphi_1$ is
again a dilation.
\end{proof}
\begin{definition}\label{def:dil.I}
For an interval $I\in\I$ the one-parameter group of dilations
associated to $I$, $s\mapsto \delta^I_s$ is defined by
$$
\delta^I_s \equiv \varphi\circ \delta_s\circ \varphi^{-1}
$$
where $\varphi \in\mob,\, \varphi(S^1_+) =I$. (By Corollary
\ref{mob.transit.onI} such a $\varphi$ exists, and by Corollary
\ref{uniq.dil.ofI} the definition is independent from the choice
of $\varphi$.)
\end{definition}
Note that up to parametrization $s\mapsto \delta^I_s$ is the
unique one-parameter group in $\mob$ preserving $I$. Note also
that by the above definition the dilational one-parameter group
$s\mapsto \delta^I_s$ is the dilational one-parameter group
associated to $S^1_+$. Moreover, as
\begin{equation}
\rho_\pi\circ\delta_s\circ\rho_\pi^{-1} = \delta_{-s}
\end{equation}
we find that $\delta^{I}_s = \delta^{I^c}_{-s}$ where $I^c$ is the
interior of the complement of $I$.
\begin{lemma}
\label{IRgenerators}
Let $I_1,I_2\in \I$ be two intervals with no common endpoints.
Then $\{\delta^K_s: K\in \{I_1,I_2\}, s\in\RR\}$ is a generating set in
$\mob$.
\end{lemma}
\begin{lemma}\label{3transit.of.vect.m}
Let $z_1,z_2,z_3\in S^1$ and $r_1,r_2,r_3\in\RR$ three real
numbers. Then there exists a unique \mo vector field
$f\in\mathfrak m$ such that $f(z_j)=r_j$ for $j=1,2,3$. Moreover,
if $r_1=r_2=0$ then this vector field is a multiple of $d^I$, that
is, the generator of dilations associated to $I$ where $I\in\I$ is
an interval with endpoints $z_1,z_2$.
\end{lemma}
\begin{lemma}\label{2diff->mob}
Let $\varphi_1,\varphi_2\in \mob$ and $z_0\in S^1$. Then
$$
\left.\begin{matrix}
\varphi_1(z_0)&=&\varphi_2(z_0)\\
\varphi_1'(z_0)&=&\varphi_2'(z_0)\\
\varphi_1''(z_0)&=&\varphi_2''(z_0)
\end{matrix}\right\} \;\;\Longleftrightarrow\;\; \phi_1 = \phi_2.
$$
\end{lemma}

\paragrf
The group $\mob$ is connected but not simply connected. Its
universal covering group $\widetilde{\mob}$ is isomorphic to the
universal covering of ${\rm PSL}(2,\RR)$; that is
$\widetilde{\mob} \simeq \widetilde{\rm SL}(2,\RR)$. Each
discussed one-parameter group lifts to a one-parameter group in
$\widetilde{\mob}$. These lifts will be denoted by simply adding a
``$\tilde{\hphantom{a}}$'', e.g.\! $\alpha \mapsto
\tilde{\rho}_\alpha$ will stand for the one-parameter group in
$\widetilde{\mob}$ obtained by lifting the one-parameter group of
rotations $\alpha\mapsto \rho_\alpha$. Note that while $\rho$ is a
periodical one-parameter group, $\tilde{\rho}$ is not. It follows
from the fact that as smooth manifolds (but not as groups) ${\rm
PSL}(2,\RR)\simeq S^1 \times \RR^2$ and $\widetilde{\rm SL}(2,\RR)
\simeq \RR^3$.

\begin{proposition}\label{iwasawa}{\bf (The Iwasawa decomposition.)} 
The map from $\RR^3$ to $\mob$ defined by the formula
$$
(\alpha,a,s) \mapsto \rho_\alpha \circ \tau_a \circ \delta_s
$$
is a smooth covering.
\end{proposition}

\noindent {\it Sketch of proof.} We shall (indicate how to) show
that for each $\varphi\in\mob$ there exists a unique $(z,s,a)\in
S^1\times \RR^2$ such that $\varphi =\rho_\alpha \circ\tau_a\circ
\delta_s $ where $e^{i\alpha}=z$ and that moreover that this
dependence is smooth.

Set $z\equiv \varphi(1)$. Then regardless what $a$ and $s$ are,
since the dilations and translations do not move the point $1\in
S^1$, we have that if $e^{i\alpha}=z$ then
\begin{equation}
(\rho_\alpha \circ \tau_a \circ \delta_s) (1) = \rho_\alpha(1) = z
= \varphi(1).
\end{equation}
Therefore $\rho_{-\alpha}\circ\varphi$ preserves the point $1$,
but not necessarily the point $-1$. In fact one sees that there
exists a unique $a\in\RR$ such that
$\rho_{-\alpha}\circ\varphi(-1)=\tau_a(-1)$. Then, as the
dilations do not move the point $-1$ we find that regardless what
$s$ is
\begin{equation}
(\rho_\alpha \circ \tau_a \circ \delta_s) (\pm 1) = \varphi(\pm
1).
\end{equation}
Therefore, by Corollary \ref{uniq.dil.of.S+} there exists a unique
$s\in\RR$ such that $\tau_{-a}\circ\rho_{-\alpha}\circ\varphi =
\delta_s$. Hence $\rho_\alpha \circ \tau_a \circ \delta_s =
\varphi$. \hfill $\Box$

\begin{corollary} \label{iwasawa.tilde}
The map from $\RR^3$ to $\widetilde{\mob}$
$$
(\alpha,a,s) \mapsto \tilde{\rho}_\alpha \tilde{\tau}_a
\tilde{\delta}_s
$$
is a diffeomorphism.
\end{corollary}

\section{Representation theory of $\widetilde{{\rm SL}(2,\RR)}$}
\label{app:mobrep}

\paragrf
In what follows $\mathbb T$ stands for the set of complex numbers
of unit length. We consider the {\bf projective unitary group}
$\U(\H)/\mathbb T$ of a Hilbert space $\H$ as a topological group
with the topology being induced by the strong operator topology on
$\U(\H)$. A {\bf strongly continuous projective representation}
$U$ of a topological group $G$ on a Hilbert space $\H$ is a
continuous homomorphism from $G$ to $\U(\H)/\mathbb T$. Sometimes
we think of $U(g)$ $(g \in G)$ as a unitary operator. Although
there are more than one ways to fix the phases, note that
expressions like Ad$(U(g))$ or $U(g) \in \M$ for a von Neumann
algebra $\M \subset {\rm B}(\H)$ are unambiguous.

Unlike with true representations, in general we cannot consider
the direct sum of projective representations. On the other hand,
we can still introduce the concept of invariant subspaces (it is
independent from the choice of phases) and thus that of
irreducible representations.

Suppose $U_1$ is obtained by fixing the phases of the projective
representation $U$ in a certain way. Then for each pair $(g,h)\in
G\times G$ there exists a unique $\omega_1(g,h) \in \mathbb T$
such that $U(g)U(h)=\omega_1(g,h) U(gh)$. Since $U(e)$, where $e$
is the identity of $G$, is a multiple of the identity, we have
$\omega_1(e,g) = \omega_1(g,e)$. Also, by the associativity one
finds that for all $g,h,k\in G$ the function $\omega$ satisfies
the following identity: $\omega_1(g,h) \omega_1(gh,k) =
\omega_1(g,hk)\omega_1(h,k)$. Thus one says that $\omega_1$ is a
$2$-cocycle on $G$. If $U_2$ is again obtained by fixing the
phases of $U$, then there exists a unique function $\eta_{2,1}: G
\rightarrow \mathbb T$ such that $U_2(g)=\eta_{2,1}(g)U_1(g)$ for
all $g\in G$. Then the $2$-cocycle $\omega_2$ associated to $U_2$
can be calculated by the formula $\omega_2(g,h)
=\eta_{2,1}(g)\eta_{2,1}(h)\overline{\eta_{2,1}(gh)}
\omega_1(g,h)$.

We can always find an open neighbourhood of the identity $O
\subset G$ and a certain fixing of the phases of $U$ in such a way
that the map $g\mapsto U(g)$ will be strongly continuous on $O$
and that the associated cocycle $\omega$ will be continuous on
$O\times O$. Indeed, take a unit vector $\Psi_0\in\H$. Whether for
a $g\in G$ the vector $U(g)\Psi_0$ is orthogonal to $\Psi_0$ is
independent from the phases. Let $O_1$ be the set where it is not
orthogonal. Then $O_1$ is an open neighbourhood of the identity.
In this neighbourhood for every $g\in G$ there is a unique way of
fixing the phase of $U(g)$ such that $\|U(g)\Psi_0 - \Psi_0\|$ is
minimal. It is clear that with this choice $g\mapsto U(g)$ is
strongly continuous on $O_1$ and that there exists an even smaller
neighbourhood of the identity $O$ such that $\omega$ is continuous
on $O\times O$.

Suppose $G$ is in particular a finite dimensional Lie group. It is
well known that if the Lie algebra of $G$ is semisimple, and $\omega$ is a 
$2$-cocycle on $G$ which is continuous in a neighbourhood of $(e,e)$, then 
$\omega$ is {\it exact} in (a certain smaller) neighbourhood of $(e,e)$. 
That is, there exists a $\eta: G \rightarrow \mathbb T$ continuous 
function such that $\eta(g)\eta(h)\overline{\eta(gh)} \omega(g,h) = 1$ for
all $(g,h)\in N \times N$ where $N$ is a certain open
neighbourhood of $e$. Thus with $\omega$ coming from $U$ where the
phases are fixed in the above explained way, we have that $g
\mapsto W(g)\equiv \eta(g) U(g)$ is strongly continuous on $N$ and
if $g,h \in N$ then $W(g)W(h)=W(gh)$. If $G$ is also connected and
simply connected, then $W$ extends to a strongly continuous
representation of $G$. These facts, as it was mentioned, are all well 
known and can be found in almost any book treating projective 
unitary representations; the original research article, establishing
the main ingredients of the theory is that \cite{bargmann} of Bargmann.

So suppose $G$ is a finite dimensional connected Lie group with
semisimple Lie algebra, and denote by $\tilde{G}$ its universal
cover. By what was said, if $U$ is a strongly continuous
projective unitary representation of $G$, then there exists a
strongly continuous unitary representation $\tilde{U}$ of
$\tilde{G}$ (on the same Hilbert space) such that in the
projective sense
\begin{equation}
\tilde{U}(\tilde{g}) = U(p(\tilde{g}))
\end{equation}
for all $\tilde{g}\in \tilde{G}$ with $p:\tilde{G}\rightarrow G$
being the covering map. If further $G$ is noncommutative with
simple Lie algebra, then such a $\tilde{U}$ is unique. Indeed, if
$\tilde{U}_1$ and $\tilde{U}_2$ are two such representations, then
there exists a continuous character $\lambda:\tilde{G}\rightarrow
\mathbb T$ such that
$\tilde{U}_2(\tilde{g})=\lambda(\tilde{g})\tilde{U}_1(\tilde{g})$
for all $\tilde{g}\in \tilde{G}$. But as $\tilde{G}$ is
noncommutative and its Lie algebra is simple, the only continuous
character it admits is the identity function.

\paragrf
Our particular interest is the (projective, positive energy)
representation theory of $\mob^{(n)}$ where $n$ is a positive
integer. The group in question is a connected Lie group with
simple Lie algebra. Moreover, there is a canonical $n$-fold
covering map from $\mob^{(n)}$ to $\mob$. Thus the universal cover
of $\mob^{(n)}$ can be canonically identified with
$\widetilde{\mob} \simeq \widetilde{{\rm SL}(2,\RR)}$, that is,
with the universal cover of $\mob$. So every strongly continuous
projective unitary representation of $\mob^{(n)}$ lifts to a
unique strongly continuous unitary representation of
$\widetilde{\mob}$. This is why we shall investigate now what we
can say about strongly continuous unitary representations of
$\widetilde{\mob}$. The following statement is demonstrated by an
adopted (and slightly modified) version of the argument used in
the proof of \cite[Prop.\! 1]{koester02}.
\begin{proposition}
\label{pos.mob^n} Let $\tilde{U}$ be a strongly continuous unitary
representation of $\widetilde{\mob}$ with $H$ and $P$ being the
self-adjoint generator of rotations and translations in
$\tilde{V}$, respectively.
Then the following four conditions are equivalent: \\
\indent $1.\;\,H$ is bounded from below, \\
\indent $2.\;\,P$ is bounded from below, \\
\indent $3.\;\,H \geq 0$, \\
\indent $4.\;\,P \geq 0$.
\end{proposition}
\begin{proof}
Recall that $\tilde{\rho}$ is the rotational one-parameter group
in $\widetilde{\mob}$ and set $P_\pi=\tilde{U}(\tilde{\rho}_\pi) P
\tilde{U}(\tilde{\rho}_\pi)^*$. Then $P_\pi$ is the self-adjoint
generator associated to the one-parameter group $a\mapsto
\widetilde{\rm Exp}(a t_\pi)$ where the vector field $t_\pi$ is
obtained by rotating the vector field of translations $t$, by
$\pi$ radian; that is $t_\pi(z)=t(-z)=
\frac{1}{2}+\frac{1}{4}(z+z^{-1})$. As $P_\pi$ is unitary
conjugate to $P$ their spectra coincide. Moreover, as $t+t_\pi =1$
on the G{\aa}rding-domain we have that $P+P_\pi=H$ which
immediately proves that if $P$ is bounded from below or positive
then so is $H$.

As for the rest of the statement, apart from the trivial
indications there remains to show that if $H$ is bounded from
below then $P$ is positive. Consider the lifted dilations $s
\mapsto \tilde{\delta}_s$. By equation (\ref{scaling}) one has
that $\tilde{U}(\tilde{\delta}_s) P \tilde{U}(\tilde{\delta}_s)^*
= e^s P$. Moreover, by direct calculation $[d,t_\pi]=-t_\pi$ so
similarly to the case of translations the dilations also ``scale''
$t_\pi$, but in the converse direction. Thus
$\tilde{U}(\tilde{\delta}_s) P_\pi \tilde{U}(\tilde{\delta}_s)^* =
e^{-s} P_\pi$. So if $H\geq r \mathbbm 1$ for some $r$ real (but
not necessarily nonnegative) number then for any vector $\xi$ in
the G{\aa}rding-domain, setting $\eta =
\tilde{U}(\tilde{\delta}_s)^*\xi$ we have that
\begin{eqnarray}
\label{ineq} r \|\xi\|^2 = r \|\eta\|^2 &\leq& \langle\eta,H\eta\rangle\, 
=\, e^s \langle\xi,P\xi\rangle + \, e^{-s} \langle\xi,P_\pi \xi\rangle
\end{eqnarray}
from which, letting $s\rightarrow \infty$ we find that $P \geq 0$.
\end{proof}
\begin{definition}
A strongly continuous unitary representation of $\widetilde{\mob}$
that satisfies any of the listed conditions of Prop.\!
\ref{pos.mob^n} is called a {\bf positive energy representation}
of $\widetilde{\mob}$.
\end{definition}
A useful observation about the G{\aa}rding-domain in connection
with the self-adjoint generator of rotations is the following
\begin{lemma}\label{GardingforH^k}
Let $\tilde{U}$ be a strongly continuous unitary representation of
$\widetilde{\mob}$ with $H$ being the self-adjoint generator of
rotations. Then the G{\aa}rding-domain of the representation is a
core for any positive power of $H$.
\end{lemma}
\begin{proof}
Let $k$ be a positive number. The Iwasawa decomposition
\begin{equation}
x,y,z \mapsto \tilde{\rho}_x  \tilde{\tau}_y \tilde{\delta}_z
\end{equation}
by Corollary \ref{iwasawa.tilde} is a diffeomorphism between
$\RR^3$ and $\widetilde{\mob}$. Hence if $\Psi\in\D(H^k)$ and
$f:\RR^3 \rightarrow \CC$ is a compactly supported smooth function
then
\begin{equation}
\Psi_f \equiv \int_{\RR^3} f(x,y,z) \tilde{U}(\tilde{\rho}_x
\tilde{\delta}_y \tilde{\tau}_z)\Psi \; dx\,dy\,dz.
\end{equation}
is in the G{\aa}rding-domain. (It does not matter that the measure
is not the Haar-measure: a replacement would only make a change in
$f$.) From here it is easy to see that if $g:\RR\rightarrow \CC$
is compactly supported and smooth then $\int_{\RR} g(x)
\tilde{U}(\tilde{\rho}_x \Psi \, dx$ is in the domain of the
operator obtained by closing the restriction of $H^k$ to the
G{\aa}rding-domain. The rest is trivial.
\end{proof}
A unitary representation $U$ of a group $G$ is called a {\bf
factor representation}, if the von Neumann algebra $U(G)''$ is a
factor. For the next statement recall the concept of direct
integrals (see e.g.\! \cite[Chapter 14.]{kadison} for direct
integral of operators, representations of $C^*$ algebras, and
direct integrals of von Neumann algebras and \cite[Remark
18.7.6]{dixmier} for the connection between direct integral of
representations of $C^*$ algebras and that of groups).

Although for any kind of direct integral construction one usually
requires separability for the Hilbert spaces in question, note
that in our case this condition can be dropped. The point is that
we can choose a countable set ${\mathcal S}\subset
\widetilde{\mob}$ which is dense. Then if $\tilde{U}$ is a
strongly continuous representation of our group then for any given
vector $\Psi$ we have that the invariant space
\begin{equation}
\overline{{\rm Span}\{\tilde{U}(\widetilde{\mob})\Psi}\} =
\overline{{\rm Span}\{\tilde{U}({\mathcal S})\Psi}\}
\end{equation}
is separable. From here it is easy to see that any strongly
continuous representation is (equivalent with) a direct sum of
representations given on separable Hilbert spaces.

 The group $\widetilde{\mob}$ is in particular second countable
and locally compact and so the $C^*$ algebras associated to it is
separable. Therefore, by the mentioned remark of the above cited
book, any strongly continuous unitary representation of
$\widetilde{\mob}$ can be decomposed as a direct integral of
factor representations. (Note that in our case the condition on
the separability of the Hilbert space can be dropped: since one
can find a countable set of points which is dense in
$\widetilde{\mob}$, by standard arguments, any strongly continuous
unitary representation of the group in question can be decomposed
into the direct sum of representations on separable Hilbert
spaces. Then on each component we can further decompose the
representation into a direct integral of factor representations.
Note also that in general for direct integral constructions there
is no need to require separability of the {\it total} Hilbert
space, what is necessary is the separability of the Hilbert spaces
that are to be ``integrated together''. Thus the integral
decompositions in the components appearing in the above direct sum
can be indeed put together to give a direct integral decomposition
of the original representation into factor representations.)
\begin{lemma}\label{mob.directint}
Every strongly continuous unitary representation of
$\widetilde{\mob}$ is (equivalent with) a direct integral of
factor representations. Moreover, a direct integral of strongly
continuous unitary representations
$$
\tilde{U}(\cdot) = {\int}_{\!\!\! X}^\oplus \tilde{U}(\cdot)(x)
d\mu(x)
$$
of $\widetilde{\mob}$ is of positive energy if and only if so is
$\tilde{U}(\cdot)(x)$ for almost every $x \in X$.
\end{lemma}
\begin{proof}
The first statement was just demonstrated (by the comments made
above). There remains to justify the affirmation about positivity.
Since positivity of energy regards only the restriction of the
representation onto the rotational subgroup, the argument will be
presented about one-parameter unitary groups, in general.

For a $t\mapsto V(t)$ strongly continuous one-parameter group of
unitaries the positivity of the self-adjoint generator is for
example equivalent with the fact that $\hat{V}(f)\equiv\int
V(t)f(t)dt = 0$ for a certain smooth, fast decreasing function $f$
whose Fourier transform is positive on $\RR^-$ and zero on
$\RR^+$. If $V$ is a direct integral of a measurable family of
strongly continuous one-parameter groups,
$V(\cdot)=\int^\oplus_{\! X} V(\cdot)(x) d\mu(x)$, then
$\hat{V}(f) =\int^\oplus_{\! X} \hat{V}(f)(x) d\mu(x)$. As
$\hat{V}(f)(x) \geq 0$ for almost every $x\in X$, the operator
$\hat{V}(f)$ is zero if and only if $\hat{V}(f)(x)=0$ for almost
every $x\in X$.
\end{proof}

Still considering a strongly continuous unitary representation
$\tilde{U}$ of $\widetilde{\mob}$, let us denote, as before, the
self-adjoint generator of rotations by $H$, the self-adjoint
generator of translations by $P$ and further, the self-adjoint
generator of dilations by $D$. Set
\begin{equation}
E_\pm \equiv (H-2T)\pm i D,
\end{equation}
then $E_\pm \subset E_\mp^*$ and by looking at the commutation
relation between the corresponding (complex) vector fields we find
that on the G{\aa}rding-domain
\begin{equation}\label{HEpm.comm.rel.}
[E_\pm, H] = \pm E_\pm,\;\; [E_-,E_+] = 2H.
\end{equation}
Further, by setting
\begin{equation}\label{Casimir}
C\equiv H^2 -H- E_+E_-,
\end{equation}
we find that $C$ is a Casimir-operator: on the G{\aa}rding-domain
it commutes with all of the $3$ generators $H,P$ and $D$. Since it
is a polinomial (quadratic) expression of the generator, and
moreover its restriction to the G{\aa}rding-domain is evidently
symmetric, it is easy to prove that the closure of this
restriction (which, by symmetricalness exists) is indeed a closed
operator commuting with both $\tilde{U}(\widetilde{\mob})$ and
$\tilde{U}(\widetilde{\mob})'$. Hence if the
representation $\tilde{U}$ is a factor representation, then there
exists a real number $\lambda \in \RR$ such that $C\Psi
 = \lambda\Psi$ for every vector $\Psi$ in
the G{\aa}rding-domain; we shall say that $\lambda$ is the ``value
of the Casimir-operator $C$ in the factor representation
$\tilde{U}$''.
\begin{proposition}
\label{factor.rep.of.mob} Let $\tilde{U}$ be a positive energy
factorial representation of $\widetilde{\mob}$, $\lambda$ the
value of the Casimir-operator defined above and $h={\rm
min\{Sp}(H)\}$. Then $\lambda = h^2-h$ and ${\rm Sp}(H)\subset
(h+\NN) \cup (-h+\NN)$.
\end{proposition}
\begin{proof}
Let $k$ be a positive integer and $q$ a polynomial with
$q([h,\infty])\subset \RR^+_0$. Then using the commutation
relations, the value of $C$ and the fact that $q(H)\geq 0$ we find
that for any $\Psi$ in the G{\aa}rding-domain
\begin{equation} \label{mobrep.eq.1}
0\leq\; \langle (E_-)^k \Psi,\,q(H) (E_-)^k \Psi\rangle\, =\,
\langle\Psi,\, p_{q,k}(H)\Psi\rangle
\end{equation}
where the polynomial $p_{q,k}$ is given by the formula
\begin{equation}
p_{q,k}(x) = q(x-k) \mathop{\Pi}_{j=0}^{k-1}
((x-j)^2-(x-j)-\lambda).
\end{equation}
Let $\Phi$ be any vector whose (nonnegative) measure $\mu$ given
by the spectral decomposition of $H$ is compactly supported. Then
on one hand we have
\begin{equation}
\langle\Phi,\,p_{q,k}(H) \Phi\rangle \; = \int p_{q,k}(x)
\,d\mu(x)
\end{equation}
but on the other hand Eq.\! (\ref{mobrep.eq.1}) and a simple use
of Lemma \ref{GardingforH^k} show that $\langle\Phi,\,p_{q,k}(H)
\Phi\rangle \,\geq 0$. This should be true for all positive
integer $k$ and polynomial $q$ with $q([h,\infty])\subset\RR^+_0$.
From here, by setting $k=1$, it is a standard, although not
trivial exercise to show that $h^2-h-\lambda = 0$. (One needs to
take vectors $\Phi_\epsilon\neq 0$ with measure supported in
$[h,h+\epsilon]$ with $\epsilon>0$ sufficiently small, and
consider for example the polynomial $q(x)=x-h$.) Moreover, if a
vector is supported in $[h,h+n]$ $(n\in\NN)$, then by setting
$k=n$, and using similar considerations one sees (by choosing the
polynomial $q$ ``well'')  that the corresponding measure must be
concentrated on the set of points at which the product
$\mathop{\Pi}_{j=0}^{k-1} ((x-j)^2-(x-j)-\lambda)$ is zero. Hence,
using that $\lambda=h^2-h$ it follows that ${\rm Sp}(H)\subset
(h+\NN) \cup (-h+\NN)$.
\end{proof}
\begin{corollary}\label{dom.of.Epm}
Let $\tilde{U}$ be a positive energy representation of
$\widetilde{\mob}$ with $H,E_\pm$ as before. Then
$\D(E_\pm)\subset\D(H)$ and
\begin{eqnarray*}
\|E_- \Psi\| &\leq& \| \sqrt{(H^2-h^2)} \Psi\| \leq  \|H \Psi \|,
\\ \|E+ \Psi\|&\leq& \|(\mathbbm 1 + H) \Psi \|.
\end{eqnarray*}
for every $\Psi \in \D(H)$, where $h$ is the minimum point of the
spectrum of $H$. So if $Q$ is any self-adjoint generator of
$\tilde{U}$ then $\D(Q)\subset \D(H)$ and there exists a constant
$r>0$ such that $\|Q \Psi\|\leq r \|(\mathbbm 1 + H) \Psi \|$ for
every $\Psi\in\D(H)$. In particular any core for $H$ is a core for
$Q$.
\end{corollary}
\begin{proof}
Once the stated energy bounds are satisfied for $E_\pm$ on the
G{\aa}rding domain, it is clear that the rest follows. If
$\tilde{U}$ is a factor representation with $\tilde{h}\geq h$
being the lowest point of ${\rm Sp}(H)$, then by our previous
proposition, taking account of the definition of $C$ we find that
on the G{\aa}rding-domain
\begin{eqnarray}\nonumber
E_+E_- &=& H^2-H-\tilde{h}^2+\tilde{h} \leq H^2-h^2,\\
E_-E_+ &=&  [E_-,E_+] + E_+E_-= H^2+H-\tilde{h}^2+\tilde{h} \leq
(\mathbbm 1 + H)^2.
\end{eqnarray}
Thus the bounds are justified by taking expectations with vectors
from the G{\aa}rding-domain. Hence the assertions for the domains
follow and the bounds hold on the whole domain of $H$.

If $\tilde{U}$ is not a factor representation, then by Lemma
\ref{mob.directint} it is a direct integral $\tilde{U}(\cdot) =
{\int}_{\!\!\! X}^\oplus \tilde{U}(\cdot)(x) d\mu(x)$ of factor
representations. It is also clear that the lowest point of the
spectrum of the self-adjoint generator of rotations $H(x)$ in the
factor representation $\tilde{U}(\cdot)(x)$, is greater or equal
to $h$ for almost every $x\in X$ (one can argue in a similar way
as in the proof of the mentioned lemma). Then the bounds hold in
almost all factor representation appearing in the direct integral
decomposition. To finish, all we need to do is to check that the
bounds on the generators can be ``integrated'' to give the right
inequalities for the direct integral decomposition; which is left
for the reader.
\end{proof}
\begin{corollary}\label{comm.on.Dinfty}
Let $\tilde{U}$ be a positive energy representation of
$\widetilde{\mob}$ with $H,E_\pm$ as before. Then $\D^\infty
\equiv \cap_{n\in\NN}$ is an invariant domain for $E_\pm$ and
hence it is an invariant core for any self-adjoint generator of
$\tilde{U}$. Moreover, not only on the G{\aa}rding-domain, but
even on $\D^\infty$ the generators satisfy the commutation
relations of the Lie algebra of $\widetilde{\mob}$.
\end{corollary}
\begin{proof}
By an easy use of the commutation relations one has that if $\Psi$
is a vector from the G{\aa}rding domain then $H^k E_- \Psi=
E_-(H-\mathbbm 1)^k \Psi$ and $H^k E_+ \Psi= E_+(H+\mathbbm
1)^k\Psi$ for all $k\in\NN^+$. So let $\Phi\in\D^\infty$. Then by
Lemma \ref{GardingforH^k}  there exists a sequence of G{\aa}rding
vectors $\Psi_n\; (n\in\NN)$ converging to $\Phi$ such that
$H^q\Psi_n \to H^q\Psi_n$ for any positive power $q\leq k+1$.
Then, since $E_\pm$ are closable operators, by the above equations
and the previous corollary we have that
\begin{eqnarray}\nonumber
E_\pm \Psi_n &\to& E_\pm \Phi \\
H^k E_\pm \Psi_n = E_\pm(H\pm\mathbbm 1)^k\Psi_n &\to&
E_\pm(H\pm\mathbbm 1)^k\Phi.
\end{eqnarray}
Hence $E_\pm \Phi \in \D(H^k)$ and $H^k E_\pm \Phi =
E_\pm(H\pm\mathbbm 1)^k\Phi$ and the invariance of $\D^\infty$
follows since $k$ was any positive integer. As for the commutation
relations, we need to consider quadratic expressions like $Q_1Q_2$
where $Q_1,Q_2$ are generators of $\tilde{U}$. On the G{\aa}rding
domain $[H,Q_2] = i Q_3$ where $Q_3$ is another generator of the
representation and by the previous corollary there exists an
$r_1>0$ and $r_3>0$ such
\begin{eqnarray} \nonumber
\|Q_1Q_2\Psi\| &\leq& r_1 \|(\mathbbm 1 + H)Q_2\Psi\| = r_1
\|(Q_2(\mathbbm 1 + H)+ iQ_3)\Psi\| \\
\nonumber &\leq&r_1^2 \|(\mathbbm 1 + H)^2\Psi\| + r_3 \|(\mathbbm
1 + H)\Psi\| \\ &\leq& (2r_1^2+2r_3) \|(\mathbbm 1 + H^2)\Psi\|
\end{eqnarray}
for any vector $\Psi$ from the G{\aa}rding domain. The rest
follows easily (again, using Lemma \ref{GardingforH^k}).
\end{proof}
\begin{corollary}\label{UmobD^infty=D^infty}
Let $\tilde{U}$ be a positive energy representation of
$\widetilde{\mob}$ with $H$ and $\D^\infty$ as before. Then
$\tilde{U}(\tilde{\gamma})\D(H^n) = \D(H^n)$ for all $n\in\NN_+$
and so $\tilde{U}(\tilde{\gamma})\D^\infty = \D^\infty$
$(\tilde{\gamma}\in\widetilde{\mob})$.
\end{corollary}
\begin{proof}
In the previous proof we have seen that an expression like
$Q_1Q_2$ admits a quadratic energy bound. In fact, it is not hard
to see that in general for a self-adjoint generator $Q$ of
$\tilde{U}$ there exists an $r_n>0$ such that $\|Q^n\|\Psi\leq
r_n\|(\mathbbm 1 + H^n)\Psi\|$ for all $\Psi\in\D^\infty$ and so
$\D(H^n)\subset \D(Q^n)$.

Suppose $\Phi\in\D(H^n)$ and $\tilde{\gamma}\in\widetilde{\mob}$.
Then $Q\equiv \tilde{U}(\tilde{\gamma})^*H
\tilde{U}(\tilde{\gamma})$ is a self-adjoint generator of
$\tilde{U}$, hence by what was established $\Phi\in\D(Q^n)$. But
it is equivalent with saying that
$\tilde{U}(\tilde{\gamma})\Phi\in\D(H^n)$. We have thus proved
that $\tilde{U}(\tilde{\gamma}) \D(H^n)\subset \D(H^n)$. Then
exchanging $\tilde{\gamma}^{-1}$ with $\tilde{\gamma}$, we see
that in fact $\tilde{U}(\tilde{\gamma}) \D(H^n)= \D(H^n)$.
\end{proof}

\paragrf
We shall now proceed in the converse direction. Namely, in stead
of considering representations of the group and then deriving
relations for the generators, we shall consider representations of
the generators of $\widetilde{\mob}$ and we shall try to integrate
them.
\begin{proposition}
\label{energybound:mob} Let $\fin\subset \H$ be a dense subspace
of a Hilbert space and $h,e_+,e_-$ three operators defined on
$\fin$ satisfying the
following properties: \\
\indent $1.$ $\fin$ is invariant for all three operators, \\
\indent $2.$ $h$ is diagonal and positive on $\fin$\\
\indent $3.$ $e_\pm \subset e_\mp^*$ \\
\indent $4.$ $[e_\pm, h] = \pm e_\pm,\;\;[e_-,e_+] = 2h$. \vspace{2mm}\\
Then $\|e_- \Psi\|\leq \|h \Psi \|$ and $\|e_+ \Psi\|\leq
\|({\mathbbm 1}+h) \Psi \|$ for every $\Psi\in\fin$. Moreover, if
$q$ is a linear combination of $h,e_+$ and $e_-$, then there
exists a $t>0$ such that
$$
\sum_{n\in\NN} t^n\frac{\|q^n \Psi\|}{n!} < \infty
$$
for every $\Psi\in\fin$.
\end{proposition}
\begin{proof}
By the second condition $\fin$ is the span of the eigenvectors of
$h$ (every vector has a unique decomposition into eigenvectors
associated to different eigenvalues of $h$) and by the commutation
relations $e_\pm$ are again ``stepping'' operators.

Since $h$ is positive, the set of eigenvalues of $h$ has a
nonnegative infinum which we shall denote by $q$. Suppose
$\Phi\in\fin$ contains only eigenvectors of $h$ with eigenvalues
in $[q,q+1)$. Then $e_-\Phi = 0$ and so $\|e_-\Phi\|\leq
\|h\Phi\|$ and $\|e_+\Phi\|^2 = \langle\Phi,\,2h\Phi\rangle \,\leq
\,\|(\mathbbm 1 + h)\Phi\|^2$. We shall proceed by induction.
Suppose these inequalities hold for all $\Phi\in\fin$ containing
eigenvectors of $h$ with eigenvalues in $[q+n-1,q+n)$ for some
$n=k$ positive integer. Suppose $\Phi$ is a combination of
eigenvectors with eigenvalues in $[q+k,q+k+1)$. Then $e_-\Phi$
belongs to $\D_{[q+k-1,q+k)}$, i.e.\! it is a combination of
eigenvectors with eigenvalues in $[q+k-1,q+k)$. Thus, setting $X_k
\equiv [q+k-1,q+k)$ we have that
\begin{eqnarray} \nonumber
\|e_-\Phi\| &=& \mathop{\rm Sup}_{\Psi\in\D_{X_k},\|\Psi\|=1}
\{|\langle\Psi,\,e_\Phi\rangle |\} = \mathop{\rm
Sup}_{\Psi\in\D_{X_k},\|\Psi\|=1}\{ |\langle e_+\Psi,\,\Phi\rangle
|\}
\\ \nonumber &\leq& \|\Phi\| \mathop{\rm Sup}_{\Psi\in\D_{X_k},\|\Psi\|=1}\{
\|e_+\Psi\|\} \\
\nonumber &\leq& \|\Phi\|\mathop{\rm
Sup}_{\Psi\in\D_{X_k},\|\Psi\|=1}\{ \|(\mathbbm 1+h)\Psi\|\} \leq
(1+(q+k))\|\Phi\| \\
&\leq& \|h\Phi\|
\end{eqnarray}
where in the third line we have used the inductive assumption.
Then, using the commutation relations, one also finds that
\begin{eqnarray}\nonumber
\|e_+\Phi\|^2 &=& \langle\Phi,\,2h\Phi\rangle + \|e_-\Phi\|^2 \\
\nonumber &\leq& 2(q+k+1)\|\Phi\|^2 + \|h\Phi\|^2 \leq
(2(q+k+1)+(q+k+1)^2)\|\Phi\|^2 \\ &\leq& \|(\mathbbm 1 +h)\Phi\|.
\end{eqnarray}
Thus those inequalities hold for all vectors from the sets
$\D_{X_n}\;(n\in\NN^+)$.

Every vector from $\fin$ can be written as a sum
$\sum_{n=1}^{\infty}\Phi_n$ where $\Phi_n\in\D_{X_n}$ and only
finite many of these vectors are non zero. These vectors are
pairwise orthogonal since $h$ is hermitian. But by the discussed
``stepping'' nature of $e_\pm$, also the vectors $e_+\Phi_n$ are
pairwise orthogonal and the same holds for the vectors
$e_-\Phi_n$. Thus by orthogonality the inequality is justified for
all vectors of $\fin$.

Let us proceed to the second part of the statement. If $q$ is a
linear combination of the operators $h,e_+$ and $e_-$ then by the
justified bounds there exists an $r>0$ such that $\|q \Psi\|\leq r
\|({\mathbbm 1}+h) \Psi \|$ for all $\Psi\in\fin$. Moreover, if
the eigenvalues appearing in the decomposition of $\Psi$ are
smaller than a certain positive integer $k$, then the eigenvalues
appearing in that of $q\Psi$ are smaller than $k+1$. Therefore, if
$n$ is an integer then
\begin{eqnarray}\nonumber
\|q^n\Psi\|&\leq& r \|(\mathbbm 1+h)q^{n-1}\Psi\|\leq r (k+n)\|q^{n-1}\Psi\| \\
\nonumber &\leq& r^2 (k+n)\|(\mathbbm 1+h)q^{n-2}\Psi\| \leq r^2
(k+n)(k+n-1)\|q^{n-2}\Psi\| \\
\nonumber && \ldots \\
&\leq& r^n \frac{(k+n)!}{k!} \|\Phi\|
\end{eqnarray}
and hence if $t< 1/r$ then
\begin{eqnarray} \nonumber
\sum_{n\in\NN} t^n\frac{\|q^n \Psi\|}{n!} &\leq& \sum_{n\in\NN}
(rt)^n \frac{(k+n)!}{k!n!} \|\Psi\| \\
&\leq& \frac{1}{k!} \|\Psi\| \sum_{n\in\NN} (rt)^n (k+n)^k <
\infty
\end{eqnarray}
which is exactly what we wanted to prove. Note that the convergence
radius depends only on $r$ (i.e.\! on the linear energy bound on the 
element $q$) but not on the vector $\Psi$.
\end{proof}
\begin{corollary}
\label{he+e-.int.repr.} If $h,e_+$ and $e_-$ are as in the
previous proposition, then they integrate to a positive energy
representation of $\widetilde{\mob}$. That it,
$$
h, \;\, t=\frac{1}{2}h-\frac{1}{4}(e_+ + e_-),\;\,
d=\frac{i}{2}(e_- -e_+)
$$
are essentially self-adjoint operators and setting $H,T,D$ for
their closure we have that the map
$$
\tilde{\rho}_\alpha \tilde{\delta}_s \tilde{\tau}_a \mapsto
e^{i\alpha H}e^{isD}e^{iaT}
$$
is a positive energy representation of $\widetilde{\mob}$.
\end{corollary}
\begin{proof}
By the previous proposition $\fin$ is a dense set of analytic
vectors for the symmetric operators $h,t,d$ (and in fact for any
linear combination of them). Thus by results of Nelson (see e.g.\!
\cite[Theorem X.39]{RSII}) they are essentially self-adjoint. Moreover,
the group $\widetilde{\mob}$ is a connected and simply connected
Lie group. Then, by the results of \cite{anal.int} 
on the integrability of the 
representations Lie algebras, the statement follows because we have a 
representation of its Lie algebra by antisymmetric operators defined on a 
common, dense, invariant, analytic set of vectors.
\end{proof}

\paragrf
Let us now return to positive energy representations of
$\widetilde{\mob}$ in which the spectrum of $H$ contains
eigenvalues only. (We have seen that this is the case if for
example the representation is factorial.) Let $\Phi$ be an
eigenvector associated to an eigenvalue $q$, then using the
commutation relations one finds that $H (E_-\Phi) = (q-1)
(E_-\Phi)$ and $H (E_+\Phi) = (q+1) (E_-\Phi)$, thus $E_\pm$
``steps up and down'' on the spectrum of $H$. In particular if
$\Phi$ is an eigenvector associated to the lowest eigenvalue
$h\geq0$, then $E_-\Phi= 0$ and using the commutation relations
$\|E_+\Phi\|^2 = 2h\|\Phi\|^2$. Thus if $h=0$ then $\Phi$ is
annihilated by all generators and it is in an invariant vector for
the representation (i.e.\! the representation contains the trivial
representation). So suppose $h>0$, $\|\Phi\|=1$ and set
$\Phi_h\equiv\Phi$ and recursively
\begin{equation} \label{recursive.phi.n}
\Phi_{h+n+1}\equiv \frac{1}{\sqrt{(2h+n)(n+1)}} E_+\Phi_{h+n}.
\end{equation}
By the mentioned ``stepping'' property
$H\Phi_{h+n}=(h+n)\Phi_{h+n}$, and by the expression for the
Casimir-operator one finds (by an inductive argument) that these
vectors have all unit length. So in fact they are a set of
orthonormal eigenvectors for $H$. Moreover, again by an inductive
argument, using the commutation relations, and the fact that
$E_-\Phi= 0$, one finds that $E_-\Phi_{h+n} =
\sqrt{(2h+n-1)n}\,\Phi_{h+n-1}$. So the span of these vectors is
an invariant subspace for the generators.

Using the above observations we shall now construct all
irreducible positive energy representations of $\widetilde{\mob}$.
Let $h>0$ be a positive number and consider the Hilbert space
$l^2(\NN)$ whose standard complete orthonormed system --- somewhat
unusually --- shall be denoted by $\{\xi_{h+n}:n\in\NN\}$.
Further, denote by $\fin$ the algebraic span of this orthonormed
system. Following the observations made, we introduce the three
operators $h,e_+$ and $e_-$ on $\fin$ by the defining them on the
base $\{\xi_{h+n}:n\in\NN\}$ as follows:
\begin{eqnarray}\nonumber
h \xi_{h+n} &\equiv& (h+n)\xi_{h+n}, \\
\nonumber e_+ \xi_{h+n} &\equiv& \sqrt{(2h+n)(n+1)}\,\xi_{h+n+1} \\
e_- \xi_{h+n} &\equiv& \sqrt{(2h+n-1)n}\, \xi_{h+n-1}.
\end{eqnarray}
(Note that the definition for $e_-$ is meaningful even when $n=0$
as in that case the coefficient is zero.) It is an elementary
exercise to check that these operators then satisfy the
requirements of Prop.\! \ref{energybound:mob} and hence by
Corollary \ref{he+e-.int.repr.} it integrates to a positive energy
representation of $\widetilde{\mob}$ which we shall denote by
$\tilde{U}_h$. Note that the spectrum of the self-adjoint
generator of rotations for this representation is exactly $h+\NN$.
\begin{proposition}
The representations $\{\tilde{U}_h: h>0\}$ are nonequivalent
irreducible representations. Moreover any nontrivial irreducible
positive energy representation of $\widetilde{\mob}$ is equivalent
with $\tilde{U}_h$ for some $h>0$.
\end{proposition}
\begin{proof}
That these representations are inequivalent is trivial, since the
lowest eigenvalue of the self-adjoint generator of rotations in
$\tilde{U}_h$ is exactly $h$. As for the irreducibility, suppose
$X\in \tilde{U}_h(\widetilde{\mob})'$. Then $X$ preserves the
eigenspaces of the self-adjoint generator of rotations, thus for
every $n\in\NN$ there exists a complex number $x_n\in\CC$ such
that $X\xi_{h+n} = x_n \xi_{h+n}$. Also, from the commutation with
the self-adjoint generator of translations and dilations it
follows that $X e_+ \xi_{h+n} = e_+ X\xi_{h+n}$. Therefore
\begin{equation}
x_{n+1}\sqrt{(2h+n)(n+1)} = x_n \sqrt{(2h+n)(n+1)}
\end{equation}
and hence $x_{n+1}=x_n$ which implies that $X$ is a multiple of
the identity. Finally, as for the last assertion, suppose
$\tilde{U}$ is an irreducible positive energy representation of
$\widetilde{\mob}$. Denote the lowest point of the spectrum of its
self-adjoint generator $H$ by $h$. As $\tilde{U}$ is in particular
a factor representation, by Prop.\! \ref{factor.rep.of.mob} there
is an eigenvector of $H$ with eigenvalue $h$. As it was already
observed, $h$ cannot be zero as then the corresponding eigenvector
would be an invariant vector for $\tilde{U}$. Let $\Phi_{h}$ be a
normalized eigenvector of $H$ associated to the eigenvalue $h>0$.
By the same observation, if define the vectors $\Phi_{h+n}$ by the
recursive formula of Eq.\! (ref{recursive.phi.n}) then the map
$\xi_{h+n}\mapsto \Phi_{h+n}$ extends to an isometric operator
$V$. Moreover, by denoting the self-adjoint generators of
translations and dilations in $\tilde{U}$ by $T$ and $D$,
respectively, we have that $h\subset VHV^*,\,t\subset VTV^*$ and
$d\subset VDV^*$. Therefore the span of $\{\Phi_{h+n}:n\in\NN\}$
is a set of invariant, analytic vectors for the self-adjoint
generators of $\tilde{U}$. Hence its closure, which is the range
of $V$ is an invariant subspace for $\tilde{U}$ and thus by
irreducibility $V$ is a unitary operator and of course it also
follows, too, that $V$ intertwines $\tilde{U}$ and $\tilde{U}_h$.
\end{proof}
\begin{proposition}\label{factor->multiple.of.irr.}
A positive energy factor representation of $\widetilde{\mob}$ is
equivalent with a multiple of an irreducible representation.
\end{proposition}
\begin{proof}
By Prop. \ref{factor.rep.of.mob} we know that the spectrum of the
self-adjoint generator of rotations $H$ in the factor
representation $\tilde{U}$ contains only eigenvalues. Let $h$ be
the lowest eigenvalue and $E_\pm$ as before. By the same argument
we employed in the proof of the previous proposition, if we pick
an eigenvector $\Phi\in {\rm Ker}(H-h\mathbbm 1)$, then the
subspace
\begin{equation}
M_\Phi \equiv \{(E_+)^n\Phi: n\in\NN\}
\end{equation}
where we have set $(E_+)^0\equiv \mathbbm 1$, is invariant and
analytic for the self-adjoint generators of $\tilde{U}$. Thus its
closure is an invariant for $\tilde{U}$. In fact by what was
explained before it is a minimal invariant closed subspace for
$\tilde{U}$. Moreover, by the fact that $E_- {\rm Ker}(H-h\mathbbm
1)$ and by the commutation relations it is clear that if
$\Phi_1,\Phi_2\in {\rm Ker}(H-h\mathbbm 1)$ then
\begin{equation}
\langle (E_+)^n\Phi_1,\,(E_+)^m \Phi_2\rangle = \delta_{n,m} a(n)
\langle \Phi_1,\, \Phi_2\rangle
\end{equation}
where $a(n)\in\CC$ is independent of the vectors $\Phi_1,\,
\Phi_2$. Thus if these vectors are perpendicular then
$\overline{M_{\Phi_1}}$ is perpendicular to
$\overline{M_{\Phi_2}}$. Hence by choosing a complete orthonormed
system $\{\Phi_\alpha: \alpha \in I\}$ in ${\rm Ker}(H-h\mathbbm
1)$ (where $I$ is some index set), we have that
$\overline{M_{\Phi_\alpha}}$ $(\alpha \in I)$ are pairwise
orthogonal minimal invariant closed subspaces and the closure of
their sum contains ${\rm Ker}(H-h\mathbbm 1)$. The orthogonal
$M^\perp$ to the subspace generated by $\{M_{\Phi_\alpha}:\alpha
\in I\}$ is non-trivial, then, by similar arguments one finds that
it contains a nontrivial irreducible representation with some
lowest eigenvalue for the self-adjoint generator of rotations
bigger then $h$. By factoriality this is not possible and thus the
statement is proved.
\end{proof}
\begin{corollary}\label{mobrep->dirint}
A positive energy representation $\tilde{U}$ of $\widetilde{\mob}$
is always equivalent with a direct integral of irreducible
representations. $\tilde{U}$ is equivalent with a direct sum of
irreducible representations if and only if the span of the
eigenvalues of the self-adjoint generator of rotations is dense.
\end{corollary}
\begin{proof}
The first statement follows directly from the previous proposition
and Lemma \ref{mob.directint}. As for the second statement: one
direction is trivial since this property is true for every
irreducible positive energy representation. For the converse
direction, suppose $\Phi$ is an eigenvector of the self-adjoint
generator of rotations $H$ and set $E_\pm$ as before. If $E_-\Phi
= 0$ then, using the notations of the proof of the previous
proposition, $\overline{M_\Phi}$ is minimal invariant closed
subspace. If $\Phi$ is orthogonal to all vectors that are
annihilated by $E_-$ and $E_-^2\Phi = 0$, then it is easy to see
that $\overline{M_{E_-\Phi}}$ is a minimal invariant closed
subspace containing $\Phi$. Similarly, if $\Phi$ is orthogonal to
all vectors that are annihilated by $E_-^2$ and $E_-^3\Phi = 0$
then we have that $\overline{M_{E^2_-\Phi}}$ is a minimal
invariant closed subspace containing $\Phi$. Since by the
``stepping'' property of $E_-$ and the positivity of $H$ the
eigenvector $\Phi$ must be annihilated by some (large enough,
positive) power of $E_-$, carrying on the sketched argument we
find that $\Phi$ is a finite sum of orthogonal vectors lying in
minimal invariant closed subspaces (in fact, the number of the
vectors needed in this decomposition is $\leq$ the eigenvalue
associated to $\Phi$ plus $1$). It is clear that together with
what was established before how this implies the statement.
\end{proof}

\chapter{Piecewise \mo transformations}
\markboth{APPENDIX B. PIECEWISE M\"OBIUS TRANSFORMATIONS}{APPENDIX
B. PIECEWISE M\"OBIUS TRANSFORMATIONS}
\section{Definitions and
properties}

We begin with some rather trivial affirmations.
Suppose that $\kappa:S^1\rightarrow S^1$ is a differentiable
bijection whose derivative $\kappa'$ is strictly positive.
Then by the inverse function theorem also $\kappa^{-1}$
is differentiable and its derivative is strictly positive.
Moreover, if $\kappa_1$ and $\kappa_2$ are differentiable bijections
with strictly positive derivative, then $\kappa_1 \circ \kappa_2$
is also a differentiable bijection with strictly positive
derivative.
\begin{proposition}
\label{pcwgroup}
If $\gamma \in \pcw$ then $\gamma'>0$ (and so it is
orientation preserving) and $\gamma^{-1} \in \pcw$.
Moreover, if $\gamma_1,\gamma_2 \in \pcw$ then
$\gamma_1\circ\gamma_2 \in \pcw$. Thus the set $\pcw$
with the composition as multiplication forms a group.
\end{proposition}
\begin{proof}
Suppose $\gamma \in \pcw$ coincides with $\phi \in \mob$ on $I \in
\I$. Then by continuity they even coincide on $\overline{I}$.
Moreover, as the derivative on $\overline{I}$ is completely
determined by the function on $\overline{I}$ (even on the
endpoints of the interval), we have that $\gamma'|_{\overline{I}}
= \phi'|_{\overline{I}}$. But a \mo transformation, being in
particular an orientation preserving diffeomorphism, has a
strictly positive derivative. Hence also the derivative of
$\gamma$ must be strictly positive at each point. The rest of the
statement is rather obvious, as also the \mo transformations, with
the composition as multiplication, form a group.
\end{proof}
\smallskip

{\it Remark.} The group $\pcw$ at least contains the \mo
transformations. But so far we have not addressed the question
whether there are nontrivial piecewise \mo transformations:
whether the group $\pcw$ is really bigger than $\mob$. Could we
have required for example twice differentiability, in stead of
once differentiability? No, since by Lemma \ref{2diff->mob} there
are no nontrivial, twice differentiable piecewise \mo
transformations: such a transformation is in fact automatically
\mo. However, ``luckily'' we have only required once
differentiability, which, as we shall soon see, gives rise to a
completely different situation.

Recall that two intervals $I_1,I_2 \in \I$ are said to be {\it distant},
if $\overline{I_1} \cap \overline{I_2} = \emptyset$. Recall also
that for each interval $I \in \I$ there is an associated
one-parameter group $s\mapsto \delta^I_s \in \mob$ which is
called the one-parameter group of dilations associated to $I$,
and that by Def.\! \ref{def:kappa}, a certain transformation 
$\kappa^{I_1,I_2}_s$ was introduced for each pair of distant intervals $I_1,I_2\in \I$ and
real number $s\in\RR$. As a dilation of a certain interval leaves fixed 
the endpoints of that interval, it is not hard to see that
the transformation $\kappa^{I_1,I_2}_s$ is at least continuous.
We shall now see that it is in fact a differentiable
bijection.
\begin{lemma}
\label{kappainpcw}
Let $I_1,I_2 \in \I$ be two distant intervals.
Then $s\mapsto \kappa^{I_1,I_2}_t$ is a one-parameter group in
$\pcw$. Moreover, if $s\neq 0$ then there are $4$ points where
$\kappa^{I_1,I_2}_s$ is not smooth, namely the endpoints of the
intervals $I_1$ and $I_2$. Thus if $s \neq 0$ then
$\delta^{I_1,I_2}_s \notin \mob$.
\end{lemma}
\begin{proof}
It is evident that $s\mapsto \kappa^{I_1,I_2}_s$ is a
one-parameter group, so in particular $\kappa^{I_1,I_2}_s$ is a
bijection for all $s\in\RR$. By its construction it is clearly a
piecewise \mo transformation if it is differentiable. Of course
the only problem with differentiability can happen at the
endpoints of the intervals appearing in the definition of the
transformation.

Recall formula \ref{formXdilations} for the dilations
$s\mapsto \delta_s$ of the interval $S^1_+$ and the convention
(identification) we use throughout this work according to which
the derivative of a $\gamma:S^1\rightarrow S^1$ transformation is
an $S^1\rightarrow \RR$ function. By direct calculation
\begin{eqnarray}\nonumber
\label{delta'} \delta_s'(z) &\equiv& -i
\overline{\delta_s(z)}\frac{d}{d\theta}\delta_s(e^{i\theta})|_{e^{i\theta}=z}\\
&=& \frac{(1+e^s)z}{(1+e^s)z+(1-e^s)}-
\frac{(1-e^s)z}{(1-e^s)z+(1+e^s)}.
\end{eqnarray}

Let $I\in \I$ be an interval with endpoints $z_+$ and $z_-$
where the $\pm$ index is given to the two endpoints in such
a way that if $p\in I$ then $z_+,p,z_-$ are three points in
the anticlockwise order. Then there exists a $\gamma \in \mob$
such that $\gamma(S^1_+)=I$. Thus $\gamma(\pm 1)=z_\pm$,
$\delta^{I}_s = \gamma\circ\delta_s\circ\gamma^{-1}$ and
using that $\delta(\pm 1)=\pm 1$ we have
\begin{eqnarray}
\label{delta'(endpoint)}
\nonumber
(\delta^{I}_s)'(z_\pm)&=&(\gamma\circ\delta_s\circ\gamma^{-1})'(z_\pm)
\\ \nonumber
&=&\gamma'(\pm 1)\, \delta_s'(\pm 1) \, (\gamma^{-1})'(\gamma(\pm 1))
\\ &=&\delta_s'(\pm 1)=e^{\pm s}
\end{eqnarray}
where the last equality we get by substituting in (\ref{delta'}).
By applying this result to our case we see that $\kappa^{I_1,I_2}_s$
is indeed differentiable.

Finally, consider the second derivative of $\kappa^{I_1,I_2}_s$ at the common
endpoint $z_{11}$ of the intervals $I_1$ and $K_1$, for example.  If it
exists then $(\delta^{I_1}_s)''(z_{11})=(\delta^{K_1}_{-s})''(z_{11})$
which by (\ref{2diff->mob}) (considering that also the transformations
and their first derivative coincide at that point) implies that
$\delta^{I_1}_s=\delta^{K_1}_{-s}$. Hence $\delta^{I_1}_s$ preserves
not only the endpoints of $I_1$ but also those of $K_1$, which by Lemma
\ref{mob.3.trans} is only possible if $s=0$.
\end{proof}

A piecewise \mo transformation is of course always smooth apart
from a finite number of points. If it is everywhere smooth then it
must be a \mo transformation. The examples of piecewise \mo
transformations which we just gave are not smooth at $4$ points.
The next lemma says that in fact in some sense these examples are
the smallest possible.
\begin{lemma}
\label{3points}
Let $\gamma \in \pcw$ be a transformation with at most $3$ points
where it is not smooth. Then $\gamma \in \mob$.
\end{lemma}
\begin{proof}
Suppose that the transformation $\gamma$ is everywhere smooth
apart from the three (different) points $z_1,z_2,z_3 \in S^1$
where the numbering of the points corresponds to the anticlockwise
direction on the circle. By Prop.\! \ref{pcwgroup} $\gamma$ is
orientation preserving, so $\gamma(z_1),\gamma(z_2),\gamma(z_3)$
are still three (different) points listed in the anticlockwise
direction. Therefore, by the $3$-transitivity of the action of the
\mo group (Lemma \ref{mob.3.trans}), there exists a $\phi \in
\mob$ such that $\phi(\gamma(z_j))=z_j$ for $j=1,2,3$.

Let $I_j \in \I$ $j=1,2,3$ be the interval that does not contain
the point $z_j$ and has the other two points (out of $z_1,z_2$ and
$z_3$) as endpoints. The piecewise \mo transformation
$\phi\circ\gamma$ preserves each of these intervals, and it is
smooth on each of these intervals. Therefore the restriction of
$\phi\circ\gamma$ on each of these intervals must correspond to a
\mo transformation; see the remarks made after Prop.\!
\ref{pcwgroup}. But any \mo transformation that preserves an
interval $I\in \I$ must be a dilation of that interval with a
certain parameter (Lemma \ref{uniq.dil.of.S+}). So there exist
three real numbers $s_1,s_2,s_3$ such that
\begin{equation}
(\phi\circ\gamma)|_{I_j}=\delta^{I_j}_{s_j}|_{I_j}
\;\;{\rm for}\;j=1,2,3.
\end{equation}
By equation (\ref{delta'(endpoint)}), the existence of the first
derivative at the point $z_1$ implies that $s_2=-s_3$. Similarly, by
the existence of the first derivative at points $z_2$ and $z_3$
we get that $s_3=-s_1$ and $s_1=-s_2$. Thus the only solution
is $s_1=s_2=s_3=0$ and so $\phi\circ\gamma = {\rm id}_{S^1}$
and $\gamma = \phi^{-1} \in \mob$.
\end{proof}
We can now present a generating set of elements in $\pcw$.
\begin{theorem}
\label{pcwgenerators}
The group $\pcw$ is generated by the set
$${\mathcal S}\equiv\mob \cup\{\rho^{I_1,I_2}_s\}_ {\{I_1, I_2 \in \I,\,
\overline{I_1}\cap\overline{I_2}=\emptyset,\, s\in \RR\}}.$$
\end{theorem}
\begin{proof}
Let $\gamma \in \pcw$. We shall show that $\gamma \in <\!{\mathcal
S}\!>$ by an inductive argument regarding the number of points at
which $\gamma$ is not smooth.

If $\gamma$ is not smooth at most three points then by Lemma
\ref{3points} it is actually a \mo transformation and so $\gamma
\in <\!{\mathcal S}\!>$. To proceed with the induction, assume
that if a piecewise \mo transformation is not smooth at most
$n\in\{3,4,..\}$ points then it indeed belongs to $<\!{\mathcal
S}\!>$. Now suppose that $\gamma$ is not smooth exactly at the
$n+1$ points $z_1,z_2,..z_n,z_{n+1}$, where the numbering
corresponds to the anticlockwise order.

Let us use the notation $I_j$ $(j=1,2,3)$ for the interval in $\I$
that does not contain any of the points where $\gamma$ is not
smooth, and has $z_j$ and $z_{j+1}$ as endpoints. Just as in the
proof of Lemma \ref{3points}, we may "adjust" $\gamma$ by a \mo
transformation in such a way that it will preserve $z_1,z_2$ and
$z_3$; i.e.\! we may assume that $\gamma(z_j)=z_j$ for $j=1,2,3$.
This means that $\gamma$ preserves the intervals $I_1$ and $I_2$
and thus, just like in the mentioned proof, we can deduce that
there exist $s_1,s_2 \in \RR$ such that
$\gamma|_{I_j}=\delta^{I_j}_{s_j}|_{I_j}$ for $j=1,2$. By equation
(\ref{delta'(endpoint)}) the differentiability of $\gamma$ at the
point $z_2$ implies that $s\equiv s_1=-s_2$. Then it is easy to
see that $\kappa^{I_1,I_3}_{-s}\circ\gamma$ acts identically on
both $I_1$ and $I_2$, and so there can be at most $n$ points at
which this composition is not smooth (this transformation can be
not smooth only at the points $z_1,..,z_{n+1}$ with the exception
of the point $z_2$). Therefore, by the inductive assumption we may
conclude that $\gamma \in <\!{\mathcal S}\!>$.
\end{proof}

We shall finish our investigation about the geometrical properties
of the piecewise \mo transformation by showing that there are
indeed ''many'' such transformation. In particular, we shall prove
that the action of the group $\pcw$ on $S^1$ is $N$-transitive for
any $N$ positive integer. In relation with local nets it is also
important that these transformations can be even required to be
{\it local}. First we shall prove the following lemma.
\begin{lemma}
\label{1localtransitivity} Let $z_0,z_0'\in \I$ be two points in
the interval $I\in \I$. Then there exists a piecewise \mo
transformation $\gamma \in \pcw$ such that $\gamma|_{I^c}={\rm
id}_{I^c}$ and $\gamma(z_0)=z_0'$
\end{lemma}
\begin{proof}
It is clear that we can choose a smaller open interval $I_1\subset
I$ which is distant from $I^c$ but still contains the given two
points. We shall denote the two distant intervals forming the
interior of the complement of $I_1\cup I^c$ by $K_+$ and $K_-$.
Note that $K_\pm \subset I$.

For a real number $s\in \RR$ set $\beta_s\equiv
\kappa^{I_1,I^c}_s\circ\delta^{I}_s \in \pcw$. As
$\delta^{I}_s=\delta^{I^c}_{-s}$, by the definition of
$\kappa^{I_1,I^c}$, the transformation $\beta_s$ acts
identically on $I^c$ for all $s\in\RR$. Moreover, as the maps
$(s,z)\mapsto \delta^{I}_s(z)$ and $(s,z)\mapsto
\kappa^{I_1,I^c}_s(z)$ are clearly continuous, so is the
the map $s\mapsto \beta_s(z_0)$ and thus its image is a
connected set in $S^1$. In fact, as $\beta$ preserves
the interval $I$ (as it acts identically on $I^c$),
it must be a connected subset of $I$; i.e.\! an interval
included in $I$.

Of course the image of $s\mapsto\delta^{I}_s(z_0)$ is the entire
interval $I$. Thus we can choose two values $s_+,s_-\in\RR$ such
that $\delta^{I^c}_{s_\pm}(z_0)\in K_\pm$. Since
$\kappa^{I_1,I^c}$ preserves the interval $K_\pm$ we have that
$\beta_{s_\pm}(z_0)\in K_\pm$. Then by what was said about the
image we can conclude that $I_1\subset \beta_{\RR}(z_0)$ and so
there exist a value $s_0\in \RR$ such that
$\beta_{s_0}(z_0)=z_0'$. Thus the piecewise \mo transformation
$\gamma\equiv \beta_{s_0}$ satisfies all what was required.
\end{proof}

Recall that for two points of the circle $p,q\in S^1$ we denoted
by $[p,q]$ the closed interval whose first point (in the
anticlockwise direction) is $p$ and last point is $q$. (We set
$[p,q]$ to be the set containing the single point $p$ if $p=q$).
\begin{theorem}
\label{localntransitivity}
Let $\{z_n\}_{n=1}^N$ and $\{z'_n\}_{n=1}^N$ be two sets of
$N$ (different) points in $S^1$ where $N$ is any positive
integer and the numbering of the points corresponds to the
anticlockwise order. Then there exists a $\gamma \in \pcw$
such that
$$
\gamma(z_n)=z'_n \;\;{\rm for}\; n=1,..,N.
$$
Moreover, if there is an interval $I\in \I$
containing both $[z_1,z_N]$ and $[z'_1,z'_N]$ then
$\gamma$ can be chosen to be localized in $I$;
i.e.\! satisfying $\gamma|_{I^c}={\rm id}|_{I^c}$.
\end{theorem}
\begin{proof}
Clearly the ``local'' version of the statement is stronger; by
applying a \mo transformation we can move the points
$z'_1,..,z_N'$ so that we will have $[z_1,z_N]=[z_1',z_N']$ in
which case of course we can always find an $I\in\I$ containing
both $[z_1,z_N]$ and $[z'_1,z'_N]$.

We shall argue by induction. For $N=1$ by Lemma
\ref{1localtransitivity} we are finished. Now suppose the
statement holds for the positive integer $N=k$ and consider the
case $N=k+1$. By our assumption there exists a $\gamma_1\in \pcw$
localized in $I$ such that $\gamma_1(z_n)=z'_n$ for $n=1,..,k$. As
$\gamma_1$ is orientation preserving and local, both $z'_{k+1}$
and $\gamma_1(z_{k+1})$ must come ``after'' (in the anticlockwise
direction) the points $z'_1,..,z'_k$ and still in $I$. Thus there
exists a smaller open interval $K\subset I$ containing both
$z'_{k+1}$ and $\gamma_1(z_{k+1})$ but not containing the points
$z'_1,..,z'_k$. By Lemma \ref{1localtransitivity} there exists a
$\gamma_2 \in\pcw$ localized in $K$ such that
$\gamma_2(\gamma_1(z_{k+1}))=z_{k+1}$. Thus we are finished as the
piecewise \mo transformation $\gamma_2\circ\gamma_1$ is still
localized in $I$ and clearly it moves points $z_1,..,z_{k+1}$ into
the points $z'_1,..,z'_{k+1}$.
\end{proof}

\section{Piecewise \mo vector fields}

Recall the notion of piecewise \mo vector fields (Definition
\ref{def:pcwm}). We begin with some elementary observation. First
of all, by its definition $\pcwm$ is a vector space: if $\lambda$
is a real number and $f_1,f_2\in \pcwm$ then $(\lambda f_1+f_2)
\in \pcwm$. Also, clearly the adjoint action of a \mo
transformation sends a piecewise \mo vector field into a piecewise
\mo vector field: if $f\in \pcwm$ and $\phi\in\mob$ then $(\phi'
f)\circ\phi^{-1} \in \pcwm$. As it was observed, a piecewise \mo
vector field is globally Liepschitz and thus by the compactness of
the circle it determines a unique flow. However, the temptation to
think of $\pcwm$ as the Lie algebra of the group $\pcw$, having in
mind the example of ${\rm Vect}(S^1)$ which is identified with the
Lie algebra of $\diff$, is misleading. First of all, the group of
piecewise \mo transformations is by no means a Lie group. In
general, the bracket of two piecewise \mo vector fields is not
differentiable, so it does not remain in $\pcwm$. Second of all,
the exponential of a piecewise \mo vector field does not need to
be a piecewise \mo transformation. On the other hand, by
differentiating a ``sufficiently regular'' one-parameter group of
$\pcw$ we indeed get an element of $\pcwm$.

\begin{lemma}
\label{aII}
\label{aI1I2z} Let $I_1,I_2 \in \I$ be two distant intervals. Then
the function $a^{I_1,I_2}$ given by the formula $z\mapsto
a^{I_1,I_2}(z)\equiv\frac{d}{ds}\kappa^{I_1,I_2}_s(z)|_{s=0}$ is
well-defined (i.e.\! the derivative exists) and it is an element
of $\pcwm$.
\end{lemma}
\begin{proof}
That this derivative indeed everywhere exists (i.e.\! even at the
endpoints of the intervals $I_1,I_2$), is trivial since the
endpoints remain fixed for all $s\in \RR$. It is also evident that
our function on each of the four open intervals determined by
these endpoints is the restriction of a \mo vector field. The only
thing to be verified is that the derivative of $a^{I_1,I_2}$ from
the left and right coincide at these four points.

Recall that the generating \mo vector field $d$ for the
one-parameter group of dilations associated to $S^1_+$ has the
expression $d(z)=\frac{-i}{2}(z^1-z^{-1})$. Thus by setting
$z=e^{i\theta}$
\begin{equation}
\label{d'(endpoint)}
d'(z)\equiv \frac{d}{d\theta} d(e^{i\theta}) =
\frac{1}{2}(z^1+z^{-1}).
\end{equation}
Let $I\in \I$ be an interval with endpoints $z_+$ and $z_-$ where
the $\pm$ index is given to the two endpoints in such a way that
if $p\in I$ then $z_+,p,z_-$ are three points in the anticlockwise
order. Then there exists a $\gamma \in \mob$ such that
$\gamma(S^1_+)=I$. Thus $\gamma(\pm 1)=z_\pm$, and the \mo vector
field $d^{I}$ generating the dilations associated to the interval
$I$ has the expression $d^{I} = (\gamma' d )\circ\gamma^{-1}$.
Thus using that $d(\pm 1)=0$ we have
\begin{eqnarray}
\label{dI'(endpoint)}
\nonumber
(d^{I})'(z_\pm)&=&(\gamma'\circ\gamma^{-1})'(z_\pm)\,d(\gamma(z_\pm)+
\gamma'(\gamma^{-1}(z_\pm))\,d'(\gamma(z_\pm))\gamma'(z_\pm)
\\ \nonumber
&=& d'(\pm 1) = \pm 1
\end{eqnarray}
where the last step is calculated by directly substituting into
equation \ref{d'(endpoint)}. Considering the result received one
can immediately see the differentiability of $a^{I_1,I_2}$.
\end{proof}

A piecewise \mo vector field is always smooth apart from a finite
set of points. If it is really smooth then it is evidently a \mo
transformation. We shall now state the counterpart of Lemma
\ref{3points} for the piecewise \mo vector fields.
\begin{lemma}
Let $f\in\pcwm$. If $f$ is not smooth at most $3$ points
then it is actually smooth and so $f\in\mathfrak m$.
\end{lemma}
\begin{proof}
Let us assume that $f$ is smooth apart from the three points
$z_1,z_2,z_3$ where the numbering corresponds to the anticlockwise
order. We shall denote by $I_j$ $(j=1,2,3)$ the element of $\I$
that does not contain the point $z_j$ and has the other two points
(out of $z_1,z_2,z_3$) as endpoints. By our assumption $f$ is the
restriction of a \mo vector field on each of the three intervals
$I_1,I_2$ and $I_3$.

By Lemma \ref{3transit.of.vect.m} there exists a unique \mo vector
field $g \in \mathfrak m$ such that $g(z_j)=f(z_j)$ for $j=1,2,3$.
Thus $(f-g)(z_j)=0$ $(j=1,2,3)$ and hence by still by the cited
Lemma, for each $j=1,2,3$ we have that $(f-g)|_{I_j}$ coincides
with the restriction of $\lambda_j d^{I_j}$ for some real multiple
$\lambda_j \in \RR$. However, by looking at equation
(\ref{dI'(endpoint)}) we see that the condition of once
differentiability at the three endpoints points implies that
$\lambda_1=-\lambda_2$, $\lambda_2=-\lambda_3$ and finally that
$\lambda_3=-\lambda_1$. Thus $\lambda_1=\lambda_2=\lambda_3=0$ and
$f=g \in \mathfrak m$.
\end{proof}

Exactly as in the case of \mo transformations, by using the
previous lemma we can prove that the examples of piecewise \mo
vector fields we have so far exhibited form a generating set.
\begin{theorem}
\label{pcwmgenerators}
$\pcwm = {\rm Span}\left\{{\mathfrak m}\cup
\{a^{I_1,I_2}\}_{\{I_1, I_2 \in \I,\,
\overline{I_1}\cap\overline{I_2}=\emptyset\}}\right\}$.
\end{theorem}
\begin{proof}
Let $f \in \pcwm$ be smooth apart from the points $z_1,..,z_N$
where $N$ is a positive integer and the numbering of the points
corresponds to the anticlockwise order. If $N\leq 3$ then by the
previous lemma $f\in \mathfrak m$ and so we are finished. We shall
proceed by induction. Assume that $f$ belongs to the Span in
question whenever $N\leq k$ where $k=3,4,..$ and suppose that
$N=k+1$. Just like in the previous lemma, we may assume that
$f(z_j)=0$ for the first three points $z_1,z_2,z_3$ (if it was not
so, we may change $f$ by adding a suitable \mo vector field to it
in such a way, that the new $f$ will have this property). Then
using equation (\ref{dI'(endpoint)}) and the same idea which was
used in the proof of the previous lemma, it is easy to see that
there exists a $\lambda \in \RR$ such that $\tilde{f}\equiv
f-\lambda a^{I_1,I_2} = 0$ on $[z_1, z_3]$ where $I_1,I_2\in \I$
are the two distant intervals with endpoints $(z_1,z_2)$ and
$(z_3,z_4)$, respectively. Then $\tilde{f}\in \pcwm$ can have at
most $k$ points where it is not smooth and thus by the assumption
of the inductive argument we are finished with our proof.
\end{proof}
As we have seen in Sect.\! \ref{sec:nonsmooth1}, from the point of
view of representation theory of $\diff$ it is useful to consider
the topology on ${\rm Vect}(S^1)$ given by the
$\|\cdot\|_{\frac{3}{2}}$ norm. It is therefore important to
investigate the relationship between this norm and the piecewise
\mo vector fields. It has been already mentioned that every
piecewise \mo vector field has a finite $\|\cdot\|_{\frac{3}{2}}$
norm (Corollary \ref{pcwmfin1.5norm}). The following lemma is the
first step in deriving a density result.
\begin{lemma}
\label{pcwmdiracdelta} Let $I\in \I$. Then there exists an
$f\in\pcwm$ such that $f\geq 0$, ${\rm Supp}(f)\subset I$ and
$\frac{1}{2\pi}\int_0^{2\pi}f(e^{i\theta})d\theta = 1$.
\end{lemma}
\begin{proof}
It is enough to find an example of a piecewise \mo vector field
which is localized {\it somewhere} and nonnegative: then by
transforming it with the adjoint action of a suitable \mo
transformation and normalizing it we get what we want.

We know that for any $I\in \I$ the function $d^I$ is strictly
positive on $I$ and negative on $I^c$. Also, if $K\subset I$ is a
nonempty, non-dense open subinterval of $I$ having one endpoint in
common with $I$, then $d^I-d^K$ is zero on the common endpoint and
strictly positive elsewhere. This is because $d^I-d^K$ must be a
positive multiple of the generator of the translations that have
the common endpoint as infinite point.

Let therefore $I_1,I_2 \in \I$ be two distant intervals with
$K_1,K_2 \in \I$ being the other two distant intervals determined
by the endpoints of $I_1,I_2$. Then $a^{I_1,I_2} - d^{I_1} \in
\pcwm$ is a nonnegative function which is zero on $I_1$ and strictly
positive on $K_1,K_2$ and $I_2$ and thus the lemma is proven.
\end{proof}

\begin{theorem}
\label{pcwmdensity}
 Let $f$ be a continuous real function on $S^1$
with finite $\|\cdot\|_{\frac{3}{2}}$ norm. Then there exists
sequence $\{f_n\}_{n=1}^\infty \in \pcwm$ converging to $f$ in the
$\|\cdot\|_{\frac{3}{2}}$ sense. Moreover, if $f$ is localized in
the interval $I\in\I$ (i.e.\! ${\rm Supp}(f)\subset I$) then the
sequence can be chosen so that every element of it is localized in
$I$.
\end{theorem}
\begin{proof}
We shall make the proof in two steps. First, let $l\in \pcwm$ any
piecewise \mo transformation. We shall prove that $(l*f)$, the
convolution defined by the formula
\begin{equation}
(f*l)(z) \equiv \frac{1}{2\pi}\int_0^{2\pi}
f(e^{i\theta})l(e^{-i\theta}z) d\theta
\end{equation}
can be approximated by elements of $\pcwm$. (Note that
$\|(f*l)\|_{\frac{3}{2}}<\infty$.) Indeed, by setting
$l_{\theta}(z)=l(e^{-i\theta}z)$ it is clear that $l_{\theta} \in
\pcwm$. Thus
\begin{equation}
h_N \equiv \frac{1}{2N\pi}\sum_{k=0}^{N-1} f(e^{i\frac{k}{2\pi}})
l_{e^{i\frac{k}{2\pi}}}\, \in \pcwm
\end{equation}
for all positive integer $N$. Then by expressing the Fourier
coefficients of $h_N$ in terms of the Fourier coefficients
$\{\hat{l}_n\}_{\{n\in\ZZ\}},\{\hat{f}_n\}_{\{n\in\ZZ\}}$ of
the functions $l$ and $f$ we find that
\begin{equation}
\|(f*l) - h_N\|_{\frac{3}{2}}= \sum^{n\in\ZZ} |\hat{l}_n|\,
|\hat{f}_n-\frac{1}{2N\pi}\sum_{k=0}^N
e^{-in\frac{k}{2\pi}}f(e^{i\frac{k}{2\pi}})|\,(1+|n|^{\frac{3}{2}}).
\end{equation}
Since $f$ is continuous, for every $n$ integer
$|\hat{f}_n-\frac{1}{2N\pi}\sum_{k=0}^N
e^{-in\frac{k}{2\pi}}f(e^{i\frac{k}{2\pi}})|$ is bounded by
$2({\rm max}\{|f|\})$ and goes to zero as $N\to\infty$. Moreover,
by Corollary\! \ref{pcwmfin1.5norm} $\sum_{n\in\ZZ} |\hat{l}_n|
(1+|n|^{\frac{3}{2}}) =\|l\|_{\frac{3}{2}}<\infty$. We can thus
conclude that $h_N \to (f*l)$ as $N\to\infty$ in the
$\|\cdot\|_{\frac{3}{2}}$ norm. Note that if ${\rm Supp}(f)\subset
I$, then if the support of $l$ is a sufficiently small interval
around the point $1 \in S^1$ then ${\rm Supp}(h_N)\subset I$ for
all positive integer $N$.

Let $I_n\in\I$ $(n=1,2,..)$ be a descending chain of intervals
contracting to the point $1\in S^1$; that is $I_j \supset I_k$ for
$j<k$ and $\cap_{n=1}^{\infty}I_n=\{1\}$. By Lemma
\ref{pcwmdiracdelta} for each $n$ positive number there exists a
nonnegative function $l_n$ such that ${\rm Supp}(l_n)\subset I_n$
and $\frac{1}{2\pi}\int^{2\pi}_0 l_n(e^{i\theta})d\theta =1$. Then
it is easy to see that as $n\to\infty$ the convolution $(f*l_n)$
converges to $f$ in the $\|\cdot\|_{\frac{3}{2}}$ norm. Moreover,
if ${\rm Supp}(f)\subset I$ then for $n$ large enough ${\rm
Supp}(f*l_n)\subset I$. Thus by further approximating the
convolution with the sum we considered before one can easily
finish the proof.
\end{proof}

\end{document}